\documentclass[aos,reqno]{imsart}

\RequirePackage{amsthm,amsmath,amsfonts,amssymb}
\RequirePackage[numbers]{natbib}
\RequirePackage[colorlinks,citecolor=blue,urlcolor=blue]{hyperref}
\RequirePackage{graphicx}
\usepackage{algorithm2e}
\usepackage{subfiles}
\usepackage{titletoc}
\usepackage{bbm}
\usepackage{multirow}
\usepackage{pdfpages}
\usepackage{mathtools}
\usepackage{algorithm2e}
\usepackage{booktabs}
\RestyleAlgo{ruled}
\SetKwComment{Comment}{/* }{ */}
\usepackage{hyperref}
\usepackage{url}
\usepackage{breqn}
\usepackage{graphicx}
\usepackage{tipa}
\usepackage{bm}
\usepackage{paralist}
\usepackage{calrsfs}
\usepackage{subfiles}
\usepackage{xcolor}
\usepackage{mwe}
\usepackage{dutchcal}
\usepackage{threeparttable}
\usepackage[labelfont=bf]{caption}

\startlocaldefs
\theoremstyle{plain}

\newtheorem{theorem}{Theorem}[section]
\newtheorem{lemma}[]{Lemma}
\newtheorem{remark}{Remark}
\newtheorem{conjecture}{Conjecture}[section]
\theoremstyle{remark}
\newtheorem{definition}[]{Definition}

\newtheorem*{fact}{Fact}

\newcommand{\bb}[1]{\mathbb{#1}}

\newcommand{\wh}[1]{\widehat{#1}}

\newcommand{\bfa}[1]{\boldsymbol{#1}}

\DeclareMathAlphabet{\pazocal}{OMS}{zplm}{m}{n}
\newcommand{\ca}[1]{\pazocal{#1}}
\newcommand{\mca}[1]{\mathcal{#1}}

\newcommand{\nnb}{\nonumber}

\newcommand{\lef}{\left}
\newcommand{\rig}{\right}

\newcommand{\tda}[1]{\tilde{\bfa{#1}}}

\newcommand{\pta}{\partial}

\newcommand{\mbbm}[1]{\mathbbm{#1}}

\newcommand{\tde}{\tilde}

\newcommand{\la}{\langle}
\newcommand{\ra}{\rangle}

\newcommand{\ub}[1]{\underbrace{#1}}

\newcommand{\bl}{{\bigg (}}
\newcommand{\br}{{\bigg )}}

\newcommand{\sm}[1]{\small{#1}\normalsize}
\newcommand{\tny}[1]{\footnotesize{#1}\normalsize}
\newcommand{\ttny}[1]{\tiny{#1}\normalsize}
\DeclareMathOperator*{\argmax}{arg\,max}
\DeclareMathOperator{\sech}{sech}
\DeclareMathOperator*{\argmin}{arg\,min}

\endlocaldefs

\begin{document}

\begin{frontmatter}
\title{Hidden Clique Inference in Random Ising Model II:
the planted Sherrington-Kirkpatrick model}
\runtitle{Hidden Clique Inference in the pSK Model}

\begin{aug}
\author[A]{\fnms{Yihan}~\snm{He}\ead[label=e1]{yihan.he@princeton.edu}},
\author[B]{\fnms{Han}~\snm{Liu}\ead[label=e2]{hanliu@northwestern.edu}}
\and
\author[A]{\fnms{Jianqing}~\snm{Fan}\ead[label=e3]{jqfan@princeton.edu}}
\address[A]{Princeton University\printead[presep={,\ }]{e1,e3}}
\address[B]{Northwestern University\printead[presep={,\ }]{e2}}
\end{aug}

\begin{abstract}
We study the problem of testing and recovering $k$-clique Ferromagnetic mean shift in the planted Sherrington-Kirkpatrick model (i.e., a type of spin glass model) with $n$ spins. The planted SK model --- a stylized mixture of an uncountable number of Ising models --- allows us to study the fundamental limits of correlation analysis for dependent random variables under misspecification. Our paper makes three major contributions: (i) We identify the phase diagrams of the testing problem by providing minimax optimal rates for multiple different parameter regimes. We also provide minimax optimal rates for exact recovery in the high/critical and low-temperature regimes. (ii) We prove a universality result implying that all the obtained rates still hold with non-Gaussian couplings. (iii) To achieve the major results, we also establish a family of novel concentration bounds and central limiting theorems for the averaging statistics in the local and global phases of the planted SK model. These technical results shed new insights into the planted spin glass models. The pSK model also exhibits close connections with a binary variant of the single-spike Gaussian sparse principle component analysis model by replacing the background identity precision matrix with a Wigner random matrix.
\end{abstract}

\begin{keyword}[class=MSC]
\kwd[Primary ]{62H22}
\kwd{62H25}
\kwd[; secondary ]{62F03}
\end{keyword}

\begin{keyword}
\kwd{Minimax Optimality}
\kwd{Exact Recovery}
\kwd{Mean Field Ising Models}
\kwd{Sparse Principle Component Analysis}
\end{keyword}

\end{frontmatter}

\section{Introduction}\label{sect1}
We study the problem of testing and recovering clique-shaped positive mean shifts in the coupling matrix of the planted SK (or pSK) model. Let $\bfa\sigma \in\Sigma_n:= \{-1, +1\}^n$ be $n$ binary random variables (or spins), 
the pSK model specifies the distribution of $\bfa\sigma$ as
\begin{align}\label{skdirect}
\bb P(\bfa\sigma|\bfa g,\bfa h)\propto {\exp\bl  \frac{\theta_1}{2k} \sum_{1\leq i< j,i,j\in S}\sigma_i\sigma_j+\frac{\theta}{\sqrt n}\sum_{1\leq i<j\leq n}g_{ij}\sigma_i\sigma_j+\sum_{i=1}^n\sigma_ih_i\br},
\end{align}
where $\bfa h=\{h_i\}_{i\in[n]}$ are realizations of arbitrary independent random variables with distribution $\mu_h$ and $S\subset[n]$ is the index set of a hidden clique of size $|S|=k$.  In \eqref{skdirect}, $\{ g_{ij}\}_{1\leq i<j\leq n}$ is a family of i.i.d. Gaussians and $\theta,\theta_1>0$ are deterministic parameters. In this work, we also study the \emph{universal} pSK model where $g_{ij}$ is i.i.d. with arbitrary symmetric distributions. We also call $\theta_1$ \emph{the inverse temperature}. There are two scaling factors $\frac{1}{k}$ and $\frac{1}{\sqrt n}$ in \eqref{skdirect}. These scaling factors guarantee the largest eigenvalue of the coupling matrix to be finite asymptotically, almost surely. In the rest of this paper, we also show that such scaling leads to interesting phase transition phenomena in high-dimensional statistical inference.

The pSK model arises from the standard SK model that first appears in \citep{sherrington1975solvable} to model spin glasses, a class of alloy exhibiting unique magnetization behavior. The characterization of the SK model is in the \emph{Gibbs Measure form} where we define $\mca H:\Sigma_n\to\bb R$ as the \emph{Hamiltonian} or energy function associated with some system and define the distribution of a spin configuration according to
\begin{align*}
    \bb P(\bfa\sigma):=\frac{\exp(-\mca H(\bfa\sigma))}{\sum_{\bfa\sigma}\exp(-\mca H(\bfa\sigma))}.
\end{align*}
Given an instance of i.i.d. Gaussian matrix $\bfa g=\{g_{ij}\}_{i,j\in[n]}$, the Hamiltonian of the SK model is given by
\begin{align}\label{skhamilt}
  \mca H^{SK}_{\theta}(\bfa \sigma) := -\frac{\theta}{\sqrt n} \sum_{1\leq i< j\leq n} g_{ij}\sigma_i\sigma_j -\sum_{i\leq n}h_i\sigma_i,\quad\theta>0
\end{align} 
with $\bfa h\in\bb R^n$ being the random vector characterizing the magnitude of the mean shift of each spin. In the following, we adopt physics terminology to denote it as the \emph{outer magnetic fields}. This random Gibbs measure is a mixture of models for $\bfa\sigma$ indexed by $\bfa h$ and $\bfa g$. The standard SK model has been rigorously studied by a large volume of works \citep{talagrand2006parisi,talagrand2010mean,guerra2001sum,guerra2003broken,aizerman1964theoretical,panchenko2012sherrington} and a phase transition phenomenon called the \emph{replica symmetry breaking} happens on a line named after Almeida and Thouless (or AT curve\footnote{In physics this is called the AT transition line, but is a curve. In this work we denote it by AT curve.} \citep{toninelli2002almeida}). The term \emph{replica} is a central concept in random Gibbs measures defined as the inner product of two sets of i.i.d. samples conditional on the randomness in the Hamiltonian. By \emph{replica symmetry}, we mean the replica inner product concentrates on a constant. When it is \emph{broken}, it does not concentrate (See \citep{panchenko2013parisi}). The AT curve is characterized by some function $\bb E[f(\theta,g,h)]=1$, and when $\bb E[f(\theta,g,h)]>1$, the replica symmetry breaking happens. In this work, we derive the corresponding AT curve for the planted SK model and study the region under this line. This region is commonly referred to as the \emph{replica symmetry} phase.

The mean-field Ferromagnetic correlation term is another important component in \eqref{skdirect}. One can also check that this ferromagnetic correlation is also analogous to the principle $k$-sparse vector in the Gaussian sparse PCA model. From the perspective of the graphical model, this corresponds to a $k$-clique. This component is akin to the Random Field Curie-Weiss (or RFCW) model, whose Hamiltonian is given by
\begin{align}\label{cwhamilt}
   \mca H^{RFCW}_{\theta}(\bfa\sigma):=- \frac{\theta}{n} \sum_{1\leq i<j\leq n}\sigma_i\sigma_j-\sum_{i\leq n}h_i\sigma_i,\quad\theta > 0,
\end{align}
with random entry-wise i.i.d. vector $\bfa h\in\bb R^n$. We also define the planted RFCW model by \eqref{cwhamilt} replacing the quadratic term by a $k$-sparse one given by $\frac{\theta_1}{k} \sum_{i,j\in S,i<j}\sigma_i\sigma_j$.

In this work, the pSK model in \eqref{skdirect} has a Hamiltonian combining \eqref{skhamilt} and \eqref{cwhamilt}. We define the Hamiltonian of the pSK model as
\begin{align}\label{hamilt}
\mca H_{\theta_1,\theta}^{pSK}(\bfa\sigma):=- \frac{\theta }{\sqrt n} \sum_{1\leq i<j\leq n}g_{ij}\sigma_i\sigma_j-\frac{\theta_1}{2k} \sum_{i,j\in S,i<j}\sigma_i\sigma_j-\sum_{i\in[n]}h_i\sigma_i,\quad\theta_1,\theta>0.
\end{align} 
In this work, we exclusively consider the case with $\bfa h$ being entry-wise i.i.d. symmetric random vector despite the fact that the non-symmetric case can be analogously derived using the same proof.
We show that the pSK model demonstrates two types of phase transitions rooted in the RFCW and SK models.  We denote $\theta_1$ to be the \emph{inverse temperature} of the pSK model. High-temperature regime corresponds to smaller $\theta_1$, and the low-temperature regime corresponds to larger $\theta_1$.
We denote $\ca G_0(\theta, n, \mu)$ to be the class of the SK models with random magnetic field $h_i$ i.i.d. sampled from a measure $\mu$ and $\ca G_1(\theta,\theta_1,k,n,\mu)$ to be the class of pSK models with the clique set $S$ unknown. 
Let $\{\bfa\sigma^{(1)}, \ldots, \bfa\sigma^{(m)}\}$ be $m$ independent observations. Our problem of interest is testing the following  hypotheses:
\begin{align}\label{mainproblem}
    \text{Null}: \bfa\sigma^{(i)}\sim \ca G_0(\theta,n,\mu)\qquad \text{vs}\qquad \text{Alternative}:\bfa\sigma^{(i)}\sim \ca G_1(\theta,\theta_1,k,n,\mu).
\end{align}
We aim to establish the statistical phase diagram for the powerful testing problem given by:
\begin{definition}[Asymptotic Power of Tests] \label{def1}
Define $\bb P_{0}$ as the probability measure under the null. Let $S$ be the index set of vertices of the hidden clique with size $|S|=k$ and define $\bb P_{S}$ as the probability measure under the alternative. Let $m$ be a sequence of values indexed by $k$ and write it as $m(k)$. Define $\bb P_{0,m}$ and $\bb P_{S,m}$ as their product of $m$ measures respectively. We say the sequence of tests $\psi\in\{0,1\}$ depending on $m$ i.i.d. samples are asymptotically powerful if
  \begin{align*}
       \lim_{k\to+\infty} \bigg[\bb P_{0,m(k)}\lef(\psi=1\rig)+\sup_{S:|S|=k}\bb P_{S,m(k)}\lef(\psi=0\rig)\bigg]=0
    .\end{align*}
Under the same notation and settings as in Definition~\ref{def1}, we say all sequences of the test are asymptotically powerless if all $\{0,1\}$-valued test statistics $\psi$ depending on $m$ i.i.d. samples $\{\bfa\sigma^{(i)}\}_{i\in[m]}$ satisfies
    \begin{align*}
        \lim_{k\to+\infty}\inf_{\psi\in\{0,1\}}\bigg[\bb P_{0,m(k)}(\psi=1)+\sup_{S:|S|=k}\bb P_{S,m(k)}(\psi=0)\bigg]=1
    .\end{align*}   
\end{definition}
In addition to testing the existence of a hidden clique, we also study the recovery problem over the set $S$.
\begin{definition}[Exact Recovery]
    Given the hidden clique index set $S$, we say an algorithm taking $m$ samples $\{\bfa\sigma^{(i)}\}_{i\in[m]}$ as input and output $\wh S$ recovers $S$ exactly if when $k\to\infty$ we have  $ \bb P(\wh S=S)=1-o(1)$.
\end{definition}

The recovery problem is more challenging than the testing problem, and the statistical gap between these two problems is the \emph{Test-Recovery Gap}. This work addresses both of them.

\subsection{Motivations} 
This work is motivated by the challenges imposed by the discrete data in high dimensional statistical analysis as a sequel over \citep{he2023hidden1}. In \citep{he2023hidden1}, we consider the problem of subset selection problem under the planted random field Curie Weiss model and establish the statistical minimax optimal rates over both the testing and the exact recovery tasks. The pRFCW model is a random Ising model suitable for characterizing the joint distribution of discrete random variables with latent random effects. Graphical models with random effects are prevalent in the statistical physics literature, with implications for real-world applications since the latent random variable can characterize the unknown feature of the data. Moreover, the problem of subset selection involves finding a subset with $k$ elements in $[n]$ with greater correlation. This problem also mathematically connected to the sparse PCA problem, which can be seen as its discrete variant. However, the pRFCW model considered in \citep{he2023hidden1} only accounts for the random effect on the mean values. To study the more complete picture where both the mean and the covariance are subject to random effects, we consider the pSK model in this work. The pSK model has a random coupling matrix with i.i.d. Wigner entries. The pSK model is also a spin glass model that can represent the unknown high dimensional covariance structure in the data, suitable for many applications in the real world. Moreover, the mathematical techniques in the analysis of the statistical minimax rates of the spin glass model are also of independent interest.

In summary, our motivations for this work are two folds: (1) To study the effect on the minimax optimal rates by the phase transitions that are intrinsic to Ising models and spin glasses; (2) To understand the effect of unknown covariance structure on the statistical procedures and in particular the minimax optimal rates of testing and inference.

\subsection{Contributions}

Our major contribution is establishing sharp minimax rates for the different parametric regimes of the pSK model. We recall that in a concurrent paper, we found that the planted RFCW model exhibits three types of phase transitions given by
\begin{enumerate}
    \item\textbf{High/Critical/Low-Temperature Regimes:} There exists a critical parameter $\theta_c$ partitioning the space of $\theta_1$ into three regimes: $\theta_1 < \theta_c$ (high-temperature regime), $\theta_1 = \theta_c$ (critical temperature regime), $\theta_1 > \theta_c$ (low-temperature regime), and the optimal sample complexity differ significantly across these regions. 
    \item\textbf{Second Transition at the Critical Temperature:} Depending on the tail heaviness of  $h$, the phase diagram at the critical temperature can vary significantly, where the optimal rate can take countable values. Moreover, the critical temperature represents an intermediate state between the high and low temperatures where we identify a co-existence of statistical phases in the high and low temperatures.
    \item \textbf{Global/Local Regimes:} We observe that a `mountain climbing' phase transition characterized by $k$ appears at all temperature regimes: when $k\gtrsim n^{\beta}$ for some $\beta\in(0,1)$, the optimal sample complexity is achieved by global tests taking all the spins as input; when $k=o(n^{\beta})$, the optimal complexity is achieved by a class of local scan tests. Moreover, $\beta$ differs across the temperature regimes and can take countable possible values at the critical temperature.
\end{enumerate}

Our first contribution is to show that all the above three regimes also exist in the pSK model, along with two new phases that appear at $k\asymp n$ that are unique to the pSK model. In \citep{he2023hidden1} we show that the phase transition of $\theta_1$ in the pRFCW model happens on $\theta_c$. In the pSK model, the transition happens on a curve characterized by a function of $\theta$ and $\theta_1$. The regimes of the pSK model hybridize those of the RFCW model and the SK model. For the pSK model, the replica symmetry-breaking regime is characterized by
\begin{align}\label{ATLINE}
        \bb E[\theta^2((1-c)\sech^4(\theta\sqrt q z+h)+c\sech^4(\theta\sqrt qz+\theta_1\mu+h))]>1,\quad 
\end{align}
with $c:=\lim_{n\to\infty}\frac{k}{n}$ and $q,\mu$ depend on $\theta,\theta_1,h$. In the standard SK model, \citep{toninelli2002almeida} proved that the phase transition happens in its AT curve, corresponding to the case $c=0$ in the pSK model. Our results stay in the regions with replica symmetry, or the complement of \eqref{ATLINE}.


Our second contribution is a series of universality results. We prove that all the results hold by replacing the i.i.d. Gaussian r.v.s. $\bfa g$ with arbitrary i.i.d. random variables (with bounded forth moment) matching the first three moments of standard Gaussian and the bounded forth moment.

Our third contribution is a family of novel concentration bounds for the \emph{average magnetization} (i.e. The average value of all spins) and replica concentration under multiple regimes. These results are crucial to obtain the sharp minimax rate of testing. We develop two novel tools to analyze the pSK model. The first tool is a modified smart path method that works for the \emph{ small clique region } that connects the upper and lower bounds with that of the pRFCW model in \citep{he2023hidden1}. The second tool is a constructive argument proving the concentration of replicas in the \emph{ large clique region }. This tool was initially invented by Talagrand \citep{talagrand2011mean2}. Then, we devise two smart paths for the cavity method (or leave-one-out method) to construct moment iterations of averaging statistics. These new results are of independent interest to analyzing other spin glass models.

In addition to the testing threshold presented above, we also present the exact recovery rate using the same algorithm as the tests. Our results imply that for small cliques, there exists no test-recovery gap.
\smallbreak
\emph{Organization.} The rest of this paper is organized as below. First, we summarize all the necessary notations. Section \ref{sect2} contains a discussion of related works; Section \ref{sect30} presents the major results in this work, including the small/large clique regimes and the recovery guarantees; Section \ref{sect6} discusses the universality results; Section \ref{sect7} presents discussions and conclusions. We delay most formal proofs and additional technical details to the supplementary material.
\smallbreak
\emph{Notations.} 
The following notations are used throughout this work. We use $:=$ as the notation for \emph{defining}. We denote $[n]:=\{1,\ldots,n\}$ and $[i:j]:=\{i,i+1,\ldots, j\}$ for $i<j$. For a vector denoted by $\bfa v=(v_1,\ldots v_n)\in\bb R^n$ we denote its $\ell_p$ norm by $\Vert \bfa v\Vert_p=\lef(\sum_{i=1}^p v_i^p\rig)^{1/p}$ for all $p\in[1,\infty)$. Denote $\Vert \bfa v\Vert_\infty=\sup_{i\in[n]}|v_i|$. For a matrix $A\in\bb R^{n\times m}$ with $m,n\in\bb N$ we denote $\Vert A\Vert_\infty = \sup_{i,j}|A_{ij}|$, $\Vert A\Vert_F=\lef(\sum_{i,j}A_{ij}^2\rig)^{1/2}$ and $\Vert A\Vert_p = \sup_{\bfa v:\Vert \bfa v\Vert_p=1}\Vert A\bfa v\Vert_p$ for all $p\in[1,\infty)$. For a vector $\bfa v\in\bb R^n$ and set $A\subset[n]$, we denote $\bfa v_{-A}$ as the vector constrained to $A^c$. We denote $\mbbm 1_{B}$ for some event $B$ as the indicator function of $B$. For some set $A\subset [n]$, we denote $\bfa v=\mbbm 1_A\in\bb R^{n}$ if $v_i=\mbbm 1_{i\in A}$. For a set $A\subset \Omega$ we denote $A^c=\Omega\setminus A$ where $\backslash$ is the notation for set minus. For another set $B\setminus \Omega$, we denote  $A \Delta B=(A\cup B)\setminus(A\cap B)$ as the symmetric difference between $A$ and $B$. Let $\bb P$ be a probability measure and $\bb P^{\otimes n}$ be the $n$-th order product measure of $\bb P$. Given two sequences $a_n$ and $b_n$, we denote $a_n\lesssim b_n$ or $a_n = O(b_n)$ if $\limsup_{n\to\infty}\lef|\frac{a_n}{b_n}\rig|<\infty$ and $a_n=o(b_n)$ if $\limsup_{n\to\infty}\lef|\frac{a_n}{b_n}\rig|=0$. Similarly, we denote $a_n\gtrsim b_n$ or $a_n =\Omega(b_n)$ if $b_n = O(a_n)$ and $a_n=\omega(b_n)$ if $b_n=o(a_n)$. We denote $a_n\asymp b_n$ or $a_n=\Theta(b_n)$ if $b_n\lesssim a_n$ and $a_n\gtrsim b_n$ both hold. For two sequence of measurable functions $f_n,g_n$ with $n\in\bb N$, we denote $f_n= O_p(g_n)$ if for all $\epsilon>0$ there exists $C>0$ such that $\limsup_n\bb P\lef (|f_n|>C|g_n|\rig )\leq \epsilon$  and $f_n = o_p(g_n)$ if for all $\delta>0$ $\limsup_{n}\bb P\lef(|f_n|>\delta |g_n|\rig) = 0$. We denote all $z$ in this work as standard Gaussian. Regarding convergence, we denote $\overset{d}{\to}$ as convergence in distribution. We denote $X\perp Y$ if two random variables are independent. If $\psi$ is a monotonic nondecreasing, convex function with $\psi(0)=0$, the Orlicz norm of an integrable random variable $X$ w.r.t. $\psi$ is given by $\Vert X\Vert_{\psi}=\lef\{u>0:\bb E\lef[\psi\lef(\frac{|X|}{u}\rig)\rig]\leq 1\rig\}$. In particular, for $\theta\in\bb R^+$ we use the notation of $\psi_{\theta}:=\exp(x^\theta)-1$. Finally, all constants are denoted by $C$ throughout this work unless specified otherwise.

For a random variable $X$, we denote $\bb E[X]$ as its mean and $\bb V[X]=\bb E[(X-\bb E[X])^2]$ as its variance.

\section{Related Works}\label{sect2}

Here, we mainly discuss the related works on spin glass models. In \citep{he2023hidden1} we extensively discuss the connections of our model's binary mathematical form with the class of Sparse PCA problems when $\theta=0$.

For the statistical analysis of random Gibbs measure, \citep{chatterjee2007estimation} studies the pseudo likelihood estimator for the inverse temperature in a few spin glass models.  Our problem differs from theirs as they focus on the estimation of $\theta$, and we focus on testing and recovering the planted structures.  Their analysis does not apply to our setting since the estimation method requires the correlation matrix among all spins to be known. As suggested in \citep{chatterjee2007estimation}, statistical inference of spin glasses is a challenging task that requires many new methodologies and frameworks.  In the past century, a rich line of work analyzed the behavior of the thermodynamical properties of the two models in \eqref{skhamilt} and \eqref{cwhamilt}.\citep{talagrand2010mean,ellis2006entropy}. However, these results are not applicable to the problems addressed in this paper. In particular, the Ferromagnetic correlation (or the planted sparse vector) in a model \eqref{hamilt} poses extra challenges to the classical analysis method for the SK model since the Ferromagnetic component is entangled with the background SK model. To handle this challenge, we need to resort to more sophisticated technical tools.  For the SK model in low temperature, \citep{talagrand2006parisi, panchenko2013parisi} focus on the limiting characterization of the partition function. \citep{chen2014mixed} conducts an asymptotic analysis of the SK model with Ferromagnetic correlation, which is a special case of the planted SK model in the regime of $k=n$. Their results are not directly usable here since no explicit convergence is derived for the averaging statistics. Also, their techniques are different from this work. In another line of work, \citep{hanen2009limit} analyzes the CLT of the standard SK model focusing on the extremely high-temperature regime and does not give a valid tail bound that is crucial in the analysis of the inference problem. Our results imply that the convergence rate for the planted SK model can be completely different from the standard SK model due to the Ferromagnetic term.

There is a rich literature applying theoretical results in the spin glass models to other fields (e.g., optimization and average computational complexity) \citep{bandeira2022franz,el2020fundamental,perry2018optimality,krzakala2016mutual}. For example, the analysis method for the cavity equation of the SK model \citep{bolthausen2014iterative} gives rise to provable algorithms in high dimensional inference, termed by the Approximate Message Passing \citep{,donoho2009message,bayati2011dynamics,javanmard2013state}.  Another example in \citep{zdeborova2016statistical} models the Stochastic Block Model with parameter $(p,q)$ using a spin glass formulation 
\begin{align*}
    \bb P(J_{ij}|\sigma_i\sigma_j) = p\delta_{1-\sigma_i\sigma_j}+q\delta_{1+\sigma_i\sigma_j},
\end{align*}
where $\delta$ is the Kronecker Delta function.  These results are 
not directly related to us, and the technical methods are largely different. Our analysis methods are based on the various interpolation and analysis methods proposed by  Talagrand \citep{talagrand2010mean,talagrand2011mean2}, Guerra \citep{guerra2001sum,guerra2003broken}, Toninelli \citep{toninelli2002almeida}, and Latala. Since most of their work focuses on the original SK model rather than planted structure and statistical inference, we develop new results and new methods. There exist other methods that can be used to analyze the concentration of the Ising model, including the concentration of measure arguments given by Chatterjee using Stein's exchangeable pairs \citep{chatterjee2005concentration,chatterjee2010applications}. However, his method is not directly applicable to the random Gibbs measures considered in this work.

\section{Main Results}\label{sect30}
In this section, we characterize a phase diagram of statistical tests by proving the minimax rates for testing in different regimes across all temperatures. Our testing results are summarized in Table \ref{table1}, and exact recovery results are summarized in Table \ref{table2}.
\begin{table}[htp!]\footnotesize
\captionsetup{labelfont=bf}
\caption{\textbf{The minimax sample complexity.} 
Presented is the sample complexity $m$ as a function of the click size $k$ and the number of spins.  Our results imply that the phase diagram of the planted SK model is similar to that of the planted RFCW model. Despite most regions being matched up to constants or logarithmic factors, our results leave two gaps: (1) the critical temperature region with $k\gtrsim n^{\frac{2\tau-1}{2\tau}}$. (2) the $\theta\gtrsim \sqrt{\frac{k}{n}}$ and $n^{1/2}\lesssim k=o(n/\log n)$ region at the high temperature regime. All of the results are stated under the maximum scaling of $\theta\lesssim 1$ unless specified otherwise. }
\makebox[0pt]{\begin{threeparttable}[t]
  \renewcommand*{\arraystretch}{2.4}
\begin{tabular}{ll|ll|ll}
\toprule
\multicolumn{2}{l|}{$\theta_1$ Regimes}                                                                                       & \multicolumn{2}{l|}{Local Tests}                                                                   & \multicolumn{2}{l}{Global Tests}                                                                                                                                                                        \\ \midrule
\multicolumn{2}{l|}{\textbf{High Temperature}}                                                                                & \multicolumn{1}{l|}{$k=o\lef(n^{\frac{1}{2}}\rig)$} & $n^{\frac{1}{2}}\lesssim k=o(n^{\frac{2}{3}})$ & \multicolumn{1}{l|}{$n^{\frac{2}{3}}\lesssim k= o\lef(\frac{n}{\log n}\rig)$}                                  & \multicolumn{1}{l}{$k\asymp n$}                                                                          \\ \hline
\multicolumn{1}{l|}{\multirow{2}{*}{$\theta_1\in(0,\frac{1}{2\bb E[\sech^2(h)]})$}}                            & UBs & \multicolumn{1}{l|}{$O(k\log n)$}                   & $O(k\log n)$                                      & \multicolumn{1}{l|}{$O\lef(\frac{n^2}{k^2}\rig)$}                                                & \multicolumn{1}{l}{$O(1)$}                                                                               \\ \cline{2-6} 
\multicolumn{1}{l|}{}                                                                                          & LBs & \multicolumn{1}{l|}{$\Omega(k\log n)$}                   & $\Omega\lef(\frac{k}{\log k}\log n\rig)$\tnote{*}                     & \multicolumn{1}{l|}{$\Omega\lef(\frac{n^2}{k^2}\rig)$\tnote{*}}                                               & \multicolumn{1}{l}{$\Omega(1)$}                                                                          \\ \hline
\multicolumn{1}{l|}{\multirow{2}{*}{$\theta_1\in\lef[\frac{1}{2\bb E[\sech^2(h)]},\frac{1}{\bb E[\sech^2(h)]}\rig)$}} & UBs & \multicolumn{1}{l|}{$O(k\log n)$}                   & $O(k\log n)$                                      & \multicolumn{1}{l|}{$O\lef(\frac{n^2}{k^2}\rig)$}                                                & \multicolumn{1}{l}{$O(1)$}                                                                               \\ \cline{2-6} 
\multicolumn{1}{l|}{}                                                                                          & LBs & \multicolumn{1}{l|}{$\Omega\lef(\frac{k}{\log k}\log n\rig)$}                   & $\Omega\lef(\frac{k}{\log k}\log n\rig)\tnote{*}$                     & \multicolumn{1}{l|}{$\Omega\lef(\frac{n^2}{k^2}\rig)$\tnote{*}}                                           & \multicolumn{1}{l}{$\Omega(1)$}                                                                          \\ \hline
\multicolumn{2}{l|}{\textbf{Critical Temperature}}                                                                            & \multicolumn{2}{l|}{$k=o\lef(n^{\frac{4\tau-2}{8\tau-5}}\rig)$}                                         & \multicolumn{1}{l|}{$n^{\frac{4\tau-2}{8\tau-5}}\lesssim k\lesssim n^{\frac{2\tau-1}{4\tau-3}}$} & $n^{\frac{2\tau-1}{4\tau-3}}\lesssim k\lesssim\frac{n^{\frac{2\tau-1}{2\tau}}}{\log^{\frac{2\tau-2}{2\tau-1}}n}$ \\ \hline
\multicolumn{2}{l|}{Upper Bounds}                                                                                             & \multicolumn{2}{l|}{$O(k^{\frac{1}{2\tau-1}}\log n)$}                                          & \multicolumn{1}{l|}{$O\lef(n^{\frac{2\tau-1}{4\tau-3}}\rig)$}                                    & $O(1)$                                                                                                    \\ \hline
\multicolumn{2}{l|}{Lower Bounds}                                                                                             & \multicolumn{2}{l|}{$\Omega\lef(\lef(\frac{k}{\log k}\rig)^{\frac{1}{2\tau-1}} \log n\rig)$}   & \multicolumn{1}{l|}{$\Omega\lef(n^{\frac{2\tau-1}{4\tau-3}}\rig)$}                               & $\Omega(1)$                                                                                               \\ \hline
\multicolumn{2}{l|}{\textbf{Low Temperature}}                                                                                 & \multicolumn{2}{l|}{$k=o\lef(n^{\frac{1}{2}}\rig)$}                                                     & \multicolumn{1}{l|}{$k=n^{\frac{1}{2}}$}                                                         & $k=\omega\lef(n^{\frac{1}{2}}\rig)$                                                                       \\ \hline
\multicolumn{2}{l|}{Upper Bounds}                                                                                             & \multicolumn{2}{l|}{$O\lef(\log n\rig)$}                                                       & \multicolumn{1}{l|}{$O(1)$}                                                                      & 1                                                                                                         \\ \hline
\multicolumn{2}{l|}{Lower Bounds}                                                                                             & \multicolumn{2}{l|}{$\Omega\lef(\log n\rig)$}                                                  & \multicolumn{1}{l|}{$\Omega(1)$}                                                                 & 1                                                                                                         \\ \bottomrule
\end{tabular}
 \begin{tablenotes}
 \footnotesize
            \item[*] Results under an additional small noise assumption of $\theta\lesssim \sqrt{\frac{k}{n}}$.
        \end{tablenotes}
\end{threeparttable}
\label{table1}}
\end{table}

\begin{table}[htp!]
\caption{\textbf{The minimax sample complexity of exact recovery } Presented is the sample complexity $m$ as a function of the click size $k$ and the number of spins. Our results leave open a few regions when the noise $\theta$ is large. All of the results are stated under the maximum scaling of $\theta\lesssim 1$ unless specified otherwise.}
  \renewcommand*{\arraystretch}{2.0}
\centering
\makebox[0pt]{\begin{threeparttable}
\begin{tabular}{l|l|l|l}
\toprule
$\theta_1$ Regimes                    & $k$ Regimes                                                                                                                                         & Upper Bounds                     & Lower Bounds                          \\ \midrule
\multirow{3}{*}{High Temperature}     & $k=O(n^{\frac{1}{2}})$                                                                                                                              & $O(k\log n)$                     & $\Omega(k\log n)$                     \\ \cline{2-4} 
                                      & $k=o(\frac{n}{\log n})$                                                                                                                             & $O(k\log n)$                     & $\Omega(k\log n)\tnote{*}$           \\ \cline{2-4} 
                                      & $\frac{n}{\log n}\lesssim k\leq n$                                                                                                                  & $O(k\log n)$                     & $\Omega(k\log n)\tnote{\P}$                               \\ \hline
\multirow{2}{*}{Low Temperature}      & $k=o(\frac{n}{\log n})$                                                                                                                             & $O(\log n)$                      & $\Omega(\log n)$                      \\ \cline{2-4} 
                                      & $\frac{n}{\log n}\lesssim k\leq n$                                                                                                                  & $O(\log n)$                      & $\Omega(\log n)\tnote{\S}$                               \\ \hline
\multirow{3}{*}{Critical Temperature} & $k=o\lef(\frac{n^{\frac{2\tau-1}{4\tau-3}}}{\log^{\frac{2}{2\tau-1}}(n-k)}\rig)$                                                                        & $O(k^{\frac{1}{2\tau-1}}\log n)$ & $\Omega(k^{\frac{1}{2\tau-1}}\log n)$ \\ \cline{2-4} 
                                      & $\frac{n^{\frac{2\tau-1}{4\tau-3}}}{\log^{\frac{2}{2\tau-1}}(n-k)}\lesssim k=o\lef(\frac{n^{\frac{2\tau-1}{2\tau}}}{\log^{\frac{2\tau-2}{2\tau-1}}(n-k)} \rig)$ & $O(k^{\frac{1}{2\tau-1}}\log n)$ & $\Omega(k^{\frac{1}{2\tau-1}}\log n)$\tnote{\ddag}                           \\
                                      \cline{2-4} 
                                      & $\frac{n^{\frac{2\tau-1}{2\tau}}}{\log^{\frac{2\tau-2}{2\tau-1}}(n-k)}\lesssim k\leq n$ & $O(k^{\frac{1}{2\tau-1}}\log n)$\tnote{\dag} & $\Omega(k^{\frac{1}{2\tau-1}}\log n)$\tnote{\ddag}\\
                                      \bottomrule
\end{tabular}
 \begin{tablenotes}
 \scriptsize
            \item[*] Results under an additional assumption of $\theta\lesssim n^{\frac{1}{2}}k^{-1}$.
            \item[\P] Results under an additional assumption of $\theta \lesssim n^{\frac{1}{2}}k^{-1}\log^{-\frac{1}{2}} n$.
            \item[\S] Results under an additional assumption of $\theta\lesssim \log^{-\frac{1}{2}} n$.
            \item[\ddag] Results under an additional assumption of $\theta\lesssim n^{\frac{1}{2}}{k^{-\frac{2\tau-1}{8\tau-6}}}{\log^{-\frac{1}{2\tau-1}}n}$.
            \item[\dag] Results under an additional assumption of $\theta\lesssim n^{\frac{1}{2}}k^{-\frac{\tau}{2\tau-1}}\log^{-\frac{\tau-1}{2\tau-1}}n$.
        \end{tablenotes}
\end{threeparttable}}
\label{table2}
\end{table}
In the pRFCW model, we divide the different regimes according to $\theta_1$, or the inverse temperatures of the pSK model. As opposed to the pRFCW model in \citep{he2023hidden1}, where the different temperature regimes are denoted by the order between $\theta_1$ and $\frac{1}{\bb E[\sech^2(h)]}$, the pSK model has temperature regimes separated by a curve on the $(\theta_1,\theta)$ plane. This curve is characterized by $\theta_1\bb E[\sech^2(\theta\sqrt{q}z+h)]=1$ with $z\sim N(0,1)$ and we denote $\theta_1<\frac{1}{\bb E[\sech^2(\theta\sqrt{q}z+h)]}$ as the high temperature regime, $\theta_1=\frac{1}{\bb E[\sech^2(\theta\sqrt{q}z+h)]}$ as the critical temperature regime and $\theta_1>\frac{1}{\bb E[\sech^2(\theta\sqrt qz+h)]}$ as the low temperature regime. Another important feature is the `replica symmetry breaking' phenomenon of the pSK model when we reach the AT curve. It is known that in the canonical pSK model \citep{talagrand2006parisi}, the replica statistics no longer concentrate on a single point as the dimension increases to infinity. Hence, we focus on the `replica symmetric' region as the replica statistics play a crucial role in the analysis of this work. 

Our separation of the high/critical and low-temperature regimes is also subject to the following constraint for $c:={k}/{n}$, which remain fixed as $n\to\infty$:
\begin{align}\label{replicasymmetrycond}
     \bb E[\theta^2((1-c)\sech^4(\theta\sqrt q z+h)+c\sech^4(\theta\sqrt qz+\theta_1\mu+h))]<1,
\end{align}
with an additional technical assumption is given by definition \ref{withintheATline}. These conditions are often referred to as the `replica symmetric region', a region where the concentration of the replica holds. We also show in this work that the concentration of replica leads to the concentration of spins, a fundamental result for the statistical analysis.
\smallbreak\emph{The Different Temperature Regimes.} We characterize the different temperature regimes according to the values of the parameter triples $(\beta,\beta_1,\mu_h)$. In particular, under the `replica symmetric condition' given by \eqref{replicasymmetrycond}, our division of the three temperature regimes is given according to the magnitude $\psi(\theta_1,\theta,\mu_h):=\theta_1\bb E[\sech^2(\theta\sqrt q z+h)]$ compared with $1$. Specifically, the temperature regime is given by $\psi(\theta_1,\theta,\mu_h)<1$, the critical temperature regime is given by $\psi(\theta_1,\theta,\mu_h)=1$ and the low temperature regime $\psi(\theta_1,\theta,\mu_h)>1$. At the critical temperature regime, we introduce a new `flatness parameter' $\tau\in\bb N$. This flatness parameter 
corresponds to the global minimum of the function $H(x)$ defined by
\begin{align}\label{taylorcond}
    H(x)=\frac{1}{2}x^2-\bb E[\log\cosh(\sqrt{\theta_1}x+\theta\sqrt qz+h)],
\end{align}
and governs the rate of convergence for the average of spins, which is formally defined by
\begin{definition}[Flatness of Local Optimum]\label{flatnessofopt}
    We call the local minimum/maximum $x^*$ of \eqref{taylorcond} is $\tau$-flat for $\tau\in\bb N\setminus\{1\}$ if $ H(x)=\frac{H^{(2\tau)}(x^*)}{(2\tau)!}(x-x^*)^{2\tau-1}+O((x-x^*)^{2\tau})$ with $H^{(2\tau)}(x_1^*)>0$ or $H^{(2\tau)}(x_1^*)<0$. 
\end{definition}
We show that the critical temperature regime only appears when $H^{(2\tau)}(0)<0$ and the scenario where $H^{(2\tau)}(0)>0$ shares the same concentration results with the low-temperature regime.
\smallbreak\emph{The Small and Large Clique Regimes.} The condition given by \eqref{replicasymmetrycond} naturally implies that the scenarios of $k/n=o(1)$ and $k/n\asymp 1$ should be treated differently, which corresponds to \emph{the large and small clique regimes}. In particular, the small clique regime is given by $c=0$, and the condition \ref{replicasymmetrycond} becomes
\begin{align}\label{atcondition}
    \theta^2\bb E[\sech^4(\theta z\sqrt q+h)]<1,\qquad \text{ where } q=\bb E[\tanh^2(\theta z\sqrt q+h)],
\end{align}
To tackle the technical difficulties arising from the two regimes we develop different proof ideas and methods.



\smallbreak\emph{Technical Challenges.} The random cross terms in the pSK model pose a few significant challenges to the statistical analysis compared with \citep{he2023hidden1}. For example, one might attempt to derive the concentration for the pSK model using the conditioning technique and argue that there exist some `high probability events' where the model demonstrates an all-positive coupling matrix to make it a Ferromagnetic model. Or one might attempt to ask if he can construct an exchangeable pair to derive the exponential inequality for the pSK model following Chatterjee's method \citep{chatterjee2005concentration}. Unfortunately, all of these methods fail to succeed in this model. In the former case, this high probability event does not exist since the Gaussian coupling is of order $\frac{1}{\sqrt n}$ in the fluctuation, whereas the positive term associated with the clique only has magnitude $\frac{1}{k}$. For the latter, one encounters a significant challenge when trying to compute the conditional expectation on all the rest of the spins $\bb E[\sigma_{k}|\bfa\sigma_{-k}]$, a common step in Chatterjee's method. One can show that conditional on the first $n-1$ spins of the pSK model, the magnetic field and the Gaussian coupling coefficients associated with the $n$-th spin undergo a distortion of the posterior distribution, which makes the conditional expectation not computable. The distribution is much more complicated to deal with, which is termed by \citep{chatterjee2010spin} as the \emph{local field} of the SK model.

Then, a fundamental question to all existing statistical analysis methods have to be re-examined: Does the law of large numbers continue to hold for the pSK model? What are the convergence rates of its statistics? As is shown in this work, many natural facts and quantities for the classical multivariate Gaussian problems and even our previous work on the pRFCW model \citep{he2023hidden1} are particularly challenging to get for the pSK model. In \citep{he2023hidden1}, we prove the convergence rate of the average magnetization for the pRFCW model, but the same method becomes invalid for the pSK model since the H-S transformation no longer decouples random correlation terms. Moreover, the dependence structure given by the random coupling terms becomes much more complicated than the pRFCW model where the out-of-clique spins just follow the i.i.d. Bernoulli distribution. These special characteristics are also not seen in traditional statistical analysis since we have: (1) The dependence between spins in $S$ (Dependent Signal). (2) The dependence between spins in $S^c$ (Dependent Noise). (3) The dependence between spins in $S$ and $S^c$ (Dependence between Signal and Noise ). We briefly discuss this idea in section \ref{sect300}. 

To overcome these technical barriers, we develop methods to decouple the different types of correlations. At the core of our analysis is the idea of the \emph{smart path method} where instead of working on the Hamiltonians of \eqref{hamilt} directly, we approximate it with that of the pRFCW model. Given an upper bound on the interpolation error, we can derive sharp results for both the upper tail bounds for magnetization and the upper bound for information divergences.
The validity behind this method is guaranteed by Stein's lemma, which argues that the difference of the above interpolation can be explicitly computed and estimated by Gaussian integration by parts. Although the idea is simple, it is often impossible to solve all the problems with a single instance of the `smart path' and we derive different interpolation methods suitable for different parameter regimes.

\smallbreak\emph{Organization.} The following sections are organized as follows: We provide the main theorems in section \ref{mainthms}; We give an overview of technical contributions in section \ref{sect300}; sections \ref{sect3}, \ref{sect4}, and \ref{sect5} provide proof outlines of the main theorems in this work.

\subsection{Main Theorems}\label{mainthms}
This section presents the main theorems of this work, reflecting the results in the two tables \ref{table1} and \ref{table2}. We provide theorems for the optimal sample complexities of testing and recovery problems under different temperature regimes. The following theorem summarizes our results given the sparsity $k\log k\lesssim n$, which corresponds to the \emph{small clique regime}.
\begin{theorem}\label{thm10}
  Assume that $k\log k=o(n)$ and the condition \eqref{atcondition} holds. Depending on the different temperature regimes, the optimal sample complexities $m$ for asymptotic powerful testing are given by :
  \begin{enumerate}
      \item At the extreme high temperature regime with $0<\theta_1<\frac{1}{2\bb E[\sech^2(\theta\sqrt q z+h)]}$,
      $$\frac{n}{k\theta^2}\wedge\frac{n^2}{k^2}\wedge k\log n\lesssim m\lesssim\frac{n^2}{k^2} \wedge k\log n.$$
      \item At the high temperature regime with $\frac{1}{2\bb E[\sech^2(\theta\sqrt q z+h)]}<\theta_1<\frac{1}{\bb E[\sech^2(\theta\sqrt qz+h)]}$,
      $$\frac{n}{k\theta^2}\wedge\frac{n^2}{k^2}\wedge k\lesssim m\lesssim\frac{n^2}{k^2} \wedge k\log n.$$
      \item At the critical temperature regime with $\theta_1=\frac{1}{\bb E[\sech^2(\theta\sqrt qz+h)]}$ and flatness parameter $\tau$ defined by \ref{flatnessofopt}, 
    $$ n^2k^{-\frac{2(4\tau-3)}{2\tau-1}}\wedge\bl\frac{k}{\log k}\br^{\frac{1}{2\tau-1}}\log n\lesssim m\lesssim k^{\frac{1}{2\tau-1}}\log n\wedge n^2k^{-\frac{2(4\tau-3)}{2\tau-1}},\text{ when }k\lesssim n^{\frac{2\tau-1}{4\tau-3}}$$ and 
    $m\asymp 1$ when $k\lesssim n^{\frac{2\tau-1}{2\tau}}\log ^{-\frac{2\tau-2}{2\tau-1}}n$.
      \item At the low temperature regime with $\theta_1>\frac{1}{\bb E[\sech^2(\theta\sqrt qz+h)]}$, when $k=o(\sqrt n)$ we have $m\asymp \log n$, when $k\asymp \sqrt n$ we have $m\asymp 1$, and when $k=\omega(\sqrt n)$ we have $m=1$.
  \end{enumerate}
\end{theorem}
An outline of the proof of the above theorem is provided in section \ref{sect3}
We shall note that lower bounds of $k=\omega(n^{\frac{2\tau-1}{2\tau}}\log^{-\frac{2\tau-2}{2\tau-1}}n)$ at the critical temperature regime remains open, which is due to a technical barrier in our proof method. Due to a similar reason, in comparison with the corresponding results given by \citep{he2023hidden1} (here we also recall that it corresponds to the special case of $\theta=0$), we observe that $\theta$ comes into the lower bounds of the high-temperature regime. This finally result in a gap from the upper bounds when $\theta$ is large whereas in \citep{he2023hidden1} our results match with the upper bounds up to a logarithmic factor.

The next result presents the main theorem for the large clique regime. Only results for the high and low-temperature regimes are provided. Before introducing the result, a few technical definitions are given as follows.
\begin{definition}[Within the AT line]\label{withintheATline}\label{hightemp}
    We define the \emph{within the AT line} region as $(\theta_1,\theta)$ satisfying
    \small{
    \begin{align*}
        \bb E\theta^2((1-c)\sech^4(\theta\sqrt q z+h)+c\sech^4(\theta\sqrt qz+\theta_1\mu+h))<1,
\end{align*}}
and
\small{
\begin{align}\label{condderi}
      \forall q^\prime >q\quad&\Rightarrow\quad\frac{\pta\Phi}{\pta m}(m,q^\prime)\bigg |_{m=1} <0,
\end{align}}
\normalsize
where $\Phi(m,q^\prime)$ is defined by (let $z,z^\prime$ being i.i.d. standard Gaussian and $\bb E^\prime$ take expectation w.r.t. $z^\prime$.)
\tny{
\begin{align}
    Y_1:&=\theta\sqrt{q}z+\theta \sqrt{q^\prime-q}z^\prime +\theta_1\mu +h,\quad\text{ and }\quad Y_2:=\theta\sqrt{q}z+\theta\sqrt{q^\prime-q}z^\prime+h,\nnb\\
\label{phidef}
    \Phi(m,q^\prime):&=\log 2+\frac{\theta^2}{4}(1-q^\prime)^2-\frac{\theta^2}{4}m(q^{\prime 2}-q^2)+\sup_{\mu\in[-1,1]}\bl\frac{k}{mn}\bb E\log \bb E^\prime\cosh^mY_1-\frac{\theta_1\mu^2}{2}\br\nnb\\
    &+\frac{n-k}{mn}\bb E\log\bb E^\prime\cosh^mY_2.
\end{align}}
\end{definition}
\normalsize
We remark that such technical definition also appears in the analysis of the replica symmetric region of SK model in \citep{talagrand2011mean2}. This definition guarantees the convergence results hold for both the average magnetization and the replicas. Given the above preparation, we are ready to state the formal results in the large clique regime.
\begin{theorem}\label{thm2}
Given $k\asymp n$, there exist asymptotically powerful testing procedures given the following:
\begin{enumerate}
    \item When $0<\theta_1<\frac{1}{\bb E[\sech^2(\theta\sqrt qz+h)]}$, $m=\omega(1)$, and the replica symmetry condition in definition \ref{withintheATline} holds with $\mu=0$ ;
    \item When $\theta_1>\frac{1}{\bb E[\sech^2(\theta\sqrt qz+h)]}$, $m=1$,  the replica symmetry condition in definition \ref{withintheATline} and the invertibility condition in the lemma \ref{cltatline} holds.
\end{enumerate}
\end{theorem}
The proof of the above theorem involves multiple steps and is outlined in section \ref{sect4}. Our results imply that the sample complexity of the small clique regime extends to the large clique regime despite the fact that a completely different analysis method is required to obtain such results. It is noted that our methods that work for both the high and low-temperature regimes do not naturally extend to the critical temperature regime, which remains open. 

The next theorem gives the sample complexity required for the exact recovery problem.
\begin{theorem}\label{exact recoveryguarantees}
    Assume that the replica symmetry condition \ref{withintheATline} holds and that the clique is positioned with index set $S$. Assume that $\wh S$ is the estimated clique set returned by algorithm \ref{alg:screen}. For the upper bounds, we have:
     \begin{enumerate}
         \item At the high temperature regime, when $m\geq Ck\log (n)$, $\bb P(\wh S=S)=1-o(1)$;
         \item At the low temperature regime, when $m\geq C\log (n)$, $\bb P(\wh S=S)=1-o(1)$;
         \item At the critical temperature regime, when $k=o( {n^{\frac{2\tau-1}{4\tau-4}}}{\log^{-\frac{2}{2\tau-1}}n})$ and $m\geq Ck^{\frac{1}{2\tau-1}}\log (n)$, $\bb P(\wh S=S)=1-o(1)$.
     \end{enumerate}
     For the lower bounds, we have:
     \begin{enumerate}
         \item At the high temperature regime, when $k=o(\sqrt n)$, $\inf_{\wh S}\sup_{\bb P\in \ca P}\bb P(S\neq\wh S)\geq 1-O\lef({m}{k^{-1}\log^{-1} (n)}\rig) $;
         \item At the high temperature regime, when $\sqrt n\lesssim k=o\lef({n}{\log^{-1} n}\rig)$, $\inf_{\wh S}\sup_{\bb P\in \ca P}\bb P(S\neq\wh S)\geq 1-O\lef({m k}{\theta^{2}n^{-1}\log^{-1} (n)}\rig) $;
         \item At the low temperature regime, when $k=o\lef({n}{\log^{-1} n}\rig)$, $\inf_{\wh S}\sup_{\bb P\in \ca P}\bb P(S\neq\wh S)\geq 1-O\lef({m }{\log^{-1} (n)}\rig) $;
         \item At the critical temperature regime, when $k=o({n^{\frac{2\tau-1}{4\tau-4}}}{\log^{-\frac{2}{2\tau-1}}n})$, $\inf_{\wh S}\sup_{\bb P\in \ca P}\bb P(S\neq\wh S)\geq 1-O({m }{k^{-\frac{1}{2\tau-1}}\log^{-1} (n)}) $.
     \end{enumerate}
\end{theorem}
The proof of the above theorem involves first showing the almost exact recovery guarantee, which is followed by a refinement procedure. In section \ref{sect5}, we outline the proof of it. It is also unknown if the dependence of $\theta^2$ in the high-temperature regime lower bounds can be removed using better technical arguments than the ones presented here.

\subsection{Overview of Techniques}\label{sect300}
This section presents an overview of the technical contributions made in this work. To prove the upper bounds, we give estimates for the moment generating functions of the average magnetizations across all the temperature regimes. To prove the lower bounds, we give sharp estimates of the information divergences between mixtures of pSK models and the null distribution, which are required by Fano's lemma and Le Cam's method. Despite the statistical framework being standard, both tasks are very challenging for the spin glass models and have little existing literature to take for reference.

Recall that in \citep{he2023hidden1}, the authors propose a sequence of methods to deal with the same problem on the pRFCW model, corresponding to $\theta=0$ in the pSK model. To deal with the upper bounds, the authors use the Hubbard–Stratonovich transformation to treat the quadratic term, the Laplace approximation method, and a carefully designed transfer principle to overcome the obstacle associated with the value of the m.g.f. for the average magnetization. Combining these methods, the authors derive the upper tail bounds. To deal with the lower bounds, the authors propose a `fake measure' method that is used together with the Hubbard–Stratonovich transformation (H-S transformation)
 to overcome the difficulties brought by the unboundedness of the chi-square divergences between a mixture of alternative and null hypotheses. However, all of the above methods given in \citep{he2023hidden1} become invalid here due to the existence of random coupling terms. For example, the H-S transformation that was originally proposed is not able to deal with the extra random quadratic terms in the pSK model. In particular, the H-S transformation is based on the Gaussian identity, which proceeds as
\begin{align*}
    \exp\bl\sum_{i,j\in S}\frac{\theta_1}{2k}\sigma_i\sigma_j\br=\int_{\bb R}\exp\bl\sum_{i\in S}\sigma_i\sqrt{\frac{\theta_1}{k}}x-\frac{x^2}{2}\br dx.
\end{align*} After the H-S transformation, the correlations between spins become decoupled and look like conditional i.i.d. ones. However, the term $\sum_{i<j\leq n}g_{ij}\sigma_i\sigma_j$ involve randomness, which makes the cross terms no longer decouple as above after the H-S transformation. Therefore, a sequence of new strategies is required in this new model to address the randomness in the coupling matrix.

To overcome the above difficulties, we develop the idea of \emph{local interpolation}. Instead of proving the necessary results like the upper bound on mgf and the information divergences in the pSK model directly, we use the ones already obtained in \citep{he2023hidden1} as intermediate results and show that the corresponding quantities for the pSK is close to the pRFCW model. So, the major technicalities lie in (1) How do we construct a valid interpolating path? (2) How do we analyze the interpolation errors between the two models? To answer the first question, we build a local Gaussian interpolation path for the function on a small portion of the spins in the pSK model. To resolve the second question, we study the local replica concentration inequalities using Latala's argument, which carefully controls the difference between the pRFCW model and the pSK model. Using these two techniques together, we achieve sharp upper bounds and lower bounds for the small clique region with $k=o(n)$. However, this method only works when our desired function is on $k$ spins where $k\ll n$. If we try to address the problem when $k\asymp n$, a new method is required, and we utilize the \emph{local cavity method} and a certificate of concentration. The local cavity method is a leave-one-out method for the local function on $k$ variables when $k\asymp n$ and both the high and low-temperature regimes. The local cavity method provides the tail bounds, finally leading to the upper bound of statistical rates. The certificate of concentration utilize the 1RSB (or 1-step replica symmetry breaking) analysis provided by Talagrand \citep{talagrand2011mean2} for the pSK model. Finally, our universality results is based on a local non-Gaussian interpolation by parts lemma between the universal variant of pSK model (where the random standard Gaussian $g$ in the pSK Hamiltonian is replaced by any forth moment bounded symmetric  distribution with unit variance ) and the pSK model, which works analogously as the Gaussian interpolation between the pSK and the pRFCW model.

\smallbreak\emph{The Local Smart Path Method.}\label{lcsmartpath} The local smart path relies on the property of the Gaussian random variables. The idea underlying this method is to `move' the random correlation terms to the random magnetic field, which makes it a pRFCW model. Specifically, we can build an interpolation path between the measures of the two models so the difference between the general functions' expectations is controlled through the derivative along the path. Then, when the derivative is of order $o(1)$, we can use the results at the easy end (the pRFCW model) to get a sharp estimate on the hard end (the pSK model). Here we provide a simple example to illustrate how it works for the pSK model (where $g$s are standard Gaussian in \eqref{hamilt}), by considering the mean of a general function on spins. Some examples of such functions include the moment generating functions and moments of the average value of the spins. In later sections, we discuss multiple variants of it that work for many different parametric regimes of the pSK model.

Consider two Gaussian processes indexed by $\Sigma_n:=\{-1,1\}^n$, define $u_{\bfa\sigma}:=\frac{1}{\sqrt n}\sum_{1\leq i<j\leq n}\theta g_{ij}\sigma_i\sigma_j$ and $v_{\bfa\sigma}:=\theta\sqrt q\sum_{i=1}^nz_i\sigma_i$ for some $q\in\bb R$ with $\{g_{ij}\}_{i,j\in[n]}$ and $\{z_i\}_{i\in[n]}$ are two independent standard i.i.d. Gaussian families (they are also independent of each other), then we can write the Hamiltonian of the pSK model with $S=[k]$ as
\begin{align*}
    -\mca H_1(\bfa\sigma):=u_{\bfa\sigma}+\sum_{i,j\in[k]}\frac{\theta_1}{2k}\sigma_i\sigma_j+\sum_{i=1}^n\sigma_ih_i.
\end{align*}
Since it is a random quadratic term, the most difficult part to analyze is $u_{\bfa\sigma}$. To resolve this, we approximate the above Hamiltonian with the following one, which replaces $u_{\bfa\sigma}$ with $v_{\bfa\sigma}$. The reader can immediately see that
\begin{align*}
    -\mca H_0(\bfa\sigma):=\sum_{i,j\in[k]}\frac{\theta_1}{2k}\sigma_i\sigma_j+\sum_{i=1}^n\sigma_i(h_i+\theta\sqrt qz_i)=v_{\bfa\sigma}+\sum_{i,j\in[k]}\frac{\theta_1}{2k}\sigma_i\sigma_j+\sum_{i=1}^n\sigma_ih_i.
\end{align*}
The reader can immediately see that $H_0$ corresponds to the Hamiltonian of an RFCW model, which we already know how to work with in \citep{he2023hidden1}. The next step is to construct an interpolating Gaussian process $u_t(\bfa\sigma):=\sqrt{t}u_{\bfa\sigma}+\sqrt{1-t}v_{\bfa\sigma}$ with $t\in[0,1]$ and its corresponding Hamiltonian,
\begin{align*}
    -\mca H_t(\bfa\sigma):=u_{t}(\bfa\sigma)+\sum_{i,j\in[k]}\frac{\theta_1}{2k}\sigma_i\sigma_j+\sum_{i=1}^n\sigma_ih_i.
\end{align*}
Then we define $\omega(\bfa\sigma)=\exp\lef(-\mca H_0(\bfa\sigma)-v_{\bfa\sigma}\rig)$. For a function $f:\Sigma_n\to\bb R$, we define $\la f\ra_t:=\frac{\sum_{\bfa\sigma}f(\bfa\sigma)\omega(\bfa\sigma)\exp(u_{\bfa\sigma}(t))}{\sum_{\bfa\sigma}\omega(\bfa\sigma)\exp(u_{\bfa\sigma}(t))}$. Then we can use Stein's lemma to get
\small{
\begin{align}\label{interpolation derivative}
    \frac{d}{dt}\bb E[\la f\ra_t]=\frac{1}{2}\sum_{i,j\in\Sigma_n}(\bb E[u_iu_j]-\bb E[v_iv_j])\bb E\bigg[\frac{\pta^2\la f\ra_t}{\pta x_i\pta x_j}\bigg]=\frac{n\theta^2}{2}\bb E[\la (R_{1,2}-q)^2(f-\la f\ra_t)\ra_t],
\end{align}
}
\normalsize
where  $R_{1,2}:=\frac{1}{n}\sum_{i=1}^n\sigma^1_i\sigma^2_i$ is the replica overlap. If we obtain an $o(1)$ upper bound on the derivative \eqref{interpolation derivative}, then we can use the expectation w.r.t. the pRFCW model to approximate that of the pSK model. The replicas of the pSK model denoted by $\bfa\sigma^1,\bfa\sigma^{2}$ are independent samples from the conditional Gibbs measure $\bb P(\bfa\sigma|\bfa g,\bfa h)$.  We note this method does not work since the replica $R_{1,2}$ has a typical convergence rate of $n^{-1/2}$ for standard spin glasses \citep{talagrand2010mean}. Therefore, this method does not yield proper interpolation if we can upper bound $\la f\ra_t$ for all $t\in[0,1]$. To overcome this obstacle, we provide a localized procedure where only the randomness of couplings within $S$ and between $S$ (recall that $S$ is the predefined clique index) and $S^c$ are approximated by the random magnetic fields. In other words, we decouple the random correlation terms as follows
\begin{align*}
    \sum_{1\leq i<j\leq n}\frac{\theta g_{ij}}{\sqrt n}\sigma_i\sigma_j=\ub{\sum_{i<j,i,j\in S}\frac{\theta g_{ij}}{\sqrt n}\sigma_i\sigma_j+\sum_{i\in S,j\in S^c}\frac{\theta g_{ij}}{\sqrt n}\sigma_i\sigma_j}_{T_1}+\ub{\sum_{i,j\in S^c,i<j}\frac{\theta g_{ij}}{\sqrt n}\sigma_i\sigma_j}_{T_2}
\end{align*}

Instead of annealing all the cross-correlation terms at once, we only approximate the crucial correlation term $T_1$ above and keep the term $T_2$ unchanged. As an extra treatment, we incorporate the following term as the $t=0$ end in replacement of the correlation term $T_1$
\begin{align*}
    \sum_{i\in S}\sigma_i\sqrt {q_1}z_i+\sum_{i\in S^c}\sigma_i\sqrt{q_2}z_i.
\end{align*}
A huge benefit of such treatment is that the spins in and out of $S$ are independent. (Of course, to make it work we still need an extra concentration of replica result from the standard spin glass model and a proper choice of $q_1$ and $q_2$) The above interpolation method yields that for any local function $f:\Sigma_S\to\bb R$ on the spins in $S$, we have
\begin{align}\label{interpolation derivative2}
    \frac{d}{dt}\bb E[\la f\ra_t]=\frac{k^2\theta^2}{2n}\bb E[\la (R^S_{1,2}-q)^2(f-\la f\ra_t)\ra_t],
\end{align}
where $R_{1,2}^S=\frac{1}{k}\sum_{i\in S}\sigma_i^1\sigma_i^2$.
Therefore, if we manage to show that $R_{1,2}^S-q$ converges with an order of roughly $k^{-1/2}$, one can use the above equation to construct a proper differential equation and give approximations for quantities, including moment generating functions. This further gives the upper tail bound through reusing our results for the pRFCW model in \citep{he2023hidden1}. The information divergences require a more delicate approximation theme but share the spirit of the above localization method. To derive the concentration of $R_{1,2}^S-q$ in the above equation, we use the next method.

\smallbreak\emph{ Latała's Argument.} To get an estimate of the quantity \eqref{interpolation derivative}, we need to study the moment bound of the term $R_{1,2}-q$ for all $t\in[0,1]$. For the standard SK model, Rafał Latała developed a method to derive the upper bound for this quantity in his unpublished manuscript \citep{talagrand2010mean}. He constructs the following quantity
\begin{align}\label{latala'squantity}
    \bb E[\la \exp(\lambda-\theta t)(R_{1,2}-q)^2\ra_t]
\end{align}
and manage to prove that its derivative \eqref{interpolation derivative} is negative. Hence, to bound the value at $t=1$ we only need to get the value at $t=0$. However, it turns out that the pSK model is more difficult to analyze at the end of $t=0$ since it correlates as opposed to the i.i.d. Bernoulli random variables for the standard SK model. Here, we use the method of asymptotic integral expansion developed in \citep{he2023hidden1} and use convex analysis to generate an upper bound for the quantity \eqref{latala'squantity}. Moreover, a special characteristic of the pSK model compared with the standard SK model is the phase transition. We show that the convergence rate of $R_{1,2}-q$ can take countable numbers of values, as opposed to the usual $\frac{1}{\sqrt n}$ rate for the standard SK model. Interestingly, despite the standard SK model, Latala's method only works for $\theta<\frac{1}{2}$, our results suggest that his method for the localized replica concentration works for all $\theta$ as long as the replica concentration holds.

For the information divergences, we work on a different form than \eqref{interpolation derivative2}, which requires us to use more complicated strategies and computations, which we discuss in section \ref{sect3.3}.

\smallbreak\emph{Guerra-Talagrand 1-step Replica Symmetry Breaking (1RSB) Analysis.} An unfortunate fact is that all the above strategies only work for the small clique regime due to the approximation error of \eqref{interpolation derivative2} being not controllable when $k\asymp n$. And it is unknown if the spins concentrate. The technique to resolve the $k\asymp n$ is more intricate, where the idea is to use the Poisson-Dirichlet Process to prove that the concentration holds. This idea is discovered independently by Guerra \citep{guerra2003broken} to give a sharper upper bound for the limiting free energy which is later generalized by Talagrand to prove the celebrated Parisi Formula \citep{talagrand2006parisi}. Although this line of work focuses on the limiting free energy of the standard SK model rather than the concentration of average magnetizations, it turns out that this interpolation strategy can also be used here to certify the concentration. 

The underlying principle is that the RSB free energy is always a lower bound of the replica symmetry free energy \citep{guerra2001sum}. Therefore if the 1-RSB upper bound is strictly smaller than the replica symmetry upper bound we can claim that the replica symmetry breaking happens.  Talagrand discovered that the gap between the two interpolating free energy being zero is directly related to the concentration of replica \citep{talagrand2011mean2}. We generalize his method to the pSK model and prove that not only the replica but also the average magnetizations of both the spins inside and outside the clique concentrate.

\smallbreak\emph{The Local Cavity Methods.}\label{localcavity} Although our 1RSB analysis implies the convergence of average magnetizations and replica, it does not imply the exact convergence rates and the moments. This is solved by an additional leave-one-out procedure, also referred to as the cavity method. The cavity method separates a single spin out combining the interpolating Hamiltonians that have been historically used to analyze the Hopfield model and the standard SK model together \citep{talagrand2010mean}. (Here, we take the simplest case where $k=n$ as an example, more complicated ones for the low-temperature regime  as well as our local variants, are presented in section \ref{sect4.3} and Appendix B.)
\tny{\begin{align}\label{cavitt}
   \mca H^{c}_t(\bfa\sigma)&=\ub{\sum_{1\leq i<j\leq n-1}\frac{\theta g_{ij}}{\sqrt n}\sigma_i\sigma_j+\sum_{i=1}^{n}h_i\sigma_i+\sqrt t\bl\frac{\theta}{\sqrt n}\sum_{i=1}^{n-1}g_{ij}\sigma_i\sigma_n\br+\sqrt{1-t}\theta\sqrt q z_i\sigma_i}_{\text{the Standard SK Model}}\nnb\\
    &+\ub{\sum_{i,j\in [k-1]}\frac{\theta_1}{2k}\sigma_i\sigma_j+\frac{t\theta_1}{k}\sum_{i\in [k]}\sigma_i\sigma_n+(1-t)\theta_1\mu\sigma_n}_{\text{the Hopfield Model}},
\end{align}}
where we denote $\nu_0$ to be the expectation over $\mca H^c_0$, $q=\nu_0(\epsilon_1\epsilon_2)$ and $\mu=\nu_0(\epsilon_1)$.
At the $t=0$ end, we note that the last spin is completely independent of the rest. Therefore, using the intuition that the function on $n-1$ spins looks like the function on $n$ spins, one can derive higher moments of replicas with lower ones. Denote $\epsilon_j:=\sigma_k^{j}$ where $j$ is the index of replica, $m_i:=\frac{1}{k}\sum_{i=1}^k\sigma_i$, and $m_i^{-}:=\frac{1}{k}\sum_{\ell=1}^{k-1}\sigma_\ell^i$, we   get the following informal derivation (The formal version is delayed to section \ref{sect4.3}) using the fact that any two spins are indistinguishable,
\tny{
\begin{align}\label{cavitysubeq1}
    \nu((m_1-\mu)^{r+1})&=\nu((\epsilon_1-\mu)(m_1-\mu)^{r})\nnb\\
        &\approx\nu((\epsilon_1-\mu)(m^-_1-\mu)^{r})+\frac{r}{k}\nu((1-\epsilon_1\mu)(m_1^--\mu)^{r-1})\nnb\\
        &\approx\nu_0((\epsilon_1-\mu)(m_1-\mu)^r)+\nu^\prime_0((\epsilon_1-\mu)(m_1-\mu)^r)+\frac{r(1-\mu^2)}{k}\nu_0((m_1-\mu)^{r-1})\nnb,\\
        &\approx\nu^\prime((\epsilon_1-\mu)(m_1-\mu)^r)+\frac{r(1-\mu^2)}{k}\nu((m_1-\mu)^{r-1}),
\end{align}}
where we already use the fact that under $\nu_0$, $\epsilon_i$ and $m_1^-$ are independent. Moreover, we use \eqref{cavitt} to derive that 
\tny{\begin{align}\label{cavitysubeq2}
    \nu_t^\prime(f)&={\theta^2\bl \nu_t(f\epsilon_1\epsilon_2(R_{1,2}-q))-2\sum_{\ell\leq 2}\nu_t(f\epsilon_{\ell}\epsilon_3(R_{\ell,3}-q))+\nu_t(f\epsilon_3\epsilon_4(R_{3,4}-q))\br}\nnb\\
    &+{\theta_1\bl\sum_{\ell\leq 2}\nu_t(f\epsilon_\ell(m_{\ell}-\mu)-2\nu_t(f\epsilon_3(m_3-\mu))\br} .
\end{align}}
We notice that the above equation implies that higher moment terms can represent the derivative of lower moment terms. 
Combining \eqref{cavitysubeq1} with \eqref{cavitysubeq2}, we get the moment iteration, representing the higher moments as a linear function of lower moment terms.

For the planted SK model, two sets of spins of $m_S:=\frac{1}{k}\sum_{i\in S}\sigma_i$ and $m_{S^c}:=\frac{1}{n-k}\sum_{i\in S^c}\sigma_i$ requires different treatments and we term this as the \emph{local} cavity method. We use three different cavity paths (one for high-temperature regime, one for low-temperature regime, and one for the out-of-clique spins) to derive the moments of the average magnetizations and the replica overlaps inductively, completing the proof of the concentration of measure together with the certificate of concentration. Using these cavity methods, we develop a linear multivariate equation that gets from the lower to the higher moments. This moment iteration and the certificate of concentration given by the 1RSB analysis recover the tail bounds of both the replicas and the average magnetization. These results imply the upper bounds of both testing and recovery.
\subsection{The Proof Outline of Theorem \ref{thm10}}\label{sect3}

This section presents a proof outline of theorem \ref{thm10}. As a road map, our proof contains the upper and lower bounds, which are further separately organized according to the different temperature regimes. For the upper bounds, our proof is built upon analyzing the testing statistics under the different temperature regimes, which requires the concentration inequalities of the spins as intermediate results. For the lower bounds, our proof involves carefully upper bounding the information divergences, which translates to the lower bounds through Le Cam's lemma.

Our presentation is organized as follows: the testing procedures, along with their theoretical guarantees and corresponding lower bounds, are presented directly according to the different temperature regimes; section \ref{sect3.2} presents a general proof outline for the upper bounds; section \ref{sect3.3} presents a general proof outline for the lower bounds.

\smallbreak\emph{High Temperature Regime.}\label{sect311}
At the high-temperature regime, we present procedure \ref{alg:two} as a combination of local and global tests. The intuitions behind the two-phase tests root in the trade-off between two contradicting factors controlling the signal-to-noise ratio. When the clique becomes larger, (1) the correlation between sites becomes less observable. (2) The proportion of the vertices in a clique gets larger and the clique itself becomes more observable. Therefore, these trade-offs finally translate to the \emph{elbow effect} on the complexity of tests. Our results are then summarized in lemma \ref{thm1}.
\begin{algorithm}[htbp]
\caption{High Temperature Test}\label{alg:two}
\KwData{$\{\bfa\sigma^{(i)}\}_{i\in[m]}$ with $\bfa\sigma\in\{-1,1\}^n$}
\eIf{$k=o(n^{\frac{2}{3}})$}{Compute empirical correlation $\wh {\bb E}[\bfa\sigma\bfa\sigma^\top]= \frac{1}{m}\sum_{i=1}^m\bfa\sigma^{(i)}\bfa\sigma^{(i)\top}$\;
Going over all subset $S\subset[n]$ with $|S|=k$.  Compute $\phi_S = \frac{1}{k}\mbbm 1_{S}^\top\wh {\bb E}[\bfa\sigma\bfa\sigma^\top]\mbbm 1_{S}$\; 
Reject Null if $\phi_1 = \sup_{S:|S|=k}\phi_S\geq \tau_\delta$ where $\tau_\delta\in\lef(0, \frac{1-\theta_1\bb E[\sech^2(\theta\sqrt qz+h)]^2}{(1-\theta_1\bb E[\sech^2(\theta\sqrt qz+h)])^2}\rig)$ \;}
{Compute empirical correlation $\phi_2= \frac{1}{mk}\sum_{i=1}^m\bfa\sigma^{(i)\top}\bfa\sigma^{(i)}-1$\;
Reject Null if $\phi_2>\tau_\delta$, with $\tau_\delta\in\lef(0,\frac{2\theta_1\bb E[\sech^2(\theta\sqrt qz+h)]-\theta_1(1+\theta_1)(\bb E[\sech^2(\theta\sqrt qz+h)])^2}{(1-\theta_1\bb E[\sech^2(\theta\sqrt qz+h)])^2}\rig)$;}
\end{algorithm}
\begin{lemma}\label{thm1}
  Assume that $k\log k=o(n)$, $0<\theta_1<\frac{1}{\bb E[\sech^2(\theta\sqrt{q}z+h)]}$, and the condition \eqref{atcondition} holds. Then algorithm \ref{alg:two} is asymptotically powerful if :
  \begin{enumerate}
      \item  $k=o(n^{\frac{2}{3}})$ and $m=\omega\lef( k\log n\rig) $;
      \item $k=\omega(n^{\frac{2}{3}})$ and $m=\omega\lef( \frac{n^2}{k^2}\rig)$.
  \end{enumerate}
    The regions of sample complexity $m$ such that all tests are asymptotic powerless are given by:
    \begin{enumerate}
        \item If $k=o( n^{\frac{1}{2}} )$ and $\theta_1<\frac{1}{2\bb E[\sech^2(\theta\sqrt qz+h)]}$, then $m=o\lef(k\log n\rig)$;
        \item If $k=\Omega( n^{\frac{1}{2}})$ and $\theta_1<\frac{1}{2\bb E[\sech^2(\theta\sqrt qz+h)]}$, then $m=o\lef(\frac{n}{k\theta^2}\wedge\frac{n^2}{k^2}\wedge k\log n\rig)$;
        \item If $k=o( n^{\frac{1}{2}})$ and $\theta_1\geq\frac{1}{2\bb E[\sech^2(\theta\sqrt qz+h)]}$, then $m=o\lef(k\frac{\log n}{\log k}\rig)$;
        \item If $k=\Omega( n^{\frac{1}{2}})$ and $\theta_1\geq\frac{1}{2\bb E[\sech^2(\theta\sqrt qz+h)]}$, then $m=o\lef(\frac{n}{k\theta^2}\wedge\frac{n^2}{k^2}\wedge k\rig)$.
    \end{enumerate}
\end{lemma}

\begin{remark}
    Recall that in the pRFCW model, this lower bound matches with the rate in the upper bound given by theorem \ref{thm1}. However, for the pSK model, a gap exists between the lower and upper bounds for the large clique region of $\Omega(n^{1/2})$ at the high-temperature regime. Despite this gap being conjectured to be improvable, and the minimax optimal rate is conjectured to be the same as the pRFCW model, it is unclear what improvement can be made in the proof to control the information divergence sharply than section \ref{sect3.2}. Moreover, if we consider a vanishing `noise to signal ratio' of $\frac{\theta}{\theta_1}\lesssim \sqrt{\frac{k}{n}}$ and $\theta_1<\frac{1}{\bb E[\sech^2(\theta\sqrt qz+h)]}$, we have matching upper and lower bounds. The reader can later see that the ratio of $\frac{k}{n}$ is a fundamental barrier to the proposed techniques in this work and any rate is subject to this restriction.
\end{remark}

\smallbreak\emph{Low Temperature Regime.}\label{sect312}
This section presents the result of low temperature regime where $\theta_1>\frac{1}{\bb E[\sech^2(\theta\sqrt{q}z+h)]}$. Our results imply that the minimax optimal rate of the low-temperature pSK model matches that of the pRFCW model. We present the test \ref{alg:three} combining the local and global procedures. 

\begin{algorithm}[htbp]
\caption{Low Temperature Test}\label{alg:three}
\KwData{$\{\bfa\sigma^{(i)}\}_{i\in[m]}$ with $\bfa\sigma\in\{-1,1\}^n$}
\uIf{$k=o(\sqrt n)$}{Going over all subset $S\subset[n]$ with $|S|=k$.  Compute $\phi_S = \frac{1}{m}\sum_{j=1}^m\lef|\frac{1}{k}\sum_{i\in S}\sigma^{(j)}_i\rig|$\; 
Reject Null if $\phi_3 = \sup_{S:|S|=k}\phi_S\geq \tau_\delta$ with $\tau_\delta\in\lef(0, x\rig)$ with $x$ defined by the positive solution to $x=\bb E[\tanh(\theta_1x+h)] $ \;}
\Else{Compute statistics $\phi_4:=\frac{1}{m}\sum_{j=1}^m\lef|\frac{1}{k}\sum_{i=1}^n\sigma^{(j)}_i\rig|$\;
\uIf{$k\asymp n$}{
Reject Null if $\phi_4 > \tau_\delta$  for $\tau_\delta\in\lef(\sqrt{\frac{2n}{\pi k^2}},\frac{\sqrt n}{k}\sqrt{\frac{2}{\pi}}\exp\lef(-\frac{x^2k^2}{2n}\rig)+x\lef[1-2\Phi\lef(-\frac{xk}{\sqrt n}\rig)\rig]\rig)$ with $\Phi$ being the cumulative distribution function of the standard Gaussian random variable\;}
\Else{Reject Null if $\phi_4>\tau_\delta$ for $\tau_\delta\in (0,x)$ with $m$ defined by the positive solution to with $x$ defined by the positive solution to $x=\bb E[\tanh(\theta_1x+h)] $\;}
}
\end{algorithm}

\begin{lemma}\label{thmlowtempfindclique}
    Assume that $k\log k=o(n)$, $\theta_1>\frac{1}{\bb E[\sech^2(\theta\sqrt{q}z+h)]}$, and the condition \eqref{atcondition} holds. Then algorithm \ref{alg:three} is asymptotically powerful if:
    \begin{enumerate}
        \item  $k=o( \sqrt n)$ and $m=\omega\lef( \log n\rig)$;
        \item  $k\asymp \sqrt n$ and $m=\omega(1)$;
        \item $k=\omega( \sqrt n)$ and $m= 1$.
    \end{enumerate}
    For the lower bound, when $k=o(\sqrt n)$ the region of sample complexity $m$ such that all tests are asymptotic powerless is given by $m\leq C\log n$.
\end{lemma}
\begin{remark}
    We notice that the approximation bounds are consistent with the minimax lower bound in the low temperature. Since under the case of $k\gtrsim\sqrt n$, we already have a constant in the upper bound, a lower bound is unnecessary. 
\end{remark}

\smallbreak\emph{Critical Temperature Regime.}\label{sect313}
This paragraph presents the results of the critical temperature regime when $\theta_1=\frac{1}{\bb E[\sech^2(\theta\sqrt qz+h)]}$. Our results imply that the statistical diagrams for the critical temperature exhibit a mixed behavior of the high and low temperatures. We present the test \ref{alg:five} as a combination of local and global tests.




\begin{algorithm}[htbp]
\caption{Critical Temperature Test }\label{alg:five}
\KwData{$\{\bfa\sigma^{(i)}\}_{i\in[m]}$ with $\bfa\sigma\in\{-1,1\}^n$}
\uIf{$k=o\lef(n^{\frac{4\tau-2}{8\tau-5}}\rig)$}{Compute scaled empirical correlation matrix $\wh{\bb E}[\bfa\sigma\bfa\sigma^\top]=\frac{1}{m}\sum_{j=1}^m\bfa\sigma^{(j)}\bfa\sigma^{(j)\top}$\;
Go over all subset $S\subset[n]$ with $|S|=k$ and  compute $\phi_{S} = k^{-(4\tau-3)/(2\tau-1)}\lef(\mbbm 1_{S}^\top\wh {\bb E}[\bfa\sigma\bfa\sigma^\top]\mbbm 1_{S}\rig)$\; 
Reject Null if $\phi_{5} = \sup_{S:|S|=k}\phi_S\geq \tau_\delta$  for $\tau_\delta\in\lef(0,\pi^{-\frac{1}{2}}(2\ca V(\tau))^{\frac{1}{2\tau-1}}\Gamma(\frac{2\tau+1}{4\tau-2})\rig)$ where $\ca V(\tau):=\frac{((2\tau)!)^2\bb V(\tanh(\theta\sqrt qz+h))(\bb E[\sech^2(\theta\sqrt qz+h)])^{4\tau-2}}{2^{2\tau-1}\lef(\bb E\lef[(1+\tanh(\theta\sqrt qz+h))\sum_{k=0}^{2\tau-1}\frac{k!}{2^k}S(2\tau-1,k)(\tanh(\theta\sqrt qz+h)-1)^{k}\rig]\rig)^2}$ and $S(n,k)$ is the second type of Stirling numbers\;}
\Else{Compute the scaled correlation $\phi_6=m^{-1}k^{-\frac{4\tau-3}{2\tau-1}}\sum_{j=1}^m\lef(\lef(\sum_{i=1}^n\sigma_i^{(j)}\rig)^2-n\rig)$\;
Reject Null if $\phi_6\geq\tau_\delta$ for $\tau_\delta\in\lef(0,\pi^{-\frac{1}{2}}(2\ca V(\tau))^{\frac{1}{2\tau-1}}\Gamma\lef(\frac{2\tau+1}{4\tau-2}\rig)\rig)$ \;}

\end{algorithm}
\begin{lemma}\label{guaranteecriticaltest}
 Assume that $k^{\frac{2\tau}{2\tau-1}}\log^{\frac{2\tau-2}{2\tau-1}}k=o(n)$, $\theta_1=\frac{1}{\bb E[\sech^2(\theta\sqrt{q}z+h)]}$, and condition \eqref{atcondition} holds.  For the upper bounds, algorithm \ref{alg:five} is asymptotically powerful if
 \begin{enumerate}
     \item  $k=o\lef(n^{\frac{4\tau-2}{8\tau-5}}\rig)$ and $m=\omega\lef( k^{1/(2\tau-1)}\log n\rig)$;
     \item   $n^{\frac{4\tau-2}{8\tau-5}}\lesssim k\lesssim n^{\frac{2\tau-1}{4\tau-3}}$ and $m=\omega\lef( n^2k^{-\frac{2(4\tau-3)}{2\tau-1}}\rig)$;
     \item  $n^{\frac{2\tau-1}{4\tau-3}}\lesssim k\lesssim n^{\frac{2\tau-1}{2\tau}}\log^{-\frac{2\tau-2}{2\tau-1}}n$ and $m=\omega(1)$.
 \end{enumerate}
    For the lower bounds, the region of sample complexity $m$ such that all tests are asymptotic powerless is if:
    \begin{enumerate}
        \item  $k=o\lef( n^{\frac{4\tau-2}{8\tau-5}} \rig)$ and $m=o\lef(\lef(\frac{k}{\log k}\rig)^{\frac{1}{2\tau-1}} \log n\rig)$; 
        \item  $n^{\frac{4\tau-2}{8\tau-5}}\lesssim k\lesssim n^{\frac{2\tau-1}{4\tau-3}}$ and $m=o\lef(n^2k^{-\frac{2(4\tau-3)}{2\tau-1}}\rig)$;
        \item $n^{\frac{2\tau-1}{4\tau-3}}\lesssim k\lesssim n^{\frac{2\tau-1}{2\tau}}\log^{-\frac{2\tau-2}{2\tau-1}}n$ and $m\leq C$ for some constant $C\geq 1$.
    \end{enumerate}
\end{lemma}

\begin{remark}
   For the critical temperature, the upper bound matches up to a logarithmic factor with the lower bound. However, compared with the results in the pRFCW model, where we obtain the minimax optimal rate for all $k<n$, the region of $k\gtrsim n^{\frac{2\tau-1}{2\tau}}$ in the pSK model remains open.
\end{remark}

\subsubsection{The Upper Bounds }\label{sect3.2}
Our upper bounds in this work are based on a precise characterization of the limiting variance and the tail bounds for the average of spins inside and outside the clique $S$. As discussed in \citep{he2023hidden1}, the pointwise convergence of mgf implies the convergence of moments but does not provide uniform tail bounds. To obtain the uniform tail bounds, one must also use the local-to-global machine given by \citep{he2023hidden1}. Our proof of the following lemma makes use of the local smart paths discussed in section \ref{sect300} to show that the m.g.f. of the average of spins in $S$ of the pSK model can be well approximated by that of the pRFCW model. At the center of our method is a novel smart path that can decouple the correlation between spins in $S$ and $S^c$.  We present our results for the upper bounds and lower bounds as follows.
\begin{lemma}[Limiting Distributions for the pSK model]\label{cltrfcr}
Assume that $h_i\sim \mu(h)$ is i.i.d. in $L_1$ and symmetric around $0$.
For the pSK model whose Hamiltonian  defined by \eqref{skhamilt} with $k$ spins in the clique and the local sum of spins in the clique denoted by $\sum_{i\in S}\sigma_i$, the following hold:
\begin{enumerate}
\item 
Assume that $\theta^2k\log k=o(n)$. In the high temperature regime with $\theta_1<\frac{1}{\bb E[\sech^2(\theta \sqrt{q}z+h)]}$ where $q$ is the root to $q=\bb E[\tanh^2(\theta\sqrt{q}z+h)]$, for $t\in\bb R$ pointwise,
\sm{\begin{align*}
    \bb E\lef[\exp\lef(\frac{t\sum_{i\in S}\sigma_i}{\sqrt k}\rig)\rig]\to \exp\lef(\frac{\ca Vt^2}{2}\rig)\quad\text{ and}\quad\bigg\Vert k^{-1/2}{\sum_{i=1}^n\sigma_i}\bigg\Vert_{\psi_2}<\infty.
\end{align*}}
with $\ca V:=\frac{1-\theta_1(\bb E[\sech^2(\theta\sqrt q z+h)])^2}{(1-\theta_1\bb E[\sech^2(\theta\sqrt q z+h)])^2}$.

 \item 
Assume that $\theta^2k\log k=o(n)$. In the low temperature regime of $\theta_1>\frac{1}{\bb E[\sech^2(\theta\sqrt qz+h)]}$, $x=\bb E[\tanh(\sqrt{\theta_1} x +h)]$ have two nonzero symmetric roots defined by $-x_1^*<0<x_1^*$. Define $q$ to be the root to $q:=\bb E[\sech^2(\theta\sqrt{q}z+\sqrt{\theta_1}x^*_1+h)]$. Define $\ca C_1=(0,\infty)$ and $\ca C_2=\ca C_1^c$. Then, for $t\in\bb R$ and $\ell\in\{1,2\}$, pointwise,
\sm{\begin{align}\label{gtr0clt}
    &\bb E\lef[\exp\lef(t\frac{\sum_{i\in S}(\sigma_i -\sqrt{\theta_1}x_\ell)}{\sqrt k}\rig)\bigg |{\frac{\sum_{i\in S}\sigma_i}{k}\in \ca C_\ell}\rig]{\to} \exp\lef(\frac{\ca V(m_1)t^2}{2}\rig)
\end{align}}
and
\sm{$
   \Vert k^{-1/2}\sum_{i=1}^k(\sigma_i-\sqrt{\theta_1}x_{\ell})|k^{-1}\sum_{i=1}^k\sigma_i\in\ca C_{\ell}\Vert_{\psi_2}<\infty
$}
with \\\sm{$\ca V(x_1^*):=\frac{(1-\theta_1(\bb E[\sech^2(\sqrt{\theta_1}x_1^*+\theta\sqrt q+h)])^2-\bb E[\tanh(\sqrt{\theta_1}x_1^*+\theta\sqrt q+h)]^2)}{(1-\theta_1\bb E[\sech^2(\sqrt{\theta_1}x_1^*+\theta\sqrt q+h)])^2}$}.
 \item 
Assume that $\theta^2k^{\frac{2\tau}{2\tau-1}}\log^{\frac{2\tau-2}{2\tau-1}}k=o(n)$. At the critical temperature $\theta_1=\frac{1}{\bb E[\sech^2(\theta\sqrt qz+h)]}$, assume that the flatness of $0$ for function $H(x)$ defined by \eqref{taylorcond} is $\tau$,  then for $t\in\bb R$, pointwise  
\sm{\begin{align}\label{supergaussian}
        \bb E\lef[\exp\lef(\frac{t\sum_{i=1}^n\sigma_i}{k^{\frac{4\tau-3}{4\tau-2}}}\rig)\rig]\to \int_{\bb R}\frac{(2\tau-1)x^{2\tau-2}}{\sqrt{2\pi v(0)}}\exp\lef(-\frac{x^{4\tau-2}}{2v(0)}+tx\rig)dx,
    \end{align}}
    and
\sm{$        \lef\Vert n^{-\frac{4\tau-3}{4\tau-2}}\sum_{i=1}^n\sigma_i\rig\Vert_{\psi_{4\tau-2}}<\infty,
$}
with \sm{\begin{align}\label{stirlingv}
    v:&=((2\tau)!)^2\bb V(\tanh(\theta\sqrt qz+h))(\bb E[\sech^2(\theta\sqrt qz+h)])^{4\tau-2}\nnb\\
    &\bb E\lef[(1+\tanh(\theta\sqrt qz+h))\sum_{k=0}^{2\tau-1}\frac{k!}{2^k}S(2\tau-1,k)(\tanh(\theta\sqrt qz+h)-1)^{k}\rig]^{-2}.
\end{align}} And if we are in the second case of \eqref{taylorcond} then \eqref{gtr0clt} holds.
\end{enumerate}
\end{lemma}


Then, we give another result on the central limit theorem of the standard SK model. This result is obtained for the extremely low value of $\theta<\frac{1}{2}$ and fixed $h$ case in \citep{hanen2009limit}. In section \ref{sect4}, we extend it towards a much larger replica symmetric region given by the condition in \ref{withintheATline}. The proof requires a technical argument by Talagrand \citep{talagrand2011mean2}. Our results are summarized as follows.
\begin{lemma}\label{cltconsk}
    We consider the model with Hamiltonian defined by \eqref{skhamilt} and i.i.d. random field $h\sim\mu$ in $L_1$.
     Define $q:=\bb E[\tanh^2(\theta\sqrt{q} z+h)]$, the mean magnetization as $m:=\frac{1}{n}\sum_{i=1}^n\sigma_i$ and limiting variance $V:=1-(\bb E[\tanh(\theta\sqrt qz+h)])^2$. Assume $\theta$ satisfies $\bb E[2\theta^2\sech^4(\theta\sqrt qz+h)]<1$ and the `within the AT line' given by definition \ref{withintheATline} holds for $c=0$. Let $z\sim N(0,1)$. Then we have for all $\eta\in\bb N$, pointwise,
     \sm{
    \begin{align*}
        \bb E\lef[\bl\sqrt{\frac{n}{V}} (m-\bb E[\tanh(\theta\sqrt qz+h)]\br^\eta\rig]\to\bb E[z^\eta]\quad\text{ and }\quad\Vert\sqrt n (m-\bb E[\tanh(\theta\sqrt qz+h)])\Vert_{\psi_2}<\infty.
    \end{align*}}
    And in particular, when $h$ is symmetric w.r.t. $0$, $V=1$ and $\bb E[\tanh(\theta\sqrt qz+h)]=0$.
\end{lemma}

\subsubsection{The Lower Bounds}\label{lbchannel}
We obtain sharp information divergence control as required from Le Cam's lemma to derive the lower bounds. However, the technical barrier exists when we hope to estimate the TV distance between discrete Gibbs measures with little existing literature. Our strategy relies on first solving the analogous problem under the pRFCW measure, whose results are given in \citep{he2023hidden1}. Then we formalize a valid smart path such that the results obtained for the pRFCW model can be used to approximate that of the pSK model. This smart path follows a different design idea from our upper bounds, introducing heterogeneous replicas. To show that the approximation is valid, we show the heterogeneous replicas (the correlation between two sets of spins that are sampled from different distributions) concentrates. 

As a starting point we illustrate how the lower bounds are connected with the information divergences through the following standard Le Cam's lemma.
\begin{lemma}[Conditional Le Cam]\label{condlecam} 
Assume that $P$, $Q$ are two probability measure dependent on random variables $\bfa h,\bfa g$ with $P\ll Q$ almost surely w.r.t. $\mu(\bfa h,\bfa g)$, we define
 $D_{\chi^2}(P,Q|\bfa h,\bfa g):=\int \lef((\frac{P(d\bfa\sigma|\bfa h,\bfa g)}{Q(d\bfa\sigma|\bfa h,\bfa g)})^2-1\rig)Q(d\bfa\sigma|\bfa h,\bfa g)$. Denote $\bb P_{S,m}=\bb P_S^{\otimes m}$ and $\bar{\bb P}_m=\frac{1}{\binom{n}{k}}\sum_{S:|S|=k}\bb P_{S,m}$. Denote $\bfa h^m=(\bfa h_1,\ldots,\bfa h_m),\bfa g^m=(\bfa g_1,\ldots,\bfa g_m)$ with $\bfa g_i$ and $\bfa h_i$ i.i.d. Then for all $\{0,1\}$-valued test statistics $\psi$ constructed by $\{\bfa\sigma^{(i)}\}_{i\in[m]}$, we have
\sm{
\begin{align*}
    \inf_{\psi}\bigg[\bb P_{0,m}(\psi=1)+\sup_{S:|S|=k}\bb P_{S,m}\lef(\psi=0\rig)\bigg]\geq 1-\frac{1}{2}\sqrt{\bb E[D_{\chi^2}(\bar{\bb P}_{m},\bb P_{0,m}|\bfa h^m,\bfa g^m)]}
.\end{align*}}
\end{lemma}
To estimate the chi-square divergence between two pSK measures, we design a novel smart path that can `decouple' the correlation between the spins within and without the clique set $S$. And our results are summarized in Lemma \ref{approximation bounds chisq}.
\begin{lemma}\label{approximation bounds chisq}
Consider two set $S$, $S^\prime$ such that $|S|=|S^\prime|=k$ and $|S\cup S^\prime |=k+r$ to be the two clique index set. Let $q$ be the root to $q=\bb E[\tanh^2(\theta\sqrt{q}z+h)]$.
Denote $\bb E_{SK,p}[R_{1,2}]$ as the mean replica overlap for the standard SK model with $p$ spins.
Let $q_{k,r}, q^{(1)}$ be the solution to the following equations
\begin{align*}
    q^{(1)}=\bb E[\tanh^2(\theta\sqrt q_{k,r}z+h)],\quad q_{k,r}=\frac{k+r}{n}q^{(1)}+\frac{n-k-r}{n}\bb E_{SK,n-k-r}[R_{1,2}].
\end{align*}
    Define $h^\prime=h+\theta\sqrt{q_{k,r}}z$ for $z\sim N(0,1)$ and $z\perp h$. Define $\bb Q_S(\cdot|\bfa h^\prime)$, $\bb Q_0(\cdot|\bfa h^\prime)$ to be the pRFCW measure with outer field measure $\mu(h^\prime)$, inverse temperature to be $\theta_1$ and the clique index set to be $S$. For arbitrary events $A_S\in\Sigma\lef(\bfa\sigma_S\rig),A_{S^\prime}\in\Sigma(\bfa\sigma_{S^\prime})$, the following holds:
    \begin{enumerate}
        \item When $\theta_1<\frac{1}{2\bb E[\sech^2(\theta\sqrt qz+h)]}$ we   have when $\theta^2k\log k=O(n)$:
        \sm{\begin{align*}
        \lef|\bb E\lef[\frac{\bb P_{S}(\cdot|\bfa h,\bfa g)\bb P_{S^\prime}(\cdot|\bfa h,\bfa g)}{\bb P_0(\cdot|\bfa h,\bfa g)}\rig]-\bb E\lef[\frac{\bb Q_{S}(\cdot|\bfa h^\prime)\bb Q_{S^\prime}(\cdot|\bfa h^\prime)}{\bb Q_0(\cdot|\bfa h^\prime)}\rig]\rig|=O\lef(\frac{\theta^2k}{n}\rig).
    \end{align*}}
     \item When $\theta_1\in\lef[\frac{1}{2\bb E[\sech^2(\theta\sqrt qz+h)]},\frac{1}{\bb E[\sech^2(\theta\sqrt qz+h)]}\rig)$ then, when $\theta^2 k\log k=O(n)$:
     \sm{
     \begin{align*}
        \bb E\lef[\frac{\mbbm 1_{A_{S}}\mbbm 1_{A_{S^\prime}}\bb P_{S}(\cdot|\bfa h,\bfa g)\bb P_{S^\prime}(\cdot|\bfa h,\bfa g)}{\bb P_0(\cdot|\bfa h,\bfa g)}\rig]\leq\exp\lef(\frac{C\theta^2k\log k}{n}\rig)\bb E\lef[\frac{\mbbm 1_{A_{S}}\mbbm 1_{A_{S^\prime}}\bb Q_{S}(\cdot|\bfa h^\prime)\bb Q_{S^\prime}(\cdot|\bfa h^\prime)}{\bb Q_0(\cdot|\bfa h^\prime)}\rig].
    \end{align*}}
    \item When $\theta_1>\frac{1}{\bb E[\sech^2(\theta\sqrt q z+h)]}$ and $\theta k=O(\sqrt n)$:
    \sm{\begin{align*}
         \bb E\lef[\frac{\mbbm 1_{A_{S}}\mbbm 1_{A_{S^\prime}}\bb P_{S}(\cdot|\bfa h,\bfa g)\bb P_{S^\prime}(\cdot|\bfa h,\bfa g)}{\bb P_0(\cdot|\bfa h,\bfa g)}\rig]\leq\exp\bl C\frac{\theta^2k^2}{n}\br\bb E\lef[\frac{\mbbm 1_{A_{S}}\mbbm 1_{A_{S^\prime}}\bb Q_{S}(\cdot|\bfa h^\prime)\bb Q_{S^\prime}(\cdot|\bfa h^\prime)}{\bb Q_0(\cdot|\bfa h^\prime)}\rig].
    \end{align*}}
    \item When $\theta_1=\frac{1}{\bb E[\sech^2(\theta\sqrt q z+h)]}$ and the global optimum $0$ of equation \eqref{taylorcond} has flatness $\tau$,  when $\theta^2k\log^{\frac{1}{2\tau-2}}k=O(n^{\frac{2\tau-1}{4\tau-4}})$:
    \tny{
    \begin{align*}
        \bb E\lef[\frac{\mbbm 1_{A_{S}}\mbbm 1_{A_{S^\prime}}\bb P_{S}(\cdot|\bfa h,\bfa g)\bb P_{S^\prime}(\cdot|\bfa h,\bfa g)}{\bb P_0(\cdot|\bfa h,\bfa g)}\rig]\leq\exp\bl\frac{C\theta^2k^{\frac{4\tau-4}{2\tau-1}}\log^{\frac{2}{2\tau-1}}k}{n}\br\bb E\bigg[\frac{\mbbm 1_{A_{S}}\mbbm 1_{A_{S^\prime}}\bb Q_{S}(\cdot|\bfa h^\prime)\bb Q_{S^\prime}(\cdot|\bfa h^\prime)}{\bb Q_0(\cdot|\bfa h^\prime)}\bigg].
    \end{align*}}
    \end{enumerate}
\end{lemma}
\begin{remark}
    We see that the approximation at the critical temperature regime has a smaller validity scope than the high and low-temperature regimes. This is due to the unique replica convergence rate in the next lemma. In addition, we also show that no log factor is needed for the high temperature regime due to a \emph{twice computation trick}.
\end{remark}
To facilitate the above results, we need to obtain a sharp upper bound on the concentration of heterogeneous replicas, which is also specified in section \ref{sect300}. This result finally translates to the upper bound on the difference between the information divergences under the pRFCW model and the pSK model.

\smallbreak\emph{The Concentration of Heterogeneous Replicas.}\label{sect3.3}
In this section, we present the exponential tail concentration inequalities for the local replicas defined by $R_{1,2}^D:=\sum_{i\in D}\frac{1}{|D|}\sigma_i^1\sigma_i^2$ for $D\subset[n]$ and $|D|=o(n)$. This result is crucial for the validity of the approximations presented in section \ref{sect3.2}, suggesting an interesting phenomenon that appears uniquely in the pSK model and not in the standard SK model: \emph{When the two local replicas are sampled from the pSK model with heterogeneous clique position, it continues to concentrate, and can have totally different rates from the classical $\frac{1}{\sqrt n}$ in the standard SK model}. Moreover, our results imply that the convergence rates depend on the temperature regimes and can take countable values in the pSK model. On the other hand, we show that the bounding of local replicas instead of global replicas (where $D=[n]$) is a necessary step for the approximation bounds to hold. The following theorem presents our main result in this section.
\begin{lemma}\label{correxponential}
    Consider $S,S^\prime\subset[n]$ and $|S|=|S^\prime|=k$. Denote $D=S\cup S^\prime$.
    Assume that $\bfa\sigma^1\sim\bb P_S(\cdot|\bfa g,\bfa h),\bfa\sigma^2\sim\bb P_{S^\prime}(\cdot|\bfa g,\bfa h)$.  Recall that $x_1^*$ in \eqref{taylorcond} has flatness $\tau$ defined by \eqref{flatnessofopt}. We define $q$ to be the solution to $q=\bb E[\tanh^2(\theta\sqrt{q}z+\theta_1\mu+h)]$, $\mu=\bb E[\tanh(\theta\sqrt qz+\theta_1\mu+h)]$. Then, the following hold:
    
    \begin{enumerate}
       \item At the high temperature regime, for $D=S\cup S^\prime$ with $|D|=k+r$, we   have $\Vert k^{1/2}(R^D_{1,2}-q)\Vert_{\psi_2}<\infty$;
       \item At the low temperature regime, and assume that $S=S^\prime$, we   have $\Vert k^{1/2}(R^S_{1,2}-q)\Vert_{\psi_2}<\infty$;
       \item At the critical temperature, when $\tde c:=\frac{r}{k}=o(n^{-\frac{\tau-2}{2\tau-1}})$, we   have $\Vert k^{\frac{\tau-1}{2\tau-1}}(R_{1,2}^D-q)\Vert_{\psi_2}<\infty$;
       \item At the critical temperature, when $\tde c=\Omega(n^{-\frac{\tau-2}{2\tau-1}})$, we  have $\Vert k^{\frac{1}{2\tau-1}}(R_{1,2}^D-q) \Vert_{\psi_2}<\infty$.
    \end{enumerate} 
\end{lemma}
To prove the above lemma, we use a unique smart path construction to approximate the result of the pSK with the pRFCW measure. The proof utilizes convex analysis and integral expansion extensively. 

\subsection{The Proof Outline of Theorem \ref{thm2}}\label{sect4}
This section presents a proof outline of theorem \ref{thm2}. As is shown in lemma \ref{cltrfcr} and lemma \ref{approximation bounds chisq}, our proof for the small clique regime does not extend to the region with $k\asymp n$ as the approximation methods present in the last subsection only hold for $k=o(n)$ due to the presence of approximation error of order $O(\frac{k\log k}{n})$. Therefore, to test the large clique that is comparable to $n$, we use an alternative 2-step method: (1) an estimation method that establishes the `weak concentration' results and (2) the cavity iterations that give the `strong concentration' results for the high and low-temperature regimes building upon the weak results given in the previous step.

This section is organized as follows: The testing procedure and its theoretical guarantees are given first, with the fundamental limits given by lemma \ref{thm4.1} and the procedure \ref{alg:seven} achieving the upper bounds; Section \ref{sect4.2} states the weak concentration results, which is further promoted to become the strong concentration results in section \ref{sect4.3}.

\begin{algorithm}[htbp]
\caption{Large Clique Test}\label{alg:seven}\label{alg:eight}
\KwData{$\{\bfa\sigma^{(i)}\}_{i\in[m]}$ with $\bfa\sigma\in\{-1,1\}^n$, $c=\frac{k}{n}$, $\mu,q$ calculated according to \eqref{meanfieldeq}}
\If{$\theta_1<\frac{1}{\bb E[\sech^2(\theta\sqrt qz+h)]}$ and $\bb E[\theta^2\sech^4(\theta\sqrt q z+h)]<1$}{Compute Empirical Correlation $\phi_7= \frac{1}{mn}\sum_{i=1}^m\bfa\sigma^{(i)\top}\bfa\sigma^{(i)}$\;
Reject $H_0$ if $\phi_7>\tau_\delta$, with $\tau_\delta\in\lef(1, c\frac{1-\theta_1(\bb E[\sech^2(\theta\sqrt qz+h+\theta_1\mu)])^2}{(1-\theta_1\bb E[\sech^2(\theta\sqrt qz+h+\theta_1\mu)])^2}+(1-c)\rig)$;}
\If{$\theta_1>\frac{1}{\bb E[\sech^2(\theta\sqrt qz+h)]}$ and $\bb E[\theta^2((1-c)\sech^4(\theta\sqrt q z+h)+c\sech^4(\theta\sqrt qz+\theta_1\mu+h))]<1$}{Compute scaled empirical correlation $\phi_8= \frac{1}{n}|\sum_{i=1}^n\sigma_i|$\; 
Reject $H_0$ if $\phi_8 \geq\tau_\delta\in\lef(0, x\rig)$ with $x$ being the positive solution to \eqref{meanfieldeq}\;}
\end{algorithm}
\begin{lemma}\label{thm4.1}
Given $k\asymp n$, the test proposed in algorithm \ref{alg:seven} is asymptotically powerful if:
\begin{enumerate}
    \item When $0<\theta_1<\frac{1}{\bb E[\sech^2(\theta\sqrt qz+h)]}$, $m=\omega(1)$, and the replica symmetry condition in definition \ref{withintheATline} holds with $\mu=0$ ;
    \item When $\theta_1>\frac{1}{\bb E[\sech^2(\theta\sqrt qz+h)]}$, $m=1$,  the replica symmetry condition in definition \ref{withintheATline} and the invertibility condition in the lemma \ref{cltatline} holds.
\end{enumerate}
\end{lemma}
\begin{remark}
     In contrast to the small clique regimes, additional conditions need to be satisfied. This is due to the proximity to the replica symmetry-breaking phase, and extra technical arguments are needed to guarantee the concentration of measure. We do not provide guarantees for the critical temperature, the result of which is closely related to the unsolved mathematical question on the concentration rate of the replica on the AT line of the SK model. (See \citep{talagrand2011mean2}) We discuss it further in section \ref{sect7}.
\end{remark}

\subsubsection{The Weak Concentration}\label{sect4.2} 
This section demonstrates the weak concentration results for the replica and the average magnetization. We recall that the only existing literature that captures the concentration of the average magnetization of the standard SK model holds in an extremely low $\theta$ region where $\theta<\frac{1}{2}$, given by \citep{hanen2009limit}. Instead, we prove a result working in almost all the replica symmetric regions, where $R_{1,2}$ concentrates around a constant value. The underlying idea of the proof is to use the Poisson Dirichlet Process to construct the smart path, an idea first proposed by Guerra in \citep{guerra2003broken}. We review the basic properties of such a process in Appendix C.  For an extensive description of Talagrand's idea in the standard SK model, we refer the reader to \citep{talagrand2011mean2} section 13. However, Talagrand's analysis does not cover the complete results here. 

Our discussion focuses on first deciding the replica-symmetry phase of the pSK model. Then, under this phase, we derive the concentration of measure results. The following results give the weak concentration result.

\begin{lemma}[Mean Field Equations]\label{mfdeq}
Assume that the condition in definition \ref{withintheATline} holds.
 Then we  have
 $R_{1,2}\overset{L_p}{\to} q$, $m:=\frac{1}{k}\sum_{i=1}^k\sigma_i$ satisfies $\lef(m|(m-\mu)\leq C\rig)\overset{L_p}{\to}\mu$ for some small $C>0$ and
the following hold for $c=\frac{r}{k}$,
\begin{align}\label{meanfieldeq}
q:&=c\bb E[\tanh^2(\theta \sqrt {q}z+\theta_1\mu+h)]+(1-c)\bb E[\tanh^2(\theta \sqrt {q}z+h)],
\nnb\\
\mu&\in\argmax_{\mu\in[-1,1]}\lef(\bb E[\log\cosh(\theta z\sqrt {q} +\theta_1\mu+h)]-\frac{\theta_1\mu^2}{2}\rig).
\end{align}
\end{lemma}
Since the precise statistical minimax rates rely on the concentration rate of the quantities given in lemma \ref{mfdeq}, we give formal results extending the above lemma in the following subsection.
\subsubsection{The Strong Concentration}\label{sect4.3}
This section presents the strong concentration results for both the replica and average magnetization in the large clique regime. We provide two main results: Lemma \ref{largcliquetailbound} gives a tail bound, which certifies the convergence, and Lemma \ref{cltatline} gives a moment convergence, which leads to the rejection regions of tests \ref{alg:seven}. Our proof idea is the \emph{cavity method}, a leave-one-out analysis. This method constructs estimates for all moments in an iterative fashion. And we recover all the limiting moments by solving linear equations. Our first tail bounds are given as follows.
\begin{lemma}[Tail Bound]\label{largcliquetailbound}
    Assume that the condition given in definition \ref{withintheATline} holds. When $\theta_1<\frac{1}{\bb E[\sech^2(\theta\sqrt qz+h)]}$ we have $\mu=0$ and
    \begin{align*}
        \Vert \sqrt n(R_{1,2}-q)\Vert_{\psi_2}<\infty,\quad \Vert\sqrt k m\Vert_{\psi_2}<\infty,\quad \Vert\sqrt{n-k}\tde m\Vert_{\psi_2}<\infty.
    \end{align*}
    When $\theta_1>\frac{1}{\bb E[\sech^2(\theta\sqrt qz+h)]}$, there are two solutions in the mean field equation \eqref{meanfieldeq}, symmetric w.r.t. $0$ which we denote by $\mu$ and $-\mu$. Then we have
    \begin{align*}
        &\Vert\sqrt n(R_{1,2}-q)\Vert_{\psi_2}<\infty,\quad\Vert\sqrt k( m-\mu)|m>0\Vert_{\psi_2}<\infty,\quad\Vert\sqrt k(m+\mu)|m<0\Vert_{\psi_2}<\infty,\\
        &\Vert\sqrt{n-k}\tde m\Vert_{\psi_2}<\infty.
    \end{align*}
\end{lemma}

To prove the above lemma, we present the following two results (lemma \ref{secondmoment} and lemma \ref{secondmomentpositive}), which give exact forms of the moment iterations using the cavity methods that we discuss in Appendix B. These two moment iterations lead to lemma \ref{largcliquetailbound}. Furthermore, together with lemma \ref{largcliquetailbound}, it further leads to our final result of lemma \ref{cltatline}.
    
\begin{lemma}[Moment Iteration for High Temperature]\label{secondmoment}
    For some $r\in\bb N$, we denote $|a-b|=\tde O(r)$ if
    \begin{align*}
        |a-b|=O\bl\frac{1}{n^{r/2}}\vee \nu(|R_{1,2}-q|^r)\vee\min_{\mu\in\ca U}\nu(|m_1-\mu|^r)\br.
    \end{align*}
    Then we consider the region where $\theta_1<\frac{1}{\bb E[\sech^2(\theta\sqrt{q}z+h)]}$.
       Define 
       \tny{\begin{gather*}
        a(2):=a(1,2)=\theta^2(1- \wh q^2),\quad a(1):=a(1,3) =\theta^2(\wh q-\wh q^2),\quad a(0):=a(3,4)=\theta^2(\wh q_4-\wh q^2),\nnb\\
        b(1):=b(1,2) = \theta^2\mu(1-\wh q),\quad b(0):=b(2,3)=\theta^2(\wh q_3-\wh q\mu),\nnb\\
        \quad d(1) =\theta_1(1-\mu^2),\quad d(0)=\theta_1(\wh q-\mu^2), \quad e(1)=\theta_1\mu(1-\wh q),\quad e(0)=\theta_1(\wh q_3 -\wh q\mu),\nnb\\
       \tde a(2):=\tde a(1,2)=\theta^2(1- \tde q^2),\quad \tde a(1):=\tde a(1,3) =\theta^2(\tde q-\tde q^2),\quad \tde a(0):=\tde a(3,4)=\theta^2(\tde q_4-\tde q^2).
    \end{gather*}}
    \tny{\begin{align*}
          \bfa A_1:&= \begin{bmatrix}
        a(2)&-4a(1)&3a(0)&2e(1)&-2e(0)\\
        a(1)&a(2)-2a(1)-3a(0)&6a(0)-3a(1)&e(1)+e(0)&e(1)-3e(0)\\
        a(0)&4a(1)-8a(0)&a(2)-8a(1)+10a(0)&2e(0)&2e(1)-4e(0)\\
        b(1)&-2b(1)-2b(0)&3b(0)&d(1)+d(0)&-2d(0)\\
        b(0)&2b(1)-6b(0)&6b(0)-3b(1)&2d(0)&d(1)-3d(0)
        \end{bmatrix}=\begin{bmatrix}
            \bfa A_{11}\in\bb R^{3\times 3}&\bfa A_{12}\\
            \bfa A_{13}&\bfa A_{14}
        \end{bmatrix},\\
        \bfa A_2:&=\begin{bmatrix}
        a(2)-2a(1)&3a(0)-2a(1)&e(1)&e(1)-2e(0)\\
        2a(1)-2a(0)&a(2)-6a(1)+6a(0)&e(0)&2e(1)-3e(0)\\
        -b(1)&b(0)&d(1)&-d(0)\\
        b(1)-2b(0)&b(0)&d(0)&d(1)-2d(0)
        \end{bmatrix}=\begin{bmatrix}
            \bfa A_{21}\in\bb R^{2\times 2}&\bfa A_{22}\\
            \bfa A_{23}&\bfa A_{24}
        \end{bmatrix},\\
        \tde{\bfa A}_1:&= \begin{bmatrix}
        a(2)&-4a(1)&3a(0)&0&0\\
        a(1)&a(2)-2a(1)-3a(0)&6a(0)-3a(1)&0&0\\
        a(0)&4a(1)-8a(0)&a(2)-8a(1)+10a(0)&0&0\\
       0&0&0&0&0\\
       0&0&0&0&0
        \end{bmatrix}=\begin{bmatrix}
            \tda A_{21}\in\bb R^{3\times 3}&\bfa 0\\
            \bfa 0&\bfa 0
        \end{bmatrix},
       \end{align*}}
    \tny{\begin{align*}
        \bfa x_{r}:&=(U_{2,r},U_{1,r},U_{0,r},V_{0,1,r},V_{0,0,r})^\top,\quad\bfa y_{r}:=(V_{1,1,r}, V_{1,0,r},W_{1,r}, W_{0,r})^\top,\\
        \bfa b_1:&=\lef(1-\wh q^2,\wh q-\wh q^2,\wh q_4-\wh q^2,\mu-\mu \wh q,\wh q_3-\mu \wh q\rig)^\top,\quad\bfa b_2:=(\mu-\wh q\mu, \wh q_3-\wh q\mu,1-\mu^2,\wh q-\mu^2)^\top,\\
        \tda b_1:&=\lef(1-\tde q^2,\tde q-\tde q^2,\tde q_4-\tde q^2,0,\tde q_3\rig)^\top,
    \end{align*}}
    \tny{\begin{gather*}
        U_{2,r}:=\nu((R_{1,2}-q)^{2r}),\quad U_{1,r}:=\nu((R_{1,2}-q)^{2r-1}(R_{1,3}-q)),\quad U_{0,r}:=\nu((R_{1,2}-q)^{2r-1}(R_{3,4}-q)),\\
        V_{0,1,r}:=\nu((R_{1,2}-q)^{2r-1}(m_1-\mu)),\quad V_{0,0,r}:=\nu((R_{1,2}-q)^{2r-1}(m_3-\mu)),\\
        V_{1,1,r}:=\nu((R_{1,2}-q)(m_1-\mu)^{2r-1}),\quad V_{1,0,r}:=\nu((R_{1,2}-q)(m_3-\mu)^{2r-1}),\\
        \quad W_{1,r}:=\nu((m_1-\mu)^{2r}),\quad W_{0,r}:=\nu((m_1-\mu)^{2r-1}(m_2-\mu)),\quad \tde W_{1,r}:=\nu(\tde m_1^{2r}),\\ 
        \quad\wh q_j:=\nu_{0,1}(\epsilon_1\cdots\epsilon_j)=\bb E[\tanh^j(\theta\sqrt q z+\theta_1\mu+h)],\\ \tde q_j:=\nu_{0,2}(\epsilon_1\ldots\epsilon_j)=\bb E[\tanh^j(\theta\sqrt qz+h)],\quad j\geq 3.
    \end{gather*}}
    Then we have
    \tny{\begin{align*}
        \bfa x_r &=(\bfa A_1+\tde{\bfa A}_1)\bfa x_r+ \frac{2r-1}{n}U_{2,r-1}(c\bfa b_1+(1-c)\tda b_1)+\tde O(2r+1),\nnb\\
        \bfa y_r &=\bfa A_2\bfa y_r+\frac{1}{k}(2r-1)W_{1,r-1}\bfa b_2+\tde O(2r+1),\qquad \tde W_{1,r}=\frac{2r-1}{n-k}\tde W_{1,r-1}+\tde O(2r+1).
    \end{align*}}
\end{lemma}
To prove the above lemma, we use the smart path constructed in lemma \ref{cavity}. Using this smart path, we construct two sets of equations that finally recover higher moments using lower moments.  The next lemma deals with the low-temperature case. Our method is the smart path given in lemma \ref{cavityII}.

\begin{lemma}[Moment Iteration at Low Temperature]\label{secondmomentpositive}
    Let $r\in\bb N$. Using similar notation as in Lemma \ref{secondmoment}, assume that $\theta_1>\frac{1}{\bb E[\sech^2(\sqrt qz+h)]}$. For $f(\bfa\sigma^1,\ldots,\bfa\sigma^r):\Sigma_{k}^r\to\bb R$, we define $\nu_t^+(f)=\nu_t(f|m_1^->0,\dots,m_r^->0)$ to be the path given by the smart path method II in Lemma \ref{zeropointdecompII}.  Define
$        \bfa y^+_{r}:=(V^+_{1,1,r}, V^+_{1,0,r},W^+_{1,r}, W^+_{0,r})^\top$ with
    \tny{\begin{gather*}
        V^+_{1,1,r}:=\nu^+((R_{1,2}-q)(m_1-\mu)^{2r-1}),\quad V^+_{1,0,r}:=\nu^+((R_{1,2}-q)(m_3-\mu)^{2r-1}),\\
        \quad W^+_{1,r}:=\nu^+((m_1-\mu)^{2r}),\quad W^+_{0,r}:=\nu^+((m_1-\mu)^{2r-1}(m_2-\mu)),\\ 
        \quad q_j:=\nu^+_0(\epsilon_1\cdots\epsilon_j)=\bb E[\tanh^j(\theta\sqrt q z+\theta_1\mu+h)],\quad j\geq 2,\quad q=q_2.
    \end{gather*}}
    Then we have
     \tny{\begin{align*}
        \bfa x_r &=(\bfa A_1+\tde{\bfa A}_1)\bfa x_r+ \frac{2r-1}{n}U_{2,r-1}(c\bfa b_1+(1-c)\tda b_1)+\tde O(2r+1),\nnb\\
        \bfa y^+_r &=\bfa A_2\bfa y^+_r+\frac{1}{k}(2r-1)W_{1,r-1}\bfa b_2+\tde O(2r+1),\qquad \tde W_{1,r}=\frac{2r-1}{n-k}\tde W_{1,r-1}+\tde O(2r+1).
    \end{align*}}
\end{lemma}

Finally, we give the limiting theorem for the large clique regime at the high and low-temperature regimes. These results finally lead to the rejection region of test statistics provided by algorithm \ref{alg:seven}.
\begin{lemma}\label{cltatline}
   Use the notations in Lemma \ref{secondmoment}. We add superscript $h$ and $l$ to denote the matrix $\bfa A$s and vector $\bfa b$s in the high and low-temperature regimes. Assume that the condition \ref{withintheATline} is satisfied. When $\theta_1<\frac{1}{\bb E[\sech^2(\theta z\sqrt q+h)]}$, for all $r\in\bb N$ we have
   \tny{
    \begin{align*}
         \bb E\bigg[\lef(\frac{n}{\ca V_r^h}\rig)^{-r/2}( R_{1,2}-q)^r\bigg],\quad\bb E\bigg[\bl\frac{k}{\ca V_m^h}\br^{-r/2}m^r\bigg],\quad\bb E[(n-k)^{-r/2}\tde m^r]\to\bb E[z^r],
    \end{align*}}
     with \tny{$\ca V_r^h=(\bfa I-c\bfa A^h_1-(1-c)\tda A^h_1)^{-1}_1(c\bfa b^h_1+(1-c)\tda b^h_1)$, $\ca V_{m}^h=(\bfa A^h_2)^{-1}_1\bfa b^h_2$}.
     
    When \tny{$\theta_1>\frac{1}{\bb E[\sech^2(\theta z\sqrt q+h)]}$},
    \tny{$\det(\bfa I-c\bfa A^l_{14}-c^2\bfa A^l_{13}\bfa B^{-1}\bfa A^l_{12})\neq 0$} for \tny{$\bfa B:=\bfa I-c\bfa A^l_{21}-(1-c)\tda A^l_{21}$}, and \tny{$\det(\bfa I-\bfa A^l_{24}-\bfa A^l_{23}(\bfa I-\bfa A^l_{21})^{-1}\bfa A^l_{22})\neq 0$},
    we have for all $r\in\bb N$,
    \tny{\begin{gather*}
         \bb E\bigg[\lef(\frac{n}{\ca V_r^l}\rig)^{-r/2}( R_{1,2}-q)^r\bigg],\quad\bb E\bigg[\bl\frac{k}{\ca V_m^l}\br^{r/2}(m-\mu)^r\bigg |m>0\bigg]\to\bb E[z^r],\\
         \bb E\bigg[\bl\frac{k}{\ca V_m^l}\br^{r/2}(m+\mu)^r\br\bigg |m<0\bigg],\quad\bb E[(n-k)^{r/2}\tde m^r]\to\bb E[z^r],
    \end{gather*}}
  with \tny{$\ca V_r^l=(\bfa I-c\bfa A^l_1-(1-c)\tda A^l_1)^{-1}_1(c\bfa b^l_1+(1-c)\tda b^l_1)$, $\ca V_{m}^l=(\bfa A^l_2)^{-1}_1\bfa b^l_2$}.
\end{lemma}

\subsection{The Proof Outline of Theorem \ref{exact recoveryguarantees}}\label{sect5}
This section provides a proof outline of the exact recovery guarantees in theorem \ref{exact recoveryguarantees}. Similar to \citep{he2023hidden1}, we first derive the almost exact recovery guarantee of the problem using the scan-based algorithms in testing. Then we provide a screening procedure to finally get the exact recovery. The following results give almost exact recovery guarantees for the scan-based algorithms.
\begin{lemma}[Almost Exact Recovery]\label{reco}
   Assume that the clique is positioned with index set $S$ and that the replica symmetry condition \ref{withintheATline} holds.
   \begin{enumerate}
       \item For all $\delta>0$, sample size $m\geq Ck\log(n)$, and $\theta_1<\frac{1}{\bb E[\sech^2(h^\prime+\theta\sqrt qz)]}$, for $\wh S\in\argmax_{S:|S|=k}\phi_S$ returned by algorithm \ref{alg:two}, we have $ \bb P\lef(|S\Delta \wh S|\geq k\delta\rig)=o(1)$.
        \item For all $\delta>0$, sample size $m\geq C\log(n)+1$, and $\theta_1>\frac{1}{\bb E[\sech^2(h^\prime+\theta\sqrt qz)]}$, for $\wh S\in\argmax_{S:|S|=k}\phi_S$ returned by algorithm \ref{alg:three}, we  have $ \bb P\lef(|S\Delta \wh S|\geq k\delta\rig)=o(1)$.
        \item For all $\delta>0$, assume that $k\lesssim\frac{n^{\frac{2\tau-1}{2\tau}}}{\log^{\frac{2\tau-2}{2\tau-1}}n}$, sample size $m\geq Ck^{\frac{1}{2\tau-1}}\log(n)$, the flatness parameter of $0$ being $\tau$, and $\theta_1=\frac{1}{\bb E[\sech^2(h^\prime+\theta\sqrt qz)]}$, for $\wh S\in\argmax_{S:|S|=k}\phi_S$ returned by algorithm \ref{alg:five} we have $ \bb P\lef(|S\Delta \wh S|\geq k\delta\rig)=o(1)$.
   \end{enumerate}
\end{lemma}
Given the above almost exact recovery guarantee, we propose the following algorithm as the screening procedure that boosts to the exact recovery guarantee. The underlying idea is to use the almost exact clique set as an anchor to the rest of the spins within the clique.
\begin{algorithm}[htbp]
\caption{Set Screening}\label{alg:screen}
\KwData{$\{\bfa\sigma^{(i)}\}_{i\in[m]}$ such that $\bfa\sigma^{(i)}\in\{-1,1\}^n$, an almost exact solution $S^\prime$ returned by algorithm \ref{alg:two}, \ref{alg:three}, or algorithm \ref{alg:five} in the high/low/critical temperature regimes respectively.}
Compute the statistics $\phi_i=\begin{cases}
     m^{-1}\sum_{\ell=1}^m\sum_{j\in S^\prime,j\neq i}\sigma^{(\ell)}_i\sigma^{(\ell)}_j&\text{ at the high temperature regime}\\
     k^{-\frac{2\tau-2}{2\tau-1}}m^{-1}\sum_{\ell=1}^m\sum_{j\in S^\prime,j\neq i}\sigma^{(\ell)}_i\sigma^{(\ell)}_j &\text{ at the critical temperature regime with flatness } \tau\\
     k^{-1}m^{-1}\sum_{\ell=1}^m\lef|\sum_{j\in S^\prime,j\neq i}\sigma^{(\ell)}_i\sigma^{(\ell)}_j\rig|&\text{ at the low temperature }
\end{cases}$\;
Then we rank $\phi_i$ and pick $S^{\prime\prime}$ by the set achieving top $k$ values of $\phi_i$.
\end{algorithm}

\section{Universality}\label{sect6}

In this section, we present the universality result for the statistical problem underlying the planted SK model. Consider a random variable $\bfa\eta$ such that $\eta_{ij}$s are i.i.d. with $\bb E[\eta^{k}_{ij}]=\bb E[g_{ij}^{k}]$ for $k\in[1:\tau]$ where $g_{ij} \sim N(0, 1)$ and $\bb E[\eta_{ij}^{\tau+1}]<\infty$, we define the universal planted SK Hamiltonian as
\begin{align}\label{upskhamilt}
   \mca H^{UpSK}_{\theta_1,\theta}(\bfa \sigma) := -\frac{\theta}{\sqrt n} \sum_{1\leq i< j\leq n} \eta_{ij}\sigma_i\sigma_j-\sum_{1\leq i<j\leq k}\frac{\theta_1}{k}\sigma_i\sigma_j -\sum_{i\leq n}h_i\sigma_i,\quad\theta,\theta_1>0.
\end{align} 
 Our main result is stated as the following theorem. To prove this theorem, our idea is to use the non-Gaussian integration by parts lemma to construct a few interpolation paths for the small and large clique regimes separately. The formal proof of this result is delayed to the supplementary material.
\begin{theorem}\label{thm6.1}
    Given that $\bb E[\eta_{ij}^k]=\bb E[g_{ij}^k]$ for $k\leq 3$ and $\bb E[\eta_{ij}^4]<\infty$, all the rates given by theorems \ref{thm10}, \ref{thm2},   and \ref{exact recoveryguarantees} hold for the universal planted SK model.
\end{theorem}

\section{Discussions}\label{sect7}
In this section, we discuss potential future works and limitations.
\smallbreak\emph{Testing.}
Out of the known statistical phases of minimax testing, the gap between the upper and lower bounds for the high-temperature large clique regime is not closed when $\theta$ is of constant order. The fundamental reason for this gap comes from that approximation error term in \eqref{interpolation derivative} is of order $\frac{\theta^2k}{n}$. And if we manage to close it to $O(\frac{\theta^2k^2}{n^2})$, the bound will be closed. However, this seems to be hopeless using current methods.  Another question for testing is the critical temperature under the large clique regime. To solve this problem, we need to develop better strategies for the cavity iterations. This problem is purely mathematical rather than statistical, which is also closely related to a long-standing open problem raised by Talagrand, asking for the convergence rate of replicas on the AT curve. We also do not touch upon the replica symmetric breaking phase since the existing results rely heavily on the replica concentration rates, a result that is proven to be false when the replica symmetry is broken.

\smallbreak\emph{Recovery.}
Seeing from the table \ref{table2}, a few gaps remain open between the minimax upper and lower bounds for recovery, given a constant order of $\theta$. For the small clique regime, instead of asking for the approximation error to be of order $O(\frac{\theta^2k^2}{n^2})$, to recover the high region with $k\log k=o(n)$ in the testing, we will need to ask for the approximation error to be of order $O(\frac{\theta^2}{k})$ to achieve the correct rate.  Similar results under the large clique regime where we cannot make any approximation is an even more difficult question.
\smallbreak\emph{Computational Constrained Inference.} All of the procedures provided in this work are not under the polynomial computational constraint. It is discussed in \citep{he2023hidden1} that the upper bounds of the same problem (using the semi-definite programs) under the computational constraints differ from the information-theoretic lower bounds for the pRFCW model. (which is the special case of $\theta=0$ in this work) A close observation of their proof suggests that the same upper bound holds for pSK model in the small clique regime (under some minor assumptions on the scaling of $\theta$). However, showing rigorous computational lower bounds might still be challenging, which is listed as following works.

\begin{appendix}
\section{The smart paths in section \ref{sect4.3}}\label{appA}
Here we present the two smart paths utilized in the proof of the moment iteration for the cavity method. The cavity method is particularly effective when analyzing the $\frac{1}{\sqrt k}$ rate since this is the correct rate when the moment iterations can be formulated. The first lemma gives a path that decouples all the correlation between the last spins from the rest of them. 
\begin{lemma}[Smart Path I]\label{cavity}
   For a function $f:\Sigma_{n}^r\to\bb R$ to be a function taking $r$ replicas as input, we define $\nu_t(f):=\bb E\lef[\frac{\sum_{\bfa\sigma} f(\bfa\sigma)\exp\lef(-H_t(\bfa\sigma)\rig)}{\sum_{\bfa\sigma}\exp\lef(-H_t(\bfa\sigma)\rig)}\rig ]$. Define $\epsilon_i:=\sigma^i_k$ and $\xi_i:=\sigma_n^i$. $m_j:=\frac{1}{k}\sum_{i=1}^k\sigma_i^j$. Consider two different smart paths for arbitrary $q$ and $\mu$.
   \begin{enumerate}
       \item For the cavity within the clique, we define
       \tny{
       \begin{align}\label{path11}
        -\mca H_{t,1}(\bfa\sigma):&= \frac{\theta}{\sqrt n}\sum_{i<j\leq n, i,j\neq k}g_{ij}\sigma_i\sigma_j+\frac{\theta\sqrt t}{\sqrt n}\sum_{i\leq n}g_{ik}\sigma_i\sigma_k+\theta\sqrt{1-t}z\sqrt q\sigma_k+\frac{\theta_1}{2 k}\sum_{i,j\leq k-1}\sigma_i\sigma_j\nnb\\
        &+\frac{\theta_1 t}{k}\sum_{i\leq k}\sigma_i\sigma_k+\theta_1(1-t)\mu\sigma_k+\sum_{i\leq n}h_i\sigma_i.
    \end{align}}
    Then, defining the expectation under $\mca H_{t,1}$ to be $\nu_{t,1}(f):=\bb E[\la f\ra_{t,1}]$,  under this smart path we have
    \tny{
    \begin{align}\label{path22}
        \nu_{t,1}^\prime(f)&=\theta^2\bl\sum_{1\leq \ell<\ell^\prime\leq r}\nu_{t,1}(f\epsilon_\ell\epsilon_{\ell^\prime}(R_{\ell,\ell^\prime}-q))\br -r\theta^2\sum_{\ell\leq r}\nu_{t,1}(f\epsilon_\ell\epsilon_{r+1}(R_{\ell,r+1}-q))\nnb\\
        &+\theta^2\frac{r(r+1)}{2}\nu_{t,1}(f\epsilon_{r+1}\epsilon_{r+2}(R_{r+1, r+2}-q))\nnb\\
        &+\theta_1\bigg(\sum_{\ell\leq r}\nu_{t,1}(f\epsilon_{\ell}(m_{\ell}-\mu))-r\nu_{t,1}(f\epsilon_{r+1}(m_{r+1}-\mu))\bigg).
    \end{align}}
    \item For the cavity without the clique, we define
    \tny{\begin{align*}
        -\mca H_{t,2}(\bfa\sigma):&= \frac{\theta}{\sqrt n}\sum_{i<j\leq n-1}g_{ij}\sigma_i\sigma_j+\frac{\theta\sqrt t}{\sqrt n}\sum_{i\leq n}g_{in}\sigma_i\sigma_n+\theta\sqrt{1-t}z\sqrt q\sigma_n+\frac{\theta_1}{2 k}\sum_{i,j\leq k}\sigma_i\sigma_j+\sum_{i\leq n}h_i\sigma_i.
    \end{align*}}
    Then, defining the expectation under $\mca H_{t,2}$ to be $\nu_{t,2}(f):=\bb E[\la f\ra_{t,2}]$, under this smart path we have
    \tny{
    \begin{align}\label{path33}
        \nu_{t,2}^\prime(f)&=\theta^2\bl\sum_{1\leq \ell<\ell^\prime\leq r}\nu_t(f\xi_\ell\xi_{\ell^\prime}(R_{\ell,\ell^\prime}-q))\br -r\theta^2\sum_{\ell\leq r}\nu_t(f\xi_\ell\xi_{r+1}(R_{\ell,r+1}-q))\nnb\\
        &+\theta^2\frac{r(r+1)}{2}\nu_t(f\xi_{r+1}\xi_{r+2}(R_{r+1, r+2}-q)).
    \end{align}}
   \end{enumerate}
    
    \end{lemma}
    The following lemma is the second smart path for studying the low temperature regime. This method separates the case between $m>0$ and $m<0$ and simultaneously decouples the correlation between the last spin and the rest. 
\begin{lemma}[Smart Path II] \label{cavityII} For a function $f:\Sigma_{n}^r\to\bb R$, we define $m^-(\bfa\sigma):=\frac{1}{k}\sum_{i=1}^{k-1}\sigma_i$, $m_i^-=m^-(\bfa\sigma^i)$, and assume that $\mu>0$. We pick
\tny{
\begin{align*}
        -\mca H_{t,3}(\bfa\sigma):&= \frac{\theta}{\sqrt n}\sum_{i<j\leq n,i,j\neq k}g_{ij}\sigma_i\sigma_j+\frac{\theta\sqrt t}{\sqrt n}\sum_{i\leq n}g_{ik}\sigma_i\sigma_k+\theta\sqrt{1-t}z\sqrt q\sigma_k+\frac{\theta_1}{2 k}\sum_{i,j\leq k-1}\sigma_i\sigma_j\\
        &+\frac{\theta_1 t}{k}\sum_{i\leq k-1}\sigma_i\sigma_k+\theta_1(1-t)\sigma_k(\mu\mbbm 1_{m^-(\bfa\sigma)>0}-\mu\mbbm 1_{m^-(\bfa\sigma)<0}) +\sum_{i\leq n}h_i\sigma_i.
    \end{align*}}
    Then, defining the expectation under $\mca H_{t,3}$ to be $\nu_{t,3}$,  following the definition in Lemma \ref{cavity}, we have
    \tny{
    \begin{align*}
        &\nu_{t,3}^\prime(f)=\theta^2\bl\sum_{1\leq \ell<\ell^\prime\leq r}\nu_{t,3}(f\epsilon_\ell\epsilon_{\ell^\prime}(R_{\ell,\ell^\prime}-q))\br -r\theta^2\sum_{\ell\leq r}\nu_{t,3}(f\epsilon_\ell\epsilon_{r+1}(R_{\ell,r+1}-q))\\
        &+\theta^2\frac{r(r+1)}{2}\nu_{t,3}(f\epsilon_{r+1}\epsilon_{r+2}(R_{r+1, r+2}-q))+\theta_1\bigg(\sum_{\ell\leq r}\nu_{t,3}(f\epsilon_{\ell}(m_{\ell}-\mu\mbbm 1_{m^-_{\ell}>0}+\mu\mbbm 1_{m^-_{\ell}<0}))\\
        &-r\nu_{t,3}(f\epsilon_{r+1}(m^-_{r+1}-\mu\mbbm 1_{m^-_{r+1}>0}+\mu\mbbm 1_{m^-_{r+1}<0}))\bigg).
    \end{align*}}
\end{lemma}

The following two lemmas imply that the correlation between the last spin and the rest can be decoupled at the $t=0$ end in the smart path Lemma \ref{cavity} and Lemma \ref{cavityII}. 
\begin{lemma}\label{zeropointdecomp}
Assume that we use the smart path defined by Lemma \ref{cavity}.
Introduce the notation $Y:=\theta\sqrt q z+\theta_1\mu+h$.
    For any function $f^-_r$ on $\Sigma_{n}^r$ and any set $I\subset\{1,\ldots,r\}$ we have
    \begin{align*}
        \nu_{0,1}\bl f^-_r\prod_{i\in I}\epsilon_i\br=\nu_{0,1}\bl \prod_{i\in I}\epsilon_i\br\nu_{0,1}(f^-_r) =\bb E[\tanh^{|I|} Y]\nu_{0,1}(f^-_r).
    \end{align*}
\end{lemma}
\begin{lemma}\label{zeropointdecompII}
    Assume we use the smart path defined by Lemma \ref{cavity}.
Introduce the notation $Y_1:=\theta\sqrt q z+\theta_1\mu+h$, $Y_2:=\theta\sqrt qz-\theta_1\mu+h$. Assume that $A_0:=[-1,0), A_1:=(0,1]$.
    For any function $f^-_r$ on $\Sigma_{n}^r$ and any set $I\subset\{1,\ldots,r\}$ , denote $A_{ij}=\{m_i^-\in A_j\}$, we have
    \begin{align*}
        \nu_{0,3}\bl f^-_r\prod_{i\in I}\epsilon_i\bigg | \cap_{i\in[I]}A_{ij}\br&=\nu_{0,3}\bl \prod_{i\in I}\epsilon_i\bigg |\cap_{i\in[I]}A_{ij}\br\nu_{0,3}(f^-_r|\cap_{i\in[I]}A_{ij})\\
        &=\bb E\bigg[\prod_{i\in[I]}\tanh^{j} Y_1\tanh^{1-j}Y_2\bigg]\nu_{0,3}(f^-_r|\cap_{i\in[I]}A_{ij}).
    \end{align*}
\end{lemma}

\end{appendix}




\begin{supplement}
\stitle{Supplement to “Hidden Clique Inference in Random Ising Model II: the planted Sherrington-Kirkpatrick model” }
\sdescription{In this supplementary material, we provide the complete proofs of the results of this work.}
\end{supplement}


\bibliographystyle{imsart-number} 
\bibliography{bib}       





\newpage

\setcounter{page}{1}
\setcounter{section}{0}
\title{Supplement to "Hidden Clique Inference in Random Ising Model II:
the planted Sherrington-Kirkpatrick model"}
\bigskip


\startcontents[supplementary]
\printcontents[supplementary]{l}{1}{\setcounter{tocdepth}{2}}

\renewcommand{\thesection}{\Roman{section}}
\renewcommand{\thesubsection}{\thesection.\roman{subsection}}
\renewcommand{\thesubsubsection}{\thesubsection.\roman{subsubsection}}
\section{The Proof Outlines of Other Lemmas and Theorems \ref{sect6} }\label{appB}
\subsection{Proof Outlines of Lemma \ref{cltrfcr}}.
    We present a heuristic description of the smart path method. At the center of the proof is the Stein's lemma, which leads to
    \begin{lemma}[Gaussian Interpolation Lemma ]\label{gil}
    Define $\bfa u,\bfa v\in\bb R^m$ to be two Gaussian process indexed by $[M]$, and define
    \begin{align*}
        u_i(t) = \sqrt t u_i+\sqrt{1-t} v_i
    .\end{align*}
    and note that $\bfa u =\bfa u(1)$ and $\bfa v =\bfa u(0)$, define function $\varphi(t)=\bb E[F(\bfa u(t))]$ with $F:\bb R^m\to\bb R$ being second order differentiable with $F$ satisfying conditions of $\lim_{x\to\infty}F(\bfa x)\exp(-C\Vert x\Vert_2)=0$ for some $C>0$, we have
    \begin{align*}
        \varphi^\prime(t)=\frac{1}{2}\sum_{i,j}(\bb E[u_iu_j]-\bb E[v_iv_j])\bb E\lef[\frac{\pta^2 F}{\pta x_i\pta x_j}(\bfa u(t))\rig]
    .\end{align*}
\end{lemma}
It is not hard to see that the Gaussian component in the pSK Hamiltonian is a Gaussian process indexed by the spin state. Then the idea is to construct another interpolating Gaussian process $\mca H_{t}(\bfa\sigma)$ such that $\mca H_{0}(\bfa\sigma)$ is some \emph{ more analyzable process } and $\mca H_{1}(\bfa\sigma)=H^{the pSK}(\bfa\sigma)$. 
Defining the \emph{Gibbs Average} of a function $f:\Sigma_n\to\bb R$ by $\la f\ra_t:=\frac{\sum_{\bfa\sigma}f\exp(-\mca H_{t}(\bfa\sigma))}{\sum_{\bfa\sigma}\exp(-\mca H_{t}(\bfa\sigma))}$. 
And one can approximate the true value of $\bb E[\la F\ra_1]$ using $\bb E[\la F\ra_0]$ and a proper upper bound on $\frac{d\bb E[\la F\ra_t]}{dt}$. Then, our smart path is defined by
\tny{\begin{align}\label{smtph1}
    &-\mca H_{t}(\bfa\sigma) := \sqrt t\bl\sum_{i,j\in [k],i\neq j}\frac{\theta g_{ij}}{\sqrt n}\sigma_i\sigma_j+\sum_{i\in[k],j\in[k+1:n]}\frac{\theta g_{ij}}{\sqrt n}\sigma_i\sigma_j\br\nnb\\
    &+\sqrt{1-t}\bl \sum_{i\in[k]}\theta\sqrt{q}z_i\sigma_i+\sum_{i=k+1}^n\theta\sqrt{q^\prime}z_i\sigma_i\br+\sum_{i,j\in[k]}\frac{\theta_1}{2k}\sigma_i\sigma_j+\sum_{k+1\leq i<j\leq n}\frac{\theta g_{ij}}{\sqrt n}\sigma_i\sigma_j+\sum_{i=1}^n\sigma_ih_i,
\end{align}}
with some $q^\prime$ and $q$ to be specified in the formal proof. One will see that this path successfully decouples the correlation between the spins within and without the clique.
\subsection{Proof Outline of Theorem \ref{thm6.1}}
The rest of this section is devoted to proving the small (section \ref{sect6.1}) and large clique (section \ref{sect6.2}) regimes respectively. At the center of our proof is the non-Gaussian integration by parts, which is given by
\begin{lemma}[Non-Gaussian integration by parts]\label{lm1.11}
    Let $\eta$ be a real random variable such that $\bb E[\eta^4]<\infty$, $\bb E[\eta]=\bb E[\eta^3]=0$ and $\bb E[\eta^2]=1$. Let $F:\bb R\to\bb R$ be three times continuously differentiable. Then there exists $\xi_1,\xi_2\in(0\wedge \eta,0\vee\eta)$ (random variables depending on $\eta$) such that
    \begin{align*}
        \bb E[\eta F(\eta)]-\bb E[\eta^2]\bb E[F^\prime(\eta)]=\bb E\lef[\frac{\eta^4}{6}F^{(3)}(\xi_1)\rig]-\bb E[\eta^2]\bb E\bigg[\frac{\eta^2}{2}F^{(3)}(\xi_2)\bigg].
    \end{align*}
\end{lemma}
\subsubsection{Small Cliques}\label{sect6.1}
For the small cliques regime, we discuss two results that imply the validity of the upper bound and lower bound respectively. One can see them correspond to Theorem \ref{cltrfcr} and \ref{approximation bounds chisq} respectively.
\begin{lemma}[Universal local limiting theorems for the planted SK model]\label{univer}
    The results in theorem \ref{cltrfcr} hold for the universal planted SK model defined by \eqref{upskhamilt}.
\end{lemma}

Then we turn to the second results regarding the lower bound analysis given by theorem \ref{approximation bounds chisq}. This is guaranteed by the following lemma.
\begin{lemma}\label{lm5.3}
    Define $\bb P^U_S(\cdot|\bfa h,\bfa \eta)$ to be the universal planted SK measure with clique planted at $S$. Define $\bb P^U_0(\cdot|\bfa h,\bfa \eta)$ to be the universal planted SK measure without the clique. Assume that we are given arbitrary events $A_S\in\Sigma\lef(\bfa\sigma_S\rig),A_{S^\prime}\in\Sigma(\bfa\sigma_{S^\prime})$. Then we have the following for sufficiently large $k$
    \begin{align*}
        \bb E\bigg[\frac{\mbbm 1_{A_S}\mbbm 1_{A_{S^\prime}}\bb P^U_S(\cdot|\bfa h,\bfa\eta)\bb P^U_{S^\prime}(\bfa h,\bfa\eta)}{\bb P_0^U(\cdot|\bfa h,\bfa\eta)}\bigg]\leq\exp\bl\frac{C\theta^2k}{n}\br\bb E\bigg[\frac{\mbbm 1_{A_S}\mbbm 1_{A_{S^\prime}}\bb P^U_S(\cdot|\bfa h,\bfa g)\bb P^U_{S^\prime}(\bfa h,\bfa g)}{\bb P_0^U(\cdot|\bfa h,\bfa g)}\bigg].
    \end{align*}
\end{lemma}
\begin{remark}
    The results for the small clique utilize a smart path constructed by approximating the Gaussian pSK model with the universal SK model. The smart path is defined by
    \sm{\begin{align*}
        -\mca H_{t}(\bfa\sigma)&=\frac{\theta}{\sqrt n}\sum_{1\leq i<j\leq k+r\text{ or }i\in[k+r],j\in[k+r+1:n]}(\sqrt{1-t}\eta_{ij}+\sqrt t g_{ij})\sigma_i\sigma_j+\frac{\theta}{\sqrt n}\sum_{1\leq i<j\leq k}\eta_{ij}\sigma_i\sigma_j\\
        &+\sum_{1\leq i<j\leq k}\frac{\theta_1}{k}\sigma_i\sigma_j+\sum_{i=1}^nh_i\sigma_i.
    \end{align*}}
    Then we go back to the smart path given by \ref{smtph1}, it is not hard to see that there are two elements that need to be considered so that we can transform this local limiting theorem from the Gaussian case to the universal case:
    \begin{enumerate}
        \item The approximation error given by the above interpolation is bounded.
        \item The smart path \eqref{smtph1} continues to hold, which is implied by the concentration of the standard universal SK model.
    \end{enumerate}
    The first element is satisfied through lemma \ref{lm5.3} and the second element is satisfied through lemma \ref{largcliquetailbound}.
\end{remark}

\subsubsection{Large Cliques}\label{sect6.2}
This section presents the results of universality for large cliques. First, we derive that the limiting free energy of the universal pSK model resembles the Gaussian pSK model. Consequently, the phase transition line will remain the same. We notice that this result has already been proved in \citep{carmona2006universality} for the standard SK model, and it is not hard to extend the results to the planted SK model. Our results are stated as follows
\begin{lemma}\label{lm5.4}
    Assume that $\eta$ satisfies $\bb E[\eta]=0,\bb E[\eta^2]=1,\bb E[|\eta|^3]<\infty$. And define 
    \begin{align*}
        Z_n(\theta,\theta_1,\eta):=\sum_{\bfa\sigma}\exp(-H^{UpSK}_{\theta_1,\theta,\bfa h}(\bfa\sigma)),\quad Z_n(\theta,\theta_1,g):=\sum_{\bfa\sigma}\exp(-H_{\theta_1,\theta,\bfa h}^{pSK}(\bfa\sigma)).
    \end{align*}
    Then we have
$        \lef|\frac{1}{n}\bb E[\log Z_n(\theta,\theta_1,\eta)]-\frac{1}{n}\bb E[\log Z_n(\theta,\theta_1,g)]\rig|=O\lef(\frac{1}{\sqrt n}\rig).
$    This immediately implies that the replica symmetric regime defined by \ref{hightemp} for the planted SK model is also the replica symmetric regime for the universal planted SK model.
\end{lemma}
Our next result implies that the concentration arguments provided in Lemma \ref{largcliquetailbound} remain unchanged for the universal planted SK model.
\begin{lemma}\label{lm5.5}
    Theorem \ref{largcliquetailbound} remains valid for the universal planted SK model.
\end{lemma}
To give an exact computation of the limiting moments, we notice that the cavity method proposed in Lemma \ref{cavity}, Lemma \ref{zeropointdecomp} will be valid, which leads to the following result.
\begin{lemma}\label{lm5.6}
    Lemma \ref{cltatline} remains valid for the universal planted SK model.
\end{lemma}
\subsection{Proof Outline of Theorem \ref{exact recoveryguarantees}}
The rest of this section is devoted to the major proof of the above results. For the almost exact recovery, the proof idea is simple and intuitive, where we construct a band of sets that overlap with $S$ with small cardinality. And with overwhelming probability, the returned set of local algorithms will be an element in this band. Hence, we need to analyze the concentration of quantity in the form of $\sum_{i\in S^\prime}\sigma_i$ when $|S^\prime\cap S|=k-r$ and $|S^\prime|=k$, assuming that the clique is planted with index set $S$. This will be resolved at the small/large clique regimes respectively in sections \ref{sect5.1} and \ref{sect5.2}. For the exact recovery, the proof idea is to consider the sum correlation between a single spin (both within and without the clique) with the spins within the set $S$, which relies on the intuition that given sufficiently large sample size, one will be able to prove that all the spins within the clique will have stronger sum correlation than anyone out of the clique, given the sum is \emph{only} taken over the spins within the clique. As opposed to the proof of the same theorem that appears in \citep{he2023hidden1}, these two results need to overcome the non-independence between 
\subsubsection{Small Cliques}\label{sect5.1}
Our results of the small clique regime are based on the following lemma, which gives the CLT for spins of a set $S^\prime$ with $|S^\prime|=k$ and $S$ overlapped with $S^\prime$. Our proof strategy is similar to the idea in section \ref{sect3} where we connect the results in the pRFCW model with the pSK model. 
\begin{lemma}[Overlapped Limiting Theorems]\label{overlapclt}
    Consider $S=[k]$ to be the index set of the hidden clique. Denote $c:=\frac{r}{k}$.  Then we have
    \begin{itemize}
        \item At high temperature, when $k\log k=o(n)$, for all $t\in\bb R$:
        \begin{align*}
    \bb E\lef[\exp\lef(\frac{t\sum_{i\in S^\prime}\sigma_i}{\sqrt k}\rig)\rig]\to \exp\lef(\frac{1}{2}((1-c)\ca V+c)t^2\rig),\qquad\Vert k^{-1/2}\sum_{i\in S^\prime}\sigma_i\Vert_{\psi_2}<\infty,
    \end{align*}
with $\ca V:=\frac{1-\theta_1(\bb E[\sech^2(\theta\sqrt q z+h)])^2}{(1-\theta_1\bb E[\sech^2(\theta\sqrt q z+h)])^2}$.
\item At low temperature, when $k\log k=o(n)$, for all $t\in\bb R$:
\tny{\begin{align*}
    &\bb E\bigg[\exp\bl\frac{t\sum_{i\in S^\prime}\sigma_i}{\sqrt k}\br\bigg|\sum_{i=r+1}^k\sigma_i>0\bigg]\to\exp\bl\frac{1}{2}((1-c)\ca V(x^*_1)+c)t^2\br,\\
    &\bb E\bigg[\exp\bl\frac{t\sum_{i\in S^\prime}\sigma_i}{\sqrt k}\br\bigg|\sum_{i=r+1}^k\sigma_i<0\bigg]\to\exp\bl\frac{1}{2}((1-c)\ca V(x^*_1)+c)t^2\br,\\
    &\bigg\Vert k^{-1/2}\bigg|\sum_{i\in S^\prime}\sigma_i\bigg|-\bb E\bigg[\bigg|\sum_{i\in S^\prime}\sigma_i\bigg|\bigg]\bigg\Vert_{\psi_2}<\infty,
\end{align*}}
with $\ca V(x_1^*)$ defined in \eqref{gtr0clt}. 
\item At critical temperature with flatness $\tau$, when $k^{\frac{2\tau}{2\tau-1}}\log^{\frac{2\tau-2}{2\tau-1}}k=o(n)$, for all $t\in\bb R$:
\scriptsize{\begin{align*}
        &\bb E\bigg[\exp\bl t\frac{\sum_{i=r+1}^{k+r}\sigma_i}{k^{\frac{4\tau-3}{4\tau-2}}}\br\bigg]\to\int_{\bb R}\frac{(2\tau-1)x^{2\tau-2}}{\sqrt{2\pi v}}\exp\bl-\frac{x^{4\tau-2}}{2v}+t(1-c)^{\frac{4\tau-3}{4\tau-2}}\sqrt{\bb E[\sech^2(\theta\sqrt qz+h)]}x\br dx,\\
        &\bigg\Vert k^{-\frac{4\tau-3}{4\tau-2}}\sum_{i=r+1}^{k+r}\sigma_i\bigg\Vert_{\psi_{4\tau-2}}<\infty,
\end{align*}}
\normalsize
where $v$ is defined in \eqref{stirlingv}.
    \end{itemize}
\end{lemma}

\subsubsection{Large Cliques}\label{sect5.2}
Then we give results for the large clique regimes. Here our results are proved using an iterative process that computes the higher order overlap in the form of $\bb E\lef[\prod_{i\in[I]}\sigma_i\rig]$ for some $I\subset[n]$ with the lower order ones using the moment iteration. Then we recover the moments for the spins in the overlap set.
\begin{lemma}\label{largecliqueweakclt}
    When $k\lesssim n$, letting $S=[k]$ be the index set of the clique, we have 
    \begin{enumerate}
        \item At high temperature, $\bb E\lef[\frac{1}{k}\lef(\sum_{i=r+1}^{k+r}\sigma_i\rig)^2\rig]=(1-c)\ca V+c+o(1)$ with\\ $\ca V:=\frac{1-\theta_1(\bb E[\sech^2(\theta\sqrt q z+h)])^2}{(1-\theta_1\bb E[\sech^2(\theta\sqrt q z+h)])^2}$ and $\lef\Vert\frac{1}{\sqrt k}\sum_{i=r}^{k+r}\sigma_i\rig\Vert_{\psi_2}<\infty$.
        \item At low temperature, $\lef|\bb E\lef[\frac{1}{k}\lef|\sum_{i=r+1}^{k+r}\sigma_i\rig|\rig]-\bb E\lef[\lef|\frac{1}{k}\sum_{i=1}^k\sigma_i\rig|\rig]\rig|\asymp 1$ and $\bigg\Vert\frac{1}{\sqrt{k}}\bl\bigg|\sum_{i=r}^{k+r}\sigma_i\bigg|-\bb E\bigg[\bigg|\sum_{i=r}^{k+r}\sigma_i\bigg|\bigg] \br\bigg\Vert_{\psi_2}<\infty.$
    \end{enumerate}
\end{lemma}
\section{Proof of Main Results}
In this section of the appendix we present the proof of major theorems.
\subsection{Proof of Lemma \ref{thm1}}
The proof of lemma \ref{thm1} contains three parts: the guarantees obtained by the local algorithm (which corresponds to the $k=o(n^{2/3})$ part in algorithm \ref{alg:two}), the guarantees obtained by the global algorithm (which corresponds to the $k=\Omega(n^{2/3})$ part in algorithm \ref{alg:two}), and the lower bounds.
\subsubsection{Proof of the Local Algorithm}
   Noticing that when $k\log k=o(n)$ we can approximate the mgf of the local spins with the results given by the RFCW model through the RFCW model. This results is given by Lemma \ref{cltrfcr}.
   Then wefollow analogously the proof of lemma 3.1 in \citep{he2023hidden1}.
\subsubsection{Proof of the Global Algorithm}
    First we notice that the alternative is given by a canonical SK model. Using Lemma \ref{cltconsk} weget that $
        \Vert\sqrt n m\Vert_{\psi_2}<\infty
$ and $\bb E[nm^2]=1+o(1)$. For the alternative we use the result given in Lemma \ref{covariance} and lemma \ref{cltrfcr} to get
\begin{align*}
    \bb E[\phi_2]&=\bb E\bigg[\frac{1}{k}\bl\sum_{i=1}^n\sigma_i\br^2\bigg]-1=\frac{1}{k}\bb E\bigg[\bl\sum_{i=1}^k\sigma_i\br^2\bigg]+\frac{1}{k}\bb E\bigg[\bl\sum_{i=k+1}^n\sigma_i\br^2\bigg]\\
    &+\frac{1}{k}\bb E\bigg[\bl\sum_{i=1}^n\sigma_i\br\bl\sum_{i=k+1}^n\sigma_i\br\bigg]-1\\
    &=\frac{1}{k}\bb E\bigg[\bl\sum_{i=1}^k\sigma_i\br^2\bigg]+\frac{1}{k}\bb E\bigg[\bl\sum_{i=k+1}^n\sigma_i\br^2\bigg]-1+O\bl\frac{1}{\sqrt k}\br\\
    &=\frac{2\theta_1\bb E[\sech^2(\theta\sqrt q z+h)]-\theta_1(1+\theta_1)(\bb E[\sech^2(\theta\sqrt qz+h)])^2}{(1-\theta_1\bb E[\sech^2(\theta\sqrt qz+h)])^2}+o(1).
\end{align*}
And the rest of the proof follows from that of the Theorem 3.2 in \citep{he2023hidden1}.
\subsubsection{Proof of the Lower Bound }
    The proof goes by a similar vein as the proof of the Theorem 3.1 in \citep{he2023hidden1} where we analyzed the limiting chi-square divergence through lemma \ref{approximation bounds chisq}. Denote $\bb P_0,\bb E_0$ to be the measure/expectation under the null hypothesis. Denote $\bb P_{S}$ to be the alternative with a clique planted at $S\subset[n]$. Denote $\bar{\bb P}_m:=\frac{1}{\binom{n}{k}}\sum_{S:|S|=k}\bb P_S$ to be the mixture measure given uniform prior on $\ca S:=\{S:|S|=k,S\subset[n]\}$. We add subscript $m$ to denote the $m$-order product measure. A meta decomposition of Chi-square is given by
    \begin{align}\label{chisqlb}
        \bb E\lef[D_{\chi^2}\lef(\bar{\bb P}_m,\bb P_{0,m}|\bfa h^m\rig)\rig]:&=\frac{1}{\binom{n}{k}^2}\sum_{S,S^\prime\in\ca S}\bb E_{0,m}\lef[\frac{\bb P_{S,m}(\bfa\sigma|\bfa g^m,\bfa h^m)\bb P_{S^\prime,m}(\bfa\sigma|\bfa g^m,\bfa h^m)}{\bb P_{0,m}^2(\bfa\sigma|\bfa g^m,\bfa h^m)}\rig]-1\nnb\\
        &=\frac{1}{\binom{n}{k}^2}\sum_{S,S^\prime\in\ca S}\bb E_{0}\lef[\frac{\bb P_{S}(\bfa\sigma|\bfa g,\bfa h)\bb P_{S^\prime}(\bfa\sigma|\bfa g,\bfa h)}{\bb P_{0}^2(\bfa\sigma|\bfa g,\bfa h)}\rig]^m-1
    .\end{align}
    Using the notation given in the statement of lemma \ref{approximation bounds chisq} we have when $\theta_1<\frac{1}{2\bb E[\sech^2(\theta\sqrt qz+h)]}$
    \begin{align*}
        \bigg|\bb E\bigg[\frac{\bb P_{S}(\cdot|\bfa g,\bfa h)\bb P_{S^\prime}(\cdot|\bfa g,\bfa h)}{\bb P_0(\cdot|\bfa g,\bfa h) }\bigg]-\bb E\bigg[\frac{\bb Q_{S}(\cdot|\bfa h^\prime)\bb Q_{S^\prime}(\cdot|\bfa h^\prime)}{\bb Q_0(\cdot|\bfa h^\prime)}\bigg]\bigg|=O\bl\frac{k}{n}\br.
    \end{align*}
    Then we have, for $S,S^\prime\subset[n]$, $|S|=|S^\prime|=k$, and $S\cup S^\prime=k+r$ we have  (Using the result in the proof of the Theorem 3.1 in \citep{he2023hidden1}):
    \begin{align*}
        \bb E\bigg[\frac{\bb P_{S}\bb P_{S^\prime}}{\bb P_0}\bigg]= O\bl\frac{k}{n}+\exp\bl\frac{(k-r)^2}{k^2}\br\br.
    \end{align*}
    And for the region where $\theta_1\in\lef[\frac{1}{2\bb E[\sech^2(\theta\sqrt qz+h)]},\frac{1}{\bb E[\sech^2(\theta\sqrt qz+h)]}\rig)$ the above procedure will encounter a unbounded chi-square. Therefore instead of starting from the second inequality of lemma \ref{condlecam} westart from the first inequality and bound it by the truncated TV distance. This method is referred to as the \emph{fake measure method} in \citep{he2023hidden1}. Define $m_{S}:=\frac{1}{k}\sum_{i\in[S]}\sigma_i$ we construct the following event and the \emph{fake measure}.
    \begin{align}\label{esgoodset}
    E_S:=\lef\{|m_S|\leq c_0\sqrt{\frac{1}{k}\log( (m\vee k)\epsilon)}\rig\},\qquad \bb P^\prime_S(\bfa\sigma)=\begin{cases}
        \bb P_S(\bfa\sigma)&\text{ if }\bfa\sigma\in E_S\\
        0&\text{ otherwise }
    \end{cases}.
\end{align}
And wecheck that by the upper bound on the mgf given in Lemma \ref{cltrfcr} for the high temperature, we have
\begin{align*}
    \Vert\bb P_{\bar S}-\bb P^\prime_{\bar S}\Vert_{TV}=\bb P(E_S^c)=O\bl\frac{1}{(m\vee k)\epsilon}\br.
\end{align*}
Then we have
    \begin{align*}
        \bb E\bigg[\frac{\bb P^\prime_{S}\bb P^\prime_{S^\prime}}{\bb P_0}\bigg]=\bb E\bigg[\frac{\mbbm 1_{E_S}\mbbm 1_{E_{S^\prime}}\bb P_S\bb P_{S^\prime}}{\bb P_0}\bigg]\leq\exp\bl C\bl\frac{k-r}{k}\br^2\log(m\vee k)\log k+C\frac{k\log k}{n}\br.
    \end{align*}
    And the rest of the proof goes through analogously.
\subsection{Proof of Lemma \ref{thmlowtempfindclique}}
Similar to the proof of lemma \ref{thm1}, the proof of lemma \ref{thmlowtempfindclique} also contains three parts: (1) The proof of the local part of algorithm \ref{alg:three} under the $k=o(\sqrt n)$ regime; (2) The proof of the global part of algorithm under the $k=\Omega(\sqrt n)$ regime; (3) The lower bounds.
\subsubsection{Proof of the Local Algorithm}
    The proof goes by noticing that Lemma \ref{cltconsk} implies that under the null,
    \begin{align*}
        \bigg\Vert\frac{1}{\sqrt k}\sum_{i\in S}\sigma_i\bigg\Vert_{\psi_2}<\infty,\qquad \bigg\Vert\frac{1}{\sqrt{m k}}\sum_{i=j}^m\sum_{i\in S}\sigma_i^{(j)}\bigg\Vert_{\psi_2}<\infty.
    \end{align*}
    And, under the alternative, by Lemma \ref{cltrfcr}, we have
    \begin{align*}
        \bb E\bigg[\bigg|\frac{1}{k}\sum_{i\in S}\sigma_i\bigg|\bigg]\asymp 1,\qquad \bigg\Vert\frac{1}{\sqrt k}\sum_{i\in S}\sigma_i-\bb E\bigg[\bigg|\frac{1}{\sqrt k}\sum_{i\in S}\sigma_i\bigg|\bigg]\bigg\Vert_{\psi_2}<\infty.
    \end{align*}
    And the rest of the proof follows similar arguments as the proof of the Theorem 3.4 in \citep{he2023hidden1}.
\subsubsection{Proof of the Global Algorithm}
    The proof goes by noticing that Lemma \ref{cltconsk} implies that under the null we have
$        \frac{1}{\sqrt n}\sum_{i=1}^n\sigma_i\overset{d}{\to} N(0,1).
$
    And under the alternative we have
    \tny{\begin{align*}
       \bl\frac{1}{\sqrt k}\sum_{i=1}^k(\sigma_i-\mu)\bigg|m_1>0\br\overset{d}{\to} N\bl 0,\frac{1-\theta_1(\bb E[\sech^2(\theta\sqrt q z+h)])^2}{(1-\theta_1\bb E[\sech^2(\theta\sqrt q z+h)])^2}\br,\thinspace \frac{1}{\sqrt{n-k}}\sum_{i=k+1}^n\sigma_i\overset{d}{\to} N(0,1).
    \end{align*}}
    And moreover, we also notice that by Lemma \ref{covariance} we have
    \begin{align*}
        \bb E\bigg[\frac{1}{\sqrt{(kn-k)}}\sum_{i=1}^k(\sigma_i-\mu)\sum_{i=k+1}^n\sigma_i\bigg]=O\bl\frac{1}{\sqrt n}\br.
    \end{align*}
    We introduce the notation $m:=\frac{1}{k}\sum_{i=1}^k\sigma_i$ and $\tde m:=\frac{1}{n-k}\sum_{i=k+1}^n\sigma_i$ to get that
    \tny{\begin{align*}
        \bl m+\frac{n-k}{k}\tde m\bigg | m>0\br\overset{d}{\to} N\bl\bb E[|m|],\frac{n}{k^2}\br,\quad\bl m+\frac{n-k}{k}\tde m\bigg | m<0\br\overset{d}{\to} N\bl-\bb E[|m|],\frac{n}{k^2}\br.
    \end{align*}}
    And this further leads  to
    \begin{align*}
        \bigg| m+\frac{n-k}{k}\tde m\bigg|\overset{d}{\to}\ca {FN}\lef(\bb E[|m|],\frac{n}{k^2}\rig),
    \end{align*}
    where we denote $\ca {FN}$ to be the folded Gaussian.
    Then the rest of the proof follows the proof of the Theorem 3.7 of \citep{he2023hidden1}.
\subsubsection{Proof of the Lower Bounds }
    The proof follows directly from the third inequality in Lemma \ref{approximation bounds chisq} and the proof of the Theorem 3.6 in \citep{he2023hidden1}. Using the same notation as in Lemma \ref{cltrfcr} we have
    \begin{align*}
        \bb E\bigg[\frac{\bb P_S(\cdot|\bfa h,\bfa g)\bb P_{S^\prime}(\cdot|\bfa h,\bfa g)}{\bb P_0(\cdot|\bfa h,\bfa g)}\bigg]\leq \bb E\lef[\frac{\bb Q_{S}(\cdot|\bfa h^\prime)\bb Q_{S^\prime}(\cdot|\bfa h^\prime)}{\bb Q_0(\cdot|\bfa h^\prime)}\rig]\exp\bl \frac{Ck^2}{n}\br\leq\exp\bl\frac{Ck^2}{n}+C(k-r)\br.
    \end{align*}
    And the rest of the proof follows from theorem 3.6 in \citep{he2023hidden1}.
\subsection{Proof of Lemma \ref{guaranteecriticaltest}}
Similar to the previous two proofs for lemmas \ref{thm1} and \ref{thmlowtempfindclique}, we also organize the proof according to the two different 
\subsubsection{Proof of the Local Algorithm}

    The proof is identical to the proof of the Theorem 3.7 in \citep{he2023hidden1}.
\subsubsection{Proof of the Global Algorithm}
     The proof is identical to the proof of the Theorem 3.8 in \citep{he2023hidden1}.
\subsubsection{Proof of the Lower Bound}
The proof will be largely identical to the proof of the Theorem 3.9 in \citep{he2023hidden1}. weuse similar method as in the lower bound proof of Theorem \ref{thm1}. Instead, here the event that we construct is
\begin{align*}
    E_S^*:=\lef\{|k^{\frac{1}{4\tau-2}}m_S|\leq C(\log(m\vee k)\log k)^{\frac{1}{4\tau-2}} \rig\},\qquad \bb P_S^*(\bfa\sigma)=\begin{cases}
        \bb P_S(\bfa\sigma)&\text{ if }\bfa\sigma\in E_S^*\\
        0&\text{ otherwise }
    \end{cases}.
\end{align*}
Then, using lemma \ref{cltrfcr} wesee that
\begin{align*}
    \Vert\bb P_{\bar S}-\bb P_{\bar S}^*\Vert_{TV}\leq\frac{1}{\binom{n}{k}}\sum_{S:|S|\leq k}\bb P_{S}(E_S^c)=O\bl\frac{1}{(m\vee k)\log k}\br.
\end{align*}
Then weuse lemma \ref{approximation bounds chisq} to get that when $k\log^{\frac{1}{2\tau-1}}k=O\lef(n^{\frac{2\tau-1}{4\tau-4}}\rig)$
\begin{align*}
    &\bb E\bigg[\frac{\bb P_S^*(\cdot|\bfa h,\bfa g)\bb P_{S^\prime}^*(\cdot|\bfa h,\bfa g)}{\bb P_0(\cdot|\bfa h,\bfa g)}\bigg]=\bb E\bigg[\frac{\mbbm 1_{E_S^*}\mbbm 1_{E_{S^\prime}^*}\bb P_S(\cdot|\bfa h,\bfa g)\bb P_{S^\prime}(\cdot|\bfa h,\bfa g)}{\bb P_0(\cdot|\bfa h,\bfa g)}\bigg]\\
    &\leq\exp\bl\frac{Ck^{\frac{4\tau-4}{2\tau-1}}\log^{\frac{2}{2\tau-1}}k}{n}+C\bl\frac{k-r}{k}\br^2k^{\frac{4\tau-4}{2\tau-1}}\wedge \bl\frac{k-r}{k}\br k^{\frac{2\tau-2}{2\tau-1}}\log^{\frac{1}{2\tau-1}}(m\vee k)\br.
\end{align*}
And wetake the cap $m\lesssim\frac{n}{k^{\frac{4\tau-4}{2\tau-1}}\log^{\frac{2}{2\tau-1}}k}$ and proceed the rest of the proof similar to that of the Theorem 3.9 in \citep{he2023hidden1}.
\subsection{Proof of Lemma \ref{cltrfcr}}
The proof will be based on a novel smart path method where we decouple the correlation of local spins in the clique $S$ from the outside. Without loss of generality we assume that $S=[k]$. Then we define the following quantities.
\begin{align*}
    u_{\bfa\sigma} &:= \sum_{i,j\in [k],i\neq j}\frac{\theta g_{ij}}{\sqrt n}\sigma_i\sigma_j+\sum_{i\in[k],j\in[k+1:n]}\frac{\theta g_{ij}}{\sqrt n}\sigma_i\sigma_j,\thinspace v_{\bfa\sigma} := \sum_{i\in[k]}\theta\sqrt{q}z_i\sigma_i+\sum_{i=k+1}^n\theta\sqrt{q^\prime}z_i\sigma_i,\\
    \omega({\bfa\sigma})&:=\exp\bl\sum_{i,j\in[k]}\frac{\theta_1}{2k}\sigma_i\sigma_j+\sum_{k+1\leq i<j\leq n}\frac{\theta g_{ij}}{\sqrt n}\sigma_i\sigma_j+\sum_{i=1}^n\sigma_ih_i\br
.\end{align*}
where $q$ and $s$ are positive numbers to be specified later. Therefore we adopt the similar notation as in Lemma \ref{gil} to define $\bfa u(t)$ indexed by $\{-1,1\}^n$ with $u_{\bfa\sigma}(t):=\sqrt{t}u_{\bfa\sigma}+\sqrt{1-t}v_{\bfa\sigma}$. Here we consider an arbitrary local function $f:\{-1,1\}^k\to\bb R$ that takes the spins within $S$. The quantity to be interpolated is then given by 
\begin{align*}
   \bb E[F(\bfa u(t))] =\bb E[\la f(\bfa\sigma)\ra_t]=\bb E\lef[ \frac{\sum_{\bfa\sigma}f(\bfa\sigma)\omega({\bfa\sigma})\exp(u_{\bfa\sigma}(t))}{\sum_{\bfa\sigma}\omega({\bfa\sigma})\exp( u_{\bfa\sigma}(t))}\rig]=:\bb E\lef[\frac{Y(\bfa u(t))}{Z(\bfa u(t))}\rig]
.\end{align*}
where $Y(\bfa x):=\sum_{\bfa\sigma}\omega(\bfa\sigma)f(\bfa\sigma)\exp(x_{\bfa\sigma})$ and $Z(\bfa x):=\sum_{\bfa\sigma}\omega(\bfa\sigma)\exp(u_{\bfa\sigma}(t))$.
We also use the notation $\la f\ra_t:=\frac{\sum_{\bfa\sigma} \omega(\bfa\sigma)f(\bfa\sigma)\exp(u_{\bfa\sigma}(t))}{\sum_{\bfa\sigma}\omega(\bfa\sigma)\exp(u_{\bfa\sigma}(t))}$ to be 'expectation' with regard to Gibbs measure for and function $\omega$. We note that since $F(\bfa x)\leq\exp(C\Vert \bfa x\Vert_\infty)$ for some universal $C$. Hence, both the function itself and its first derivative satisfies the \emph{moderate growth condition}. Then we introduce $R_{1,2}^S:=\frac{1}{k}\sum_{i\in S}\sigma_i^1\sigma_i^2$ and $R_{1,2}^{S^c}:=\frac{1}{n-k}\sum_{i\in S^c}\sigma_i^1\sigma_i^2$ for $\bfa\sigma^1,\bfa\sigma^2\in\{-1,1\}^n$. We also define the following quantity that is used in Gaussian interpolation in Lemma \ref{gil}, with the introduction of $q^{(1)}, q^{(2)}$ to be specified later.
\tny{\begin{align}\label{usigma1sigma2}
    &U(\bfa\sigma,\bfa\tau) :=\bb E[u_{\bfa\sigma^1}(t)u_{\bfa\sigma^2}(t)-v_{\bfa\sigma^1}(t)v_{\bfa\sigma^2}(t)]\nnb\\
    &=\frac{\theta^2}{2n} \lef(k^2R_{1,2}^{S,2}-1+2k(n-k)R_{1,2}^SR_{1,2}^{S^c}\rig)-\theta^2qkR_{1,2}^S-\theta^2q^\prime(n-k)R_{1,2}^{S^c}\nnb\\
    &=\frac{\theta^2}{2n}k^2\lef (R_{1,2}^{S,2}-2q^{(1)}R_{1,2}^S-\frac{1}{k^2}\rig )+\frac{\theta^2}{n}k(n-k)(R_{1,2}^S-q^{(1)})(R_{1,2}^{S^c}-q^{(2)})-\frac{\theta^2}{n}k(n-k)q^{(1)}q^{(2)}.
    \end{align}}
    And to make the equation above consistent with $u_{\bfa\sigma}$ and $v_{\bfa\sigma}$, we set the values of $q$s as follows:
\tny{\begin{align*}
    q^{(1)}&=\bb E[\tanh^2(\sqrt{\theta_1}x_1^*+\theta\sqrt q z+h)],\thinspace
    q^{(2)}=\nu_0(R_{1,2}^{S^c}),\thinspace q=\frac{k}{n}q^{(1)}+\frac{n-k}{n}q^{(2)},\thinspace q^\prime=\frac{k+r}{n}q^{(1)}.
\end{align*}}
Note that here $\nu_0(R_{1,2}^{S^c})$ corresponds to the 
Then we define $\varphi(t):=\bb E[F(\bfa u(t))]$ and by Gaussian integration by parts, the following holds
\tny{\begin{align*}
    \varphi^\prime(t)=\frac{k^2\theta^2}{2n}\bb E[\la(R_{1,2}^S-q^{(1)})^2(f-\la f\ra_t)\ra_t]+\frac{\theta^2k(n-k)}{n}\bb E\lef[\la(R_{1,2}^S-q^{(1)})(R_{1,2}^{S^c}-q^{(2)})(f-\la f\ra_t)\ra_t \rig].
\end{align*}}
The goal is then to upper bound the first and second term separately. First, recall that using the result in Lemma \ref{correxponential} and the exponential inequality given by the theorem 13.7.1 in \citep{talagrand2011mean2} we have at the high/low temperature $\theta_1\neq\frac{1}{\bb E[\sech^2(\theta\sqrt qz+h)]}$ there exists $\lambda>0$ such that
\begin{align}\label{exponentialbounds}
    &\bb E[\la\exp(\lambda k(R^S_{1,2}-q^{(1)})^2)\ra_{t}]\leq C,\qquad
    \bb E[\la\exp(\lambda(n-k)(R_{1,2}^{S^c}-q^{(2)})^2)\ra_0]\leq C.
\end{align}
Then we check that for all $\tau>0$
\tny{\begin{align*}
    \bb P(|R_{1,2}^S-q^{(1)}|\geq \tau)\leq\exp(-k C\tau^2)\thinspace\text{ for all }t\in[0,1],\thinspace\bb P(|R_{1,2}^{S^c}-q^{(2)}|\geq \tau)\leq\exp(-k C\tau^2)\text{ for }t=0.
\end{align*}}
Therefore we can set the truncation threshold at $\sqrt{\frac{\log k}{k}}$, further noticing that $|R_{1,2}-q|\leq 2$, the first term can be analyzed as
\tny{\begin{align}\label{truncate1}
    \frac{k^2\theta^2}{2n}&\bb E[\la(R_{1,2}^S-q^{(1)})^2(f-\la f\ra_t)\ra_t]=\frac{k^2\theta^2}{2n}\bb E\lef[\bigg\la\bl\mbbm 1_{|R_{1,2}^S-q^{(1)}|\leq{\sqrt \frac{\log k}{k}}}\br (R_{1,2}^S-q^{(1)})^2(f-\la f\ra_t)\bigg\ra_t\rig]\nnb\\
    &+\frac{k^2\theta^2}{2n}\bb E\lef[\bigg\la\bl\mbbm 1_{|R_{1,2}^S-q^{(1)}|>{\sqrt \frac{\log k}{k}}}\br (R_{1,2}^S-q^{(1)})^2(f-\la f\ra_t)\bigg\ra_t\rig]\lesssim \lef(\frac{k\log k}{n}+\frac{k}{n}\rig)\varphi(t).
\end{align}}
And analogously we analyze the second term through Taylor expansion, recalling that $f$ is a function on the first $k$ spins and under $\la\ra_0$ we have $\bfa\sigma_{S}\perp\bfa\sigma_{S^c}$ and $\bb E[\la(R_{1,2}^{S^c}-q^{(2)})^{2\tau} \ra_0]=O\lef(\frac{1}{n^{\tau}}\rig)$,
\tny{\begin{align}\label{crossterms1}
    &\bb E[\la(R_{1,2}^S-q^{(1)})(R_{1,2}^{S^c}-q^{(2)})(f-\la f\ra_t)\ra_{t}]=\bb E[\la(R_{1,2}^S-q^{(1)})(R_{1,2}^{S^c}-q^{(2)})(f-\la f\ra_0)\ra_{0}]\nnb\\
    &+\frac{\theta^2(k+r)(n-k-r)t}{n}
   \bb E[\la(R_{1,2}^S-q^{(1)})^2(R_{1,2}^{D^c}-q^{(2)})^2(f-\la f\ra_0)\ra_{0}]\bl 1+ O\bl\sqrt{\frac{k\log k}{n}}\br\br\nnb\\
   &\lesssim \frac{k}{n}\bb E[\la(R_{1,2}^D-q^{(1)})^2(f-\la f\ra_0))\ra_{0}]\bl 1+O\bl\sqrt{\frac{k\log k}{n}}\br\br\lesssim \frac{\log k}{n}\varphi(0).
\end{align}}
Collecting the above pieces weget:
\begin{align*}
    \varphi^\prime(t)\leq\frac{Ck\log k}{n}\lef(\varphi(t)+\varphi(0)\rig)\qquad\Rightarrow\qquad\varphi(t)\leq\exp\lef(\frac{Ctk\log k}{n}\rig)\varphi(0).
\end{align*}
In particular, here we also derive that when $\varphi(0)$ is bounded and $\frac{k\log k}{n}=O(1)$, $\varphi(t)$ is uniformly bounded. Then we can apply H\"older's inequality to achieve that
\begin{align*}
    |\varphi(t)-\varphi(0)|=O\lef(\frac{k}{n}\rig).
\end{align*}
Then we discuss over the critical temperature case when $\theta_1=\frac{1}{\bb E[\sech^2(\theta\sqrt qz+h)]}$ and $x_1^*$ has flatness $\tau$. Similar to the high temperature, here we have
\begin{align}\label{exponentialbounds}
    \bb E[\la\exp(\lambda k^{\frac{2\tau-2}{2\tau-1}}(R^S_{1,2}-q^{(1)})^2)\ra_{t}]\leq C,\qquad
    \bb E[\la\exp(\lambda(n-k)(R_{1,2}^{S^c}-q^{(2)})^2)\ra_0]\leq C.
\end{align}
Then wecheck that for all $\tau>0$:
\begin{align*}
    \bb P(|R_{1,2}^S-q^{(1)}|&\geq \tau)\leq\exp(-k^{\frac{2\tau-2}{2\tau-1}} C\tau^2)\thinspace\text{ for all }t\in[0,1],\\
    \bb P(|R_{1,2}^{S^c}-q^{(2)}|&\geq \tau)\leq\exp(-k C\tau^2)\text{ for }t=0.
\end{align*}
Then, using the same truncation method at the threshold of $\lef(\frac{\log k}{k}\rig)^{\frac{\tau-1}{2\tau-1}}$, we have
\begin{align*}
    \varphi(t)\leq\exp\bl\frac{Ck^{\frac{2\tau}{2\tau-1}}\log^{\frac{2\tau-2}{2\tau-1}}k}{n}\br\varphi(0).
\end{align*}
Then the next step is to set $f$ be the mgf of average magnezation in the local spins. This results appeared in \citep{he2023hidden1}, where we proved the following theorem for the RFCW model, which is put here for completeness.
\begin{theorem}[Limiting Theorem for the Random Field Curie-Weiss Model with Symmetric $h$]\label{cltrfcrrfcw}
Assume that $h_i\sim \mu(h)$ is i.i.d. in $L_1$.
For a random field Curie-Weiss model whose Hamiltonian,
\begin{enumerate}
   \item 
In the high temperature regime with $\theta_1<\frac{1}{\bb E[\sech^2(h)]}$,   for $t\in\bb R$,  
\begin{align*}
    \bb E\bigg[\exp\bl n^{-1/2}{t\sum_{i=1}^n\sigma_i}\br\bigg]\to \exp\lef({\ca Vt^2}/2\rig),\qquad\bigg\Vert n^{-1/2}{\sum_{i=1}^n\sigma_i}\bigg\Vert_{\psi_2}<\infty.
\end{align*}
with $\ca V:=\frac{1-\theta_1(\bb E[\sech^2(h)])^2}{(1-\theta_1\bb E[\sech^2(h)])^2}$.

 \item 
In the low temperature regime of $\theta_1>\frac{1}{\bb E[\sech^2(h)]}$, $x=\bb E[\tanh(\sqrt{\theta_1} m +h)]$ will have two nonzero solutions defined by $x_1<0<x_2$. Define $\ca C_1=(0,\infty)$ and $\ca C_2=\ca C_1^c$. Then we have for $t\in\bb R$,  for $\ell\in\{1,2\}$
\begin{align}\label{gtr0clt2}
    &\bb E\lef[\exp\lef(t\frac{\sum^n_{i=1}(\sigma_i -\sqrt{\theta_1}x_\ell)}{\sqrt n}\rig)\bigg |{\frac{\sum_{i=1}^n\sigma_i}{n}\in \ca C_\ell}\rig]{\to} \exp\lef(\frac{\ca V(m_1)t^2}{2}\rig),
\end{align}
and
$
   \lef\Vert n^{-1/2}\sum_{i=1}^n(\sigma_i-\sqrt{\theta_1}x_{\ell})\bigg|n^{-1}\sum_{i=1}^n\sigma_i\in\ca C_{\ell}\rig\Vert_{\psi_2}<\infty,
$ with\\ $\ca V(x):=\frac{(1-\theta_1(\bb E[\sech^2(\sqrt{\theta_1}x+h)])^2-\bb E[\tanh(\sqrt{\theta_1}x+h)]^2)}{(1-\theta_1\bb E[\sech^2(\sqrt{\theta_1}x+h)])^2}$.
 \item 
 At the critical temperature $\theta_1=\frac{1}{\bb E[\sech^2(h)]}$, assume that the critical value is of flatness is $\tau$, then for $t\in\bb R$  ,  
\begin{align}\label{supergaussian2}
        \bb E\lef[\exp\lef(\frac{t\sum_{i=1}^n\sigma_i}{n^{\frac{4\tau-3}{4\tau-2}}}\rig)\rig]\to \int_{\bb R}\frac{(2\tau-1)x^{2\tau-2}}{\sqrt{2\pi v(0)}}\exp\lef(-\frac{x^{4\tau-2}}{2v(0)}+tx\rig)dx,
    \end{align}
    and
$        \lef\Vert n^{-\frac{4\tau-3}{4\tau-2}}\sum_{i=1}^n\sigma_i\rig\Vert_{\psi_{4\tau-2}}<\infty,
$ with \begin{align}\label{stirlingv2}
    v(x):=&((2\tau)!)^2\bb V(\tanh(\sqrt{\theta_1}x+h))(\bb E[\sech^2(\sqrt{\theta_1}x+h)])^{4\tau-2}\nnb\\
    &\cdot\bb E\bigg[(1+\tanh(\sqrt{\theta_1}x+h))\sum_{k=0}^{2\tau-1}\frac{k!}{2^k}S(2\tau-1,k)(\tanh(\sqrt{\theta_1}x+h)-1)^{k}\bigg]^{-2}.
\end{align}
     and $H$ is the function defined in \eqref{taylorcond}. And if we are at the second case of \eqref{taylorcond},   \eqref{gtr0clt2} holds.
\end{enumerate}
\end{theorem}
And using the above interpolation results we finalize the proof.
\subsection{Proof of Lemma \ref{approximation bounds chisq}}
Our proof strategy is to consider conditional Chi-square as lower bound and use the following lemma

The proof relies on a novel construction of the smart path that accounts heterogeneous distribution. Consider two measure of the form $\bb P_S(\cdot|\bfa h,\bfa g),\bb P_S^\prime(\cdot|\bfa h,\bfa g)$ with different position of cliques where $S:=[k]$ and $S^\prime:=[r+1:k+r]$ and Hamiltonian
\begin{align*}
    \mca H_A(\bfa\sigma):=\theta\sum_{1\leq i<j\leq n}\frac{g_{ij}}{\sqrt n}\sigma_i\sigma_j+\sum_{i,j\in A}\frac{\theta_1}{2k}\sigma_i\sigma_j+\sum_{i=1}^n\sigma_ih_i.\thinspace\text{ for }A\in\{S,S^\prime\}.
\end{align*}
And Similarly, we consider
\tny{\begin{align}\label{chisqlb}
        \bb E\lef[D_{\chi^2}\lef(\bar{\bb P}_m,\bb P_{0,m}|\bfa h^m,\bfa g^m\rig)\rig]:=\frac{1}{\binom{n}{k}^2}\sum_{S,S^\prime\in\ca S}\bb E_{0,m}\lef[\frac{\bb P_{S,m}(\bfa\sigma|\bfa h^m,\bfa g^m)\bb P_{S^\prime,m}(\bfa\sigma|\bfa h^m,\bfa g^m)}{\bb P_{0,m}^2(\bfa\sigma|\bfa h^m,\bfa g^m)}\rig]-1
.\end{align}}
Written in detail and define $m_S:=\frac{1}{k}\sum_{i\in S}\sigma_i$ and $m_S^\prime:=\frac{1}{k}\sum_{i\in S^\prime}\sigma_i$, westudy the following quantity.
\begin{align}\label{crossterm}
    \bb E_{0,m}&\lef[\frac{\bb P_{S,m}(\bfa\sigma|\bfa h^m,\bfa g^m)\bb P_{S^\prime,m}(\bfa\sigma|\bfa h^m,\bfa g^m)}{\bb P_{0,m}^2(\bfa\sigma|\bfa h^m,\bfa g^m)}\rig]=\bb E_0\lef[\frac{\bb P_S(\bfa\sigma|\bfa h,\bfa g)\bb P_{S^\prime}(\bfa\sigma|\bfa h,\bfa g)}{\bb P^2_0(\bfa\sigma|\bfa h,\bfa g)}\rig]^m\nnb\\
    &=\bb E\bigg[\frac{\sum_{\bfa\sigma}\exp\lef(\sum_{1\leq i<j\leq n}\frac{g_{ij}}{\sqrt n}\sigma_i\sigma_j+\frac{\theta_1k}{2}\lef(m_S^2+m_{S^\prime}^2\rig)+\sum_{i=1}^nh_i\sigma_i\rig)}{\sum_{\bfa\sigma}\exp(\sum_{1\leq i<j\leq n}\frac{g_{ij}}{\sqrt n}\sigma_i\sigma_j+\frac{\theta_1k}{2}m_S^2+\sum_{i=1}^nh_i\sigma_i)}\nnb\\
    &\cdot\frac{\sum_{\bfa\sigma}\exp\lef(\sum_{1\leq i<j\leq n}\frac{\theta g_{ij}}{\sqrt n}\sigma_i\sigma_j+\sum_{i=1}^nh_i\sigma_i\rig)}{\sum_{\bfa\sigma}\exp(\sum_{1\leq i<j\leq n}\frac{g_{ij}}{\sqrt n}\sigma_i\sigma_j+\frac{\theta_1k}{2}m_{S^\prime}^2+\sum_{i=1}^nh_i\sigma_i)}\bigg]^m\nnb\\
     &=\bb E\bigg[\bigg\la\exp\bigg(\frac{\theta_1}{2k}\sum_{i,j\in S^\prime}\sigma_i\sigma_j\bigg)\bigg\ra_S\bigg\la\exp\bigg(-\frac{\theta_1}{2k}\sum_{i,j\in S^\prime           }\sigma_i\sigma_j\bigg)\bigg\ra_{S^\prime}\bigg]^m
.\end{align}
From here on we assume that $S=[k]$ and $S^\prime=[r+1:r+k]$ without loss of generality. Hence one will notice that their overlap is $|S\cap S^\prime|=k-r$. Notice that the above expectation cannot be computed explicitly. To analyze it, we construct a smart path to transform the two expectation simultaneously. wedefine the intermediate Hamiltonian w.r.t. the measure $\bb P_S(\cdot|\bfa h,\bfa g)$ as follows for $t\in[0,1]$
\tiny{\begin{align}\label{interpolating_hamilt}
    \mca H_{S}&(\bfa\sigma, t):=-\bl\sqrt{1-t}\bl\sum_{1\leq i< j\leq k+r}\frac{\theta g_{ij}}{\sqrt n}\sigma_i\sigma_j+\sum_{i\in[k+r],j\in[k+r+1:n]}\frac{\theta g_{ij}}{\sqrt n}\sigma_i\sigma_j\br\nnb\\
    &+\sqrt t\bl\sum_{i=1}^{k+r}\theta\sqrt{q}z_i\sigma_i+\sum_{i=k+r+1}^n\theta\sqrt{q^\prime}z_i\sigma_i\br+\sum_{i,j\in[k]}\frac{\theta_1}{2k}\sigma_i\sigma_j+\sum_{k+r+1\leq i<j\leq n}\frac{\theta g_{ij}}{\sqrt n}\sigma_i\sigma_j+\sum_{i=1}^n\sigma_ih_i\br,
\end{align}}\normalsize
where we introduce $q,q^\prime\in\bb R$ to be specified later. 
Analogously, wedefine $u_{\bfa\sigma}(t):=\sum_{1\leq i< j\leq k+r}\frac{\theta g_{ij}}{\sqrt n}\sigma_i\sigma_j+\sum_{i\in[k+r],j\in[k+r+1:n]}\frac{\theta g_{ij}}{\sqrt n}\sigma_i\sigma_j$ and $v_{\bfa\sigma}(t):=\sum_{i=1}^{k+r}\theta \sqrt qz_i\sigma_i$. Then, for $\bfa\sigma^1,\bfa\sigma^2$ being two replicas, define $R_{1,2}^{D}:=\frac{1}{k+r}\sum_{i=1}^{k+r}\sigma^1_i\sigma^2_i$ and $R_{1,2}^{D^c}:=\frac{1}{n-k-r}\sum_{i=k+r+1}^n\sigma_i^1\sigma_i^2$ and we have by similar argument as \eqref{usigma1sigma2}
\tny{\begin{align*}
    &U(\bfa\sigma^1,\bfa\sigma^2)=\frac{1}{2}\bb E[u_{\bfa\sigma^1}(t)u_{\bfa\sigma^2}(t)-v_{\bfa\sigma^1}(t)v_{\bfa\sigma^2}(t)]\\
    &=\frac{\theta^2}{2n} \lef((k+r)^2R_{1,2}^{D,2}-1+2(k+r)(n-k-r)R_{1,2}^DR_{1,2}^{D^c}\rig)-\theta^2q(k+r)R_{1,2}^D-\theta^2q^\prime(n-k-r)R_{1,2}^{D^c}\\
    &=\frac{\theta^2}{2n}(k+r)^2\lef (R_{1,2}^{D,2}-2q^{(1)}R_{1,2}^D-\frac{1}{(k+r)^2}\rig )+\frac{\theta^2}{n}(k+r)(n-k-r)(R_{1,2}^D-q^{(1)})(R_{1,2}^{D^c}-q^{(2)})\\
    &-\frac{\theta^2}{n}(k+r)(n-k-r)q^{(1)}q^{(2)}.
\end{align*}}
Here we introduce two new values $q^{(1)},q^{(2)}$ and the relationship between $q,q^\prime,q^{(1)},q^{(2)}$ is given by
\tny{\begin{align*}
    q^{(1)}&=\frac{k-r}{k+r}\bb E[\tanh^2(\sqrt{\theta_1}x^*_1+\theta\sqrt q z+h)]+\frac{2r}{k+r}\bb E[\tanh(\sqrt{\theta_1}x_1^*+\theta\sqrt qz+h)\tanh(\theta \sqrt{q}z+h)],\\
    q^{(2)}&=\nu_0(R_{1,2}^{D^c}),\qquad q=\frac{k+r}{n}q^{(1)}+\frac{n-k-r}{n}q^{(2)},\qquad q^\prime=\frac{k+r}{n}q^{(1)}.
\end{align*}}
And in particular when $x_1^*=0$ we have
$    q^{(1)}=\bb E[\tanh^2(\theta\sqrt{q}z+h)]
$. By the standard theory of SK model $\lim_{n\to\infty }q^{(2)}=\lim_{n\to\infty}\bb E\lef[\tanh^2\lef(\theta\sqrt{\frac{n-k}{n}}\sqrt{q^{(2)}}z+h\rig)\rig]$.
Analogously, we can define $\mca H_{S^\prime}(\bfa\sigma,t)$. And we extend the previous definition by letting $\la\ra_{A,t}$ to be the Gibbs average w.r.t. $\mca H_{A,t}(\bfa\sigma)$ for $A\in\{S,S^\prime\}$. A new notation $\la \ra_{S,S^\prime,t}$ is also introduced such that for all $F\in\{+1,-1\}^{n}\times\{+1,-1\}^{n}\to\bb R$:
\begin{align*}
\la F(\bfa\sigma_1,\bfa\sigma_2) \ra_{S,S^\prime, t}:= \frac{\sum_{\bfa\sigma_1}\sum_{\bfa\sigma_2}F(\bfa\sigma_1,\bfa\sigma_2)\exp(-\mca H_{S}(\bfa\sigma_1)-\mca H_{S^\prime}(\bfa\sigma_2))}{\sum_{\bfa\sigma_1}\sum_{\bfa\sigma_2}\exp(-\mca H_{S}(\bfa\sigma_1)-\mca H_{S^\prime}(\bfa\sigma_2))}.
\end{align*}
The next step is to study the error terms given by interpolation. wenotice that for $
    \ca A(t):= \la f_{S^\prime}\ra_{S,t}$ and $\ca B(t):=\la f_{S^\prime}^{-1}\ra_{S^\prime,t}$. We can then write $\varphi(t):=\bb E\lef[\la f_{S^\prime}\ra_{S,t}\la f_{S^\prime}^{-1}\ra_{S^\prime, t}\rig] =\bb E\lef[\ca A(t)\ca B(t)\rig]$ we have by Lemma \ref{gil}:
\tny{\begin{align}\label{derivativebound}
    &\varphi^\prime(t)=\bb E\bigg[\sum_{\bfa\sigma,\bfa\tau\in\{-1,1\}^{k+r}}U(\bfa\sigma,\bfa\tau)\frac{\pta^2 \ca A(t)\ca B(t)}{\pta x_{\bfa\sigma}\pta x_{\bfa\tau}}\bigg]\nnb\\
    &=\bb E\bigg[\sum_{\bfa\sigma,\bfa\tau\in\{-1,1\}^{k+r}}U(\bfa\sigma,\bfa\tau)\bigg(\frac{\pta^2 \ca A(t)}{\pta x_{\bfa\sigma}\pta x_{\bfa\tau}}\ca B(t)+2\frac{\pta\ca A(t)}{\pta x_{\bfa\sigma}}\frac{\pta \ca B(t)}{\pta x_{\bfa\tau}}+\frac{\pta^2 \ca B(t)}{\pta x_{\bfa\sigma}\pta x_{\bfa\tau}}\ca A(t)\bigg)\bigg]\nnb\\
    &=\frac{(k+r)^2\theta^2}{2n}\bigg(\bb E[\la (R^D_{1,2}-q^{(1)})^2(f_{S^\prime}-\la f_{S^\prime}\ra_{S,t})\ra_{S,t}\la f_{S^\prime}^{-1}\ra_{S^\prime,t}]\nnb\\
    &+\bb E[\la (R^{D}_{1,2}-q^{(1)})^2(f_{S^\prime}^{-1}-\la f_{S^\prime}^{-1}\ra_{S^\prime,t})\ra_{S^\prime,t}\la f_{S^\prime}\ra_{S,t}]\nnb\\
    &+2\bb E\lef[\la (R^{D}_{1,2}-q^{(1)})^2(f_{S^\prime}(\bfa\sigma_1)-\la f_{S^\prime}(\bfa\sigma_1)\ra_{S,t})( f_{S^\prime}^{-1}(\bfa\sigma_2)-\la f^{-1}_{S^\prime}(\bfa\sigma_2)\ra_{S^\prime,t})\ra_{S,S^\prime,t} \rig]\bigg)\nnb\\
    &+\frac{\theta^2(k+r)(n-k-r)}{n}\bl\bb E[\la(R_{1,2}^D-q^{(1)})(R_{1,2}^{D^c}-q^{(2)})(f_{S^\prime}-\la f_{S^\prime}\ra_{S,t})\ra_{S,t}\la f^{-1}_{S^\prime}\ra_{S^\prime,t}]\nnb\\
    &+\bb E[\la(R_{1,2}^D-q^{(1)})(R_{1,2}^{D^c}-q^{(2)})(f^{-1}_{S^\prime}-\la f^{-1}_{S^\prime}\ra_{S^\prime,t})\ra_{S^\prime,t}\la f_{S^\prime}\ra_{S^\prime, t}]\nnb\\
    &+2\bb E\lef[\la (R^{D}_{1,2}-q^{(1)})(R_{1,2}^{D^c}-q^{(2)}) (f_{S^\prime}(\bfa\sigma_1)-\la f_{S^\prime}(\bfa\sigma_1)\ra_{S,t})( f_{S^\prime}^{-1}(\bfa\sigma_2)-\la f^{-1}_{S^\prime}(\bfa\sigma_2)\ra_{S^\prime,t})\ra_{S,S^\prime,t}\rig]\br
.\end{align}}
To upperbound the above quantity, we observe that at the high temperature $x_1^*=0$, there exists $\lambda>0$ such that 
\begin{align}\label{exponentialbounds}
    &\bb E[\la\exp(\lambda(k+r)(R^D_{1,2}-q^{(1)})^2)\ra_{A,t}]\leq C,\text{ for }A\in\{S,S^\prime\}\\
    &\bb E[\la\exp(\lambda(k+r)(R^D_{1,2}-q^{(1)})^2)\ra_{S,S^\prime,t}]\leq C,\thinspace
    \bb E[\exp(\lambda(n-k-r)(R_{1,2}^{D^c}-q^{(2)})^2)]\leq C.
\end{align}

Then we have the following tail bound by the property of $\psi_2$ Orlicz norm being bounded.
\tny{\begin{align*}
    \bb P_{S,S^\prime, t}(|R_{1,2}^D-q^{(1)}|\geq t)\leq\exp(-C kt^2),\thinspace\bb P_{A,t}(|R_{1,2}^D-q^{(1)}|\geq t)\leq\exp(-C kt^2)\text{ for } A\in\{S,S^\prime\}.
\end{align*}}
And also noticing that $|R^D_{1,2}-q|\vee|R^{D^c}_{1,2}-q|\leq 2$, we collect the above pieces \eqref{derivativebound}, \eqref{crossterms} to get that
Then, we first do Taylor expansion on the cross terms, considering a general function $g(\bfa\sigma)$ on the spins in $D=S\cap S^\prime$ and recall that under $\nu_0$ the spins in $D$ and $D^c$ are independent, $q^{(2)}=\nu_0(R_{1,2}^{D^c})$,
\begin{align}\label{crossterms}
    &\bb E[\la(R_{1,2}^D-q^{(1)})(R_{1,2}^{D^c}-q^{(2)})g(\bfa\sigma)\ra_{S,S^\prime,t}]=\bb E[\la(R_{1,2}^D-q^{(1)})(R_{1,2}^{D^c}-q^{(2)})g(\bfa\sigma)\ra_{S,S^\prime,0}]\nnb\\
    &+\frac{\theta^2(k+r)(n-k-r)}{n}
   \bb E[\la(R_{1,2}^D-q^{(1)})^2(R_{1,2}^{D^c}-q^{(2)})^2g(\bfa\sigma)\ra_{S,S^\prime,0}]\bl 1+ O\bl\sqrt{\frac{k\log k}{n}}\br\br\nnb\\
   &\lesssim \frac{k}{n}\bb E[\la(R_{1,2}^D-q^{(1)})^2g(\bfa\sigma)\ra_{S,S^\prime,0}]\bl 1+O\bl\sqrt{\frac{k\log k}{n}}\br\br.
\end{align}
Then, weconclude that, uniformly across all $t\in[0,1]$:
\begin{align*}
    \varphi^\prime(t)&\lesssim \frac{k\log k}{n}\bb E\lef[\mbbm 1_{|R_{1,2}^D-q^{(1)}|\leq{\sqrt \frac{\log k}{k}}}\ca A(t)\ca B(t)\rig]+\bb E\lef[\mbbm 1_{|R^D_{1,2}-q^{(1)}|\geq\sqrt{\frac{\log k}{k}}}\ca A(t)\ca B(t)\rig]\\
    &\lesssim \frac{k\log k}{n}\varphi(t)+4\bb P\bl|D_{1,2}^D-q^{(1)}|\geq\sqrt{\frac{\log k}{k}}\br\varphi(t)\lesssim \frac{k\log k}{n}\varphi(t).
\end{align*}
Hence we solve the above differential equation to get
\begin{align*}
    \varphi(t)\leq\exp\bl\frac{Ck\log kt}{n}\br\varphi(0).
\end{align*}
However, this is not the best possible rate one can get, wethen utilize the boundedness of $\varphi(t)$ for all $\theta_1<\frac{1}{2}\theta_c$, $\theta<\theta_d$ defined by the AT line and H\"older's inequality to get uniformly across all $t\in[0,1]$
\tny{\begin{align*}
    \varphi^\prime(t)\lesssim \frac{k^2}{n}\bb E[(R_{1,2}^D-q^{(1)})^{1+\delta_1}]^{\frac{1}{1+\delta_1}}\bb E[\ca A^{1+\delta_2}(t)\ca B^{1+\delta_2}(t)]^{\frac{1}{1+\delta_2}}\lesssim \frac{k}{n}\thinspace\Rightarrow\thinspace|\varphi(1)-\varphi(0)|=O\lef( \frac{k}{n}\rig),
\end{align*}}
where $\delta_1,\delta_2>0$ and satisfies $\frac{1}{1+\delta_1}+\frac{1}{1+\delta_2}=1$. Analogously, notice that at the critical temperature $\theta_1=\theta_c$ we have
\begin{align*}
    &\bb E[\la\exp(\lambda(k+r)^{\frac{2}{2\tau-1}}(R^D_{1,2}-q^{(1)})^2)\ra_{A,t}]\leq C,\text{ for }A\in\{S,S^\prime\}\\
    &\bb E[\la\exp(\lambda(k+r)^{\frac{2}{2\tau-1}}(R^D_{1,2}-q^{(1)})^2)\ra_{S,S^\prime,t}]\leq C,\thinspace
    \bb E[\exp(\lambda(n-k-r)(R_{1,2}^{D^c}-q^{(2)})^2)]\leq C,
\end{align*}
which leads to, for all $t\in[0,1]$:
\begin{align*}
    \bb P_{S,S^\prime, t}(|R_{1,2}^D-q^{(1)}|\geq t)&\leq\exp(-C k^{\frac{2}{2\tau-1}}t^2),\\
    \bb P_{A,t}(|R_{1,2}^D-q^{(1)}|\geq t)&\leq\exp(-C k^{\frac{2}{2\tau-1}}t^2)\text{ for } A\in\{S,S^\prime\}.
\end{align*}
And following similar path wearrive at
\begin{align*}
    \varphi(t)\leq\exp\bl\frac{Ctk^{\frac{4\tau-4}{2\tau-1}}\log^{\frac{2}{2\tau-1}}k}{n}\br\varphi(0).
\end{align*}
Then we study the low temperature case, it is worth noticing that in this case wenot have \eqref{exponentialbounds} and we can only study the $k=O(\sqrt n)$ case as
\begin{align*}
    \varphi^\prime(t)\lesssim\frac{k^2}{n}\varphi(t)\thinspace\Rightarrow\thinspace\varphi(t)\leq\exp\lef(Ct\frac{k^2}{n}\rig)\varphi(0).
\end{align*}
Given the above interpolation bound, wethen proceed to the analysis of $\varphi(0)$. We notice that at the $0$ point the Hamiltonian is decoupled by
\tny{\begin{align*}
    -\mca H_{S}(\bfa\sigma,0 )&=\sum_{i,j\in[k]}\frac{\theta_1}{2k}\sigma_i\sigma_j+\ub{\sum_{i=1}^{k+r}(\theta\sqrt qz_i+h_i)\sigma_i}_{-\mca H_0(\bfa\sigma_{[k+r]})}+\ub{\sum_{k+r+1\leq i<j\leq n}\frac{g_{ij}}{\sqrt n}\sigma_i\sigma_j+\sum_{i=k+r+1}^n(\theta\sqrt{q^\prime}z_i+h_i)\sigma_i}_{-\mca H_1(\bfa\sigma_{[k+r+1:n]})}.
\end{align*}}
Then wecheck that $\varphi(0)$ can be rewritten as
\tiny{\begin{align*}
    &\bb E\lef[\frac{\sum_{\bfa\sigma}\exp\lef(-\frac{k\theta_1}{2}(m_S^2+m_{S^\prime}^2)-\mca H_0(\bfa\sigma_{[k+r]})-\mca H_1(\bfa\sigma_{[k+r+1:n]})\rig)\sum_{\bfa\sigma}\exp\lef(-\mca H_0(\bfa\sigma_{[k+r]})-\mca H_1(\bfa\sigma_{[k+r+1:n]})\rig)}{\sum_{\bfa\sigma}\exp(-\frac{k\theta_1}{2}m_S^2-\mca H_0(\bfa\sigma_{[k+r]})-\mca H_1(\bfa\sigma_{[k+r+1:n]}))\sum_{\bfa\sigma}\exp(-\frac{k\theta_1}{2}m_S^2-\mca H_0(\bfa\sigma_{[k+r]})-\mca H_1(\bfa\sigma_{[k+r+1:n]}))}\rig]\\
    &=\bb E\lef[\frac{\sum_{\bfa\sigma_{[k+r]}}\exp\lef(-\frac{k}{2}\theta_1(m_S^2+m_{S^\prime}^2)-\mca H_0(\bfa\sigma_{[k+r]})\rig)\sum_{\bfa\sigma_{[k+r]}}\exp\lef(-\mca H_0(\bfa\sigma_{[k+r]})\rig)}{\sum_{\bfa\sigma_{[k+r]}}\exp\lef(-\frac{k}{2}\theta_1m_S^2-\mca H_0(\bfa\sigma_{[k+r]})\rig)\sum_{\bfa\sigma_{[k+r]}}\exp\lef(-\frac{k}{2}\theta_1m_{S^\prime}^2-\mca H_0(\bfa\sigma_{[k+r]})\rig)}\rig].
\end{align*}}\normalsize
This gives the standard chi-square divergence of the RFCW model.
\subsection{Proof of Lemma \ref{correxponential}}

For a set of replicas $\bfa\sigma_1\sim\bb P_S(\cdot|\bfa h,\bfa g),\bfa\sigma_2\sim\bb P_{S^\prime}(\cdot|\bfa h,\bfa g)$, weconstruct similar smart path using the Hamiltonian given by \eqref{interpolating_hamilt}. 
Recall the notation $\nu_t(f):=\bb E[\la f(\bfa\sigma_1,\bfa\sigma_2)\ra_{S,S^\prime,t}]$ to denote the expectation of $ f:\{-1,1\}^{2(k+r)}\to \bb R$ that is a function on the local spins in $S\cup S^\prime$. Here we introduce the notation
    \tny{\begin{align*}
        Z_A(\bfa x_1) &:= \sum_{\bfa\sigma}\omega_A(\bfa\sigma)\exp(x_{\bfa\sigma})\thinspace\text{ for }A\in\{S,S^\prime\},
        \thinspace F_1(\bfa x) :=\sum_{\bfa\sigma^1,\bfa\sigma^2}\omega_S(\bfa\sigma^1)\omega_{S^\prime}(\bfa\sigma^2)f(\bfa\sigma^1,\bfa\sigma^2)\exp(x_{\bfa\sigma^1}+x_{\bfa\sigma^2}),\\
        F(\bfa x) 
        &:=Z_S^{-1}(\bfa x)Z_{S^\prime}^{-1}(\bfa x) F_1(\bfa x),\quad \omega_A({\bfa\sigma}):=\exp\bl\sum_{i,j\in A}\frac{\theta_1}{2k}\sigma_i\sigma_j+\sum_{k+r+1\leq i<j\leq n}\frac{\theta g_{ij}}{\sqrt n}\sigma_i\sigma_j+\sum_{i=1}^n\sigma_ih_i\br
    .\end{align*}}
    And we can see that
    \tny{\begin{align*}
    \frac{\pta F(\bfa x)}{\pta x_{\bfa\sigma}}= Z_S^{-1}(\bfa x)Z_{S^\prime}^{-1}(\bfa x)\frac{\pta F_1(\bfa x)}{\pta x_{\bfa\sigma}}-Z_S^{-2}(\bfa x)Z_{S^\prime}^{-1}(\bfa x)\frac{\pta Z_S(\bfa x)}{\pta x_{\bfa\sigma}} F_1(\bfa x)-Z_S^{-1}(\bfa x)Z_{S^\prime}^{-2}(\bfa x)\frac{\pta Z_{S^\prime}(\bfa x)}{\pta x_{\bfa\sigma}} F_1(\bfa x)
    .\end{align*}}
    Therefore we have
    \ttny{\begin{align}\label{complicate}
        \frac{\pta F(\bfa x)}{\pta x_{\bfa\sigma}\pta x_{\bfa\tau}} &=Z_S^{-1}(\bfa x)Z_{S^\prime}^{-1}(\bfa x)\frac{\pta^2 F_1(\bfa x)}{\pta x_{\bfa\sigma}\pta x_{\bfa\tau}}-Z_S^{-2}(\bfa x)Z_{S^\prime}^{-1}(\bfa x)\frac{\pta^2 Z_S(\bfa x)}{\pta x_{\bfa\sigma}\pta x_{\bfa\tau}} F_1(\bfa x)-Z_S^{-1}(\bfa x)Z_{S^\prime}^{-2}(\bfa x)\frac{\pta^2 Z_{S^\prime}(\bfa x)}{\pta x_{\bfa\sigma}\pta x_{\bfa\tau}} F_1(\bfa x)\nnb\\
        &-Z_S^{-2}(\bfa x)Z_{S^\prime}^{-1}(\bfa x)\frac{\pta F_1(\bfa x)}{\pta x_{\bfa\sigma}}\frac{\pta Z_S(\bfa x)}{\pta x_{\bfa\tau}}-Z_{S^\prime}^{-2}(\bfa x)Z_{S}^{-1}(\bfa x)\frac{\pta F_1(\bfa x)}{\pta x_{\bfa\sigma}}\frac{\pta Z_{S^\prime}(\bfa x)}{\pta x_{\bfa\tau}}+Z_S^{-2}(\bfa x)Z_{S^\prime}^{-2}(\bfa x)\frac{\pta Z_S(\bfa x)}{\pta x_{\bfa\sigma}}\frac{\pta Z_{S^\prime}(\bfa x)}{\pta x_{\bfa\tau}} F_1(\bfa x)\nnb\\
        &-Z_S^{-2}(\bfa x)Z_{S^\prime}^{-1}(\bfa x)\frac{\pta F_1(\bfa x)}{\pta x_{\bfa\tau}}\frac{\pta Z_S(\bfa x)}{\pta x_{\bfa\sigma}}-Z_{S^\prime}^{-2}(\bfa x)Z_{S}^{-1}(\bfa x)\frac{\pta F_1(\bfa x)}{\pta x_{\bfa\tau}}\frac{\pta Z_{S^\prime}(\bfa x)}{\pta x_{\bfa\sigma}}+Z_S^{-2}(\bfa x)Z_{S^\prime}^{-2}(\bfa x)\frac{\pta Z_S(\bfa x)}{\pta x_{\bfa\tau}}\frac{\pta Z_{S^\prime}(\bfa x)}{\pta x_{\bfa\sigma}} F_1(\bfa x)\nnb\\
        &+2Z_{S^\prime}^{-3}(\bfa x)Z_{S}^{-1}(\bfa x)\frac{\pta Z_{S^\prime}(\bfa x)}{\pta x_{\bfa\sigma}}\frac{\pta Z_{S^\prime}(\bfa x)}{\pta x_{\bfa\tau}} F_1(\bfa x)+2Z_S^{-3}(\bfa x)Z_{S^\prime}^{-1}(\bfa x)\frac{\pta Z_S(\bfa x)}{\pta x_{\bfa\sigma}}\frac{\pta Z_S(\bfa x)}{\pta x_{\bfa\tau}} F_1(\bfa x)
    .\end{align}}
    Then we observe that
    \begin{align*}
        \frac{\pta Z_A(\bfa x)}{\pta x_{\bfa\sigma}}=\omega_{A}(\bfa\sigma)\exp(x_{\bfa\sigma}),\qquad\text{ for } A\in\{S,S^\prime\}
    \end{align*}
    which implies that the term involving \emph{only} the first order derivatives of $Z_A$s in \eqref{complicate} can be written as:
    \tny{\begin{align*}
        C(\bfa\sigma,\bfa\tau,\bfa x) :&=Z_{S^\prime}^{-1}(\bfa x)Z_{S}^{-1}(\bfa x)\exp\lef(x_{\bfa\sigma}+x_{\bfa\tau}\rig) F_1(\bfa x)\big(2\omega_{S^\prime}(\bfa\sigma)\omega_{S^\prime}(\bfa\tau)Z_{S^\prime}^{-2}(\bfa x)+2\omega_{S}(\bfa\sigma)\omega_{S}(\bfa\tau)Z_{S}^{-2}(\bfa x)\\
        &+(\omega_{S}(\bfa\sigma)\omega_{S^\prime}(\bfa\tau)+\omega_{S^\prime}(\bfa\sigma)\omega_{S}(\bfa \tau))Z_{S^\prime}^{-1}(\bfa x)Z_{S}^{-1}(\bfa x)\big) F_1(\bfa x)
    .\end{align*}}
    And consequently 
    \tny{\begin{align*}
\sum_{\bfa\sigma,\bfa\tau}&U(\bfa\sigma,\bfa\tau)C(\bfa\sigma,\bfa\tau,\bfa u(t))= 2Z_{S^\prime}^{-1}(\bfa u(t))Z_{S}^{-1}(\bfa u(t))\sum_{\bfa\sigma,\bfa\tau}U(\bfa\sigma,\bfa\tau)\omega_S(\bfa\sigma)\omega_{S^\prime}(\bfa\tau)\exp\lef(u_{\bfa\sigma}(t)+u_{\bfa\tau}(t)\rig)F(\bfa u(t))\\
&+\sum_{A\in\{S,S^\prime\}}2 Z_{A}^{-2}(\bfa u(t))\sum_{\bfa\sigma,\bfa\tau} U(\bfa\sigma,\bfa\tau)\omega_A(\bfa\sigma)\omega_A(\bfa\tau)\exp( u_{\bfa\sigma}(t)+u_{\bfa\tau}(t))F(\bfa u(t))\\
&=2\la U(\bfa\sigma^{3},\bfa\sigma^{4})f(\bfa\sigma^1,\bfa\sigma^2)\ra_{S,S^\prime}+2\la U(\bfa\sigma^{3},\bfa\sigma^{5})f(\bfa\sigma^1,\bfa\sigma^2)\ra_{S,S^\prime}+2\la U(\bfa\sigma^{4},\bfa\sigma^{6})f(\bfa\sigma^1,\bfa\sigma^2)\ra_{S,S^\prime}
.\end{align*}}
where we let $\bfa\sigma^{2k}\sim \bb P_S$ and $\bfa \sigma^{2k+1}\sim \bb P_{S^\prime}$ for all  $k\in\bb N$. And the rest of terms in \eqref{complicate} is analyzed similarly, note that
\tny{\begin{align*}
    \frac{\pta^2 F_1(\bfa x)}{\pta x_{\bfa\sigma}\pta x_{\bfa\tau}}=\sum_{\ell,\ell^\prime\leq 2}C_{\ell,\ell^\prime}(\bfa\sigma,\bfa\tau,\bfa x) = \sum_{\ell,\ell^\prime\leq 2}\sum_{\bfa\sigma^1,\bfa\sigma^2}\mbbm 1_{\bfa\sigma^{\ell}=\bfa\sigma}\mbbm 1_{\bfa\sigma^{\ell^\prime}=\bfa\tau}\omega_{S}(\bfa\sigma^1)\omega_{S^\prime}(\bfa\sigma^2)f(\bfa\sigma^1,\bfa\sigma^2)\exp\bigg(\sum_{i\in[2]} x_{\bfa\sigma^{i}}\bigg)
.\end{align*}}
Therefore the first term is written as
\begin{align*}
    \sum_{\ell,\ell^\prime\leq 2}\sum_{\bfa\sigma,\bfa\tau}U(\bfa\sigma,\bfa\tau)C_{\ell,\ell^\prime}(\bfa\sigma,\bfa\tau,\bfa x)=\sum_{\ell,\ell^\prime\leq 2}\la U(\bfa\sigma^{\ell},\bfa\sigma^{\ell^\prime})f(\bfa\sigma^1,\bfa\sigma^2) \ra_{S,S^\prime}
.\end{align*}
And the other terms can be analogously checked.
\tny{\begin{align*}
    \frac{\pta F_1(\bfa x)}{\pta x_{\bfa\sigma}}\frac{\pta Z_A(\bfa x)}{\pta x_{\bfa\tau}}=\sum_{\ell\in\{1,2\}}\sum_{\bfa\sigma^1,\bfa\sigma^2}\mbbm 1_{\bfa\sigma^\ell=\bfa\sigma}\omega_S(\bfa\sigma^1)\omega_{S^\prime}(\bfa\sigma^2)f(\bfa\sigma^1,\bfa\sigma^2)\exp\bigg(\sum_{i\in[2]}x_{\bfa\sigma^i}\bigg)\omega_A(\tau)\exp(x_{\bfa\tau})
.\end{align*}}
And consequently 
\begin{align*}
    \sum_{\bfa\sigma,\bfa\tau}U(\bfa\sigma,\bfa\tau)Z_S^{-2}Z_{S^\prime}^{-1}(\bfa u(t))\frac{\pta F_1(\bfa u(t))}{\pta x_{\bfa\sigma}}\frac{\pta Z_{S}(\bfa u(t))}{\pta x_{\bfa\tau}}&=\sum_{\ell\in\{1,2\}}\la U(\bfa\sigma^\ell,\bfa\sigma^3)f(\bfa\sigma^1,\bfa\sigma^2)\ra_{S,S^\prime,t}
.\end{align*}
Dropping terms concerning the second derivatives of $Z$ and given that $U(\bfa\sigma,\bfa\sigma)$ is a constant we arrive at
\begin{align*}
    \nu_t(f)^\prime &= 2\nu_t(U(\bfa\sigma^1,\bfa\sigma^{2})f(\bfa\sigma^1,\bfa\sigma^2))-4\nu_t(U(\bfa\sigma^{3},\bfa\sigma^{2})f(\bfa\sigma^1,\bfa\sigma^2))-4\nu_t(U(\bfa\sigma^{4},\bfa\sigma^{1})f(\bfa\sigma^1,\bfa\sigma^2))\\
    &+2\nu_t(U(\bfa\sigma^3,\bfa\sigma^4)f(\bfa\sigma^1,\bfa\sigma^2))+ 2\nu_t(U(\bfa\sigma^3,\bfa\sigma^5)f(\bfa\sigma^1,\bfa\sigma^2)) +2\nu_t(U(\bfa\sigma^4,\bfa\sigma^6)f(\bfa\sigma^1,\bfa\sigma^2))
.\end{align*}
However, observing the fact that the when $\bfa\sigma^1,\bfa\sigma^2$ are exchangeable in function $f$ and the fact that $S$ and $S^\prime$ are exchangeable upon taking the expectation of Gaussian part we conclude that
\begin{align*}
    \nu_t(f)^\prime&=2\nu_t(U(\bfa\sigma^1,\bfa\sigma^2)f(\bfa\sigma^1,\bfa\sigma^2))-8\nu_t(U(\bfa\sigma^3,\bfa\sigma^2)f(\bfa\sigma^1,\bfa\sigma^2))+2\nu_t(U(\bfa\sigma^3,\bfa\sigma^4)f(\bfa\sigma^1,\bfa\sigma^2))\\&+4\nu_t(U(\bfa\sigma^3,\bfa\sigma^5)f(\bfa\sigma^1,\bfa\sigma^2))
.\end{align*}
Here wereplace $U(\bfa\sigma,\bfa\tau)$ with \eqref{usigma1sigma2} and recalling that $S=[k]$ and $S=[r+1:k+r]$ we have
\ttny{\begin{align}\label{derivative-corr}
    \nu_t^\prime(f)&=\frac{(k+r)^2\theta^2}{n}\bigg(\nu_t((R^D_{1,2}-q^{(1)})^2f)-4\nu_t((R^D_{1,3}-q^{(1)})^2f)+\nu_t((R^D_{3,4}-q^{(1)})^2f)+2\nu_t((R^D_{3,5}-q^{(1)})^2f)\bigg)\nnb\\
    &+\frac{2\theta^2k(n-k)}{n}\bl\nu_t((R_{1,2}^D-q^{(1)})(R_{1,2}^{D^c}-q^{(2)})f)-4\nu_t((R_{1,3}^D-q^{(1)})(R_{1,3}^{D^c}-q^{(2)})f)\nnb\\
    &+\nu_t((R_{3,4}^D-q^{(1)})(R_{3,4}^{D^c}-q^{(2)})f)+2\nu_t((R_{3,5}^D-q^{(1)})(R_{3,5}^{D^c}-q^{(2)})f)\br
.\end{align}}
To deal with the second term, assume that $\nu_0(f^{1+\delta})<\infty$ for some $\delta>0$, then we use the result of $t=0$ in Lemma \ref{zeropointrpk} and Lemma \ref{cltconsk} to get that the cross terms satisfies (here we take $R_{1,2}$ as an example and other terms can be analogously derived)
\begin{align*}
    \frac{2\theta^2k(n-k)}{n}&\nu_t((R_{1,2}^D-q^{(1)})(R_{1,2}^{D^c}-q^{(2)})f)\\
    &\leq \frac{Ck^2(n-k)^2}{n^2}\nu_0((R_{1,2}^D-q^{(1)})^2(R_{1,2}^{D^c}-q^{(2)})^2f)\bl 1+O\bl\sqrt{\frac{k}{n}}\br\br\\
    &=O\bl\frac{k}{n}\br\nu_0(f^{1+\delta})^{\frac{1}{1+\delta}}=O\bl\frac{k}{n}\br,
\end{align*} when the temperature regime is high/low.
And when the temperature regime is critical we have when $c=o(k^{-\frac{\tau-2}{2\tau-1}})$
\begin{align*}
    \frac{2\theta^2k(n-k)}{n}\nu_t((R_{1,2}^D-q^{(1)})(R_{1,2}^{D^c}-q^{(2)})f)=O\bl\frac{k^{\frac{2\tau}{2\tau-1}}}{n}\br.
\end{align*}
And when $c=\Omega(k^{-\frac{\tau-2}{2\tau-1}})$, we have
\begin{align*}
    \frac{2\theta^2k(n-k)}{n}\nu_t((R_{1,2}^D-q^{(1)})(R_{1,2}^{D^c}-q^{(2)})f)=O\bl\frac{k^{\frac{4\tau-4}{2\tau-1}}}{n}\br.
\end{align*}
And wearrive at (Noticing that $f$ is positive and the negative terms can be dropped)
\tny{\begin{align}\label{derivativehetero}
    \nu_t^\prime(f)\leq\frac{C(k+r)^2}{n}\bigg(\nu_t((R^D_{1,2}-q^{(1)})^2f)+\nu_t((R^D_{3,4}-q^{(1)})^2f)+2\nu_t((R^D_{3,5}-q^{(1)})^2f)\bigg)+O\bl\frac{k^{2-\beta}}{n}\br.
\end{align}}
where $\beta$ is taking different values for different temperature regimes.
And then our next step is to construct a special $f$ such that the above inequality can be addressed. Here we can replace $q^{(1)}$ with $q$ as their difference is sufficiently small. A natural choice is to let $f$ be in the form of $\exp((k+r)^{\beta_1}(\lambda-\alpha t)(R_{1,2}^D-q)^2)$ where $\beta_1,\alpha$ are to be specified later. However, we observe that R.H.S. has two terms of different nature (1) $R_{3,4}$ is the overlap between two replicas coming from  $\bb P_S(\cdot|\bfa g,\bfa h) ,\bb P_{S^\prime}(\cdot|\bfa g,\bfa h)$ respectively. (2) $R_{3,5}$ is the overlap between two replicas coming from the same $\bb P_{S}(\cdot|\bfa g,\bfa h)$. Hence our strategy is the two step procedure: $\emph{(I)}$: We derive the concentration for the case of $S=S^\prime$. $\emph{(II)}$: We derive the concentration for the common case with $S\neq S^\prime$ using the result obtained in the first step.

\begin{center}
    \textbf{Step I: $S=S^\prime$}
\end{center}
In step I weutilize the following lemma 
\begin{lemma}\label{interpolate-corr1}
    Consider any number $\lambda>0$. Then for all $\beta>0$, $D=[k+r]$, we have
    \begin{align*}
        \nu_t((R^D_{3,4}-q)^2\exp\lambda (k+r)^{\beta}(R^D_{1,2}-q)^2)\leq\nu_t((R^D_{1,2}-q)^2\exp\lambda (k+r)^{\beta}(R^D_{1,2}-q)^2)
    .\end{align*}
\end{lemma}
Using the above lemma, wepick $\alpha=\frac{Ck^{2-\beta}}{n}$ (Recall that when $S=S^\prime=D$, $r=0$) and $\beta$ according to the different temperature regimes as stated in Lemma \ref{zeropointrpk}. Then, when $\alpha\leq C$ for some $C>0$ we have
\begin{align}\label{derivativelessthan0}
    \frac{d}{dt}\nu_t\lef(\exp\lef(\lambda-\alpha t\rig)k^{\beta}(R^S_{1,2}-q)^2\rig)&=\nu^\prime_t\lef(\exp\lef(\lambda-\alpha t \rig)k^{\beta}(R^S_{1,2}-q)^2\rig)\nnb\\
    &- \alpha k^{\beta}\nu_t\lef((R^S_{1,2}-q)^2\exp\lef(\lambda-\alpha t\rig)k^{\beta}(R^S_{1,2}-q)^2\rig)<0
.\end{align}
Henceforth, when $k\leq Cn^{\frac{1}{2-\beta}}$, there exists $C$ such that for all $\lambda<C$:
\begin{align*}
    \nu_1(\exp(\lambda-\alpha)k^{\beta}(R_{1,2}^S-q)^2)\leq\nu_0(\lambda k^{\beta}(R_{1,2}^S-q)^2)<\infty.
\end{align*}
Then we can extend the results from $R_{1,2}^S$ to $R_{1,2}^D$ with $D=[k+r]$, using the following lemma.
\begin{lemma}[\citep{talagrand2010mean}]\label{secondterm}
For $i\in[k+1:k+r]$ when $S=[k]$, pick $q$ satisfies $q=\bb E[\tanh^2(\theta z\sqrt {q} +h)]$. We have when $k$ is sufficiently large for all $u\in\bb R$:
    $\nu_0\lef(\exp u\lef(\sigma_i^1\sigma_i^2-q\rig)\rig)\leq\exp\lef(\frac{u^2}{2}\rig)$.
\end{lemma}
First wedenote $R^{D\setminus S}_{1,2}:=\frac{1}{r}\sum_{i\in D\setminus S}\sigma_i^1\sigma_i^2$. It is checked by the above lemma, \eqref{gstrick}, \eqref{gstrick1}, and \eqref{gstrick2} we have there exists $C>0$ such that for all $\lambda<C$ we have
\begin{align*}
    \nu_0(\lambda(k+r)^{\beta}(R_{1,2}^D-q)^2)<\infty,
\end{align*}
where $\beta$ is different across temperature but identical to the exponential bound of $R_{1,2}^S-q$. Then we proceed similarly as the argument given by \eqref{derivativelessthan0}, there exists $C>0$ such that when $k\leq C_1n^{\frac{1}{2-\beta}}$ for some $C_1>0$ and $\lambda<C$ we have
\begin{align}\label{uniformr35}
    \nu_t(\lambda(k+r)^\beta(R_{1,2}^D-q)^2)<\infty,\qquad\text{ uniformly over all }t\in[0,1].
\end{align}
This further leads to the following lemma.
\begin{lemma}\label{lmoverlappedreplicacon}
    Let $q$ be the solution to $q_1=\bb E[\tanh^2(\theta\sqrt {q_1}z+\theta_1\mu+h)]$ for $\mu=\bb E[\tanh(\theta\sqrt{q_1}z+h)]$ and $q_2=\bb E[\tanh^2(\theta\sqrt {q_2}x+h)]$. Let $|D|=k+r$ and $r=ck$, wethen have
    \begin{enumerate}
        \item At the high temperature, $\Vert k^{1/2}(R_{1,2}^D-q_2) \Vert_{\psi_2}<\infty$.
        \item At the low temperature, $\Vert k^{1/2}(R_{1,2}^D-(1-c)q_1-cq_2)\Vert_{\psi_2}<\infty$.
        \item At the critical temperature, $\Vert k^{\frac{\tau-1}{2\tau-1}}(R_{1,2}^D-q)\Vert_{\psi_2}<\infty$.
    \end{enumerate}
\end{lemma}
\begin{center}
    \textbf{Step II: $S\neq S^\prime$}
\end{center}
Then we can start our analysis over the $S\neq S^\prime$ case. Here wediscuss only the high and critical temperature case as at the low temperature case $R_{1,2}^D$ and $R_{1,3}^D$ does not converge to the same limit. We go back to \eqref{derivativehetero} and use lemma \ref{interpolate-corr1}, \eqref{uniformr35} to get
\begin{align*}
    \nu_t^\prime(\exp(&\lambda(k+r)^{\beta}(R^D_{1,2}-q)^2))\leq \frac{Ck^2}{n}\bl\nu_t((R_{1,2}-q)^2\exp(\lambda(k+r)^{\beta}(R^D_{1,2}-q)^2))\\
    &+\nu_t((R^D_{3,5}-q)^2\exp(\lambda(k+r)^{\beta}(R^D_{1,2}-q)^2))\br+O\lef(\frac{k^{2-\beta}}{n}\rig)\\
    &\leq \frac{Ck^2}{n}\bl\nu_t((R_{1,2}-q)^2\exp(\lambda(k+r)^{\beta}(R_{1,2}^D-q)^2)\br+O\bl\frac{k^{2-\beta}}{n}\br.
\end{align*}
Then wecheck that through a similar argument as \eqref{derivativelessthan0} we have when $k\leq Cn^{\frac{1}{2-\beta}}$, there exists $C$ such that for all $\lambda<C$ we have
\begin{align*}
    \nu_1(\exp(\lambda (k+r)^{\beta}(R_{1,2}^D-q)^2))<\infty.
\end{align*}
And we finish the proof by replacing $\beta$ with their corresponding value for the different temperature regimes.

In the last, we give the following results for the concentration at the $t=0$ end
\begin{lemma}\label{zeropointrpk}
Let $q$ be the solution to $q=\bb E[\tanh^2(h+\theta\sqrt{q}z)]$, then for the pRFCW model with the random field $h^\prime=h+\theta\sqrt qz$ i.i.d. we have
\begin{enumerate}
    \item In the high temperature regime, there exists $C>0$ such that for all $\lambda<C$,\\
$     \nu_0\lef(\exp\lambda (k+r)(R^D_{1,2}-q)^2\rig)<\infty
$.
    \item In the low temperature regime, let $c:=\frac{r}{k}$ and $q=\frac{1-c}{1+c}\bb E[\tanh^2(\sqrt{\theta_1}x_1^*+h+\theta\sqrt{q}z)]+\frac{2c}{1+c}\bb E[\tanh(\sqrt{\theta_1}x_1^*+h+\theta\sqrt{q}z)\tanh(h+\theta\sqrt{q}z)].$ Then, there exists $C>0$ such that for all $\lambda<C$ we have
$    \nu_0(\exp\lambda(k+r)(R_{1,2}^D-q)^2)<\infty.
$    
\item At the critical temperature regime with flatness $\tau\in\bb N\setminus\{1\}$, assume that $q=\bb E[\tanh^2(\theta\sqrt qz+h)]$. When $c=o(k^{-\frac{\tau-2}{2\tau-1}})$ there exists $C>0$ such that for all $\lambda<C$,
$    \nu_0(\exp\lambda (k+r)^{\frac{2\tau-2}{2\tau-1}}(R_{1,2}^D-q)^2)<\infty.
$    
\item When $c=\Omega(k^{-\frac{\tau-2}{2\tau-1}})$, there exists $C>0$ such that for all $\lambda<C$,
$    \nu_0(\exp\lambda (k+r)^{\frac{2}{2\tau-1}}(R_{1,2}^D-q)^2)<\infty.
$\end{enumerate}
\begin{remark}
    An interesting phenomenon is that at the critical temperature, weobserve a two phase behavior on the convergence rate of the local replica. Furthermore, this rate will be weaker than the $\frac{1}{\sqrt k}$ that is commonly observed in the standard SK model. This rate is closely connected to the unknown critical temperature region that we know nothing about. We elaborate more on this fact in section \ref{sect7}.
\end{remark}
\end{lemma}

\subsection{Proof of Lemma \ref{thm4.1}}
The proof of lemma \ref{thm4.1} is separately organized according to the different temperature regimes.
\subsubsection{Proof of the high temperature regime}
    First we consider the null. Using Lemma \ref{cltconsk} we have
    \begin{align*}
        \bb E\bigg[\frac{1}{n}\bl\sum_{i=1}^n\sigma_i\br^2\bigg]\to 1,\qquad \bigg\Vert\frac{1}{\sqrt n}\sum_{i=1}^n\sigma_i\bigg\Vert_{\psi_2}<\infty\thinspace\Rightarrow\thinspace\bigg\Vert\frac{1}{n}\bl\sum_{i=1}^n\sigma_i\br^2\bigg\Vert_{\psi_1}<\infty.
    \end{align*}
    Therefore one will have by Bernstein inequality,
    \begin{align*}
        \bb P_0(\phi_7-\bb E[\phi_7]>t)=\bb P_0\bl\frac{1}{n}\bl\sum_{i=1}^n\sigma_i\br^2-1\geq t \br\leq\exp(-Cmt^2\wedge mt).
    \end{align*}
    Consider the alternative. Using Lemma \ref{cltatline} we have
    \begin{align*}
        \frac{1}{\sqrt n}\sum_{i=1}^n\sigma_i=\ub{\frac{\sqrt{k}}{\sqrt n}\frac{1}{\sqrt k}\sum_{i\in S}\sigma_i}_{T_1}+\ub{\frac{\sqrt{n-k}}{\sqrt n}\frac{1}{\sqrt{n-k}}\sum_{i\in S^c}\sigma_i}_{T_2}.
    \end{align*}
    wethen notice that by Lemma \ref{covariance} we have $\bb E[T_1T_2]=O\lef(\frac{1}{\sqrt n}\rig)$ and here we can compute $\ca V_m^h=\frac{1-\theta_1(\bb E[\sech^2(\theta\sqrt qz+h+\theta_1\mu)])^2}{(1-\theta_1\bb E[\sech^2(\theta\sqrt qz+h+\theta_1\mu)])^2}$ to get
    \begin{align*}
        \bb E[(T_1+T_2)^2]&=\bb E[T_1^2]+\bb E[T_2^2]+2\bb E[T_1T_2]\\
        &=c\frac{1-\theta_1(\bb E[\sech^2(\theta\sqrt qz+h+\theta_1\mu)])^2}{(1-\theta_1\bb E[\sech^2(\theta\sqrt qz+h+\theta_1\mu)])^2}+(1-c)+O\bl\frac{1}{\sqrt n}\br.
    \end{align*}
    And moreover weuse the sub-additivity of Orlicz norm and theorem \ref{largcliquetailbound} to get
    \begin{align*}
        \lef\Vert T_1+T_2\rig\Vert_{\psi_2}\leq \Vert T_1\Vert_{\psi_2}+\Vert T_2\Vert_{\psi_2}<\infty\thinspace\Rightarrow\thinspace\Vert (T_1+T_2)^2\Vert_{\psi_1}<\infty.
    \end{align*}
    Then one will use Bernstein inequality again to get
    \begin{align*}
        \bb P_S(\phi_4-\bb E[\phi_4]<t)&=\bb P_S\bl\frac{1}{nm}\sum_{j=1}^m\bl\sum_{i=1}^n\sigma^{(j)}_i\br^2 -\bb E[(T_1+T_2)^2]\leq -t\br\\
        &\leq\exp(-Cmt^2\wedge mt).
    \end{align*}
    And we pick large $m\asymp 1$ to satisfy the criteria of asymptotically powerful tests.
\subsubsection{Proof of the low temperature regime}
    The proof goes by noticing that under the null, using Lemma \ref{cltconsk} we have
    \begin{align*}
        \frac{1}{n}\bb E\bigg[\bigg|\sum_{i=1}^n\sigma_i\bigg|\bigg]\leq\frac{1}{n}\bb E\bigg[\bl\sum_{i=1}^n\sigma_i\br^2\bigg]^{1/2}=o(1).
    \end{align*}
    And then wealso have by the Chebyshev inequality
    \begin{align*}
        \bb P(\phi_8\geq t)\leq\frac{\bb E[\phi_8^2]}{t^2}=O\bl\frac{1}{nt^2}\br.
    \end{align*}
    Then we consider the alternative and notice that for $\mu\in\ca U$ being the positive optimality point of \eqref{meanfieldeq}, we have
    \tny{\begin{align*}
        \bb E[T_1]-\bb E[T_2]=\bb E\bigg[\frac{1}{n}\bigg|\sum_{i=1}^k\sigma_i+\sum_{i=k+1}^n\sigma_i\bigg|\bigg]=\bb E\bigg[\bigg|\frac{1}{n}\sum_{i=1}^k\sigma_i\bigg|\bigg]+O\bl\bb E\bigg[\bigg|\frac{1}{n}\sum_{i=k+1}^n\sigma_i\bigg|\bigg]\br=\frac{n}{k}\mu+o(1).
    \end{align*}}
    And then weanalyze the two parts $T_1:=|\frac{1}{n}\sum_{i=1}^k\sigma_i|$ and $T_2:=|\frac{1}{n}\sum_{i=k+1}^{n}\sigma_i|$ separately. Using theorem \ref{largcliquetailbound}, we have $\Vert \sqrt n(T_1-\bb E[T_1])\Vert_{\psi_2}<\infty$ and $\Vert\sqrt nT_2\Vert_{\psi_2}<\infty$. This further leads to $\Vert\sqrt n(T_2-\bb E[T_2]) \Vert_{\psi_2}<\infty$ and we have for all $t>0$:
    \begin{align*}
        \bb P(|T_1-\bb E[|T_1|]|\geq t)\leq 2\exp(-Cnt^2),\qquad \bb P(|T_2-\bb E[T_2]|\geq t)\leq 2\exp(-Cnt^2).
    \end{align*}
    Then we further notice that $ \phi_8\geq T_1-T_2$ and we have for all $t>0$,
    \begin{align*}
        \bb P(\phi_8-\bb E[T_1]+\bb E[T_2]\leq -t)&\leq\exp(T_1-\bb E[T_1]+\bb E[T_2]-T_2\leq -t)\\
        &\leq \exp\bl T_1- \bb E[T_1]<-\frac{t}{2}\br+\exp\bl \bb E[T_2]-T_2\leq-\frac{t}{2}\br\\
        &\leq C\exp(-Cnt^2)=o(1).
    \end{align*}
    Therefore we complete the proof that the test given by \ref{alg:eight} is asymptotically powerful.
\subsection{Proof of Lemma \ref{mfdeq}}
   The proof of lemma \ref{mfdeq} involves a few technical steps: \textbf{(1)} The $1$-step replica symmetry breaking bound on the free energy. \textbf{(2)} An alternative description of the replica symmetry regime. \textbf{(3)} Bounds on the single copy and coupled copies, which complete the proof when used together with the Gaussian concentration inequality. The proof idea of (1) and (2) is a comparison strategy, where one can derive two upper bounds (The two smart paths are the Poisson Dirichlet process versus the ordinary one appears in \ref{sect3}.) for the same free energy of the pSK model and compare them to get that at the replica symmetry region they match with each other. A similar strategy also helps us quantitatively prove the part {(3)}, despite requiring a few more algebraic arguments. And the results of step {(3)} finally give rise to the following lemma.
   \begin{lemma}[Large Deviation]
When the condition in definition \ref{withintheATline} is satisfied. Define the replica overlap to be
$    R_{1,2}:=\frac{1}{n}\sum_{i\leq n}\sigma_i^1\sigma_i^2$ and $m_1:=\frac{1}{k}\sum_{i=1}^k\sigma_i$. Define $\ca U:=\{\mu:\mu\text{ is a maximum of \eqref{meanfieldeq}}\}$. Then, we have for all $\epsilon>0$,
$    \bb P\lef(|R_{1,2}-q|\geq \epsilon\rig)\leq C(\epsilon)\exp\lef(-C(\epsilon)n\rig),
$when $|\ca U|=1$, we have
$   \bb P(|m_1|\geq\epsilon)\leq C(\epsilon)\exp\lef(-C(\epsilon)n\rig), 
$. Otherwise we have for $\mu$ being the positive element in $\ca U$,
$    \bb P(||m_1|-\mu|\geq\epsilon)\leq C(\epsilon)\exp(-C(\epsilon)n),
$ with $C(\epsilon)$ not dependent on $n,k, u,t$.
\end{lemma}

Here we introduce the following notation for the free energy, 
\begin{align}\label{pkn}
    p_{k,n}(\theta_1,\theta) =\frac{1}{n}\bb E[\log Z^{pSK}_{\theta_1,\theta}]=\frac{1}{n}\bb E\bigg[\log\sum_{\bfa\sigma}\exp\lef(-H^{pSK}_{\theta_1,\theta}(\bfa\sigma)\rig)\bigg].
\end{align}
Here we introduce the upper bounds of 1-step Replica Symmetry Breaking, which is proved using the interpolation given by the Poisson Dirichlet Process. This results is important as it gives a tighter upper bound for the limiting free energy. weshow that this tighter upper bound is crucial in proving that the averaging statistics concentrates. 
\begin{lemma}[1RSB Bound]\label{1rsbbound} Consider $0\leq  q,q^\prime\leq 1$ and $0<m<1$. Set
\begin{align}\label{yprime}
    Y_1:=\theta\sqrt{q}z+\theta \sqrt{q^\prime-q}z^\prime +\theta_1\mu +h,\thinspace\text{ and }\thinspace Y_2:=\theta\sqrt{q}z+\theta\sqrt{q^\prime-q}z^\prime+h,
\end{align}
where $z,z^\prime$ are independent standard Gaussian r.v.s. independent of $h$. Denote $\bb E^\prime$ to be the expectation w.r.t. $z^\prime$ only.  we have the following upper bound for $p_{k,n}(\theta_1,\theta)$:
\sm{\begin{align*}
    p_{k,n}(\theta_1,\theta)&\leq\log 2+\frac{\theta^2}{4}(1-q^\prime)^2-\frac{\theta^2}{4}m(q^{\prime 2}-q^2)+c\sup_{\mu\in[-1,1]}\lef(\frac{1}{m}\bb E\log\bb E^\prime\cosh^m Y_1 -\frac{\theta_1\mu^2}{2}\rig)\\
    &+\frac{\log(k)}{n}+\frac{1-c}{m}\bb E\log\bb E^\prime\cosh^mY_2.
\end{align*}}
\end{lemma}
Then, we introduce the notation
\sm{\begin{align}\label{phidef}
    \Phi(m,q^\prime):&=\log 2+\frac{\theta^2}{4}(1-q^\prime)^2-\frac{\theta^2}{4}m(q^{\prime 2}-q^2)+\sup_{\mu\in[-1,1]}\bl\frac{k}{mn}\bb E\log \bb E^\prime\cosh^mY_1-\frac{\theta_1\mu^2}{2}\br\nnb\\
    &+\frac{n-k}{mn}\bb E\log\bb E^\prime\cosh^mY_2,\\
    \ca P_1(\theta_1,\theta):&=\inf_{m,q^\prime}\Phi(m,q^\prime).\nnb
\end{align}}
And weimmediately see that for sufficiently large $k$ we have (noticing that $pSK(\theta_1,\theta)=\Phi(1,q)$) and weimmediately have
\begin{align*}
    \ca P_1(\theta_1,\theta)\leq pSK(\theta_1,\theta).
\end{align*}

Essentially lemma \ref{1rsbbound} implies that for large $k$ we have
\begin{align*}
    \lim_{k\to\infty} p_{n,k}(\theta_1,\theta)\leq\ca P_1(\theta_1,\theta) \leq pSK(\theta_1, \theta).
\end{align*}
Therefore, when $\ca P_1(\theta_1,\theta)<pSK(\theta_1,\theta)$, we are out of the replica symmetry phase. In \citep{talagrand2011mean2}, the following lemma alternatively states the exact temperature regime for the standard SK model.
\begin{lemma}\label{alternativersb1}
    The replica symmetry regime of planted SK model is the region where $pSK(\theta_1,\theta)=\ca P_1(\theta_1,\theta)$. And the region defined in definition \ref{withintheATline} satisfies it.
\end{lemma}
A conjecture for the replica symmetric regime of the SK model is stated in \citep{talagrand2011mean2} Chapter 13. Alternatively, the following will be the corresponding conjecture for the planted SK model
\begin{conjecture}
    The replica symmetric regime of planted SK model is the region where
    \begin{align*}
        \bb E\theta^2((1-c)\sech^4(\theta\sqrt q z+h)+c\sech^4(\theta\sqrt qz+\theta_1\mu+h))<1,
\end{align*}
\end{conjecture}
This conjecture comes from \citep{toninelli2002almeida} who proved that replica symmetry is broken when the $\theta$ satisfies
\begin{align*}
        \bb E\theta^2((1-c)\sech^4(\theta\sqrt q z+h)+c\sech^4(\theta\sqrt qz+\theta_1\mu+h))>1.
\end{align*}

To prove the concentration, a center argument is an upper bound on the coupled copies. Before we state the results, we introduce the following Hamiltonian interpolating path. One will notice that at $t=1$ this interpolating Hamiltonian corresponds to the pSK measure.
\begin{align*}
\mca H_t(\bfa\sigma)&:=\frac{\theta\sqrt t}{\sqrt n}\sum_{i<j\leq n}g_{ij}\sigma_i\sigma_j+\theta\sqrt{1-t}\sum_{i\leq n}\sigma_iz_i\sqrt q+\frac{\theta_1k}{2}m^2+\sum_{i\leq n}h_i\sigma_i,\nnb\\
    \psi(t,u)&:=\frac{1}{n}
        \bb E\log\sum_{\bfa\sigma^1,\bfa\sigma^2}\mbbm 1_{R_{1,2}=u}\exp\lef(-\mca H_t(\bfa\sigma^1)-\mca H_t(\bfa\sigma^2)\rig),\\
        \eta(t,u):&=\frac{1}{n}\bb E\log\sum_{\bfa\sigma}\mbbm 1_{m=u}\exp\lef(-\mca H_t(\bfa\sigma)\rig)
\end{align*}
Then we introduce the following two quantities, the first one is the coupled copies which is the partial free energy given $R_{1,2}=u$. The second one is the partial free energy given $m=u$. These two quantities play the crucial parts in proving the concentration arguments.
\begin{align*}
    r_n(\theta_1,\theta,h,u):&=\frac{1}{k}\bb E\log\sum_{\bfa\sigma}\mbbm 1_{R_{1,2}=u}\exp\lef(-\mca H_1(\bfa\sigma^1)-\mca H_1(\bfa\sigma^2)\rig),\\
    t_{k}(\theta_1,\theta,h,u):&=\frac{1}{k}\bb E\log\sum_{\bfa\sigma}\mbbm 1_{m=u}\exp\lef(-\mca H_1(\bfa\sigma^1)\rig).
\end{align*}
The underlying logic is that when we obtain a \emph{proper} upper bound for the above two quantities and a \emph{proper} lower bound for $p_{k,n}$ which is given by lemma \ref{alternativersb1}, one will notice that by Gaussian concentration inequality, the following holds
\begin{align}\label{heuristicsub}
    \lef|\frac{1}{n}\log\la \mbbm 1_{R_{1,2}=u}\ra-(\psi(1,u)-2p_{k,n})\rig|&=O_p\lef(\frac{1}{\sqrt n}\rig),\nnb\\ 
    \lef|\frac{1}{n}\log\la \mbbm 1_{m=u}\ra-(\eta(1,u)-p_{k,n})\rig|&=O_p\lef(\frac{1}{\sqrt n}\rig).
\end{align}
Then applying the truncation argument, one will get the upper bound for $\bb P(R_{1,2}=u)$ and $\bb P(m=u)$. Then in the following, we achieve the two goals separately. The proof of lemma \ref{alternativersb1} utilized lemma \ref{lm4.3}, which is a direct result of lemma \ref{intermediatelmpsk}.And lemma \ref{intermediatelmpsk} is a direct result of lemma \ref{rkpkub}. Finally, we give lemma \ref{exponentialeq} as the consequence of all of them.
\begin{lemma}\label{lm4.3}
    Assume that $pSK(\theta_1,\theta) = \ca P_1(\theta_1,\theta)$. Define $m:=\frac{1}{k}\sum_{i=1}^k\sigma_i$ and $\tde m:=\frac{1}{n-k}\sum_{i=k+1}^n\sigma_i$. The function defined by
    \begin{align}\label{definitions}
        \psi(t)&:=\log 2+\frac{\theta^2}{4}t(1-q^2)+\sup_{\mu\in[-1,1]}\lef\{c\bb E\lef[\log\cosh\lef(\theta z \sqrt q+\theta_1\mu+h\rig)\rig]-\frac{t\theta_1\mu^2}{2}\rig\}\nnb\\
        &+(1-c)\bb E[\log\cosh(\theta z\sqrt q+h)].
    \end{align}
    satisfies that for all $t\leq t_0$ for some $t_0<1$ there exists constant $C$ independent of $n$ such that
    \sm{\begin{align*}
        \psi(t,u)&\leq 2\psi(t)-C(u-q)^2,\thinspace
        \eta(t,u)\leq \psi(t)-C\begin{cases}
        \min_{\mu\in\ca U}\{(u-\mu)^2\}&\text{ at high/low temperature }\\
        u^{2\tau}&\text{ at critical temperature}
        \end{cases}.
    \end{align*}}
\end{lemma}
The above lemma is a direct result of the next lemma that we present, where some change of variable will immediately lead to the conclusion.
\begin{lemma}\label{intermediatelmpsk}
    Assume that the condition in definition \ref{withintheATline} holds.
    Then there exists $C>0$ such that the following holds
    \begin{align}\label{bddcopiesub}
         r_n(\theta_1,\theta,h, u)&\leq 2pSK(\theta_1,\theta)-C(u-q)^2 +\frac{2\log k}{n},\\
        t_{k}(\theta_1,\theta,h,u)&\leq pSK(\theta_1,\theta,h)-C\begin{cases}
        \min_{\mu\in\ca U}\{(u-\mu)^2\}&\text{ at high/low temperature }\\
        u^{2\tau}&\text{ at critical temperature}
        \end{cases}.
    \end{align}
\end{lemma}
To prove the above lemma, one will have to use another Poisson Dirichlet process approximation. The result of this approximation, is given by the following lemma. And weuse the detailed discussion of Proposition 13.6.6 in \citep{talagrand2011mean2} (which involves some algebraic manipulations) to complete the proof of the above technical lemma.
    \begin{lemma}\label{rkpkub}
        Assume that there exist standard normal r.v.s. $z_1,z_2,z_1^\prime,z_2^\prime$ such that $(z_1,z_2)\perp(z_1^\prime,z_2^\prime)$. We define
        \begin{align*}
            c_0=\bb E[z_1z_2],\thinspace c_1=\bb E[z_1^\prime z_2^\prime].
        \end{align*}
        Then for $0\leq q_1\leq q_2\leq 1$, assume that $u=c_0q_1+c_1(q_2-q_1)$. For $j\in\{1,2\}$, we set for $u\in[-1,1]$:
        \begin{align*}
            Y_j(u)&=\theta z_j\sqrt{q_1}+\theta z_j^\prime\sqrt{q_2-q_1}+\theta_1u+h.
        \end{align*} Denote $\bb E^\prime$ as expectation w.r.t. $(z_1^\prime,z_2^\prime)$. Then if $0\leq m\leq 1$, for any $\lambda$,
        \tny{\begin{align*}
            r_n(&\theta_1,\theta,h,u)\leq\frac{2\log k}{n}+2\log2+\frac{\theta^2}{2}((1-q_2)^2-m(q_2^2-q_1^2)-m(u^2-c_0^2q_1^2))\\
            &+c\sup_{\mu\in[-1,1]}\lef\{\frac{1}{m}\bb E\log\bb E^\prime(\cosh Y_1(\mu)\cosh Y_2(\mu)\cosh\lambda+\sinh Y_1(\mu)\sinh Y_2(\mu)\sinh \lambda)^m-\theta_1\mu^2\rig\}-\lambda u\\
            &+\frac{1-c}{m}\bb E\log\bb E^\prime(\cosh Y_1(0)\cosh Y_2(0)\cosh\lambda+\sinh Y_1(0)\sinh Y_2(0)\sinh \lambda)^m,\\
             t_n(&\theta_1,\theta,h,u)\leq \frac{\theta^2}{4}((1-q_2)^2-m(q_2^2-q_1^2))+\log 2-\frac{\theta_1u^2}{2}+\frac{c}{m}\bb E\log\bb E^\prime\cosh^m(Y_i(u))\\
            &+\frac{1-c}{m}\bb E\log\bb E^\prime\cosh^m(Y_i(0)).
        \end{align*}}
    \end{lemma}
    And finally, wearrive at the desired results, which served as a key component to the limiting theorems. The following lemma is directly proved using the heuristic discussion given in \eqref{heuristicsub}.
\begin{lemma}[Large Deviation]\label{exponentialeq}
When the condition in definition \ref{withintheATline} is satisfied. Define the replica overlap to be
$    R_{1,2}:=\frac{1}{n}\sum_{i\leq n}\sigma_i^1\sigma_i^2$ and $m_1:=\frac{1}{k}\sum_{i=1}^k\sigma_i$. Define $\ca U:=\{\mu:\mu\text{ is a maximum of \eqref{meanfieldeq}}\}$. Then, we have for all $\epsilon>0$,
\begin{align*}
    \bb P\lef(|R_{1,2}-q|\geq \epsilon\rig)\leq C(\epsilon)\exp\lef(-C(\epsilon)n\rig),
\end{align*}
when $|\ca U|=1$, we have
\begin{align*}
   \bb P(|m_1|\geq\epsilon)\leq C(\epsilon)\exp\lef(-C(\epsilon)n\rig), 
\end{align*}
and otherwise we have for $\mu$ being the positive element in $\ca U$.
\begin{align*}
    \bb P(||m_1|-\mu|\geq\epsilon)\leq C(\epsilon)\exp(-C(\epsilon)n),
\end{align*}
with $C(\epsilon)$ not dependent on $n,k, u,t$.
\end{lemma}
\subsection{Proof of Lemma \ref{largcliquetailbound}}
We first state the following lemma giving the iterations of moments.
\begin{lemma}\label{inductivehypo}
    If it is true that for $ r\in\bb N$ there exists $C>$ such that $\forall \ell\leq r$
    \begin{align*}
        \min_{\mu\in\ca U}\{\nu((m-\mu)^{2\ell})\}\leq\lef(\frac{C\ell}{k}\rig)^\ell.
    \end{align*}
    Recall the notation $m^-=\frac{1}{k}\sum_{i\leq k-1}\sigma_i$. Then we have $\forall j\leq 2r$,
    \begin{align*}
        \min_{\mu\in\ca U}\{\nu(|m-\mu|^j)\}\leq\lef(\frac{C(j+1)}{2k}\rig)^{j/2},\thinspace \min_{\mu\in\ca U}\{\nu((m^--\mu)^{2r})\}\leq 3\lef(\frac{C(r+1)}{k}\rig)^r.
    \end{align*}
\end{lemma}
    We prove by induction and assume that for all $r^\prime\leq r$ there exists universal $C>0$ such that
    \begin{align*}
        \nu((R_{1,2}-q)^{2r^\prime})\vee\nu((m_1-\mu)^{2r^\prime})\leq C^r\lef(\frac{r^\prime+1}{k}\rig)^{r^\prime}.
    \end{align*} 
    The above trivially holds when $r=0$.
    Then, by Lemma \ref{secondmoment} and lemma \ref{secondmomentpositive} the following holds for some $C_1>0$:
    \begin{align}\label{momentiter}
        \nu((R_{1,2}-q)^{2r+2})&\leq \frac{Cr}{n}\nu((R_{1,2}-q)^{2r})+C_1\nu((R_{1,2}-q)^{2r+3})\vee\nu((|m_1|-\mu)^{2r+3}),\nnb\\
        \nu((|m_1|-\mu)^{2r+2})&\leq\frac{Cr}{k}\nu((|m_1|-\mu)^{2r})+C_1\nu((R_{1,2}-q)^{2r+3})\vee\nu((|m_1|-\mu)^{2r+3}),\nnb\\
        \nu((\tde m_1)^{2r+2})&\leq\frac{Cr}{n-k}\nu((\tde m_1)^{2r})+C_1\nu((R_{1,2}-q)^{2r+3}).
    \end{align}
    Note that by Lemma \ref{exponentialeq} it is easily checked that we have for all $\epsilon>0$,
    \begin{align*}
        \nu(|R_{1,2}-q|^{2r+3})&\leq \epsilon\nu(|R_{1,2}-q|^{2r+2})+2^{2r+3}\bb P(|R_{1,2}-q|\geq\epsilon)\\
        &\leq\epsilon\nu(|R_{1,2}-q|^{2r+2})+2^{2r+3}k\exp\lef(-Ck\rig),\\
        \nu(||m_1|-\mu|^{2r+3})&\leq \epsilon\nu(||m_1|-\mu|^{2r+2})+2^{2r+3}\bb P(||m_1|-\mu|\geq\epsilon)\\&\leq\epsilon\nu(||m_1|-\mu|^{2r+2})+2^{2r+3}k\exp\lef(-Ck\rig).
    \end{align*}
    Therefore picking small $\epsilon>0$, we have for large $k$ the exponential term is sufficiently small and
    \tny{\begin{align}\label{importanttailproof}
        \nu(|R_{1,2}-q|^{2r+2})\vee\nu(||m_1|-\mu|^{2r+2})\leq\frac{1}{2}\nu(|R_{1,2}-q|^{2r+2})\vee\nu(||m_1|-\mu|^{2r+2})+\frac{C^{r+1}r}{k}\lef(\frac{r+1}{k}\rig)^{r}.
    \end{align}}
    And we complete the induction. We note that the inequality of all  moments (the even ones can be obtained through H\"older's) implies the final result. And after we get $\Vert\sqrt n(R_{1,2}-q)\Vert_{\psi_2}<\infty$, weget $\Vert\sqrt{n-k}\tde m\Vert_{\psi_2}<\infty$ using \eqref{momentiter}.
    
\subsection{Proof of Lemma \ref{secondmoment}}
First we define $\bfa\sigma^1,\ldots\bfa\sigma^r$ to be $r$ replicas. And $\epsilon_j:=\sigma_k^j$, $\xi_j:=\sigma_n^j$ with the clique assumed to be indexed by $[k]$.
    We define $ R_{1,2}^-:=R_{1,2}-\epsilon_1\epsilon_2$, $\ca R^-_{1,2}:=R_{1,2}-\xi_1\xi_2$, $m_j:=\frac{1}{k}\sum_{i=1}^k\sigma^j_i$, and $\tde m_j:=\frac{1}{n-k}\sum_{i=k+1}^n\sigma_i^j$.
    Note that by the notation given in lemma \ref{cavity} and \ref{cavityII}, considering a function $f:\Sigma_n^r\to\bb R$, we have for $\tau_1,\tau_2>0$ and $\frac{1}{\tau_1}+\frac{1}{\tau_2}=1$ the following holds for some $C_1,C_2>0$:
    \tny{\begin{align}\label{cauchyschdire}
        |\nu(f)-\nu_{0,i}(f)|&\leq C_1(k,\theta,\theta_1)\nu(|f|^{\tau_1})^{1/\tau_1}\lef(\nu(|R_{1,2}^--q|^{\tau_2})^{1/\tau_2}\vee\nu(|m_1^--q|^{\tau_2})^{1/\tau_2}\rig),\nnb\\
        |\nu(f)-\nu_{0,i}(f)-\nu_{0,i}^\prime(f)|&\leq C_2(k,\theta,\theta_1)\nu(|f|^{\tau_1})^{1/\tau_1}\lef(\nu(|R_{1,2}^--q|^{2\tau_2})^{1/\tau_2}\vee\nu(|m_1^--q|^{2\tau_2})^{1/\tau_2}\rig).
    \end{align}}
    We define $\wh q:=\bb E[\tanh^2(\theta\sqrt qz+\theta_1\mu+h)]=\nu_{0,1}(\epsilon_1\epsilon_2)$ and $\tde q:=\bb E[\tanh^2(\theta\sqrt qz+h)]=\nu_{0,2}(\xi_1\xi_2)$.  It is checked that $q=c\wh q+(1-c)\tde q$. Then we can expand as follows
    \tny{\begin{align*}
        \nu((R_{1,2}-&q)^{r+1})=c\nu\big((\epsilon_1\epsilon_2-\wh q)(R_{1,2}-q)^{r}\big)+(1-c)\nu\big((\xi_1\xi_2-\tde q)(R_{1,2}-q)^r\big)\\
        &=c\nu\big((\epsilon_1\epsilon_2-\wh q)(R_{1,2}^--q)^{r}\big)+(1-c)\nu\big((\xi_1\xi_2-\tde q)(\ca R_{1,2}^--q)^r\big)\\
        &+\frac{cr}{n}\nu((1-\epsilon_1\epsilon_2\wh q)(R_{1,2}^--q)^{r-1}) 
        +\frac{(1-c)r}{n}\nu((1-\xi_1\xi_2\tde q)(\ca R_{1,2}^--q)^{r-1})+\tde O(r+2).
        \end{align*}}
        And again we define $m_j:=\frac{1}{k}\sum_{i=1}^k\sigma^j_i$, and notice that
        \begin{align*}
        \nu((m_1-\mu)^{r+1})&=\nu((\epsilon_1-\mu)(m_1-\mu)^{r})\\
        &=\nu((\epsilon_1-\mu)(m^-_1-\mu)^{r})+\frac{r}{k}\nu((1-\epsilon_1\mu)(m_1^--\mu)^{r-1})+\tde O(r+2),\\
        \nu(\tde m_1^{r+1})&=\nu(\xi_1\tde m_1^{-,r})+\frac{r}{n-k}\nu(\tde m_1^{-,r-1})+\tde O(r+2).
    \end{align*}
    Then we can make use of Lemma \ref{zeropointdecomp} and \eqref{cauchyschdire} by setting $\tau_1=\infty$ and $\tau_2=1$ :
    \begin{align*}
        \nu((1-\epsilon_1\epsilon_2\wh q)(R^-_{1,2}-q)^{r-1})&=\nu_{0,1}(1-\epsilon_1\epsilon_2 \wh q)\nu_{0,1}((R_{1,2}-q)^{r-1})+ \tde O(r)\\
        &=(1-\wh q^2)\nu((R_{1,2}-q)^{r-1})+\tde O(r),\\
        \nu((1-\epsilon_1\mu)(R^-_{1,2}-q)^{r-1})&=\nu_{0,1}(1-\epsilon_1\mu)\nu((R_{1,2}-q)^{r-1})+\tde O(r)\\
        &=(1-\mu^2)\nu((R_{1,2}-q)^{r-1})+\tde O(r),\\
        \nu((1-\xi_1\xi_2\tde q)(\ca R_{1,2}^--q)^{r-1})&=\nu_{0,2}(1-\tde q^2)\nu((R_{1,2}-q)^{r-1})+\tde O(r).
    \end{align*}
    Note that by lemma \ref{zeropointdecomp} we have $\nu_{0,1}\big((\epsilon_1\epsilon_2-\wh q)(R_{1,2}^--q)^{r}\big),\nu_{0,2}\big((\xi_1\xi_2-\tde q)(\ca R_{1,2}^--q)^r\big)=0$. Then, using \eqref{cauchyschdire} we have
    \begin{align}\label{firsttermsecond}
        \nu\big((\epsilon_1\epsilon_2-\wh q)(R_{1,2}^--q)^{r}\big)&=\nu_{0,1}^\prime\big((\epsilon_1\epsilon_2-\wh q)(R_{1,2}^--q)^{r}\big)+\tde O(r+2),\nnb\\
        \nu\big((\xi_1\xi_2-\tde q)(\ca R_{1,2}^--q)^r\big)&=\nu^\prime_{0,2}\big((\xi_1\xi_2-\tde q)(\ca R_{1,2}^--q)^r\big)+\tde O(r+2).
    \end{align}
    Then, change $r$ to $2r$, we can rewrite the above equations as:
    \ttny{\begin{align}\label{concentrations}
        \nu((R_{1,2}-q)^{2(r+1)})&=\frac{2r+1}{n}(c(1-\wh q^2)+(1-c)(1-\tde q^2))\nu((R_{1,2}-q)^{2r}) +c\nu_{0,1}^\prime((\epsilon_1\epsilon_2-\wh q)(R_{1,2}^--q)^{2r+1})\nnb\\
        &+(1-c)\nu^\prime_{0,2}((\xi_1\xi_2-\tde q)(\ca R_{1,2}^--q)^{2r+1})+\tde O(2r+2),\nnb\\
        \nu((m_1-\mu)^{2(r+1)})&=\frac{(2r+1)(1-\mu^2)}{k}\nu((m_1-\mu)^{2r})+\nu_{0,1}^\prime((\epsilon_1-\mu)(m_1^--q)^{2r+1})+\tde O(2r+2),\nnb\\
        \nu(\tde m_1^{2(r+1)})&=\frac{2r+1}{n-k}\nu(\tde m^{2r})+\nu^\prime_{0,2}(\xi_1(\tde m_1^-)^{2r+1})+\tde O(2r+2).
    \end{align}}
    Using lemma \ref{cavity} to analyze the derivative term from above, we note that the following quantities always appear:
    \begin{align}\label{aellellpdef}
   a(\ell,\ell^\prime):=\theta^2\nu_{0,1}(\epsilon_\ell\epsilon_{\ell^\prime}(\epsilon_1\epsilon_2-\wh q)),&\quad b(\ell,\ell^\prime):=\theta^2\nu_{0,1}((\epsilon_\ell\epsilon_{\ell^\prime}-\wh q)\epsilon_1),\nnb\\
   c(\ell):=\theta_1\nu_{0,1}((\epsilon_{3-\ell}-\mu)\epsilon_1\epsilon_2),&\quad d(\ell):=\theta_1\nu_{0,1}(\epsilon_{2-\ell}(\epsilon_1-\mu)),\nnb\\
   \tde a(\ell,\ell^\prime):=\theta^2\nu_{0,2}(\xi_\ell\xi_{\ell^\prime}(\xi_1\xi_2-\tde q)),
    \end{align}
   and
    \tny{\begin{gather*}
        U_{2,r}:=\nu((R_{1,2}-q)^{2r}),\quad U_{1,r}:=\nu((R_{1,2}-q)^{2r-1}(R_{1,3}-q)),\quad U_{0,r}:=\nu((R_{1,2}-q)^{2r-1}(R_{3,4}-q)),\\
        V_{0,1,r}:=\nu((R_{1,2}-q)^{2r-1}(m_1-\mu)),\quad V_{0,0,r}:=\nu((R_{1,2}-q)^{2r-1}(m_3-\mu)),\\
        V_{1,1,r}:=\nu((R_{1,2}-q)(m_1-\mu)^{2r-1}),\quad V_{1,0,r}:=\nu((R_{1,2}-q)(m_3-\mu)^{2r-1}),\\
        \quad W_{1,r}:=\nu((m_1-\mu)^{2r}),\quad W_{0,r}:=\nu((m_1-\mu)^{2r-1}(m_2-\mu)),\quad \tde W_{1,r}:=\nu(\tde m_1^{2r}),\\ 
        \quad\wh q_j:=\nu_{0,1}(\epsilon_1\cdots\epsilon_j)=\bb E[\tanh^j(\theta\sqrt q z+\theta_1\mu+h)],\quad \tde q_j:=\nu_{0,2}(\epsilon_1\ldots\epsilon_j)=\bb E[\tanh^j(\theta\sqrt qz+h)],\quad j\geq 3.
    \end{gather*}}
    Then we can write 
    \tny{\begin{gather}\label{definiabcd}
        a(2):=a(1,2)=\theta^2(1- \wh q^2),\quad a(1):=a(1,3) =\theta^2(\wh q-\wh q^2),\quad a(0):=a(3,4)=\theta^2(\wh q_4-\wh q^2),\nnb\\
        b(1):=b(1,2) = \theta^2\mu(1-\wh q),\quad b(0):=b(2,3)=\theta^2(\wh q_3-\wh q\mu),\nnb\\
        \quad d(1) =\theta_1(1-\mu^2),\quad d(0)=\theta_1(\wh q-\mu^2), \quad e(1)=\theta_1\mu(1-\wh q),\quad e(0)=\theta_1(\wh q_3 -\wh q\mu),\nnb\\
       \tde a(2):=\tde a(1,2)=\theta^2(1- \tde q^2),\quad \tde a(1):=\tde a(1,3) =\theta^2(\tde q-\tde q^2),\quad \tde a(0):=\tde a(3,4)=\theta^2(\tde q_4-\tde q^2).
    \end{gather}}
     Then we define
    \ttny{\begin{gather*}
        a(\ell,\ell^\prime;x,y)=a(|\{x,y\}\cap\{\ell,\ell^\prime\}|),\quad b(\ell,\ell^\prime;x)=b(\mbbm 1_{x\in\{\ell,\ell^\prime\}}),\quad e(\ell;x,y) =e(\mbbm 1_{\ell\in\{x,y\}}),\quad d(\ell,x)=d(\mbbm 1_{\ell=x}),\\
        \tde a(\ell,\ell^\prime;x,y)=\tde a(|\{x,y\}\cap\{\ell,\ell^\prime\}|).
    \end{gather*}}
    Using lemma \ref{cavity}, for a function $f$ on $\Sigma_{n}^r$ and two integers $x,y\leq n$
    \tny{\begin{align*}
        \nu_{0,1}^\prime&((\epsilon_x\epsilon_y-\wh q)f)=\sum_{1\leq\ell<\ell^\prime\leq r}a(\ell,\ell^\prime;x,y)\nu_{0,1}(f(R^-_{\ell,\ell^\prime}-q))-r\sum_{\ell\leq r}a(\ell,r+1;x,y)\nu_{0,1}(f(R_{\ell,r+1}^--q))\\
        &+\frac{r(r+1)}{2}a(0)\nu_{0,1}(f(R_{r+1,r+2}^--q))+\sum_{\ell\leq r}e(\ell;x,y)\nu_{0,1}(f(m_\ell^--\mu))-re(0)\nu_{0,1}(f(m_{r+1}^--\mu)).
    \end{align*}}
    And for the second path we have
     \tny{\begin{align*}
        \nu_{0,2}^\prime((\xi_x\xi_y-\wh q)f)&=\sum_{1\leq\ell<\ell^\prime\leq r}\tde a(\ell,\ell^\prime;x,y)\nu_{0,2}(f(\ca R^-_{\ell,\ell^\prime}-q))-r\sum_{\ell\leq r}\tde a(\ell,r+1;x,y)\nu_{0,2}(f(\ca R_{\ell,r+1}^--q))\\
        &+\frac{r(r+1)}{2}\tde a(0)\nu_{0,2}(f(\ca R_{r+1,r+2}^--q)).
    \end{align*}}
    And we can proceed similarly to get
    \ttny{\begin{align*}
        \nu_{0,1}^\prime((\epsilon_x-\mu)f)&=\sum_{1\leq\ell<\ell^\prime\leq r}b(\ell,\ell^\prime;x)\nu_{0,1}(f(R^-_{\ell,\ell^\prime}-q))-r\sum_{\ell\leq r}b(\ell,r+1;x)\nu_0(f(R^-_{\ell,r+1}-q))\\
        &+\frac{r(r+1)}{2}b(0)\nu_{0,1}(f(R_{r+1,r+2}^--q))+\sum_{\ell\leq r}d(\ell,x)\nu_{0,1}(f(m_\ell^--\mu))-rd(0)\nu_{0,1}(f(m_{r+1}^--\mu)),
    \end{align*}}
    and
    \tny{\begin{align*}
        \nu_{0,2}^\prime((\xi_x-\mu)f)&=\sum_{1\leq\ell<\ell^\prime\leq r}\tde b(\ell,\ell^\prime;x)\nu_{0,2}(f(\ca R^-_{\ell,\ell^\prime}-q))-r\sum_{\ell\leq r}\tde b(\ell,r+1;x)\nu_{0,2}(f(\ca R^-_{\ell,r+1}-q))\\
        &+\frac{r(r+1)}{2}\tde b(0)\nu_{0,2}(f(\ca R_{r+1,r+2}^--q)).
    \end{align*}}
    Therefore, using lemma \ref{cavity}, we can further write \eqref{concentrations} as:
    \tny{\begin{align}\label{firstder}
        U_{2,r} &= \frac{2r-1}{n}(c(1-\wh q^2)+(1-c)(1-\tde q^2))U_{2,r-1}+(ca(2)+(1-c)\tde a(2))U_{2,r}-4(ca(1)+(1-c)\tde a(1))U_{1,r}\nnb\\
        &+3(c a(0)+(1-c)\tde a(0))U_{0,r}+2ce(1)V_{0,1,r}-2ce(0)V_{0,0,r} +\tde O(2r+1),\nnb\\
        W_{1,r} &=\frac{2r-1}{n}(c(1-\mu^2)+(1-c))W_{1,r-1}-(c b(1)+(1-c)\tde b(1))V_{1,1,r} +(cb(0)+(1-c)\tde b(0))V_{1,0,r}\nnb\\
        &+cd(1)W_{1,r}-cd(0)W_{0,r} +\tde O(2r+1).
    \end{align}}
    Proceed similarly to other terms, we see that
    \tny{\begin{align*}
        U_{1,r}
        &=\frac{2r-1}{n}(c(\wh q-\wh q^2)+(1-c)(\tde q-\tde q^2))U_{2,r-1}+(ca(1)+(1-c)\tde a(1))U_{2,r}\\
        &+(c(a(2)-2a(1)-3a(0))+(1-c)(\tde a(2)+2\tde a(1)-3\tde a(0)))U_{1,r}+(c(6a(0)-3a(1))\\
        &+(1-c)(6\tde a(0)-3\tde a(1)))U_{0,r}+c(e(1)+e(0))V_{0,1,r}+c(e(1)-3e(0))V_{0,0,r}+\tde O(2r+1),\\
         U_{0,r}
        &=\frac{2r-1}{n}(c(\wh q_4-\wh q^2)+(1-c)(\tde q_4-\tde q^2))U_{2,r-1}+(ca(0)+(1-c)\tde a(0))U_{2,r}\\
        &+(c(4a(1)-8a(0))+(1-c)(4\tde a(1)-8\tde a(0)))U_{1,r}+c(a(2)-8a(1)+10a(0))U_{0,r}\\
        &+(1-c)(\tde a(2)-8\tde a(1)+10\tde a(0))U_{0,r}+2e(0)V_{0,1,r}+(2e(1)-4e(0))V_{0,0,r}+\tde O(2r+1),\\
         V_{0,1,r}
        &=\frac{2r-1}{n}c(\mu-\wh q\mu)U_{2,r-1}+cb(1)U_{2,r}-2c(b(1)+b(0))U_{1,r}+3cb(0)U_{0,r}\\
        &+c(d(1)+d(0))V_{0,1,r}-2cd(0)V_{0,0,r}+\tde O(2r+1),\\
         V_{0,0,r}
        &=\frac{2r-1}{n}(c({\wh q_3-\wh q\mu})+(1-c)({\tde q_3-\tde q\mu}))U_{2,r-1}+cb(0)U_{2,r}+c(2b(1)-6b(0))U_{1,r}\\
        &+(6b(0)-3b(1))U_{0,r}+2d(0)V_{0,1,r}+(d(1)-3d(0))V_{0,0,r}+\tde O(2r+1),\\
        V_{1,1,r}
        &=(2r-1)\frac{\mu-q\mu}{k}W_{1,r-1}+ (a(2)-2a(1))V_{1,1,r}+(3a(0)-2a(1))V_{1,0,r}\\&+e(1)W_{1,r}+(e(1)-2e(0))W_{0,r}+\tde O(2r+1).\\
        V_{1,0,r}
        &=(2r-1)\frac{q_3-q\mu}{k}W_{1,r-1}+(2a(1)-3a(0))V_{1,1,r}+(a(2)-6a(1)+6a(0))V_{1,0,r}\\
        &+e(0)W_{1,r}+(2e(1)-3e(0))W_{0,r}+\tde O(2r+1).\\
         W_{0,r}
        &=(2r-1)\frac{q-\mu^2}{k}W_{1,r-1}+(b(1)-2b(0))V_{1,1,r}+b(0)V_{1,0,r}+d(0)W_{1,r}+(d(1)-2d(0))W_{0,r}\\&+\tde O(2r+1).
    \end{align*}}
    And then considering the second path we have
    \begin{align*}
        \tde W_{1,r}&=(2r-1)\frac{1}{n-k}\tde W_{1,r-1}+\tde O(2r+1).
    \end{align*}
    And we also define $U_{1,0}=U_{0,0}=1$.
    Then, we define 
    \begin{align*}
          \bfa A_1:&= \begin{bmatrix}
        a(2)&-4a(1)&3a(0)&2e(1)&-2e(0)\\
        a(1)&a(2)-2a(1)-3a(0)&6a(0)-3a(1)&e(1)+e(0)&e(1)-3e(0)\\
        a(0)&4a(1)-8a(0)&a(2)-8a(1)+10a(0)&2e(0)&2e(1)-4e(0)\\
        b(1)&-2b(1)-2b(0)&3b(0)&d(1)+d(0)&-2d(0)\\
        b(0)&2b(1)-6b(0)&6b(0)-3b(1)&2d(0)&d(1)-3d(0)
        \end{bmatrix}\\
        &=\begin{bmatrix}
            \bfa A_{11}&\bfa A_{12}\\
            \bfa A_{13}&\bfa A_{14}
        \end{bmatrix},\\
        \bfa A_2:&=\begin{bmatrix}
        a(2)-2a(1)&3a(0)-2a(1)&e(1)&e(1)-2e(0)\\
        2a(1)-2a(0)&a(2)-6a(1)+6a(0)&e(0)&2e(1)-3e(0)\\
        -b(1)&b(0)&d(1)&-d(0)\\
        b(1)-2b(0)&b(0)&d(0)&d(1)-2d(0)
        \end{bmatrix}=\begin{bmatrix}
            \bfa A_{21}&\bfa A_{22}\\
            \bfa A_{23}&\bfa A_{24}
        \end{bmatrix},
       \end{align*}
       with $\bfa A_{11}\in\bb R^{3\times 3}$ and $\bfa A_{21}\in\bb R^{2\times 2}$ and other matrices defined similarly. Then wealso define
        \begin{align*}
          \tde{\bfa A}_1:&= \begin{bmatrix}
        a(2)&-4a(1)&3a(0)&0&0\\
        a(1)&a(2)-2a(1)-3a(0)&6a(0)-3a(1)&0&0\\
        a(0)&4a(1)-8a(0)&a(2)-8a(1)+10a(0)&0&0\\
       0&0&0&0&0\\
       0&0&0&0&0
        \end{bmatrix}.
       \end{align*}
     We also define
    \tny{\begin{align*}
        \bfa x_{r}:&=(U_{2,r},U_{1,r},U_{0,r},V_{0,1,r},V_{0,0,r})^\top,\quad\bfa y_{r}:=(V_{1,1,r}, V_{1,0,r},W_{1,r}, W_{0,r})^\top,\\
        \bfa b_1:&=\lef(1-\wh q^2,\wh q-\wh q^2,\wh q_4-\wh q^2,\mu-\mu \wh q,\wh q_3-\mu \wh q\rig)^\top,\quad\bfa b_2:=(\mu-\wh q\mu, \wh q_3-\wh q\mu,1-\mu^2,\wh q-\mu^2)^\top,\\
        \tda b_1:&=\lef(1-\tde q^2,\tde q-\tde q^2,\tde q_4-\tde q^2,0,\tde q_3\rig)^\top.
    \end{align*}}
    Written in matrix form , we have
    \begin{align}\label{cavityiteration}
        \bfa x_r &=(\bfa A_1+\tde{\bfa A}_1)\bfa x_r+ \frac{2r-1}{n}U_{2,r-1}(c\bfa b_1+(1-c)\tda b_1)+\tde O(2r+1),\nnb\\
        \bfa y_r &=\bfa A_2\bfa y_r+\frac{1}{k}(2r-1)W_{1,r-1}\bfa b_2+\tde O(2r+1).
    \end{align}
\subsection{Proof of Lemma \ref{secondmomentpositive}}
    The proof will largely follows from lemma \ref{secondmoment} upon observing that in the smart path method II given by lemma \ref{cavityII}. Then we have
    \begin{align*}
        \nu_{0,3}^+\bl\prod_{i\in I}\epsilon_i f_r^-\br=\nu_{0,3}^+\bl\prod_{i\in I}\epsilon_i\br\nu_{0,3}^+( f_r^-)=\bb E[\tanh^{|I|}(\theta_1\mu+\theta\sqrt{q}z +h)]\nu^+_{0,3}(f_r^-),
    \end{align*}
    which essentially recovers the analogue of Lemma \ref{cavity}.
    And moreover, the smart path method given by \eqref{cavityII} implies that
    \ttny{\begin{align*}
        \nu_{t,3}^{+\prime}(f) &=\theta^2\bl\sum_{1\leq \ell<\ell^\prime\leq r}\nu^+_{t,3}(f\epsilon_\ell\epsilon_{\ell^\prime}(R_{\ell,\ell^\prime}-q))\br -r\theta^2\sum_{\ell\leq r}\nu^+_{t,3}(f\epsilon_{\ell}\epsilon_{r+1}(R_{\ell,r+1}-q))\\
        &+\theta^2\frac{r(r+1)}{2}\nu^+_{t,3}(f\epsilon_{r+1}\epsilon_{r+2}(R_{r+1, r+2}-q))+\theta_1\bigg(\sum_{\ell\leq r}\nu_{t,3}^+(f\epsilon_{\ell}(m_{\ell}-\mu))-r\nu_{t,3}^+(f\epsilon_{r+1}(m_{r+1}-\mu))\bigg).
    \end{align*}}
    And one will then the result of the results in an identical way as in the proof of Lemma \ref{secondmoment}. For the iteration of $\tde W_{1,r}$ we note that the proof does not change. And, we notice that the iteration of $U_{2,r}$ is invariant when flipping the sign of $\mu$.

\subsection{Proof of Lemma \ref{cltatline}}
We first prove the case when $\mu$, the solution of \eqref{meanfieldeq} is unique and is $0$. One will first check that $c\bfa A_1+(1-c)\tda A_1$ has eigenvalues 
\tny{\begin{align}
    &\theta^2\lef(c(1-2\wh q+\wh q_4)+(1-c)(1-2\tde q+\tde q_4)\rig),\quad\theta^2(c(1-4\wh q+3\wh q_4)+(1-c)(1-4\tde q+4\tde q_4)),\quad \theta_1q.
\end{align}}
And therefore with the condition given by definition \ref{withintheATline} we have $\bfa I-c\bfa A_1-(1-c)\tda A_1$ is invertible. Similarly, $\bfa A_2$ has the same eigenvalues and $\tda A_2$ has the eigenvalues given by $q\theta_1$. Then, under the high temperature, all of the inverse matrices are well defined under the non-critical assumption. For the low temperature the results are more complicated, where wenotice that the invertibility relies on the the following fact
\tny{\begin{align*}
\bfa B:&=\bfa I-c\bfa A_{21}-(1-c)\tda A_{21},\thinspace\det\begin{bmatrix}
        \bfa B&c\bfa A_{12}\\
        c\bfa A_{13} &\bfa I -c\bfa A_{14}
    \end{bmatrix}=\det(\bfa B)\det(\bfa I-c\bfa A_{14}-c^2\bfa A_{13}\bfa B^{-1}\bfa A_{12}).
\end{align*}}
And the second determinate is undecided. Similarly, the invertibility of the second vector equation is decided by $\bfa I-\bfa A_{24}-\bfa A_{23}(\bfa I-\bfa A_{21})^{-1}\bfa A_{22}$. 
Note that by Lemma \ref{secondmoment} and \eqref{inductivehypo},
\tny{\begin{align*}
    U_{2,r}&=\nu((R_{1,2}-q)^{2r})=\frac{(2r-1)!!}{n^r}\lef((\bfa I-c\bfa A_{1}-(1-c)\tda A_1)_1^{-1}(c\bfa b_{1}+(1-c)\tda b_1)\rig)^{2r}+\tde O(2r+1),\\
    W_{1,r}&=\nu((m_1-\mu)^{2r})=\frac{(2r-1)!!}{k^r}((\bfa I-\bfa A_{2})_1^{-1}\bfa b_{2})^{2r}+\tde O(2r+1),\\
    \tde W_{1,r}&=\frac{(2r-1)!!}{(n-k)^r}((\bfa I-\tda A_2)_1^{-1}\tda b_2)^{2r}+\tde O(2r+1).
\end{align*}}
Then we move toward the case conditional on positive $m_1$. We note that the negative case holds analogously with the positive since $m$ has symmetric distribution w.r.t. $0$. Consider the two set
\begin{align*}
    \ca A^-:=\{m_1^-> 0\},\qquad \ca A:=\{m_1> 0\}.
\end{align*}
Then we have
\begin{align*}
    \ca A^-\Delta\ca A\subset \ca B:= \lef\{m_1\in\lef[-\frac{1}{k},\frac{1}{k}\rig]\rig\}.
\end{align*}
And therefore, by the concentration inequality given by Lemma \ref{exponentialeq}, $\nu(\ca A^-\Delta\ca A)\leq\exp(-Ck)$ for some $\mu$. Then weconclude that
\begin{align*}
   \bb E\bigg[\bl\frac{k}{\ca V_m^l}\br^{r/2}(m-\mu)^r\bigg |m>0\bigg]=\bb \bb E\bigg[\bl\frac{k}{\ca V_m^l}\br^{r/2}(m-\mu)^r\bigg |m^->0\bigg]+o(1)\to\bb E[z^r],\\
        \bb E\bigg[\bl\frac{k}{\ca V_m^l}\br^{r/2}(m-\mu)^r\bigg |m<0\bigg]=\bb \bb E\bigg[\bl\frac{k}{\ca V_m^l}\br^{r/2}(m-\mu)^r\bigg |m^-<0\bigg]+o(1)\to\bb E[z^r].
\end{align*}

\subsection{Proof of Theorem \ref{exact recoveryguarantees}}
    The proof follows from utilizing lemma \ref{overlapclt} , lemma \ref{largecliqueweakclt}, and the proof of Corollary 3.2.1, 3.4.1, and 3.7.1 in \citep{he2023hidden1}, and following the proof of theorem 3.9 in \citep{he2023hidden1} we complete the whole proof.
\subsection{Proof of Theorem \ref{thm6.1}}
    The proof is a direct result of lemma \ref{univer}, \ref{lm5.3}, \ref{lm5.4}, \ref{lm5.5}, \ref{lm5.6}. We proceed the same proof as the Gaussian case to get the desired result.
\section{Proof of Results in the Appendix}
\subsection{Proof of Lemma \ref{cavityII}}
    Note that one can divide $\mca H_{t,3}(\bfa\sigma):=\mca H_{t,3,1}(\bfa\sigma)+\mca H_{t,3,2}(\bfa\sigma)$. We define
    \begin{align}\label{decompose1}
        \mca H_{t,3,2}(\bfa\sigma) =\frac{\theta_1t}{k}\sum_{i\leq k}\sigma_i\sigma_k+\theta_1(1-t)\sigma_k(\mu\mbbm 1_{m>0}-\mu\mbbm 1_{m<0}).
    \end{align} to be the contribution of the Ferromagnetic correlation in the Hamiltonian. It is also checked that $\mca H_{1,t}(\bfa\sigma)$ is the standard interpolation path for SK model. Then wesee that
    \tny{\begin{align*}
        \nu_{t,3}^\prime(f) = \ub{\sum_{i\in[r]}\sum_{\bfa\sigma^i\in\{-1,1\}^{n}}\bb E\bigg[\frac{\pta \la f\ra_{t,3}}{\pta \mca H_{t,3,1}(\bfa\sigma^i)}\frac{\pta \mca H_{t,3,1}(\bfa\sigma^i)}{\pta t}\bigg]}_{T_1}+\ub{\sum_{i\in[r]}\sum_{\bfa\sigma^i\in\{-1,1\}^n}\bb E\bigg[\frac{\pta \la f\ra_{t,3}}{\pta \mca H_{t,3,2}(\bfa\sigma^i)}\frac{\pta \mca H_{t,3,2}(\bfa\sigma^i)}{\pta t}\bigg]}_{T_2}.
    \end{align*}}
    Therefore only the second part needs to be studied since the first one follows from the results in Lemma 1.6.3. of \citep{talagrand2010mean}. 
    We then proceed by
    \begin{align*}
        \sum_{\bfa\sigma^i}\bb E\bigg[\frac{\pta \la f\ra_{t,3}}{\pta \mca H_{t,3,2}(\bfa\sigma^i)}\frac{\pta \mca H_{t,3,2}(\bfa\sigma^i)}{\pta t}\bigg]&=\bb E\lef[\theta_1\lef(\la \epsilon_i(f-\la f\ra_{t,3})(m_i^--\mu\mbbm 1_{m^-_i>0}+\mu\mbbm 1_{m^-_i<0})\ra_{t,3}\rig)\rig].
    \end{align*}
    And noticing that $\la f\ra_{t,3}\la m_i^--\mu\mbbm 1_{m_i>0}+\mu\mbbm 1_{m_i<0}\ra_{t,3}=\la f\epsilon_{r+1}(m_i^--\mu\mbbm 1_{m_i>0}+\mu\mbbm 1_{m_i<0})\ra_{t,3}$, we can sum over all replicas to arrive at the following,
    \tny{\begin{align*}
        T_2&= \theta_1\bl\sum_{\ell\leq r}\nu_{t,3}(f\epsilon_\ell(m_i^--\mu\mbbm 1_{m_i^->0}+\mu\mbbm 1_{m_i^-<0})-r\nu_{t,3}(f\epsilon_{r+1}(m_{r+1}^--\mu\mbbm 1_{m_i^->0}+\mu\mbbm 1_{m_i^-<0})))\br\\
        &=\theta_1\bl\sum_{\ell\leq r}\nu_{t,3}(f\epsilon_\ell(m_i-\mu\mbbm 1_{m_i^->0}+\mu\mbbm 1_{m_i^-<0}))-r\nu_{t,3}(f\epsilon_{r+1}(m_{r+1}-\mu\mbbm 1_{m_i^->0}+\mu\mbbm 1_{m_i^-<0}))\br.
    \end{align*}}
\subsection{Proof of Lemma \ref{zeropointdecompII}}
    The proof follows immediately upon observing that conditional on $A_{ij}$ we have decomposibility in lemma \ref{zeropointdecomp} when $t=0$.
\subsection{Proof of Lemma \ref{univer}}
Similar to the cavity method, here we consider the following smart path:
\begin{align*}
    -\mca H_t(\bfa\sigma):&=\theta\sqrt t\sum_{1\leq i<j\leq k}\frac{ g_{ij}}{\sqrt n}\sigma_i\sigma_j+\theta\sqrt t\sum_{i\in[k],j\in[k+1:n]}\frac{g_{ij}}{\sqrt n}\sigma_i\sigma_j+\theta\sqrt{1-t}\sum_{1\leq i<j\leq k}\frac{\eta_{ij}}{\sqrt n}\sigma_i\sigma_j\\
    &+\theta\sqrt{1-t}\sum_{i\in[k],j\in[k+1:n]}\frac{\eta_{ij}}{\sqrt n}\sigma_i\sigma_j+\sum_{i,j\in[k]}\frac{\theta_1}{2 k}\sigma_i\sigma_j+\sum_{i=1}^n\sigma_ih_i.
\end{align*}
Therefore one will get
\tny{\begin{align*}
    -\frac{\pta \mca H_t(\bfa\sigma) }{\pta t}&=\frac{\theta}{2\sqrt t}\sum_{1\leq i<j\leq k}\frac{g_{ij}}{\sqrt n}\sigma_i\sigma_j+\frac{\theta}{2\sqrt t}\sum_{i\in[k],j\in[k+1:n]}\frac{g_{ij}}{\sqrt n}\sigma_i\sigma_j-\frac{\theta}{2\sqrt{1-t}}\sum_{1\leq i<j\leq k}\frac{\eta_{ij}}{\sqrt n}\sigma_i\sigma_j\\
    &-\frac{\theta}{2\sqrt{1-t}}\sum_{i\in[k],j\in[k+1:n]}\frac{\eta_{ij}}{\sqrt n}\sigma_i\sigma_j.
\end{align*}}
And for $\bfa\tau\in\Sigma_n$, we have for $f:\Sigma_k\to\bb R^+$:
\ttny{\begin{align*}
    \frac{1}{\sqrt n}\frac{d\la f\sigma_i\sigma_j\ra_t}{d g_{ij}} &=\frac{\sum_{\bfa\tau}f(\bfa\tau)\frac{\theta\sqrt t}{ n}\exp(-\mca H_t(\bfa\tau))}{\sum_{\bfa\sigma}\exp(-\mca H_t(\bfa\sigma))}-\frac{\sum_{\bfa\tau}f(\bfa\tau)\sigma_i\sigma_j\frac{\theta\sqrt t}{ n}\exp(-\mca H_t(\bfa\tau))\sum_{\bfa\sigma}\frac{\theta\sqrt t\sigma_i\sigma_j}{\sqrt n}\exp(-\mca H_t(\bfa\sigma))}{\lef(\sum_{\bfa\sigma}\exp(-\mca H_t(\bfa\sigma))\rig)^2}\\
    &=\frac{\theta\sqrt t}{n}[\la f\tau_i^2\tau_j^2\ra_t-\la f\sigma_i\sigma_j\ra_t\la\sigma_i\sigma_j\ra_t]. \\
    \frac{1}{\sqrt n}\frac{d\la f\sigma_i\sigma_j\ra_t}{d \eta_{ij}} &=\frac{\theta\sqrt {1-t}}{n}[\la f\tau_i^2\tau_j^2\ra_t-\la f\sigma_i\sigma_j\ra_t\la\sigma_i\sigma_j\ra_t].
\end{align*}}

And for the second order derivatives we have for all bounded functions
\begin{align*}
    \frac{1}{\sqrt n}\frac{d^2\la f\sigma_i\sigma_j\ra_t}{d \eta^2_{ij}} &=\frac{\theta^2(1-t)}{n^{3/2}}[2\la f\sigma_i\sigma_j\ra_t\la\sigma_i\sigma_j\ra_t^2-2\la f\ra_t\la\sigma_i\sigma_j\ra_t],\\
    \bigg |\frac{1}{\sqrt n}\frac{d^3\la f\sigma_i\sigma_j\ra_t}{d \eta^3_{ij}}\bigg| &=\bigg|\frac{2\theta^3(1-t)}{n^{2}}(\la f\ra_t-\la f\sigma_i\sigma_j\ra_t\la\sigma_i\sigma_j\ra_t)(1-3\la\sigma_i\sigma_j\ra^2)\bigg|=O\bl\frac{1}{n^2}\br\la f\ra_t.
\end{align*}
Analogously, we have
\begin{align*}
    \bigg |\frac{1}{\sqrt n}\frac{d^3\la f\ra_t\la\sigma_i\sigma_j\ra_t}{d\eta^3_{ij}}\bigg |=O\bl\frac{1}{n^2}\br\la f\ra_t.
\end{align*}
And then one consider the function $\varphi(t):=\bb E[\la f\ra_t]:=\bb E\bigg[\frac{\sum_{\bfa\sigma}f(\bfa\sigma)\exp(-\mca H_t(\bfa\sigma))}{\sum_{\bfa\sigma}\exp(-\mca H_t(\bfa\sigma))}\bigg]$, 
\begin{align*}
    \varphi^\prime(t)=\bb E\bigg[\bigg\la -f\frac{\pta \mca H_t}{\pta t}\bigg\ra_t\bigg]+\bb E\bigg[\la f\ra_t\bigg\la\frac{\pta \mca H_t}{\pta t}\bigg\ra_t\bigg].
\end{align*}
Therefore we have by Gaussian integration by parts in lemma \ref{gtrbt} and Non-Gaussian integration by parts in lemma \ref{lm1.1} to get that there exists $\xi_1\in(0\wedge \eta,0\vee\eta)$ depending on $\eta_{ij}$ such that
\begin{align}\label{pieces1}
    \bb E\bigg[\frac{1}{\sqrt n}\sigma_i\sigma_j\frac{\exp(-\mca H_t(\bfa\sigma))}{\sum_{\bfa\sigma}\exp(-\mca H_t(\bfa\sigma))}\bl\frac{g_{ij}}{\sqrt t}-\frac{\eta_{ij}}{\sqrt{1-t}}\br\bigg]=O\bl\frac{1}{n^2}\bb E\bigg[\eta^4\la f|\eta_{ij}=\xi_1\ra_t\bigg]\br.
\end{align}
Further we notice that by Taylor expansion at $\eta$ and using the fact that $f>0$ we have for all $\xi_1\in(0\wedge \eta,0\vee\eta)$ the following holds (notice that $\xi_1$'s order does not depend on $n$)
\begin{align}\label{pieces2}
    \la f|\eta_{ij}=\xi_1\ra_t&=\la f|\eta_{ij}=0\ra_t+\sum_{i=1}^\infty\frac{d^i\la f|\eta_{ij}=0\ra_t}{d\eta^i_{ij}}\frac{1}{i!}\xi_1^i\nnb\\
    &=\la f|\eta_{ij}=0\ra_t+O\bl\frac{1}{\sqrt n}|\xi_1|\la f|\eta_{ij}=0\ra_t\br=O\lef(\la f|\eta_{ij}=0\ra_t\rig).
\end{align}
Using lemma \ref{lm1.1} and notice that there exists only $O(kn)$ i.i.d. $\eta_{ij}$ to be replaced by Gaussians, combining pieces from \eqref{pieces1} and \eqref{pieces2}, we have by symmetry, for all $t\in[0,1]$,
\begin{align}\label{derivativeub}
    |\varphi^\prime(t)|\leq C\frac{k}{n}\bb E[\eta^4\la f|\eta_{ij}=0\ra_t]=C\frac{k}{n}\bb E[\eta^4]\bb E[\la f|\eta_{ij}=0\ra_t]=O\bl\frac{k}{n}\br\bb E[\eta^4]\varphi(t).
\end{align}
Therefore weconclude that for some $C_1,C_2>0$:
\begin{align*}
    \exp\bl\frac{-C_1kt}{n}\br\bb E[\la f\ra_0]\leq\bb E[\la f\ra_t]\leq\exp\bl\frac{C_2kt}{n}\br\bb E[\la f\ra_0].
\end{align*}
Therefore we let $f$ be the mgf of the local spins, the CLT in theorem \ref{cltrfcr} then immediately holds.
\begin{lemma}[Non-Gaussian integration by parts]\label{lm1.1}
    Let $\eta$ be a real random variable such that $\bb E[\eta^4]<\infty$, $\bb E[\eta]=\bb E[\eta^3]=0$ and $\bb E[\eta^2]=1$. Let $F:\bb R\to\bb R$ be thrice continuously differentiable. Then there exists $\xi_1,\xi_2\in(0\wedge \eta,0\vee\eta)$ such that
    \begin{align*}
        \bb E[\eta F(\eta)]-\bb E[\eta^2]\bb E[F^\prime(\eta)]=\bb E\lef[\frac{\eta^4}{6}F^{(3)}(\xi_1)\rig]-\bb E[\eta^2]\bb E\bigg[\frac{\eta^2}{2}F^{(3)}(\xi_2)\bigg].
    \end{align*}
\end{lemma}
The proof follows by elementary Taylor expansion, noticing that there exists $\xi_1,\xi_2\in(0\wedge \eta,0\vee\eta)$ such that
\sm{\begin{align*}
    &F(\eta)-F(0)-\eta F^{\prime}(0)-\frac{\eta^2}{2}F^{(2)}(0)=\frac{\eta^3}{6} F^{(3)}(\xi_1),\quad F^\prime(\eta)-F^\prime(0)-\eta F^\prime(0)= \frac{\eta^2}{2} F^{(2)}(\xi_2).
\end{align*}}
Therefore,
\tny{\begin{align*}
    \bb E[\eta F(\eta)]&-\bb E[\eta^2]\bb E[F^\prime(\eta)]=\bb E[\eta F(\eta)]-\bb E[\eta^2]\bb E[F^\prime(\eta)]-F(0)\bb E[\eta]\\
    &=\bb E\bigg[\eta \bl F(\eta)-F(0)-\eta F^\prime(0)-\frac{\eta^2}{2}F^{(2)}(0)\br\bigg]-\bb E[\eta^2]\bb E\bigg[F^\prime(\eta)- F^\prime(0)-\eta F^{(2)}(0)\bigg]\\
    &=\bb E\lef[\frac{\eta^4}{6}F^{(3)}(\xi_1)\rig]-\bb E[\eta^2]\bb E\bigg[\frac{\eta^2}{2}F^{(3)}(\xi_2)\bigg].
\end{align*}}
And one can prove the second inequality analogously by treating each coordinates separately.

\subsection{Proof of Lemma \ref{lm5.3}}
    The proof will follow the similar idea as lemma \ref{univer}. Consider two sets $S:=[k]$ and $S^\prime:=[r+1:k+r]$ weuse the following smart path
    \ttny{\begin{align*}
        -\mca H_{t,S}(\bfa\sigma)&=\theta\sqrt t\sum_{1\leq i<j\leq k+r}\frac{ g_{ij}}{\sqrt n}\sigma_i\sigma_j+\theta\sqrt t\sum_{i\in[k+r],j\in[k+r+1:n]}\frac{g_{ij}}{\sqrt n}\sigma_i\sigma_j+\theta\sqrt{1-t}\sum_{1\leq i<j\leq k+r}\frac{\eta_{ij}}{\sqrt n}\sigma_i\sigma_j\\
    &+\theta\sqrt{1-t}\sum_{i\in[k+r],j\in[k+r+1:n]}\frac{\eta_{ij}}{\sqrt n}\sigma_i\sigma_j+\sum_{i,j\in S}\frac{\theta_1}{2 k}\sigma_i\sigma_j+\sum_{i=1}^n\sigma_ih_i.
    \end{align*}}
    Using the definition of $\la f\ra_{S,t}:=\frac{\sum_{\bfa\sigma} f(\bfa\sigma)\exp(-\mca H_{t,S}(\bfa\sigma))}{\sum_{\bfa\sigma}\exp(-\mca H_{t,S}(\bfa\sigma))}$. Then wecheck that for $f_1,f_2:\Sigma_{k+r}\to\bb R$,
    \begin{align*}
        \varphi(t):=\bb E[\la f_1\ra_{S,t}\la f_2\ra_{S^\prime,t}].
    \end{align*}
    And we have
    \ttny{\begin{align*}
        \varphi^\prime(t)=\bb E\bigg[\bl\bigg\la\frac{\pta \mca H_t}{\pta t}\bigg\ra_{S,t}+\bigg\la\frac{\pta \mca H_t}{\pta t}\bigg\ra_{S^\prime,t}\br\la f_1\ra_{S,t}\la f_2\ra_{S,t}\bigg]-\bb E\bigg[\bigg\la f_1\frac{\pta \mca H_t}{\pta t}\bigg\ra_{S,t}\la f_2\ra_{S^\prime,t}\bigg]-\bb E\bigg[\bigg\la f_2\frac{\pta \mca H_t}{\pta t}\bigg\ra_{S^\prime,t}\la f_1\ra_{S,t}\bigg].
    \end{align*}}
    And we also notice that by Gaussian integration by parts it is checked that
    \ttny{\begin{align*}
        \frac{1}{\sqrt {nt}}\frac{d}{dg_{ij}}\la f_1\sigma_i\sigma_j \ra_{S,t}\la f_2\ra_{S^\prime,t}&=\frac{\theta}{n}\lef(\la f_1\ra_{S,t}\la f_2\ra_{S^\prime,t}-\la f_1\sigma_i\sigma_j \ra_{S,t}\la f_2\ra_{S^\prime,t}\la\sigma_i\sigma_j\ra_{S,t}+\la f_1\sigma_i\sigma_j \ra_{S,t}\la f_2\sigma_i\sigma_j \ra_{S^\prime,t}\rig. \\
        &\lef. -\la f_1\sigma_i\sigma_j\ra_{S,t}\la f_2\ra_{S^\prime,t}\la\sigma_i\sigma_j\ra_{S^\prime,t} \rig)=\frac{1}{\sqrt{n(1-t)}}\frac{d}{d\eta_{ij}}\la f_1\sigma_i\sigma_j\ra_{S,t}\la f_2\ra_{S,t}.
    \end{align*}}
    Similarly wecheck that
    \begin{align*}
        \frac{1}{\sqrt {nt}}\frac{d}{dg_{ij}}\la\sigma_i\sigma_j\ra_{S,t}\la f_1 \ra_{S,t}\la f_2\ra_{S^\prime,t}=\frac{1}{\sqrt {n(1-t)}}\frac{d}{d\eta_{ij}}\la f_1\ra_{S,t}\la f_2\ra_{S^\prime,t}.
    \end{align*}
    Therefore the two principal terms will cancel out in the derivative.
    And for the third order derivative we have
    \begin{align*}
        \bigg|\frac{1}{\sqrt n}\frac{d^3}{d\eta_{ij}^3}\la f_1\sigma_i\sigma_j \ra_{S,t}\la f_2\ra_{S^\prime,t}\bigg|=O\bl\frac{1}{n^2}\br\la f_1\ra_{S,t}\la f_2\ra_{S^\prime,t}.
    \end{align*}
    And therefore we have by similar vein of \eqref{derivativeub}:
    \begin{align*}
       \forall t\in[0,1],\quad |\varphi^\prime(t)|=O\bl\frac{k}{n}\br\varphi(t)\quad\Rightarrow\quad\varphi(0)\leq\exp\lef(C\frac{k}{n}\rig)\varphi(1).
    \end{align*}
    Then we let $f_1=\mbbm 1_{A_S}\mbbm 1_{A_{S^\prime}}\exp\lef(\frac{\theta_1}{2k}\sum_{i,j\in S}\sigma_i\sigma_j\rig)$, $f_2:=\exp\lef(-\frac{\theta_1}{2k}\sum_{i,j\in S^\prime}\sigma_i\sigma_j\rig)$ to complete the proof.
\subsection{Proof of Lemma \ref{lm5.4}}
    The proof follows identically to the \citep{carmona2006universality} and weomit it here.

\subsection{Proof of Lemma \ref{lm5.5}}
    The proof follows from the smart path of 
    \begin{align*}
    -\mca H_t(\bfa\sigma):&=\theta\sqrt t\sum_{1\leq i<j\leq n}\frac{ g_{ij}}{\sqrt n}\sigma_i\sigma_j+\theta\sqrt{1-t}\sum_{1\leq i<j\leq n}\frac{\eta_{ij}}{\sqrt n}\sigma_i\sigma_j+\sum_{i,j\in[k]}\frac{\theta_1}{2 k}\sigma_i\sigma_j+\sum_{i=1}^n\sigma_ih_i.
\end{align*}
And wenotice that for $\varphi(t):=\bb E[\la f\ra_t]$, by similar vein of \eqref{derivativeub} we have for all $t\in[0,1]$
\begin{align*}
    |\varphi^\prime(t)|\leq C\varphi(t),\thinspace\Rightarrow\thinspace\varphi(1)<\infty\thinspace\Leftrightarrow\thinspace\varphi(0)<\infty.
\end{align*}
Then we replace $f$ with $\exp(nt(R_{1,2}-q)^2)$, $\exp(tk(m-\mu)^2)\mbbm 1_{m>0}$, $\exp(tk(m-\mu)^2)\mbbm 1_{m<0}$, and $\exp(t(n-k)\tde m^2)$ for some constant $t$ and use theorem \ref{largcliquetailbound} to complete the proof.

\subsection{Proof of Lemma \ref{lm5.6}}
    The proof goes by ensuring the cavity method in lemma \ref{cavity}, \ref{cavityII} is valid. weconsider the following smart path as the example and the other two smart paths in lemma \ref{cavity} and \ref{cavityII} are stated analogously.
    \ttny{\begin{align*}
        -\mca H_{t,1}(\bfa\sigma):&= \frac{\theta}{\sqrt n}\sum_{i<j\leq n, i,j\neq k}\eta_{ij}\sigma_i\sigma_j+\frac{\theta\sqrt t}{\sqrt n}\sum_{i\leq n}\eta_{ik}\sigma_i\sigma_k+\theta\sqrt{1-t}z\sqrt q\sigma_k+\frac{\theta_1}{2 k}\sum_{i,j\leq k-1}\sigma_i\sigma_j+\frac{\theta_1 t}{k}\sum_{i\leq k}\sigma_i\sigma_k\\
        &+\theta_1(1-t)\mu\sigma_k+\sum_{i\leq n}h_i\sigma_i.
    \end{align*}}
    And one will use lemma \ref{lm1.1} to get (Using similar derivation as \eqref{pieces2}:
    \ttny{\begin{align*}
        \nu_{t,1}^\prime(f)&=\theta^2\bl\sum_{1\leq \ell<\ell^\prime\leq r}\nu_{t,1}(f\epsilon_\ell\epsilon_{\ell^\prime}(R_{\ell,\ell^\prime}-q))\br -r\theta^2\sum_{\ell\leq r}\nu_{t,1}(f\epsilon_\ell\epsilon_{r+1}(R_{\ell,r+1}-q))\\
        &+\theta^2\frac{r(r+1)}{2}\nu_{t,1}(f\epsilon_{r+1}\epsilon_{r+2}(R_{r+1, r+2}-q))+\theta_1\bigg(\sum_{\ell\leq r}\nu_{t,1}(f\epsilon_{\ell}(m_{\ell}-\mu))-r\nu_{t,1}(f\epsilon_{r+1}(m_{r+1}-\mu))\bigg)\\
        &+O\bl\frac{1}{n}\br\nu_{t,1}(f).
    \end{align*}}
    And we notice that the rest of the two paths in \eqref{path22} and \eqref{path33} can be derived analogously. 
    Then one will see that the moment iteration in lemma \ref{secondmoment} and lemma \ref{secondmomentpositive} also holds. And apply similar arguments as the derivation of theorem \ref{cltatline} one will complete the proof.

\subsection{Proof of Lemma \ref{overlapclt}}
Before we begin the proof of the overlapped limiting theorem, the following lemma will be noticed:
\begin{lemma}[Covariance]\label{covariance}
    The covariance between the spins in the clique and out of the clique satisfies $       \nu(m_1\tde m_1)=O\lef(\frac{1}{n\sqrt k}\rig).
$\end{lemma}
    The proof is based on the following meta lemma
    \begin{lemma}
    For a function $f:\Sigma_{k+r}\to\bb R$, define the interpolating Hamiltonian as \eqref{interpolating_hamilt}, we have that under the high/ low temperature, when $k\log k=O(n)$ there exists $C>0$ such that
    \begin{align*}
        \exp\lef(-\frac{C k\log k}{n}\rig)\bb E[\la f\ra_0]\leq\bb E[\la f\ra_1]\leq\exp\lef(\frac{Ck\log k}{n}\rig)\bb E[\la f\ra_0].
    \end{align*}
    And at the critical temperature when $k^{\frac{2\tau}{2\tau-1}}\log^{\frac{2\tau-2}{2\tau-1}}k=O(n)$:
    \begin{align*}
         \exp\lef(-\frac{C k\log k}{n}\rig)\bb E[\la f\ra_0]\leq\bb E[\la f\ra_1]\leq\exp\lef(\frac{Ck\log k}{n}\rig)\bb E[\la f\ra_0].
    \end{align*}
\end{lemma}
We consider the set $S^\prime=[r+1:k+r]$ and the original index set of clique to be $S:=[k]$. Then, we consider the smart path with Hamiltonian given by \eqref{interpolating_hamilt}. Define $\varphi(t):=\bb E[\la f\ra_t]$ for some function $f:\Sigma_{k+r}\to\bb R$. First wenotice that the local replica concentration has already been derived in lemma \ref{lmoverlappedreplicacon}. Then, weuse similar derivation as \eqref{truncate1} and \eqref{crossterms1} to get that at high temperature and low temperature, for all $t\in[0,1]$:
\begin{align*}
    \varphi^\prime(t)\leq\frac{k\log k}{n}(\varphi(t)+\varphi(0))\quad\Rightarrow \quad\varphi(t)\leq\exp\bl\frac{k\log k}{n}t\br\varphi(0).
\end{align*}
And at the critical temperature with flatness $\tau$ we have for all $t\in[0,1]$:
\begin{align*}
    \varphi^\prime(t)\leq\frac{k^{\frac{2\tau}{2\tau-1}}\log^{\frac{2\tau-2}{2\tau-1}}k}{n}\varphi(t)\quad\Rightarrow\quad\varphi(t)\leq\exp\bl\frac{k^{\frac{2\tau}{2\tau-1}}\log^{\frac{2\tau-2}{2\tau-1}}k}{n}t\br\varphi(0).
\end{align*}
And we replace $f$ with the mgf of $\sum_{i=r+1}^{k+r}\sigma_i$ and notice that the results at $t=0$ corresponds to the RFCW model which we already proved in the proof of Corollary 3.2.1, 3.4.1 and 3.7.1. in \citep{he}.

\subsection{Proof of Lemma \ref{largecliqueweakclt}}
    We first prove the first part of the lemma regarding the high temperature regime. It is noticed that by lemma \ref{covariance}, 
    \begin{align*}
        \bb E\bigg[\frac{1}{k}\bl\sum_{i=r+1}^{k+r}\sigma_i\br^2\bigg]&=\bb E\bigg[\frac{1}{k}\bl\sum_{i=r+1}^k\sigma_i\br^2\bigg]+\bb E\bigg[\frac{1}{k}\bl\sum_{i=k+1}^{k+r}\sigma_i\br^2\bigg]+O\lef(\frac{1}{\sqrt n}\rig)\\
        &=1+k(1-c)\lef(1-c-\frac{1}{k}\rig)\bb E[\sigma_1\sigma_2]+c\lef(c-\frac{1}{k}\rig)\bb E[\sigma_{k+1}\sigma_{k+2}].
    \end{align*}
    wealso notice that by the convergence of moments given by theorem \ref{cltatline}, one will get
    \begin{align*}
        1+(k-1)\bb E[\sigma_1\sigma_2]=\ca V+o(1),\qquad 1+(n-k-1)\bb E[\sigma_{k+1}\sigma_{k+2}]=1+o(1).
    \end{align*}
    Therefore we have
    \begin{align*}
        \bb E\bigg[\frac{1}{k}\bl\sum_{i=r+1}^{k+r}\sigma_i\br^2\bigg]=(1-c)\ca V+c+o(1).
    \end{align*}
    Furthermore we can apply this procedure to higher moment inductively achieve that $\forall r\in\bb N$ not growing with $k$, the following holds
    \begin{align}\label{convergenceofpartialmoment}
        \bb E\bigg[\bl\frac{1}{\sqrt k\sqrt{(1-c)\ca V}}\sum_{i=r+1}^{k}\sigma_i\br^{2r}\bigg]\to\bb E[z^r]. 
    \end{align}
    Then, we have for all $r\in\bb N$, 
    \begin{align*}
        \bb E\bigg[\bl\sum_{i=r+1}^k\sigma_i\br^{2r}\bigg]\leq \bb E\bigg[\bl\sum_{i=1}^k\sigma_i\br^{2r}\bigg].
    \end{align*}
    Hence wefinally obtain that for some $t>0$:
    \begin{align*}
        \bb E\bigg[\exp\bl t\bl\frac{1}{\sqrt k}\sum_{i=r+1}^k\sigma_i\br^2\br\bigg]\leq \bb E\bigg[\exp\bl t\bl\frac{1}{\sqrt k}\sum_{i=1}^k\sigma_i\br^2\br\bigg]\quad\Rightarrow\quad \bigg\Vert\frac{1}{\sqrt k}\sum_{i=r+1}^k\sigma_i\bigg\Vert_{\psi_2}<\infty.
    \end{align*}
    And analogously we can derive that $\lef\Vert\frac{1}{\sqrt k}\sum_{i=k+1}^{k+r}\sigma_i\rig\Vert_{\psi_2}<\infty $. Finally we have by the sub-additivity of the Orlicz norm
    \begin{align*}
        \bigg\Vert\frac{1}{\sqrt k}\sum_{i=r+1}^{r+k}\sigma_i \bigg\Vert_{\psi_2}<\infty.
    \end{align*}
    For the low temperature regime, first wenotice that
    \tny{\begin{align*}
        \bb E\bigg[\frac{1}{k^2}\bl\sum_{i=1}^k\sigma_i\br^2\bigg]=\frac{1}{k}+\frac{k-1}{k}\bb E[\sigma_1\sigma_2],\thinspace \bb E\bigg[\frac{1}{k^2}\bl\sum_{i=r+1}^k\sigma_i\br^2 \bigg]=\frac{1}{k}+(1-c)\lef(1-c-\frac{1}{k}\rig)\bb E[\sigma_1\sigma_2].
    \end{align*}}
    We notice $\frac{1}{k}\sum_{i=1}^k\sigma_i\leq 1$, and
    \begin{align*}
        \bb E\bigg[\frac{1}{k}\bigg|\sum_{i=1}^k\sigma_i\bigg|\bigg]-\bb E\bigg[\frac{1}{k}\bigg|\sum_{i=r+1}^{k}\sigma_i\bigg|\bigg]\geq \frac{1}{2}\bl\bb E\bigg[\frac{1}{k}\bigg(\sum_{i=1}^k\sigma_i\bigg)^2\bigg]-\bb E\bigg[\frac{1}{k}\bigg(\sum_{i=r+1}^{k}\sigma_i\bigg)^2\bigg] \br\asymp 1.
    \end{align*}
    Further notice that by theorem \ref{cltatline}
    \sm{\begin{align*}
        \bb E\bigg [\frac{1}{k}\bigg |\sum_{i=r+1}^k\sigma_i\bigg|\bigg]\leq\frac{1}{k}\bb E\bigg[\bl\sum_{i=r+1}^k\sigma_i\br^2\bigg]^{\frac{1}{2}}\quad\Rightarrow\quad \bb E\bigg[\frac{1}{k}\bigg|\sum_{i=r+1}^{k+r}\sigma_i\bigg|\bigg]=\bb E\bigg[\frac{1}{k}\bigg|\sum_{i=r+1}^{k}\sigma_i\bigg|\bigg]+O\bl\frac{1}{k}\br.
    \end{align*}}
    And we finish the proof of the moment inequality. Then weanalyze the tail bound. Notice that in a similar vein as \eqref{convergenceofpartialmoment} we have for all $r\in\bb N$ not growing with $k$:
    \tny{\begin{align*}
        \bb E\bigg[\frac{1}{k^r}\bl\sum_{i=r+1}^k\sigma_i-(k-r)\bb E[\sigma_i|m_1>0]\br^{2r}\bigg|m_1>0\bigg]\lesssim\bb E\bigg[\frac{1}{k^r}\bl \sum_{i=1}^k\sigma_i-k\bb E[\sigma_i|m_1>0]\br^{2r}\bigg|m_1>0\bigg].
    \end{align*}}
    And wefinally get 
    \begin{align*}
        \bigg\Vert\frac{1}{\sqrt k}\bl\sum_{i=r+1}^k\sigma_i-(k-r)\bb E[\sigma_i|m_1>0]\br\bigg|m_1>0\bigg\Vert_{\psi_2}<\infty,\\
        \bigg\Vert\frac{1}{\sqrt k}\bl\sum_{i=r+1}^k\sigma_i-(k-r)\bb E[\sigma_i|m_1<0]\br\bigg|m_1<0\bigg \Vert_{\psi_2}<\infty.
    \end{align*}
    Therefore we finally complete the proof.
\section{Proof of Auxiliary Lemmas}
In this section of the appendix we provide proof of major technical lemmas.

\subsection{Proof of Lemma \ref{zeropointrpk}}
Without loss of generality we can define $D=S\cup S^\prime=[k+r]$, $S=[k]$ and $S^\prime=[r+1:k+r]$.
 We introduce the notation of $h_i^\prime := h_i+\theta\sqrt {q} z_i$ for $z_i\sim N(0,1)$ i.i.d. and checked that the following holds
\ttny{\begin{align*}
    &\nu_0(\exp(u(k+r)(R^D_{1,2}-q)))=\nu_0\bigg(\exp\bigg(u\sum_{i\in[k+r]}(\sigma_i^1\sigma_i^2-q)\bigg)\bigg)\\
    &=\bb E\lef[\frac{\sum_{\bfa\sigma^1,\bfa\sigma^2}\exp\lef(u\sum_{i\in[r+k]}(\sigma_i^1\sigma_i^2-q)+\sum_{i\in[k]}\frac{\theta_1}{2k}\sigma^1_i\sigma^1_j+\sum_{i\in[r+1:k+r]}\frac{\theta_1}{2k}\sigma^2_i\sigma^2_j+\sum_{i\in[k+r]}h_i^\prime(\sigma^1_i+\sigma_i^2)\rig)}{\sum_{\bfa\sigma^1,\bfa\sigma^2}\exp\lef(\sum_{i\in[k]}\frac{\theta_1}{2k}\sigma^1_i\sigma^1_j+\sum_{i\in[r+1:k+r]}\frac{\theta_1}{2k}\sigma^2_i\sigma^2_j+\sum_{i\in[k+r]}h_i^\prime(\sigma^1_i+\sigma_i^2)\rig)}\rig]\\
    &=\bb E\lef[\frac{\exp(-u(k+r)q)\int\int\exp(-k\mca H_{0,k}(u,x,y,\bfa h))dxdy}{2^{2(r+k)}\int\exp(-k\mca H_{1,k}(x,\bfa h))dx\int\exp(-k\mca H_{2,k}(x,\bfa h))dx}\rig]
.\end{align*}}
And we give reason going from the second line to the third line. This comes from the Laplace approximation. We note that the denominator can be written as
\ttny{\begin{align*}
&\sum_{\bfa\sigma^2_{[r+1:k+r]}}\exp\bigg(\frac{\theta_1}{2k}\bigg(\sum_{i=r+1}^{k+r}\sigma_i^2\bigg)^2+\sum_{i=r+1}^{k+r}h^\prime_i\sigma_i^2\bigg)\sum_{\bfa\sigma^1_{[k]}}\exp\bigg(\frac{\theta_1}{2k}\bigg(\sum_{i\in[k]}\sigma_i^1\bigg)^2+\sum_{i\in[k]}\sigma_i^1h^\prime_i\bigg)\prod_{i\in[r]\cup[k+1:k+r]}\cosh(h^\prime_i)\\
&=\sum_{\bfa\sigma^2_{[r+1:k+r]}}\frac{1}{2\pi}\int\exp\bl-\frac{x^2}{2}+x\sqrt{\frac{\theta_1}{k}}\sum_{i=r+1}^{k+r}\sigma_i^2+\sum_{i=r+1}^{k+r}h^\prime_i\sigma_i^2 \br dx\\
&\cdot\sum_{\bfa\sigma^1_{[k]}}\int\exp\bl -\frac{x^2}{2}+x\sqrt{\frac{\theta_1}{k}}\sum_{i\in[k]}\sigma_i^1+\sum_{i\in[k]}\sigma_i^1h^\prime_i\br dx
\prod_{i\in[r]\cup[k+1:k+r]}\cosh(h^\prime_i)\\
&=\frac{2^{2(r+k)}k}{2\pi}\int\exp( -k \mca H_{1,k}(x,\bfa h))dx\int\exp(-k\mca H_{2,k}(x,\bfa h))dx
.\end{align*}}
where 
\begin{align*}
\mca H_{1,k}(x,\bfa h)&:=\frac{x^2}{2}-\frac{1}{k}\sum_{i=r+1}^{k+r}\log\cosh(\sqrt{\theta_1} x+h^\prime_i)-\frac{1}{k}\sum_{i\in[r]}\cosh(h^\prime_i),\\
\mca H_{2,k}(x,\bfa h)&:=\frac{x^2}{2}-\frac{1}{k}\sum_{i\in[k]}\log\cosh(\sqrt{\theta_1} x+h^\prime_i)-\frac{1}{k}\sum_{i=k+1}^{k+r}\cosh(h^\prime_i)
.\end{align*} For the nominator, we can check that
\ttny{\begin{align*}
    &\sum_{\bfa\sigma^1,\bfa\sigma^2}\exp\bl u\sum_{i\in[k+r]}\sigma_i^1\sigma_i^2+\frac{\theta_1}{k}\bl\sum_{i\in[k]}\sigma_i^1\br^2+\frac{\theta_1}{k}\bl \sum_{i=r+1}^{k+r}\sigma_i^2\br^2+\sum_{i\in[k+r]}h^\prime_i(\sigma_i^1+\sigma_i^2)\br\\
    &=\frac{1}{2\pi}\sum_{\bfa\sigma^1,\bfa\sigma^2}\int\int\exp\bl -\frac{x^2+y^2}{2}+u\sum_{i\in[k+r]}\sigma_i^1\sigma_i^2+\sqrt{\frac{\theta_1}{k}}\sum_{i\in[k]}\sigma_i^1 x +\sqrt{\frac{\theta_1}{k}}\sum_{i=r+1}^{k+r}\sigma_i^2y+\sum_{i\in[k+r]}h^\prime_i(\sigma_i^1+\sigma_i^2)\br dxdy\\
    &=\frac{2^{k+r}}{2\pi}\int\int\exp\bl -\frac{x^2+y^2}{2}\br\sum_{\bfa\sigma^1}\prod_{i\in[r+1:k+r]}\cosh\lef(u\sigma_i^1+\sqrt{\theta}y+h^\prime_i\rig)\\
 &\cdot\prod_{i\in[r]}\cosh\lef(u\sigma_i^1+h^\prime_i\rig)\exp\bl\sqrt{\frac{\theta_1}{k}}\sum_{i\in[k]}\sigma_i^1x+\sum_{i\in[k+r]}h^\prime_i\sigma_i^1\br dydx\\
 &=\frac{k}{2\pi}\int\int\exp\lef(-k\frac{x^2+y^2}{2}\rig)\prod_{i\in[r]}f(u,h^\prime_i,\sqrt{\theta_1} x+h^\prime_i)\\
 &\cdot\prod_{i\in[r+1:k]}f(u,\sqrt{\theta_1} y+h^\prime_i,\sqrt{\theta_1} x+h^\prime_i)\prod_{i\in[k+1:k+r]}f(u,\sqrt{\theta_1} y+h^\prime_i,h^\prime_i)dxdy\\
 &=\frac{k}{2\pi}\int\int \exp(-k\mca H_{0,k}(u,x,y,\bfa h))dxdy
.\end{align*}}
where we define 
\begin{align*}
    f(a,b,d):&=2\exp(a)\cosh(b+d)+2\exp(-a)\cosh(b-d)\\
    &=4\cosh(b)\cosh(d)\cosh(a)+4\sinh(b)\sinh(d)\sinh(a),
\end{align*} and
\ttny{\begin{align*}
    \mca H_{0,k}(u,x,y,\bfa h)&=\frac{x^2+y^2}{2}-\frac{1}{k}\sum_{i=r+1}^{k}\log f(u,\sqrt{\theta_1} y+h^\prime_i,\sqrt{\theta_1} x+h^\prime_i)-\frac{1}{k}\sum_{i\in[r]}\log f(u,h^\prime_i,\sqrt{\theta_1}x+h^\prime_i)\\
    &-\frac{1}{k}\sum_{i=k+1}^{k+r}\log f(u,\sqrt{\theta_1}y+h^\prime_i,h^\prime_i)
.\end{align*}}
Introduce $c:=\frac{r}{k}$ and we defining the population variants of $\mca H_{0,k}$ and $\mca H_{1,k}, \mca H_{2,k}$ by
\begin{align*}
    \mca H_0(u,x,y) &= \frac{x^2+y^2}{2}-c\bb E[(\log f(u,h^\prime,\sqrt{\theta_1} x+h^\prime)+\log f(u,\sqrt{\theta_1} y+h^\prime, h^\prime))]\\
    &-(1-c)\bb E[\log f(u,\sqrt{\theta_1} y+h^\prime,\sqrt{\theta_1} x+h^\prime)],\\
    \mca H_1(x)&=\frac{x^2}{2}-(1-c)\bb E[\log\cosh(\sqrt{\theta_1} x+h^\prime)]-c\bb E[\log\cosh(h^\prime)]=\mca H_2(x)
.\end{align*}
The following lemma gives uniform convergence guarantees for $\mca H_{0,k}, \mca H_{1,k}$, and $\mca H_{2,k}$ to their respective limit. Assume that $\mca H_{0}(u,x,y)$, $\mca H_1$, $\mca H_2$ have a sequence of stationary point $(x^{(i)},y^{(i)}), x_1^{(i)}, x_2^{(i)}$ for $i\in[2]$  given $u$ and . Then, $\mca H_{0,k}$ has a sequence of stationary point $(x_{k}^{(i)},y_{k}^{(i)})_{k\in\bb N}$ of $\mca H_{0,k}$ satisfying $(x_{k}^{(i)},y_{k}^{(i)})\to(x^{(i)},y^{(i)})$. Similar also holds when assuming $\mca H_{1,k}(x)$,$\mca H_{2,k}(x)$  have series of fixed point $x^{(i)}_{1,k}\to x^{(i)}_1$ and  $x_{2,k}^{(i)}\to x^{(i)}_2$ respectively. And moreover we note that as $u\to 0$ we will have $(x^{(i)},y^{(i)})\to (x_1^{(i)},x_2^{(i)})$ by uniform convergence in $u$. Therefore in what follows we omit the super script $(i)$ since the results hold for all minimum point (since the maximum point is trivial solution) and replace instead with a $*$. 
\begin{lemma}[Regularity Conditions]\label{convergeas}
 Almost surely in $\mu(\bfa h^\prime)$ and uniformly on $(u,x,y)$ we will have 
        \begin{align}\label{laplace1}
            \mca H_{0,k}^{(j_1,j_2,j_3)}(s_1,s_2,s_3,\bfa h^\prime):=\frac{\pta^{\sum_ij_i}\mca H_{0,k}( s_1,s_2,s_3,\bfa h^\prime)}{\prod_{i\in[3]}\pta s_i^{j_i} }\to \mca H_0^{(j_1,j_2,j_3)}(\bfa s)
        .\end{align} with $\mca H_{0,k}^{(0,0,0)}=\mca H_{0,k}$
         and 
         \begin{align*}
             \mca H_{1,k}^{(i)}(s,\bfa h^\prime)\to \mca H_1^{(i)}(s),\qquad\qquad \mca H_{2,k}^{(i)}(s,\bfa h^\prime)\to \mca H_2^{(i)}(s)
         .\end{align*}
    The condition \ref{laplace2} and \ref{laplace3} holds for in lemma \ref{laplace} holds for 
    \begin{align*}
        \mca H_{0,k}\to \mca H_0\qquad \mca H_{1,k}\to \mca H_1\qquad \mca H_{2,k}\to \mca H_2=\mca H_1
    .\end{align*} with parameter set $\bfa s=(x,y)$ and $\bfa s=(x)$ respectively. Moreover we have $\mca H_{0,k}\to \mca H_0$ uniformly over $u\in\bb R$.
\end{lemma}
For the first condition, we define
\tny{\begin{align*}
    \varphi_k(u,x,y,\bfa h^\prime)&:=-\frac{1}{k}\bl\sum_{i\in[r]}\log f(u,h^\prime_i,\sqrt{\theta_1} x+h^\prime_i) +\sum_{i\in[r+1:k]}\log f(u,\sqrt{\theta_1} y+h^\prime_i,\sqrt{\theta_1} x+h^\prime_i) \\
    &+\sum_{i\in[k+1:k+r]}\log f(u,\sqrt{\theta_1} y+h^\prime_i,h^\prime_i)\br,
\end{align*}}
and
\begin{align*}
    \varphi(u,x,y)&:=-c\bb E[(\log f(u,h,\sqrt{\theta_1} x+h^\prime)+\log f(u,\sqrt{\theta_1} y+h^\prime, h^\prime))]\\
    &-(1-c)\bb E[\log f(u,\sqrt{\theta_1} y+h^\prime,\sqrt{\theta_1} x+h^\prime)]
.\end{align*}
It is not hard to see that by SLLN almost surely we have $\varphi_k(u,x,y,\bfa h^\prime)\to\varphi(u,x,y)$ point-wise. 
Then we can check that for $(u_1,x_1,y_1)$ and $(u_2,x_2,y_2)\in\bb R^3$ we will have
\begin{align*}
    \varphi_k(u_1,x_1,y_1,\bfa h^\prime)-\varphi_k(u_2,x_2,y_2,\bfa h^\prime)\leq 2\lef(|u_1-u_2|+\sqrt{\theta_1}|x_1-x_2|+\sqrt{\theta_1}|y_1-y_2|\rig),\forall n
.\end{align*}
implies that $\varphi_k$ form an uniformly equicontinuous sequence. Since countable intersection of sets with measure $1$ has also measure $1$ we conclude that it is possible to choose $A\subset\Omega$ such that $\mu(A)=1$ such that $\forall \bfa h^\prime\in A$,$\varphi_k(u,x,y,\bfa h^\prime)\to\varphi(u,x,y)$. This implies that $\mca H_{0,k}\to \mca H_0$ uniformly almost surely (A simple exercise using Arzelà–Ascoli theorem). Similar argument can be easily verified to hold for $\mca H_{1,k}$ and we omit it here.

Then we move toward the discussion over the derivatives. Since we verified that the derivatives of $\varphi^{(i,j,n)}(u,x,y,\bfa h^\prime):=\frac{\pta^{i+j+n}\varphi(u,x,y,\bfa h^\prime)}{\pta u^i\pta x^j\pta y^n}$ is bounded according to \ref{circderivative} and seeing that $|h|,|l|,|g|\leq |f|$. Therefore, we conclude that $\varphi^{(i,j,n)}(u,x,y,\bfa h^\prime):=\frac{\pta^{i+j+n}\varphi(u,x,y,\bfa h^\prime)}{\pta u^i\pta x^j\pta y^n}$ is equicontinuous and hence uniformly almost surely converging to $\varphi^{(i,j,n)}(u,x,y):=\frac{\pta^{i+j+n}\varphi(u,x,y)}{\pta u^i\pta x^j\pta y^n}$. This implies that the derivatives also converges unformly almost surely. Similar arguments can be analogously applied to $\mca H_1$.

For the second condition, noticing that $\log f(a,b,c)\leq 2\log2 +|a|+|b|+|c|$ we see that:
\begin{align*}
-\varphi_k(u,x,y,\bfa h^\prime)&\leq \frac{1}{k}\sum_{i=1}^{k+r}\lef(|u|+2|h^\prime_i|\rig)+\sqrt{\theta_1}|x|+\sqrt{\theta_1}|y|+4\log 2\\
&\leq 2|u|+2|x|+2|y|+4\log 2+\frac{2}{k+r}\sum_{i\in[k+r]}|h^\prime_i|
.\end{align*}
which consequently shows that
\begin{align*}
    \mca H_{0,k}(u,x,y,\bfa h^\prime)\geq\frac{x^2+y^2}{2}-2|u|-2|x|-2|y|-4\log 2-\frac{2}{k+r}\sum_{i\in[k+r]}|h^\prime_i|
.\end{align*}
We denote  $\tau =2$ and $C(\bfa h^\prime)=16\exp\lef(\frac{2}{k+r}\sum_{i\in[k+r]}|h^\prime_i|\rig)$ it is checked that by dominated convergence theorem and $h$ is in $L_1$ that for all $u<\infty$ we will have
\tny{\begin{align*}
    \int\exp(-\mca H_0(u,x,y))dxdy&=\lim_{k\to\infty}\int\exp\lef(-\mca H_{0,k}(u,x,y,\bfa h^\prime)\rig)dxdy\\
    &\leq \exp\lef(\int_{\bb R}2|h|d\mu(h)\rig)\int_C\exp\lef(-\frac{x^2+y^2}{2}+2|u|+2|x|+2|y|\rig)dxdy\\
    &\leq A\exp(2|u|)<\infty
.\end{align*}}
for some constant $A$ not dependent on $x,y,u$
Similar argument also holds for $\mca H_{1}$ and we complete the proof.

In what follows we will omit the dependence between  $\mca H_{0,k},\mca H_{1,k}$ and $\bfa h^\prime$ just for notation simplicity. We will rewrite the target as
\ttny{\begin{align}\label{corrsubgaus}
    &\nu_0(\exp(u(k+r)(R^D_{1,2}-q)))\nnb\\
    &=\bb E\lef[\frac{\exp(-k\mca H_{0,k}(u,x_k,y_k)+k\mca H_{1,k}(x_{1,k})+k\mca H_{2,k}(x_{2,k}))\int\exp(k(\mca H_{0,k}(u,x_k,y_k)-\mca H_{0,k}(x,y)))dxdy}{2^{2(r+k)}\int\exp(k(\mca H_{1,k}(x_{1,k})-\mca H_{1,k}(x)))dx\int\exp(k(\mca H_{2,k}(x_{2,k})-\mca H_{2,k}(x)))dx}\rig]
.\end{align}}
Then we need to estimate the difference between $\mca H_{0,k}(u,x_k,y_k)$ and $\mca H_{0,k}(0,x_{1,k},x_{2,k})$. We then introduce $g(a,b,d):=\frac{\pta f}{\pta b} = 2\exp(a)\sinh(b+d)+2\exp(-a)\sinh(b-d)$ and $(x^*(u),y^*(u))\in\argmin_{x,y}\mca H_0(u,x,y)$, $h(a,b,d):=2\exp(a)\cosh(b+d)-2\exp(-a)\cosh(b-d)$ and $l(a,b,d):= 2\exp(a)\sinh(b+d)-2\exp(-a)\sinh(b-d)$. Finally recall that $f(a,b,d):=2\exp(a)\cosh(b+d)+2\exp(-a)\cosh(b-d)$. For some expression denoted by $f_0(f,g,h,l)$ we write
$$\tau(p_1,p_2,p_3):=f_0(f(p_1,p_2,p_3),g(p_1,p_2,p_3),h(p_1,p_2,p_3),l(p_1,p_2,p_3)).$$
with the further simplification of
\begin{align*}
\bar\tau_{1}&=\frac{1}{k}\sum_{i=r+1}^k\tau(u,\sqrt{\theta_1}y+h^\prime_i,\sqrt{\theta_1}x+h^\prime_i),\bar\tau_{2}=\frac{1}{k}\sum_{i=1}^{r}\tau(u,\sqrt{\theta_1}x+h^\prime_i,h^\prime_i),\\
\bar\tau_{3}&=\frac{1}{k}\sum_{i=k+1}^{k+r}\tau(u,h^\prime_i,\sqrt{\theta_1}y+h^\prime_i)
.\end{align*}
We will have the following relations regarding the above quantities:
\begin{align}\label{circderivative}
    &\frac{\pta f}{\pta a} =h,\quad \frac{\pta f}{\pta b}= g,\quad \frac{\pta f}{\pta d}= l,\quad\frac{\pta h}{\pta a}=f,\quad\frac{\pta h}{\pta b} =l,\quad\frac{\pta h}{\pta d}=g,\nnb\\&\frac{\pta g}{\pta a} =l,\quad\frac{\pta g}{\pta b}=f,\quad\frac{\pta g}{\pta d}=h,\quad\frac{\pta l}{\pta a}=g,\quad\frac{\pta l}{\pta b}=h,\quad\frac{\pta l}{\pta d} =f
.\end{align}
And we can write the Hessian of $\mca H_{0,k}$ as
\ttny{\begin{align*}
\begin{bmatrix}
-1-c+\overline{\lef(\frac{h}{f}\rig)_3^2}+\overline{\lef(\frac{h}{f}\rig)_2^2}+\overline{\lef(\frac{h}{f}\rig)_1^2}&\sqrt{\theta_1}\lef(-\overline{\lef(\frac{gf-hl}{f^2}\rig)_3}-\overline{\lef(\frac{gf-hl}{f^2}\rig)_1}\rig)&\sqrt{\theta_1}\lef(-\overline{\lef(\frac{lf-gh}{f^2}\rig)_2}-\overline{\lef(\frac{lf-gh}{f^2}\rig)_1}\rig)\\
\sqrt{\theta_1}\lef(-\overline{\lef(\frac{gf-hl}{f^2}\rig)_3}-\overline{\lef(\frac{gf-hl}{f^2}\rig)_1}\rig)&1-\theta_1+\theta_1\lef(\overline{\lef(\frac{l}{f}\rig)^2_3}+\overline{\lef(\frac{l}{f}\rig)_1^2}\rig)&-\theta_1\overline{\frac{hf-gl}{f^2}_1}\\
\sqrt{\theta_1}\lef(-\overline{\lef(\frac{lf-gh}{f^2}\rig)_2}-\overline{\lef(\frac{lf-gh}{f^2}\rig)_1}\rig)&-\theta_1\overline{\frac{hf-gl}{f^2}_1}&1-\theta_1+\theta_1\lef(\overline{\lef(\frac{g}{f}\rig)^2_2}+\overline{\lef(\frac{g}{f}\rig)_1^2}\rig)
\end{bmatrix}
.\end{align*}}

Upon observing that $\lef(\frac{h}{f}\rig)^2\leq 1$ we will have
$
    -1-c\leq\frac{\pta^2 \mca H_{0,k}}{\pta u^2}\leq 0.
$ We further define
\begin{align}\label{etah}
    &\eta(\bfa h):=-\frac{1}{1+c}\frac{d \mca H_{0,k}(0,x_{k}(0),y_{k}(0))}{d u} =-\frac{1}{1+c}\frac{\pta \mca H_{0,k}(0,x_{1,k},x_{2,k})}{\pta u} \nnb\\
    &=\frac{1}{k+r}\bl \sum_{i=r+1}^k\tanh(h_i^\prime+\sqrt{\theta_1}x_{1,k})\tanh(h^\prime_i+\sqrt{\theta_1}x_{2,k})\nnb\\
    &+\sum_{i=1}^r\tanh(h_i^\prime)\tanh(h_i^\prime+\sqrt{\theta_1}x_{1,k})+\sum_{i=k+1}^{k+r}\tanh(h_i^\prime)\tanh(h^\prime_i+\sqrt{\theta_1}x_{2,k})\br
.\end{align}
Using the fact that at all $(u,x_k,y_k)$ we will have by Fermat's condition:
\begin{align*}
    \frac{\pta \mca H_0(u,x_k,y_k)}{\pta x}=0=\frac{\pta \mca H_0(u,x_k,y_k)}{\pta y}=0
.\end{align*}
And we also have
\begin{align*}
    \frac{d^2\mca H_{0,k}(u,x_k,y_k)}{du^2}&=\frac{\pta^2 \mca H_{0,k}(u,x_k,y_k)}{\pta u^2}+\frac{\pta^2\mca H_{0,k}(u,x_k,y_k)}{\pta u\pta x}\frac{dx}{du} +\frac{\pta^2\mca H_{0,k}(u,x_k,y_k)}{\pta u\pta y}\frac{dy}{du}\\&=\frac{\pta^2\mca H_{0,k}(u,x_k,y_k)}{\pta u^2}
.\end{align*}
which implies that 
\begin{align*}
    \mca H_{0,k}(u, x_k,y_k)-\mca H_{0,k}(0,x_{1,k},x_{2,k})&=\int_0^u\lef(\frac{\pta \mca H_{0,k}}{\pta u}+\frac{\pta \mca H_{0,k}}{\pta x}\frac{dx}{du}+\frac{\pta \mca H_{0,k}}{\pta y}\frac{dy}{du}\rig)du \\&\geq (-1-c)\lef(\frac{1}{2}u^2+\eta(\bfa h) u\rig)
.\end{align*}
Hence we will have
\sm{\begin{align*}
    \exp\lef(-(k+r)uq-k(\mca H_{0,k}(u,x_k,y_k)-\mca H_{0,k}(0,0,0))\rig)\leq\exp\lef((k+r)\lef(\frac{1}{2}u^2+(\eta(\bfa h)-q)u\rig)\rig)
.\end{align*}}
For the second term in \ref{corrsubgaus} we notice that $
    \mca H_{1,k}^{(2)}(x_{1,k})=1-\theta_1+\theta_1\frac{1}{k}\sum_{i\in[k]}\tanh(h_i^\prime)\tanh(h^\prime_i+x_{1,k})
$ and $
    \mca H_{2,k}^{(2)}(x_{2,k})=1-\theta_1+\theta_1\frac{1}{k}\sum_{i\in[r+1:k]}\tanh(h_i^\prime)\tanh(h^\prime_i+x_{2,k})
$. Using lemma \ref{convergeas} we will notice that almost surely :
\begin{align*}
\mca H_{1,k}^{(2)}(x_{1,k})\to \mca H_{1}(0)=1-\theta_1+\theta_1\bb E[\tanh^2(\sqrt{\theta_1}x_1^*+h^\prime)],\\
 \mca H_{2,k}^{(2)}(x_{1,k})\to \mca H_{1}(0)=1-\theta_1+\theta_1\bb E[\tanh^2(\sqrt{\theta_1}x_2^*+h^\prime)]
.\end{align*}
Using the fact that by Taylor expansion we will have
\tny{\begin{align*}
     (hf-gl)(a,b,d) &= 4\exp(2a)\lef(\cosh^2(b+d)-\sinh^2(b+d)\rig) -4\exp(-2a)\lef(\cosh^2(b-d)-\sinh^2(b-d) \rig)\\
    &=8\sinh(2a)=16 a+O(a^2)
.\end{align*}}
Then we see that Taylor expand around $0$ for $u$ we will have
\begin{align*}
  \bl\lef(\overline{\frac{hf-gl}{f}}\rig)_1+\lef(\overline{\frac{hf-gl}{f}}\rig)_3\br^2 = O(u^2).
\end{align*}
And, note that for $\lef(\frac{l}{f}\rig)^2,\lef(\frac{g}{f}\rig)^2$ we will have when $u=O(1)$:
\begin{align*}
    &\overline{\lef(\frac{g}{f}\rig)^2_2}+\overline{\lef(\frac{g}{f}\rig)_1^2}=\frac{1}{k}\sum_{i=1}^{k}\tanh^2(\sqrt{\theta_1}x_k+h^\prime_i)+O(u),\\
    &\overline{\lef(\frac{l}{f}\rig)^2_3}+\overline{\lef(\frac{l}{f}\rig)_1^2}=\frac{1}{k}\sum_{i=r+1}^{k+r}\tanh^2(\sqrt{\theta_1}y_k+h^\prime_i)+O(u).
\end{align*}
To analyze the second term in \eqref{corrsubgaus}, our analysis will be divided into two cases and treated individually. (1) High/Low Temperature Regime with $\mca H_{1,k}^{(2)}(x_1^*)>0$. (2) Critical Temperature Regime with $\mca H_{1,k}^{(i)}(x_1^*)=0$ for all $i\leq 2\tau-1$ and $H^{(2\tau)}_{1,k}(x_1^*)>0$.
\begin{center}
    \textbf{I: High/Low Temperature Regime }
\end{center}
First we consider the high/low temperature regime, notice that the \eqref{corrsubgaus} can be rewritten as according to lemma \ref{laplace}:
\ttny{\begin{align*}
    \bb E\bigg[\exp\lef(-(k+r)q-k\lef(\mca H_{0,k}(u,x_k,y_k)-\mca H_{0,k}(0,x_{1,k},x_{2,k})\rig)\rig)\frac{(\mca H_{1,k}^{(2)}(x_{1,k})^{1/2}\mca H_{2,k}^{(2)}(x_{2,k})^{1/2})}{ \det\lef(\nabla^2_{xy}\mca H_0(u,x_k,y_k)\rig)^{1/2}}\bigg]\lef(1+O\lef(\frac{1}{k}\rig)\rig)
\end{align*}} where we use $\nabla^2_{xy}H:=\begin{bmatrix}
    \frac{\pta^2 H}{\pta x^2}&\frac{\pta^2 H}{\pta x\pta y}\\
    \frac{\pta^2 H}{\pta x\pta y}&\frac{\pta^2 H}{\pta y^2}
\end{bmatrix}$ to denote the sub matrix of Hessian containing only derivatives of $x$ and $y$.
And we have
\tny{\begin{align*}
    \sqrt{k}\lef(x_{1,k}-x_1^*\rig)=\frac{1}{\sqrt k\mca H_{1,k}^{(2)}(x_1^*,\bfa h)}\sum_{i=1}^k\lef(\bb E[\tanh(\sqrt{\theta_1}x_1^*+h_i^\prime)]-\tanh(\sqrt{\theta_1}x_1^*+h_i^\prime)\rig)+o_{\psi_2}(1),\\
    \sqrt{k}(x_{2,k}-x_2^*)=\frac{1}{\sqrt k\mca H_{2,k}^{(2)}(x_1^*,\bfa h)}\sum_{i=r+1}^{k+r}\lef(\bb E[\tanh(\sqrt{\theta_1}x_2^*+h_i^\prime)]-\tanh(\sqrt{\theta_1}x_2^*+h_i^\prime)\rig)+o_{\psi_2}(1).
\end{align*}}
Then we will also have
\ttny{\begin{align*}
    \det\lef(\nabla^2_{xy}\mca H_0(u,x_k,y_k)\rig) = \lef(1-\theta_1\sum_{i=1}^k\frac{1}{k}\sech^2(\sqrt{\theta_1}x_k+h^\prime_i)\rig)\lef(1-\theta_1\sum_{i=r+1}^{k+r}\frac{1}{k}\sech^2(\sqrt{\theta_1}y_k+h^\prime_i)\rig)+O\lef(u\rig)
.\end{align*}} 
This, together with,  $x^*-x_1^*=O(u), y^*-x_1^*=O(u)$ implies that
\begin{align}\label{secondorderpre}
    \frac{\mca H_{1,k}^{(2)}(x_{1,k})^{1/2}\mca H_{2,k}^{(2)}(x_{2,k})^{1/2}}{ \det\lef(\nabla^2_{xy}\mca H_0(u,x_k,y_k)\rig)^{1/2}}=1+O(u+u^2).
\end{align}
And consequently, using \ref{convergeas}, we will have that for all $u\in\bb R$ the following holds:
\tny{\begin{align}\label{mgfub}
    \nu_0(\exp(u(k+r)(R^D_{1,2}-q)))\leq \bl 1+O(u)+O(u^2)+O\bl \frac{1}{k}\br\br\bb E\lef[\exp\lef((k+r)\lef(\frac{1}{2}u^2+(\eta(\bfa h)-q^\prime)u\rig)\rig)\rig]
.\end{align}}
And we can take $U\sim N(0,\frac{2\lambda}{k+r})$ and use the following fact to conclude that
\begin{fact}\label{gaussianfact}
        When $Z\sim N(0,\sigma^2)$, then for $2a\leq\frac{1}{\sigma^2}$ and any $b\in\bb R$ we will have $\bb E[\exp\lef(aZ^2+bZ)\rig)]=\frac{1}{\sqrt{1-2a\sigma^2}}\exp\lef(\frac{\sigma^2b^2}{2(1-2a\sigma^2)}\rig)$
\end{fact}
\begin{align}\label{gstrick}
    \nu_0&(\exp(\lambda (k+r)(R^D_{1,2}-q)^2))= \nu_0\lef(\bb E_U[\exp(U (k+r)(R^D_{1,2}-q))]\rig)\nnb\\
    &\leq\bb E_{\bfa h}\lef[\bb E_{U}\lef[(1+O(U)+O(U^2))\exp\lef((k+r)\lef(\frac{1}{2}U^2+(\eta(\bfa h)-q^\prime)U\rig)\rig)\rig]\rig]\lef(1+o(1)\rig)\nnb\\
    &\leq\frac{1}{\sqrt{1-2\lambda}}\bb E_{\bfa h}\lef[\exp\lef(\frac{\lambda (k+r)}{(1-2\lambda)}\lef(\eta(\bfa h)-q\rig)^2\rig)\rig]\lef(1+o(1)\rig),
\end{align}
where the $U$ terms are evaluated by the H\"older's inequality.
And we can move back to \ref{etah}, note that by Taylor expansion at $x_1^*$ we will have
\tny{\begin{align}\label{withop}
   \sqrt{k+r}\eta(\bfa h)&=\frac{1}{\sqrt{k+r}}\bl\sum_{i=1}^{r}\tanh(h_i^\prime)\tanh(\sqrt{\theta_1}x^*_{1}+h_i^\prime)+\sum_{i=k+1}^{k+r}\tanh(h_i^\prime)\tanh(\sqrt{\theta_1}x_2^*+h_i^\prime)\nnb\\
    &+\sum_{i=r+1}^k\tanh(\sqrt{\theta_1}x_1^*+h_i^\prime)\tanh(\sqrt{\theta_1}x_2^*+h_i^\prime)\br+O_{\psi_2}\lef(\frac{1}{k}\rig)
.\end{align}}
Then we consider the low temperature $(\theta_1> \frac{1}{\bb E[\sech^2(h^\prime)]})$ and the high temperature $(\theta_1<\frac{1}{\bb E[\sech^2(h^\prime)]})$ respectively. It is checked that for non vanishing $r$ in the low temperature the replica symmetric is essentially breaking (This result is also not necessary from the perspective of utility.). And, for vanishing $r$ in the low temperature, as well as for the whole $r\leq k$ region in the high temperature case, we will have the following weak convergence result
\begin{align*}
    \sqrt{k+r}&\eta(\bfa h)\overset{d}{\to}N\lef(\frac{(1-c)E_1+2cE_2}{1+c},\frac{(1-c)V_1+2cV_2}{1+c}\rig),
\end{align*}
with $E_1:=\bb E[\tanh^2(\sqrt{\theta_1}x_1^*+h^\prime)], V_1:=\bb  V[\tanh^2(\sqrt{\theta_1}x_1^*+h^\prime)]$, $E_2:=\bb E[\tanh(\sqrt{\theta_1}x_1^*+h^\prime)\tanh(h^\prime)]$, and $V_2:=\bb V[\tanh(\sqrt{\theta_1}x_1^*+h^\prime)\tanh(h^\prime)]$.
Hence, we will pick $q=\frac{(1-c)E_1+2cE_2}{1+c}$ and evaluate the term $O_{\psi_2}(\frac{1}{k})$ in \eqref{withop} can be evaluated with H\"older's inequality. Then for arbitrarily small $\delta>0$ we will have
\begin{align*}
    \bb E\lef[\exp\lef(\frac{\lambda}{1-2\lambda}k\lef(\eta(\bfa h)-q\rig)^2\rig)\rig]\leq \frac{1+o(1)}{\sqrt{1-\frac{2\lambda(1+\delta)}{1-2\lambda}\frac{(1-c)V_1+2cV_2}{1+c}}}
.\end{align*}
Therefore we will check that for all $\lambda<\frac{1}{2\lef(1+\frac{(1-c)V_1+2cV_2}{1+c}\rig)}$,
\begin{align*}
    \nu_0(\exp(\lambda (k+r)(R_{1,2}^D-q)^2))<\infty.
\end{align*}
Notice that in the high temperature case we will have $x_1^*=x_2^*=0$ and $\frac{(1-c)V_1+2cV_2}{1+c}=\bb V[\tanh^2(h^\prime)]$.
\begin{center}
    \textbf{II: Critical Temperature with $2\tau$-th Derivative}
\end{center}
Under the Critical Temperature Case with $2\tau$ positive derivative, then we will have
\sm{\begin{align}\label{2taux}
    \sqrt{k}x_{1,k}^{2\tau-1}&=\frac{(2\tau-1)!}{\sqrt k \mca H_{1,k}^{(2\tau)}(0,\bfa h)}\sum_{i=1}^k\lef(\bb E[\tanh(\sqrt{\theta_1}x_1^*+h_i^\prime)]-\tanh(\sqrt{\theta_1}x_1^*+h_i^\prime)\rig)+o_{\psi_2}(1),\\
    \sqrt{k}x_{2,k}^{2\tau-1}&=\frac{(2\tau-1)!}{\sqrt k \mca H_{2,k}^{(2\tau)}(0,\bfa h)}\sum_{i=r+1}^{k+r}\lef(\bb E[\tanh(\sqrt{\theta_1}x_2^*+h_i^\prime)]-\tanh(\sqrt{\theta_1}x_2^*+h_i^\prime)\rig)+o_{\psi_2}(1).
\end{align}}
And for the $\mca H_{0,k}$, defining $\bfa\delta:=(x_k-x^*,y_k-y^*)^\top$, then we will have
\begin{align*}
    \nabla \mca H_{0,k}(u,x_k,y_k)&=0=\nabla \mca H_{0,k}(u,x^*,y^*)\times\bfa\delta+\nabla^2 \mca H_{0,k}(u,x^*,y^*)\times\bfa\delta^2+\ldots\\&+\nabla^{2\tau}\mca H_{0,k}(u,x^*,y^*)\times\bfa\delta^{2\tau-1}+O(\Vert\bfa\delta\Vert_2^{2\tau+1}).
\end{align*}
Similar to the derivation of the high/low temperature regime, for $i\in\bb N$ we will have $\nabla^i_x \mca H_{0,k}(u,x_k,y_k)-\nabla^i_x\mca H_{1,k}(x_{1,k})=O(u)$, $\nabla^i_y \mca H_{0,k}(u,x_k,y_k)-\nabla^i_x\mca H_{2,k}(x_{2,k})=O(u)$, and
$    \nabla^2_{xy}\mca H_{0,k}(u,x^*,y^*)=O(u^2)
$. 
To analyze the second term in \ref{corrsubgaus} we notice that by Taylor expansion,
\ttny{\begin{align*}
    -k(&\mca H_{0,k}(u,x,y)-\mca H_{0,k}(u,x_k,y_k))=k\nabla^2_{xy}\mca H_{0,k}(u,x_k,y_k)(x-x_k)(y-y_k)+k\nabla_x^{(2\tau)}\mca H_{0,k}(u,x_k,y_k)(x-x_k)^{2\tau}\\
    &+k\nabla^{(2\tau)}_{y}\mca H_{0,k}(u,x_k,y_k)(y-y_k)^{2\tau}+O(k((x-x_k)^{2\tau+1}\vee(y-y_k)^{2\tau+1}))+O(u^2(x-x_k)^2(y-y_k)).
\end{align*}}
And we will use Taylor expansion with the cross term and use the Laplace method \ref{laplace} to deal with the principle terms to get that when $u=O(1)$ the following holds
\tny{\begin{align*}
    &\frac{\int_{\bb R^2}\exp\lef(-k(\mca H_{0,k}(u,x,y)-\mca H_{0,k}(u,x_k,y_k))\rig)dxdy}{\int_{\bb R}\exp\lef(-k(\mca H_{1,k}(x)-\mca H_{1,k}(x_{1,k}))\rig)dx\int_{\bb R}\exp\lef(-k(\mca H_{2,k}(x)-\mca H_{2,k}(x_{2,k}))\rig)dy}\\
    &=\frac{\int_{\bb R}\exp(-k\sum_{i=1}^\infty\nabla^{i}_x\mca H_{0,k}(x-x_k)^i)dx\int_{\bb R}\exp(-k\sum_{i=1}^\infty\nabla^{i}_y\mca H_{0,k}(y-y_k)^i)dy}{\int_{\bb R}\exp\lef(-k(\mca H_{1,k}(x)-\mca H_{1,k}(x_{1,k}))\rig)dx\int_{\bb R}\exp\lef(-k(\mca H_{2,k}(x)-\mca H_{2,k}(x_{2,k}))\rig)dy}(1+O(u^2)k^{\frac{\tau-1}{\tau}})\\
    &=1+O(u)+O(u^2k^{\frac{\tau-1}{\tau}})+O\lef(\frac{1}{k}\rig).
\end{align*}}
Then we analyze the term $\eta(\bfa h)$, noticing that by Taylor expansion we will have (recall that $x_1^*=0=x_2^*$).
\tny{\begin{align}\label{etatcritical}
   \eta(\bfa h)&=\frac{1}{k+r}\bl \sum_{i=r+1}^k\tanh(h_i^\prime+\sqrt{\theta_1}x_{1,k})\tanh(h^\prime_i+\sqrt{\theta_1}x_{2,k})+\sum_{i=1}^r\tanh(h_i^\prime)\tanh(h_i^\prime+\sqrt{\theta_1}x_{1,k})\nnb\\
   &+\sum_{i=k+1}^{k+r}\tanh(h_i^\prime)\tanh(h_i^\prime+\sqrt{\theta_1}x_{2,k})\br.
\end{align}}
We notice that \eqref{etatcritical} might have different convergence rate that depends on $r$. First we consider $c\gtrsim k^{-\frac{\tau-2}{2\tau-1}}$ and we will have
\tny{\begin{align*}
    \eta(\bfa h)
    &=\frac{1}{k+r}\sum_{i=1}^{k+r}\tanh^2(h_i^\prime)+\frac{\sqrt{\theta_1}}{k+r}\sum_{i=1}^k\tanh(h_i^\prime)\sech^2(h_i^\prime)x_1-\frac{\theta_1}{k+r}\sum_{i=1}^k\tanh^2(h_i^\prime)\sech^2(h_i^\prime)x^2_1\nnb\\
    &-\frac{\theta_1}{k+r}\sum_{i=r+1}^{k+r}\tanh^2(h_i^\prime)\sech^2(h_i^\prime)x_2^2+\frac{\theta_1}{k+r}\sum_{i=r+1}^k\sech^4(h_i^\prime)x_1x_2+O(x_1^3\vee x_2^3).
\end{align*}}
Therefore, we use the fact that $\Vert k^{\frac{1}{4\tau-2}}x_1\Vert_{\psi_{4\tau-2}}\asymp \Vert k^{\frac{1}{4\tau-2}}x_2\Vert_{\psi_{4\tau-2}}\asymp 1$ and $\Vert k^{\frac{1}{2\tau-1}}x_1x_2\Vert_{2}\asymp 1$ , pick $q=\bb E[\tanh^2(\sqrt{\theta_1}x_1^*+h^\prime)]$ and $U\sim N\lef(0,{2\lambda}{(k+r)^{-\frac{4\tau-4}{2\tau-1}}}\rig)$. Similar to \eqref{gstrick} we will use H\"older's inequality to get
\begin{align}\label{gstrick1}
    \nu_0(\exp(\lambda& (k+r)^{\frac{2}{2\tau-1}}(R^D_{1,2}-q)^2))=\nu_0\lef(\bb E_U[\exp(U (k+r)(R^D_{1,2}-q))]\rig)\nnb\\
    &\leq\bb E_{\bfa h}\lef[\bb E_{U}\lef[(1+O(U)+O(U^2k^{\frac{\tau-1}{\tau}}))\exp\lef((k+r)\lef(\frac{1}{2}U^2+(\eta(\bfa h)-q)U\rig)\rig)\rig]\rig]\nnb\\
    &\leq\bb E_{\bfa h}\lef[\exp\lef(\lambda(1+o(1)) (k+r)^{\frac{2}{2\tau-1}}\lef(\eta(\bfa h)-q\rig)^2\rig)\rig](1+o(1)).
\end{align}
Then we will check that by \eqref{etatcritical}, there exists $C>0$ such that for $\lambda<C<\frac{1}{2}$ we will have 
\begin{align*}
    \nu_0&(\exp(\lambda (k+r)^{\frac{2}{2\tau-1}}(R^D_{1,2}-q)^2))<\infty.
\end{align*}
And then we consider the case of $c=o(k^{-\frac{\tau-2}{2\tau-1}})$,  by \eqref{2taux} we will have $|x_{2,k}-x_{1,k}|=O\lef(cx_{1,k}^{-(2\tau-1)}\wedge c^{\frac{1}{2\tau-1}}\rig)$. We also define
$    \mca H_{3,k}(x):=\frac{1}{k+r}\sum_{i=r+1}^k\log\cosh(\sqrt{\theta_1}x+h_i^\prime)
$. Then we notice that $\mca H_{3,k}^{(i)}\to 0$ for all $i\in[2\tau-1]$ and hence,
\ttny{\begin{align*}
    \eta(\bfa h)&=\frac{1}{k+r}\bl\sum_{i=1}^{k+r}\tanh^2(h_i^\prime)-2\sqrt{\theta_1}\sum_{i=r+1}^k\tanh(h_i^\prime)\sech^2(h_i^\prime)x_{1,k}+\ldots+\frac{H^{(2\tau)}_{3,k}(0)}{\theta_1(2\tau-2)!}x_{1,k}^{2\tau-2}+O(x_{1,k}^{2\tau-1})\br\\
    &+O\lef(c x^2_{1,k}\rig).
\end{align*}}
And we pick $q=\bb E[\tanh^2(h^\prime)]$ and $U\sim N(0,2\lambda (k+r)^{-\frac{2\tau}{2\tau-1}})$ and noticing that $\Vert k^{\frac{\tau-1}{2\tau-1}}(\eta(\bfa h)-q)\Vert_{\psi_2}<\infty$ there exists $C>0$ such that for all $\lambda<C$:
\tny{\begin{align}\label{gstrick2}
    \nu_0(\exp(\lambda(k+r)^{\frac{2\tau-2}{2\tau-1}}(R_{1,2}^D-q)^2))
    &\leq\bb E_{\bfa h}\lef[\exp\lef(\lambda(1+o(1)) (k+r)^{\frac{2\tau-2}{2\tau-1}}\lef(\eta(\bfa h)-q\rig)^2\rig)\rig](1+o(1))<\infty.
\end{align}}
\subsection{Proof of Lemma \ref{1rsbbound}}
The proof relies on the previous result in lemma \citep{guerra2001sum}, which utilizes the interpolation of $\{v_{\alpha}\}_{\alpha\geq 1} \sim PD(m,0) $ process. Here we define 
\begin{align*}
    \mca H_t(\bfa\sigma,\alpha,\bfa h):=\frac{\theta  \sqrt t}{\sqrt n}\sum_{i<j\leq n}g_{ij}\sigma_i\sigma_j +\theta \sqrt{1-t}\sum_{i\leq n}\sigma_iz_{i,\alpha}\sqrt{q}+\sum_{i=1}^kh_i\sigma_i.
\end{align*}

Using $\Theta_k=\{-1,-\frac{k-2}{k},\ldots,\frac{k-2}{k},1\}$, introducing $\bfa h^\prime_i=\bfa h$ for $i\in[k+1:n]$ and $\bfa h^\prime_i=\bfa h+\theta_1\mu$ for $i\in[k]$, notice that the following holds
\begin{align}\label{varphirelationship}
    \varphi_0(t)&=\frac{1}{n}\bb E\log\sum_{\alpha\geq 1}v_{\alpha}\sum_{\mu\in\Theta_k}\exp\bl-\frac{k\theta_1\mu^2}{2}\br\sum_{\bfa\sigma}\mbbm 1_{m_1=\mu}\exp\bl-\mca H_t(\bfa\sigma,\alpha,\bfa h)+\sum_{i\leq k}\sigma_i\mu\br\nnb\\
    &\leq\frac{1}{n}\bb E\log\sum_{\alpha\geq 1}v_{\alpha}\sum_{\mu\in\Theta_k}\exp\bl-\frac{k\theta_1\mu^2}{2}\br\sum_{\bfa\sigma}\exp\bl-\mca H_t(\bfa\sigma,\alpha,\bfa h^\prime)\br\nnb\\
    &\leq\frac{\log k}{n}+\ub{\sup_{\mu\in [-1,1]}\bl-\frac{c\theta_1\mu^2}{2}+\bb E\log\sum_{\alpha\geq 1}v_{\alpha}\sum_{\bfa\sigma}\exp\bl-\mca H_t(\bfa\sigma,\alpha,\bfa h^\prime)\br\br}_{=:\varphi(t)}.
\end{align}
Checking utilizing the standard results for Poission Dirichlet process that is given by (See, for example, Lemma 3.1 in \citep{panchenko2012sherrington})
\begin{align*}
    \bb E\bigg[\log\sum_{\alpha\geq 1}v_{\alpha}\sum_{\bfa\sigma}\exp(F(\bfa\sigma,\alpha))\bigg]=\bb E\bigg[\sum_{\bfa\sigma}\log\sum_{\alpha\geq 1}v_{\alpha}\exp(F(\bfa\sigma,\alpha))\bigg].
\end{align*}
Then, using lemma \ref{logtompower}, we will see that
\tny{\begin{align*}
    \varphi(0)=\log 2+c\sup_{\mu\in[-1,1]}\bl\frac{1}{m}\bb E\log\bb E^\prime\cosh^m Y_1-\frac{\theta_1\mu^2}{2}\br+\frac{1-c}{m}\bb E\log\bb E^\prime\cosh^mY_2.
\end{align*}}
Given a function $f$ on $\Sigma_n\times\bb N$, define
$    \la f\ra_t=\frac{1}{Z(t)}\sum_{\bfa\sigma,\alpha}v_{\alpha}f(\bfa\sigma,\alpha)\exp(-\mca H_t(\bfa\sigma,\alpha)).
$ We will check that
\ttny{\begin{align*}
    2\varphi^\prime(t)&=\frac{\theta}{n\sqrt t}\bb E\bigg\la \frac{1}{\sqrt n}\sum_{1\leq i<j\leq n}g_{ij}\sigma_i\sigma_j\bigg\ra_t-\frac{\theta}{n\sqrt{1-t}}\bb E\bigg\la\sum_{i\leq n}\sigma_iz_i\sqrt {q}\bigg\ra_t
    -\frac{\theta}{n\sqrt{1-t}}\bb E\bigg\la\sum_{i=1}^n\sigma_iz^\prime_{i,\alpha}\sqrt{q^\prime-q}\bigg\ra_t.
\end{align*}}
Then we will analyze the derivative, by integration by parts, the following holds (See also (13.27) of \citep{talagrand2011mean2})
\begin{align*}
    \varphi^\prime(t)\leq\frac{\theta^2}{4}\big((1-q^\prime)^2-m(q^{\prime 2}-q^2)\big).
\end{align*}
Then we will conclude that
\begin{align*}
    \varphi_0(1)\leq\frac{\log k}{n}+\varphi(1)\leq\frac{\log k}{n}+\frac{\theta^2}{4}(1-q^\prime)^2-\frac{\theta^2}{4}m(q^{\prime 2}-q^2)+\varphi(0).
\end{align*}
\subsection{Proof of Lemma \ref{lm4.3}}
For $z\sim N(0,1)$ independent of $h$ we define
\begin{align*}
    \theta_t:=\sqrt{t}\theta,\quad h_t=h+\sqrt{1-t}\theta\sqrt q.
\end{align*}
And the proof goes by observing the fact that
\begin{align*}
    \theta_tz\sqrt q+\theta_{1}\mu+h_t\overset{d}{=}\theta z\sqrt q+\theta_1\mu+h,\quad \theta_tz\sqrt q+h_t\overset{d}{=}\theta z\sqrt q+h,
\end{align*}
which implies that (we recover the implicit dependence of $pSK(\theta,\theta_1)$ with $h$ by writing it as $pSK(\theta,\theta_1,h)$)
\begin{align*}
    \psi(t)=pSK(\theta_{1},\theta_t,h_t),\quad\psi(t,u) = r_n(\theta_{1},\theta_t,h_t,u),\quad \eta(t,u)=t_{k}(\theta_{1},\theta_t,h_t,u).
\end{align*}
Then we use lemma \ref{intermediatelmpsk} to complete the proof.

\subsection{Proof of Lemma \ref{intermediatelmpsk}}
    The proof goes by first estimating the two quantity through the smart path method given by lemma \ref{rkpkub}. The bound of $r_n$ follows analogous to the proof of \citep{talagrand2011mean2} Proposition 13.6.6. And for the second inequality we note that by letting $m=1$, defining
    \tny{\begin{align*}
       \ca P_2(\theta_1,\theta,u_1,u_2,m):&=\log 2+c\bb E\lef[\log\cosh\lef(\theta z\sqrt {q}+\theta_1u_1+h\rig)\rig]+(1-c)\bb E[\log\cosh(\theta z\sqrt {q}+u_2+h)]\\
       &+\frac{\theta^2}{4}(1-q)^2-\frac{c\theta_1 u_1^2}{2},
    \end{align*}}
    we will immediately see that
    \begin{align*}
       pSK(\theta_1,\theta,h)=\sup_{u_1,u_2\in[-1,1]}\ca P_2(\theta_1,\theta,u_1,u_2,1).
    \end{align*}
    And
    \begin{align*}
        t_{1,k}(\theta_1,\theta,h,u)\leq\ca P_2(\theta_1,\theta,u,0,1),\quad t_{2,k}(\theta_1,\theta,h,u)\leq\sup_{\mu}\ca P_2(\theta_1,\theta,\mu,u,1).
    \end{align*}
    Therefore, we will see that at the high temperature
    \begin{align*}
        \frac{\pta^2\ca P_2}{\pta u_1^2}\bigg |_{u_1=0}=c\theta_1(-1+\theta_1\bb E[\sech^2(\theta z\sqrt q+h)])<0.
    \end{align*}
    At the low temperature,
    \begin{align*}
        \frac{\pta^2\ca P_2}{\pta u_1^2}\bigg |_{u_1=\mu}=\frac{\pta^2\ca P_2}{\pta u_1^2}\bigg |_{u_1=-\mu}=c\theta_1(-1+\theta_1\bb E[\sech^2(\theta z\sqrt q+\theta_1\mu+h)])<0.
    \end{align*}
    And at the critical temperature,
    \begin{align*}
        \frac{\pta^{2\tau}\ca P_2}{\pta u_1^{2\tau}}\bigg |_{u_1=0}<0.
    \end{align*}
\subsection{Proof of Lemma \ref{rkpkub}}
        The proof goes by noticing that defining $(z_{i,1},z_{i,2})$ be the independent copies of $(z_1,z_2)$ and independent copies $(z_{i,1}^\prime,z_{i,2}^\prime)$ and $(z_{i,\alpha,1}^\prime,z_{i,\alpha,2}^\prime)$ of the pair $(z_1^\prime,z_2^\prime)$. Then we define
        \tny{\begin{align*}
            -\mca H_{1,t}(\bfa\sigma^1,\bfa\sigma^2,\alpha,\bfa h):&=\theta\sqrt{\frac{t}{n}}\sum_{i<j}g_{ij}(\sigma_i^1\sigma_j^1+\sigma_i^2\sigma_j^2)+\theta\sqrt{1-t}\sum_{i\leq k}\sum_{j=1,2}\sigma_i^j(z_{i,j}\sqrt {q_1}+z_{i,\alpha,j}^\prime\sqrt{q_2-q_1})\\
            &+\sum_{i\leq k}h_i(\sigma_i^1+\sigma_i^2).
        \end{align*}}
        Then we will check that defining $\Theta=\{-1,-\frac{k-2}{k},\ldots,\frac{k-2}{k},1\}$, using similar derivation as in \eqref{varphirelationship} we will have 
        \ttny{\begin{align*}
            \varphi(t)&=\frac{1}{n}\bb E\log\sum_{\alpha\geq 1}v_{\alpha}\sum_{R_{1,2}=u}\exp\lef(-\mca H_t(\bfa\sigma^1,\bfa\sigma^2,\alpha,\bfa h)+\frac{\theta_1k}{2}(m_1^2+m_2^2)\rig)\\
            &=\frac{1}{n}\bb E\log\sum_{\alpha\geq 1}v_\alpha\sum_{\mu_1,\mu_2\in\Theta}\sum_{\bfa\sigma^1,\bfa\sigma^2}\mbbm 1_{m_1=\mu_1}\mbbm 1_{m_2=\mu_2}\mbbm 1_{R_{1,2}=u}\exp\lef(-\mca H_t(\bfa\sigma^1,\bfa\sigma^2,\alpha,\bfa h)+\frac{\theta_1k}{2}(m_1^2+m_2^2)\rig)\\
            &\leq\frac{2\log k}{n}+\sup_{\mu_1,\mu_2\in[-1,1]}\bigg\{-\theta_1c\frac{\mu_1^2+\mu_2^2}{2}+\frac{1}{n}\bb E\log\sum_{\alpha\geq 1}v_{\alpha}\sum_{\bfa\sigma}\exp\bl-\mca H_t(\bfa\sigma^1,\bfa\sigma^2,\alpha,\bfa h)+\sum_{i\leq k}\theta_1(\mu_1\sigma_i^1+\mu_2\sigma_i^2)\br\bigg\}\\
            &=\frac{2\log k}{n}+\sup_{\mu\in[-1,1]}\bigg\{-\theta_1c\mu^2+\frac{1}{n}\bb E\log\sum_{\alpha\geq 1}v_{\alpha}\sum_{\bfa\sigma}\exp\bl- \mca H_t(\bfa\sigma^1,\bfa\sigma^2,\alpha,\bfa h^\prime)\br\bigg\}.
        \end{align*}}
        where we define $\bfa h^\prime_{i}=\begin{cases}
            h_i+\theta_1(\mu_1\sigma_i^1+\mu_2\sigma_i^2)&\text{ if }i\in[k]\\
            h_i&\text{ if }i\in[k+1:n]
        \end{cases}$.
        And the rest of the proof for the upper bound on $r_n$ will follows from Theorem 13.5.1 in \citep{talagrand2011mean2}.
        Then we consider the first set of spins $\bfa\sigma$ and we denote $z_{j}:=z_{1,j}$,  $z_{\alpha,j}:=z_{1,\alpha,j}$ to be i.i.d. Gaussians. Define
        \begin{align*}
            Y_{\alpha,i}(u):=\theta z_i\sqrt{q_1}+\theta z_{\alpha,i}^\prime\sqrt{q_2-q_1}+\theta_1u+h.
        \end{align*}
        We prove the second bound on $p_k$ through the following definition:
        \sm{\begin{align*}
            \mca H_t(\bfa\sigma^1,\alpha,\bfa h):=\theta\sqrt{\frac{t}{n}}\sum_{i<j}g_{ij}\sigma_i^1\sigma_j^1+\theta\sqrt{1-t}\sum_{i\leq n}\sigma^1_i(z\sqrt {q_1}+z_{i,\alpha}\sqrt{q_2-q_1})+\sum_{i\leq n}h_i\sigma_i^1+\frac{k\theta_1}{2}m_1^2.
        \end{align*}}
        Then we define
        \begin{align*}
            \varphi(t):&=\frac{1}{n}\bb E\log \sum_{\alpha\geq 1}v_{\alpha}\sum_{\bfa\sigma^1}\mbbm 1_{m_1=u}\exp\lef(\mca H_t(\bfa\sigma^1,\alpha,\bfa h)\rig).
        \end{align*}
        One will immediately see that under the above definition $\varphi(1)=p_k(\theta,\theta_1, h,u )$.
        Then, using lemma \ref{logtompower}, we get 
        \begin{align*}
            \varphi(0)&=\frac{1}{n}\bb E\log\sum_{\alpha\in\bb N}v_{\alpha}\bl\sum_{\bfa\sigma}\mbbm 1_{m_1=u}\exp\bl\sum_{i\leq k}\sigma_i^1Y_{\alpha,i}(u)+\sum_{i=k+1}^n\sigma_i^1Y_{\alpha,i}(0)-\frac{k\theta_1}{2}m_1^2\br\br\\
            &=-\frac{\theta_1u^2}{2}+\frac{1}{nm}\bb E\log\bb E^\prime\bl\sum_{\bfa\sigma}\mbbm 1_{m_1= u}\exp\bl\sum_{i\leq k}\sigma_i^1Y^\prime_{i}(u)+\sum_{i=k+1}^n\sigma_i^1Y_i(0)\br\br^m\\
            &\leq-\frac{\theta_1 u^2}{2}+\frac{1}{nm}\bb E\log\bb E^\prime\bl\sum_{\bfa\sigma}\exp\bl\sum_{i\leq k}\sigma_i^1Y^\prime_i(u)+\sum_{i=k+1}^n\sigma_i^1Y_i(0)\br\br^m\\
            &\leq \log 2 -\frac{\theta_1u^2}{2}+\frac{c}{m}\bb E\log\bb E^\prime\cosh^m(Y^\prime_i(u))+\frac{1-c}{m}\bb E\log\bb E^\prime\cosh^m(Y_i(0)).
        \end{align*}
        Defining $\la f(\bfa\sigma,\alpha) \ra_t:=\frac{1}{Z}\sum_{\alpha\geq 1}v_{\alpha}\sum_{\bfa\sigma}\mbbm 1_{m_1=u}f(\bfa\sigma,\alpha)\exp\lef(-\mca H_t(\bfa\sigma,\alpha,\bfa h)\rig)$.
        And by lemma \ref{pdprop}, the derivative satisfies:
        \sm{\begin{align*}
            2\varphi^\prime(t) = \frac{\theta}{\sqrt {tn}}\bb E\bigg\la\sum_{i<j}g_{ij}\sigma^1_i\sigma_j^1\bigg\ra_t-\frac{\theta}{\sqrt{1-t}}\bb E\bigg\la\sum_{i\leq n}\sigma_i^1z_i\sqrt{q_1}\bigg\ra_t-\frac{\theta}{\sqrt{1-t}}\bb E\bigg\la\sum_{i\leq n}\sigma_i^1z_{\alpha,i}\sqrt{q_2-q_1}\bigg\ra_t.
        \end{align*}}
        Using lemma \ref{pdprop} 
        $U_\alpha=W_{\alpha}=V_{\alpha}=\exp\lef(-\mca H_t(\bfa\sigma,\alpha)\rig)$, one immediately get for $\gamma$ a replica of $\alpha$,
        \begin{align*}
            \bb E\la\mbbm 1_{\alpha=\gamma}\ra_t=1-m.
        \end{align*}
        And by Gaussian integration by parts we will arrive at
        \tny{\begin{align*}
            2\varphi^\prime(t)&\leq\frac{\theta^2}{2}(1-\bb E\la R^2_{1,2}\ra_t)-\theta^2(q_2-q_1\bb E\la R_{1,2}\ra_t-(q_2-q_1)\bb E\la R_{1,2}\mbbm 1_{\alpha=\gamma}\ra_t)\\
            &=\frac{\theta^2}{2}\lef((1-2q_2)+q_1^2\bb E\la\mbbm 1_{\alpha=\gamma}\ra_t+q_2^2\bb E\la \mbbm 1_{\alpha=\gamma}\ra_t-\bb E\la(R_{1,2}-q)^2\mbbm 1_{\alpha\neq \gamma}\ra-\bb E\la(R_{1,2}-q_2)^2\mbbm 1_{\alpha=\gamma}\ra_t\rig)\\
            &\leq \frac{\theta^2}{2}\lef((1-q_2)^2-m(q_2^2-q_1^2)\rig).
        \end{align*}}
        Then, collecting pieces, we will have
        \begin{align*}
            t_n(\theta_1,\theta,h,u)&\leq \frac{\theta^2}{4}((1-q_2)^2-m(q_2^2-q_1^2))+\log 2-\frac{\theta_1u^2}{2}+\frac{c}{m}\bb E\log\bb E^\prime\cosh^m(Y_i(u))\\
            &+\frac{1-c}{m}\bb E\log\bb E^\prime\cosh^m(Y_i(0)).
        \end{align*}
\subsection{Proof of Lemma \ref{exponentialeq}}
        The proof goes by noticing that defining $(z_{i,1},z_{i,2})$ be the independent copies of $(z_1,z_2)$ and independent copies $(z_{i,1}^\prime,z_{i,2}^\prime)$ and $(z_{i,\alpha,1}^\prime,z_{i,\alpha,2}^\prime)$ of the pair $(z_1^\prime,z_2^\prime)$. Then we define
        \tny{\begin{align*}
            -\mca H_{1,t}(\bfa\sigma^1,\bfa\sigma^2,\alpha,\bfa h):&=\theta\sqrt{\frac{t}{n}}\sum_{i<j}g_{ij}(\sigma_i^1\sigma_j^1+\sigma_i^2\sigma_j^2)+\theta\sqrt{1-t}\sum_{i\leq k}\sum_{j=1,2}\sigma_i^j(z_{i,j}\sqrt {q_1}+z_{i,\alpha,j}^\prime\sqrt{q_2-q_1})\\
            &+\sum_{i\leq k}h_i(\sigma_i^1+\sigma_i^2).
        \end{align*}}
        Then we will check that defining $\Theta=\{-1,-\frac{k-2}{k},\ldots,\frac{k-2}{k},1\}$, using similar derivation as in \eqref{varphirelationship} we will have 
        \ttny{\begin{align*}
            \varphi(t)&=\frac{1}{n}\bb E\log\sum_{\alpha\geq 1}v_{\alpha}\sum_{R_{1,2}=u}\exp\lef(-\mca H_t(\bfa\sigma^1,\bfa\sigma^2,\alpha,\bfa h)+\frac{\theta_1k}{2}(m_1^2+m_2^2)\rig)\\
            &=\frac{1}{n}\bb E\log\sum_{\alpha\geq 1}v_\alpha\sum_{\mu_1,\mu_2\in\Theta}\sum_{\bfa\sigma^1,\bfa\sigma^2}\mbbm 1_{m_1=\mu_1}\mbbm 1_{m_2=\mu_2}\mbbm 1_{R_{1,2}=u}\exp\lef(-\mca H_t(\bfa\sigma^1,\bfa\sigma^2,\alpha,\bfa h)+\frac{\theta_1k}{2}(m_1^2+m_2^2)\rig)\\
            &\leq\frac{2\log k}{n}+\sup_{\mu_1,\mu_2\in[-1,1]}\bigg\{-\theta_1c\frac{\mu_1^2+\mu_2^2}{2}+\frac{1}{n}\bb E\log\sum_{\alpha\geq 1}v_{\alpha}\sum_{\bfa\sigma}\exp\bl-\mca H_t(\bfa\sigma^1,\bfa\sigma^2,\alpha,\bfa h)+\sum_{i\leq k}\theta_1(\mu_1\sigma_i^1+\mu_2\sigma_i^2)\br\bigg\}\\
            &=\frac{2\log k}{n}+\sup_{\mu\in[-1,1]}\bigg\{-\theta_1c\mu^2+\frac{1}{n}\bb E\log\sum_{\alpha\geq 1}v_{\alpha}\sum_{\bfa\sigma}\exp\bl- \mca H_t(\bfa\sigma^1,\bfa\sigma^2,\alpha,\bfa h^\prime)\br\bigg\}.
        \end{align*}}
        where we define $\bfa h^\prime_{i}=\begin{cases}
            h_i+\theta_1(\mu_1\sigma_i^1+\mu_2\sigma_i^2)&\text{ if }i\in[k]\\
            h_i&\text{ if }i\in[k+1:n]
        \end{cases}$.
        And the rest of the proof for the upper bound on $r_n$ will follows from Theorem 13.5.1 in \citep{talagrand2011mean2}.
        Then we consider the first set of spins $\bfa\sigma$ and we denote $z_{j}:=z_{1,j}$,  $z_{\alpha,j}:=z_{1,\alpha,j}$ to be i.i.d. Gaussians. Define
        \begin{align*}
            Y_{\alpha,i}(u):=\theta z_i\sqrt{q_1}+\theta z_{\alpha,i}^\prime\sqrt{q_2-q_1}+\theta_1u+h.
        \end{align*}
        We prove the second bound on $p_k$ through the following definition:
        \sm{\begin{align*}
            \mca H_t(\bfa\sigma^1,\alpha,\bfa h):=\theta\sqrt{\frac{t}{n}}\sum_{i<j}g_{ij}\sigma_i^1\sigma_j^1+\theta\sqrt{1-t}\sum_{i\leq n}\sigma^1_i(z\sqrt {q_1}+z_{i,\alpha}\sqrt{q_2-q_1})+\sum_{i\leq n}h_i\sigma_i^1+\frac{k\theta_1}{2}m_1^2.
        \end{align*}}
        Then we define
        \begin{align*}
            \varphi(t):&=\frac{1}{n}\bb E\log \sum_{\alpha\geq 1}v_{\alpha}\sum_{\bfa\sigma^1}\mbbm 1_{m_1=u}\exp\lef(\mca H_t(\bfa\sigma^1,\alpha,\bfa h)\rig).
        \end{align*}
        One will immediately see that under the above definition $\varphi(1)=p_k(\theta,\theta_1, h,u )$.
        Then, using lemma \ref{logtompower}, we get 
        \begin{align*}
            \varphi(0)&=\frac{1}{n}\bb E\log\sum_{\alpha\in\bb N}v_{\alpha}\bl\sum_{\bfa\sigma}\mbbm 1_{m_1=u}\exp\bl\sum_{i\leq k}\sigma_i^1Y_{\alpha,i}(u)+\sum_{i=k+1}^n\sigma_i^1Y_{\alpha,i}(0)-\frac{k\theta_1}{2}m_1^2\br\br\\
            &=-\frac{\theta_1u^2}{2}+\frac{1}{nm}\bb E\log\bb E^\prime\bl\sum_{\bfa\sigma}\mbbm 1_{m_1= u}\exp\bl\sum_{i\leq k}\sigma_i^1Y^\prime_{i}(u)+\sum_{i=k+1}^n\sigma_i^1Y_i(0)\br\br^m\\
            &\leq-\frac{\theta_1 u^2}{2}+\frac{1}{nm}\bb E\log\bb E^\prime\bl\sum_{\bfa\sigma}\exp\bl\sum_{i\leq k}\sigma_i^1Y^\prime_i(u)+\sum_{i=k+1}^n\sigma_i^1Y_i(0)\br\br^m\\
            &\leq \log 2 -\frac{\theta_1u^2}{2}+\frac{c}{m}\bb E\log\bb E^\prime\cosh^m(Y^\prime_i(u))+\frac{1-c}{m}\bb E\log\bb E^\prime\cosh^m(Y_i(0)).
        \end{align*}
        Defining $\la f(\bfa\sigma,\alpha) \ra_t:=\frac{1}{Z}\sum_{\alpha\geq 1}v_{\alpha}\sum_{\bfa\sigma}\mbbm 1_{m_1=u}f(\bfa\sigma,\alpha)\exp\lef(-\mca H_t(\bfa\sigma,\alpha,\bfa h)\rig)$.
        And by lemma \ref{pdprop}, the derivative satisfies:
        \sm{\begin{align*}
            2\varphi^\prime(t) = \frac{\theta}{\sqrt {tn}}\bb E\bigg\la\sum_{i<j}g_{ij}\sigma^1_i\sigma_j^1\bigg\ra_t-\frac{\theta}{\sqrt{1-t}}\bb E\bigg\la\sum_{i\leq n}\sigma_i^1z_i\sqrt{q_1}\bigg\ra_t-\frac{\theta}{\sqrt{1-t}}\bb E\bigg\la\sum_{i\leq n}\sigma_i^1z_{\alpha,i}\sqrt{q_2-q_1}\bigg\ra_t.
        \end{align*}}
        Using lemma \ref{pdprop} 
        $U_\alpha=W_{\alpha}=V_{\alpha}=\exp\lef(-\mca H_t(\bfa\sigma,\alpha)\rig)$, one immediately get for $\gamma$ a replica of $\alpha$,
        \begin{align*}
            \bb E\la\mbbm 1_{\alpha=\gamma}\ra_t=1-m.
        \end{align*}
        And by Gaussian integration by parts we will arrive at
        \tny{\begin{align*}
            2\varphi^\prime(t)&\leq\frac{\theta^2}{2}(1-\bb E\la R^2_{1,2}\ra_t)-\theta^2(q_2-q_1\bb E\la R_{1,2}\ra_t-(q_2-q_1)\bb E\la R_{1,2}\mbbm 1_{\alpha=\gamma}\ra_t)\\
            &=\frac{\theta^2}{2}\lef((1-2q_2)+q_1^2\bb E\la\mbbm 1_{\alpha=\gamma}\ra_t+q_2^2\bb E\la \mbbm 1_{\alpha=\gamma}\ra_t-\bb E\la(R_{1,2}-q)^2\mbbm 1_{\alpha\neq \gamma}\ra-\bb E\la(R_{1,2}-q_2)^2\mbbm 1_{\alpha=\gamma}\ra_t\rig)\\
            &\leq \frac{\theta^2}{2}\lef((1-q_2)^2-m(q_2^2-q_1^2)\rig).
        \end{align*}}
        Then, collecting pieces, we will have
        \begin{align*}
            t_n(\theta_1,\theta,h,u)&\leq \frac{\theta^2}{4}((1-q_2)^2-m(q_2^2-q_1^2))+\log 2-\frac{\theta_1u^2}{2}+\frac{c}{m}\bb E\log\bb E^\prime\cosh^m(Y_i(u))\\
            &+\frac{1-c}{m}\bb E\log\bb E^\prime\cosh^m(Y_i(0)).
        \end{align*}

\subsection{Proof of Lemma \ref{inductivehypo}}
        We check that only even $j$ needs to be considered. Picking $\ell$ such that $2\ell=j+1$. Then using H\"older's inequality we will have
        \begin{align*}
            \min_{\mu\in\ca U}\{\nu(|m-\mu|^j)\}\leq\min_{\mu\in\ca U}\{(\nu(m-\mu)^{2\ell})^{j/2\ell}\}\leq \lef(\frac{C\ell}{k}\rig)^{j/2} = \lef(\frac{C(j+1)}{2k}\rig)^{j/2}.
        \end{align*}
        Then we consider the second inequality regarding the cavity quantity, note that by the upper bound, we will have
        \ttny{\begin{align*}
            \min_{\mu\in\ca U}\{\nu((m^-_k-\mu)^{2r})\}&\leq\sum_{0\leq j\leq 2k}\binom{2r}{j}\frac{1}{k^{2r-j}}\min_{\mu\in\ca U}\{\nu(|m-\mu|^j)\}\leq\sum_{0\leq j\leq 2k}\binom{2r}{j}\frac{1}{k^{2r-j}}\lef(\frac{C(j+1)}{2k}\rig)^{j/2}\\
            &\leq\sum_{0\leq j\leq 2k}\binom{2r}{j}\frac{1}{k^{2r-j}}\lef(\frac{C(2k+1)}{2k}\rig)^{j/2}=\lef(\frac{C(2r+1)}{2k}\rig)^{r}\lef(1+\sqrt{\frac{2}{rC(2r+1)}}\rig)^{2r}.
        \end{align*}}
        And using $1+x\leq\exp x$ we complete the proof with $C\geq 4$.
\subsection{Proof of Lemma \ref{covariance}}
    The proof follows from the cavity construction. We notice that by the exponential inequality in lemma \ref{largcliquetailbound}, we will have
    \sm{\begin{align*}
        \nu(m_1\tde m_1)&=\nu(m_1\xi_1)=\nu_{0,2}(m_1\xi_1)+\nu_{0,2}^\prime(m_1\xi_1)+O\bl\frac{1}{n\sqrt k}\br=\nu_{0,2}^\prime(m_1\xi_1)+O\bl\frac{1}{n\sqrt k}\br\\
        &=\theta^2\nu_{0,2}^\prime(m_1\xi_1\xi_2\xi_3(\ca R^-_{2,3}-q))-\theta^2\nu_{0,2}^\prime(m_1\xi_2(\ca R^-_{1,2}-q))+O\bl\frac{1}{n\sqrt k}\br=O\bl\frac{1}{n\sqrt k}\br.
    \end{align*}}



\section{Auxilliary Lemmas}
\begin{lemma}[Gaussian Integration by Parts]\label{gtrbt}
        Given i.i.d. standard Gaussian $g, z_1,\ldots, z_n$ and a smooth function $F:\bb R^n\to \bb R$ that satisfies the \emph{moderate growth condition} defined by
        \begin{align*}
            \lim_{\Vert \bfa x\Vert_2\to\infty}|F(\bfa x)|\exp(-a\Vert \bfa x\Vert_2^2)=0
        .\end{align*}
        for all $a>0$ we have
        \begin{align*}
            \bb E[gF(z_1,\ldots, z_n)]=\sum_{\ell\leq n}\bb E[gz_{\ell}]\bb E\frac{\pta F}{\pta z_{\ell}}(z_1,\ldots,z_n)
        .\end{align*}
\end{lemma}
\begin{lemma}[Gaussian Interpolation Lemma ]\label{gil}
    Define $\bfa u,\bfa v\in\bb R^m$ to be two Gaussian process indexed by $[M]$, and define
    \begin{align*}
        u_i(t) = \sqrt t u_i+\sqrt{1-t} v_i
    .\end{align*}
    and note that $\bfa u =\bfa u(1)$ and $\bfa v =\bfa u(0)$, define function $\varphi(t)=\bb E[F(\bfa u(t))]$ with $F:\bb R^M\to\bb R$ being second order differentiable with $F$ satisfying conditions in \ref{gtrbt} we have
    \begin{align*}
        \varphi^\prime(t)=\frac{1}{2}\sum_{i,j}(\bb E[u_iu_j]-\bb E[v_iv_j])\bb E\lef[\frac{\pta^2 F}{\pta x_i\pta x_j}(\bfa u(t))\rig]
    .\end{align*}
\end{lemma}
\begin{lemma}[Multivariate Laplace Method]\label{laplace}
    Suppose we are given r.v.s. $\bfa h\in\bb R^d$, parameters $\bfa s\in\bb S\subset\bb R^d$ and $\{\Gamma_n( \bfa s,\bfa h)\}$ is a family of random variables in $\Omega$ with $\Gamma_n$ infinitely differentiable w.r.t. $\bfa s$. Furthermore, let us assume that $\Gamma_n$ has unique global minimum almost surely for all $n\in\bb N$ within $\bb S$, and the following are satisfied:
    \begin{enumerate}
        \item  There exists $C(\bfa h)>0$, independent of $n$ and real $\tau$ such that almost surely
        \begin{align}\label{laplace2}
            \exp\lef(-\Gamma(\bfa s,\bfa h)\rig)\leq C(\bfa h)\exp\lef(-\Vert \bfa s\Vert_2^2/2+\tau\Vert\bfa s\Vert_1\rig)
        \end{align} uniformly on compact sets in $\bb R$.
        \item We have  almost surely:
        \begin{align}\label{laplace3}
            \int_{\bb S}\exp(-\Gamma(\bfa s,\bfa h))d\bfa s:=\int_{S_1}\cdots\int_{S_n} \exp(-\Gamma(\bfa s,\bfa h))\prod_{i\in[d]}ds_i<\infty
        .\end{align}
    \end{enumerate}
    Then, we will have almost surely there exists random variables $a_1(\bfa h),\ldots a_M(\bfa h)$ for all  $M\in\bb N$ such that
    \tny{\begin{align*}
    \int_{\bb S}\exp(-n\Gamma_n(\bfa s,\bfa h))d\bfa s\sim\exp(-n\Gamma_n(\bfa s_n^*,\bfa h))\det\lef(\frac{n\nabla^2\Gamma_n(\bfa s^*_n,\bfa h)}{2\pi}\rig)^{-1/2}\lef(1+\frac{a_1(\bfa h)}{n}+\ldots+\frac{a_M(\bfa h)}{n^M}\rig) 
    ,\end{align*}}
    where $\nabla$ only take derivative w.r.t. $\bfa s$.
\end{lemma}
    The proof goes by first slicing the integral into two parts denoted by $\bfa V_n(\delta):=\{\bfa s:\Vert \bfa s-\bfa s_n^*\Vert_2\leq\delta\}$ that contains $\bfa s_n^*:=\argmin_{\bfa s\in\bb R^d}\Gamma_n(\bfa s,\bfa h)$ and let $\bfa V^c(\delta)$ be its complement. Note that there exists $\epsilon>0$ such that
    \begin{align*}
        \inf_{\bfa s\in\bfa V^c(\delta)}\Gamma_n(\bfa s,\bfa h)-\inf_{\bfa s\in\bb R^d}\Gamma_n(\bfa s,\bfa h)\geq\epsilon
    .\end{align*}
Hence, using \eqref{laplace2} and \eqref{laplace3} we note that
\sm{\begin{align*}
        \exp(n\Gamma_n(\bfa s_n^*,\bfa h))&\int_{\bfa V^c(\delta)}\exp(-n\Gamma_n(\bfa s,\bfa h))d\bfa s\\
        &=\exp(n\Gamma_n(\bfa s_n^*,\bfa h))\int_{\bfa V^c(\delta)}\exp\lef(-(n-1)\Gamma_n(\bfa s,\bfa h)\rig)\exp\lef(-\Gamma_n(\bfa s,\bfa h)\rig)d\bfa s\\
        &\leq\exp\lef(n\Gamma_n(\bfa s_n^*,\bfa h)-(n-1)\inf_{\bfa s\in\bfa V^c(\delta)}\Gamma_n(\bfa s^*,\bfa h)\rig)\int_{\bfa V^c(\delta)}\exp\lef(-\Gamma_n(\bfa s,\bfa h) \rig)d\bfa s\\
        &\leq O(\exp(-n\epsilon))
    .\end{align*}}

    Then we review in the following an important fact and the divergence theorem in vector calculus.
    \begin{fact}\label{fact2}
        Let $\bfa 0$ lie in the interior of $D\subset \bb R^d$. Then as $\lambda\to\infty$ we will have
        \begin{align*}
            \int_{D}\exp\lef(-\frac{\lambda}{2}\bfa\xi^\top\bfa\xi\rig)d\bfa\xi =\lef(\frac{2\pi}{\lambda}\rig)^{d/2} +o(\lambda^{-m})
        \end{align*}
        for all  $m\in\bb N$.
    \end{fact}
    \begin{theorem}[Divergence Theorem]\label{diverge} Suppose $\bfa D$ is a subset of  $\bb R^d$ with  $\bfa D$ a compact space with piecewise smooth boundary  $\bfa S=\partial \bfa D$. If $\bfa F$ is a continuously differentiable vector field defined on a neighborhood of $\bfa D$ then
    \begin{align*}
        \int_{\bfa D}(\nabla\cdot \bfa F)dV=\oint_{\bfa S}(\bfa F\cdot\bfa n)dS
    .\end{align*}
    where $\bfa n$ is the unit outward normal vector to $\bfa S$ and $d S$ is the differential element on the hypersurface $\bfa S$.
    
    Changing $\bfa F$ to $\bfa Fg$ for some smooth scalar function  $g$ we will have
    \begin{align*}
        \int_{\bfa D}\lef(\bfa F\cdot\nabla g+g\nabla\cdot\bfa F\rig)dV=\oint_{\bfa S}g\bfa F\cdot \bfa ndS
    .\end{align*}
    \end{theorem}
    
    The next step is to consider what will lie in $\bfa V_n(\delta_1)$.  The proof strategy follows from \citep{barndorff1989asymptotic} and \citep{bleistein1975asymptotic}. By Taylor expansion there exists $\delta_2>0$ sufficiently small such that for all  $\bfa s\in\bfa V_n(\delta_2)$ we will have
    \begin{align*}
        \Gamma_{n}(\bfa s,\bfa h) - \Gamma_n(\bfa s_n^*,\bfa h)&=\frac{1}{2}(\bfa s-\bfa s_n^*)^\top\nabla^2 \Gamma_n(\bfa s_n^*,\bfa h)(\bfa s-\bfa s_n^*)+o\lef(\Vert \bfa s-\bfa s^*_n\Vert_2^2\rig)\\
        &=\frac{1}{2}\bfa z^\top\bfa z +o(\Vert\bfa z\Vert_2^2)
    .\end{align*}
    where  $\bfa z:=\lef(\nabla^2\Gamma_n(\bfa s_n^*,\bfa h)\rig)^{1/2}(\bfa s-\bfa s_n^*)$.
    Then we can introduce $\bfa m:\bb R^d\to\bb R^d$ such that $m_i(\bfa z)= z_i+o(z_i)$ as $z_i\to 0$ and satisfying
    \begin{align*}
        \Gamma_n(\bfa s,\bfa h)-\Gamma_n(\bfa s^*_n,\bfa h)&=\frac{1}{2}\bfa m^\top(\bfa s)\bfa m(\bfa s)
    .\end{align*}
    Defining the function $ G_0(\bfa m):=\ca J(\bfa m)=\frac{\pta (s_1,\ldots,s_d)}{\pta (m_1,\ldots,m_d)}$ to be the Jacobian at $\bfa s$ and we note that
$  \ca J(\bfa 0)=\lef|\det(\nabla^2\Gamma_n(\bfa s_n^*,\bfa h))\rig|^{-1/2}
$   Introducing  $\bfa D$ to be the image of $\bfa V_n(\delta)$ under the two round of change of variables and $\bfa S=\pta\bfa D$. Therefore with the above preparation we can write the integral as:
\begin{align*}
        &\int_{\bfa V_n(\delta_1)}\exp\lef(-n\Gamma_n(\bfa s,\bfa h)\rig)d\bfa s = \exp(-n\Gamma_n(\bfa s_n^*,\bfa h))\int_{\bfa D} G_0(\bfa m)\exp\lef(-\frac{n}{2}\bfa m^\top\bfa m\rig)d\bfa m
.\end{align*}
Note that there exists a function $\bfa H_0:\bb R^d\to\bb R^d$ such that $ G_0(\bfa m)= G_0(\bfa 0)+\bfa m^\top\bfa H_0(\bfa m)$.
We will then use theorem \ref{diverge} to get 
\begin{align*}
    &\int_{\bfa V_n(\delta_1)}\exp\lef(-n\Gamma_n(\bfa s,\bfa h)\rig)d\bfa s=\exp(-n\Gamma_n(\bfa s_n^*,\bfa h))\bigg[\int_{\bfa D} G_0(\bfa 0)\exp\lef(-\frac{n}{2}\bfa m^\top\bfa m\rig)d\bfa m\\
    &-\frac{1}{n}\int_{\bfa S}(\bfa H_0(\bfa m)\cdot\bfa n)\exp\lef(-\frac{n}{2}\bfa m^\top\bfa m\rig)dS+\frac{1}{n}\int_{\bfa D} G_1(\bfa m)\exp\lef(-\frac{n}{2}\bfa m^\top\bfa m\rig)d\bfa m\bigg]
.\end{align*}
And we can do the above process recursively and get
\begin{align*}
    &\int_{\bfa V_n(\delta_1)}\exp\lef(-n\Gamma_n(\bfa s,\bfa h)\rig)d\bfa s=\exp(-n\Gamma_n(\bfa s_n^*,\bfa h))\bigg[\sum_{j=0}^{M} G_j(\bfa m)\int_{\bfa D}\exp\lef(-\frac{n}{2}\bfa m^\top\bfa m\rig)d\bfa m\\
    &-\frac{1}{n^M}\int_{\bfa D}G_M(\bfa s)\exp\lef(-\frac{n}{2}\bfa m^\top\bfa m\rig)d\bfa s\bigg]
.\end{align*}
since we note that the boundary integral is exponentially small almost surely according to \ref{laplace2} as $n\to\infty$ and could be ignored here. Note that $ G_j$ is defined recursively as
\begin{align*}
    G_j(\bfa m):=G_j(\bfa 0)+\bfa m^\top\cdot \bfa H_j(\bfa m),\\
    G_{j+1}(\bfa m):=\nabla\cdot \bfa H_j(\bfa m)
.\end{align*}
Further notice that by \ref{laplace2} and \ref{laplace3} we can check that almost surely:
\begin{align*}
    \lef|\frac{1}{n^M}\int_{D}\exp\lef(-\frac{n}{2}\bfa m^\top\bfa m\rig)G_M(\bfa m)d\bfa m\rig|=O\lef(\frac{1}{n^M}\rig)
.\end{align*}
Together with the fact \ref{fact2} we can see that almost surely:
\sm{\begin{align*}
    \int_{\bfa V_n(\delta_1)}\exp(-n\Gamma_n(\bfa s,\bfa h))d\bfa s=\exp\lef(-n\Gamma_n(\bfa s_n^*,\bfa h)\rig)\lef(\frac{2\pi}{n}\rig)^{d/2}\bl\sum_{j\in[M-1]}\frac{G_j(\bfa 0)}{n^j}+O\lef(\frac{1}{n^M}\rig)\br
.\end{align*}}
Note that $G_j(\bfa 0)$ are functions of $\bfa h$ we complete the proof by defining $a_k(\bfa h)=\frac{G_k(\bfa 0)}{G_0(\bfa 0)}$.
\begin{align*}
.\end{align*}

\section{Poisson Dirichlet Process}
Here we review a few properties of the Poisson Dirichlet Process used in the proof of section \ref{sect4.2}

\begin{definition}[Poisson Point Process]
Consider a positive measure on $\Omega$ of finite total mass $|\mu|$ defined by the Lebesgue integral over $\Omega$. Assume $\mu$ has no atoms.
    We call a random finite subset $\Pi$ a Poisson point process with intensity measure $\mu$ if it has the following properties:
    \begin{enumerate}
        \item $|\Pi|$ is a Poisson r.v. with expectation $|\mu|$.
        \item Given $|\Pi|=k$, $\Pi\overset{d}{=}\{X_1,\ldots, X_k\}$ with $X_i$ being i.i.d. r.v.s. with law $\frac{\mu}{|\mu|}$.
    \end{enumerate}
\end{definition}
Then we are ready to define the Poisson Dirichlet process.
\begin{definition}[Poisson Dirichlet Process \citep{pitman1997two}]
    For $\alpha\in[0,1)$ and $\theta>-\alpha$, suppose that a probability measure $\bb P_{\alpha,\theta}$ governs independent random variables $\tde Y_n$ such that $\tde Y_n$ has $Beta(1-\alpha,\theta+n\alpha)$ distribution. For $n\geq 2$, let
    \begin{align*}
        \tde V_1=\tde Y_1,\quad \tde V_n=(1-\tde Y_1)\ldots(1-\tde Y_{n-1})\tde Y_n,
    \end{align*}
    and let $V_1\geq V_2\geq\cdots$ be the non-increasing ordered values of $\tde V_n$. Define the Poisson Dirichlet distribution with parameters $(\alpha, \theta)$ ($PD(\alpha,\theta)$) to be the distribution of $\{V_n\}_{n\geq 1}$.
\end{definition}
And one can simply checked that a $PD(m,0)$ can be alternatively constructed as follows: Define a Poisson point process $\Pi$ with intensity measure $\mu_m:=x^{-m-1}$ for $m\in(0,1)$ on $\Omega=\bb R^+$. Then order the elements in $\Pi$ by a decreasing order as $\{u_{\alpha}\}_{\alpha\geq 1}$. The process $\{v_{\alpha}\}_{\alpha\geq 1}\sim PD(m,0)$ can be alternatively defined as
\begin{align*}
    v_{\alpha}:=\frac{u_{\alpha}}{\sum_{\gamma\geq 1}u_{\gamma}}.
\end{align*}
A useful property of this property is stated as follows:
\begin{lemma}[\citep{talagrand2010mean}]\label{logtompower}
    Consider $0<m<1$. Consider i.i.d. copies $(V_{\alpha})_{\alpha\geq 1}$ of a r.v. $V>0$ with $\bb E[V^m]<\infty$, that are independent of a sequence $(v_{\alpha})_{\alpha\geq 1}$ of distribution $PD(m,0)$. Then
    \begin{align*}
        \bb E\bigg[\log\sum_{\alpha\geq 1}v_{\alpha}V_{\alpha}\bigg]=\frac{1}{m}\log\bb E[V^m].
    \end{align*}
\end{lemma}
and
\begin{lemma}[\citep{talagrand2010mean}]\label{pdprop}
    Consider $0<m<1$, a triple $(U,V,W)$ of r.v.s. with $V\geq 1$ and assume that $\bb E[V^m]<\infty$, $\bb E[U^2]+\bb E[W^2]<\infty$. Consider independent copies $(U_\alpha, V_\alpha, W_{\alpha})$ of this triple, which are independent of a sequence $(v_{\alpha})_{\alpha\geq 1}$ with distribution $PD(m,0)$. Then we have
    \begin{align*}
        \bb E\frac{\sum_{\alpha\geq 1}v_{\alpha}U_{\alpha}}{\sum_{\alpha\geq 1}v_{\alpha}V_{\alpha}}&=\frac{\bb E UV^{m-1}}{\bb E V^m},\\
        \bb E\frac{\sum_{\alpha\geq 1}v_{\alpha}^2U_{\alpha}W_{\alpha}}{(\sum_{\alpha\geq 1}v_{\alpha} V_{\alpha})^2}&=(1-m)\frac{\bb EUWV^{m-2}}{\bb E V^m},\\
        \bb E\frac{\sum_{\alpha\neq\gamma}v_{\alpha}v_{\gamma}U_{\alpha} W_{\gamma}}{(\sum_{\alpha\geq 1}v_{\alpha}V_{\alpha})^2}&=m\frac{\bb EUV^{m-1}\bb EWV^{m-1}}{(\bb E V^m)^2}.
    \end{align*}
\end{lemma}
    
       

\end{document}